\documentclass[12pt,twoside,leqno]{amsart}
\usepackage{amssymb,amsbsy,amsmath,amsfonts,amssymb,amscd,times,graphics,color,wasysym}
\newcommand{\C}                 {\mathbb{C}}

\newcommand{\K}                 {\mathbb{K}}
\newcommand{\R}                 {\mathbb{R}}

\newtheorem{theorem}{Theorem}
\newtheorem{lemma}{Lemma}
\newtheorem{proposition}{Proposition}
\newtheorem{corollary}{Corollary}

\newtheorem{openproblem}{Open Problem}

\begin{document}





\title[Characterization of Newtonian free
particles]{Characterization of the Newtonian 
\\
free particle system 
\\
in $m\geqslant 2$ dependent variables}

\author{Jo\"el Merker}

\address{
CNRS, Universit\'e de Provence, LATP, UMR 6632, CMI, 
39 rue Joliot-Curie, F-13453 Marseille Cedex 13, France. \ \ 
\\
{\it Internet}:
{\tt http://www.cmi.univ-mrs.fr/$\sim$merker/index.html}}

\email{merker@cmi.univ-mrs.fr} 


\date{\number\year-\number\month-\number\day}

\begin{abstract}
We treat the problem of linearizability of a system of second order
ordinary differential equations. 
The criterion we provide has applications to nonlinear
Newtonian mechanics, especially in three-dimensional space.

Let $\K=\R$ or $\C$, let $x \in \K$, let $m \geqslant 2$, let
$y:=(y^1, \dots, y^m) \in \K^m$ and let
$$
y_{xx}^1
=
F^1 
\left( 
x, y, y_x 
\right), 
\ldots\dots, 
y_{xx}^m 
= 
F^m\left( x,y,y_x \right), 
$$ 
be a collection of $m$ analytic second order ordinary differential
equations, in general nonlinear. We obtain a new and applicable
necessary and sufficient condition in order that this system is
equivalent, under a point transformation
$$
(x, y^1,
\dots, y^m) \mapsto \left( X(x,y), Y^1(x,y),\dots, Y^m(x, y)
\right), 
$$
to the Newtonian free particle system $Y_{XX}^1 =\dots =Y_{XX}^m
=0$.

Strikingly, the explicit differential system that we obtain is of
first order in the case $m\geqslant 2$, whereas according to a classical
result due to Lie, it is of second order the case of a single equation
$(m=1)$.

\bigskip

\noindent
CNRS, Universit\'e de Provence, LATP, UMR 6632, CMI, 
39 rue Joliot-Curie, F-13453 Marseille Cedex 13, France.
\\
{\tt merker@cmi.univ-mrs.fr}
\\
{\tt http://www.cmi.univ-mrs.fr/$\sim$merker/index.html}

\end{abstract}

\maketitle

\begin{center}
\begin{minipage}[t]{13cm}
\baselineskip =0.35cm
{\scriptsize

\centerline{\bf Table of contents}

\smallskip

{\bf 1.~Introduction \dotfill 1.}

{\bf 2.~Proof of Lie's theorem \dotfill 7.}

{\bf 3.~Systems of second order ordinary
differential equations equivalent to free particles \dotfill 17.}

{\bf 4.~Compatibility conditions for the second auxiliary system
\dotfill 47.}

{\bf 5.~General point transformation of the free particle system
\dotfill 62.}

}\end{minipage}
\end{center}

\bigskip

\section*{\S1.~Introduction}

Several physically meaningful systems of ordinary differential
equations are of second order, as for instance the free particle in
$m$-dimensional space, the damped or undamped harmonic oscillator,
coupled or not, having constant time-dependent frequency or not, {\it
etc.} Such systems are ubiquitous in Newtonian Mechanics, in
Hamiltonian Dynamics and in General Relativity.

Two classical major problems are to classify these systems up to point
or contact equivalence (Lie's Grail) and to recognize when they
coincide with the Euler equations associated to a Lagrangian
(inverse variational problem). In small dimensions, complete results
hold (Lie, Tresse, Cartan; Darboux, Douglas). However, in arbitrary
dimension, both tasks quickly exceed the human as well as the digital
computer scale, due to the intrinsic complexity of the underlying
symbolic computations (explosion, swelling) and to the exponentially
increasing number of cases to be treated. We refer to Olver's
monograph~\cite{ ol1995} for a panorama of problems, methods and
results.

At least, as a first step in classification, with respect to
applications, there are both a mathematical and a physical interest in
determining concrete, explicit and applicable (``ready-made'')
criteria for a system of ordinary differential equation to be
equivalent, via a local point transformation, to a linear equation.

\subsection*{1.1.~Scalar equation}
In this respect, we remind the celebrated linearizability criterion
for a single equation, due to Lie.  Let $\K=\R$ of $\C$. Let $x\in \K$
and $y\in \K$. Consider a local second order ordinary differential
equation $y_{xx } = F( x, y, y_x)$, possibly nonlinear, with a locally
$\K$-analytic right-hand side\footnote{By borrowing techniques
developed in~\cite{ ma2003}, this theorem as well as the next both hold
under weaker smoothness assumptions, namely with a $\mathcal{ C}^2$ or
a $W_{\rm loc}^{ 1, \infty}$ right-hand side.}.

\def\thetheorem{1.2}\begin{theorem}
{\rm (\cite{lie1883}, pp.~362--365; \cite{ gtw1989}; \cite{ ol1995},
p.~406)} The following four conditions are equivalent{\rm :}

\begin{itemize}

\smallskip\item[{\bf (1)}]
$y_{ xx} = F (x, y, y_x)$ is equivalent under a local
point transformation $(x, y) \mapsto (X, Y)$ to the free particle
equation $Y_{ XX} = 0${\rm ;}

\smallskip\item[{\bf (2)}]
$y_{ xx} = F (x, y, y_x)$ is equivalent to some linear equation 
$Y_{ XX} = G_0 (X) +  G_1 (X) \, Y + H (X) \, Y_X${\rm ;}

\smallskip\item[{\bf (3)}]
the local Lie symmetry group of $y_{ xx} = F (x, y, y_x)$ is
eight-dimensional, locally isomorphic to the group ${\rm PGL}(3, \K)$
of all projective transformations of $P_2 (\K)$.

\smallskip\item[{\bf (4)}]
$F_{y_x y_x y_x y_x }=0$, or equivalently $F=G +y_x\,H+ (y_x)^2 \,
L+(y_x)^3\, M$, where $G, H, L, M$ are functions of $(x,y)$ that
satisfy{\rm :}
\begin{small}
$$
\aligned
0
& 
=
-
2\,G_{yy} 
+
\frac{4}{3}\,H_{xy}
-
\frac{1}{3}\,L_{xx}
+
2\,(G\,L)y
-
2\,G_x\,M
-
4\,G\,M_x
+
\frac{2}{3}
\, 
H\,L_x
-
\frac{4}{3}\,H\,H_y, 
\\
0 
&
= 
-
\frac{2}{3}\,H_{yy}
+
\frac{4}{3}\,L_{xy} 
-
2\,
M_{xx}
+
2\,G\,M_y
+
4\,G_y\,M
-
2\,(H\,M)_x
-
\frac{2}{3}\, 
H_y\,L
+
\frac{ 4}{3}\,L\,L_x.
\endaligned
$$
\end{small}
\end{itemize}
\end{theorem}

Section~2 of this paper is devoted to a detailed exposition of the
original proof of the equivalence between {\bf (1)} and {\bf (4)},
following~\cite{ lie1883}. In the contemporary literature, to the
author's knowledge, there is no modern restitution of Lie's elegant
proof, whereas the description of an alternative proof of Theorem~1.2
as a byproduct of \'E.~Cartan's equivalence algorithm appears in the
references~\cite{ tr1896}, \cite{ ca1924}, \cite{ gtw1989}, \cite{
ol1995}, \cite{ ns2003}\footnote{ We note that in these references,
the already substantial computations are stopped just after the
reduction to an $\{e\}$-structure on an eight-dimensional (local)
principal bundle over the three-dimensional first order jet space. The
vanishing of two (among four) fundamental tensors in the structure
equations of the obtained $\{e\}$-structure yields two partial
differential equations satisfied by the right-hand side $F(x,y,y_x)$,
which are equivalent to {\bf (4)} of Theorem~1.2. We mention that with
the help of Maple programming, the complete reduction to an
$\{e\}$-structure {\it on the base}\, (not only on the principal
bundle) is achieved in~\cite{ hk1989}, in the simpler case of
so-called {\sl fiber-preserving}\, transformations, namely point
transformations leaving invariant the ``vertical'' foliation $\{x={\rm
ct}.\}$. To the author's knowledge, the complete confirmation of
Tresse's results by means of \'E.~Cartan's method has never been
achieved.}.

\smallskip

Lie's Grail would comprise:

\begin{itemize}

\smallskip\item[$\bullet$]
a complete classification of all Lie algebras of local vector fields;
Lie achieved this task in dimension $2$ over $\C$ (real case: \cite{
gko1992}); however, as soon as the dimension is $\geqslant 3$, the complete
classification is unknown, due to the intrinsic richness of
imprimitive Lie algebras of vector fields;

\smallskip\item[$\bullet$]
a list of all the possible Lie algebras that can be realized as
infinitesimal Lie symmetry algebras of partial differential equations,
together with their Levi-Mal$\check{\text{\rm c}}$ev decomposition;

\smallskip\item[$\bullet$]
an explicit Gr\"obner basis of the (noncommutative)
algebra of all differential invariants of each equation in
the list.

\end{itemize}\smallskip

However, complete results hold only for the scalar second order
ordinary differential equation. Some tables extracted from Lie's {\em
Gesammelte Abhandlungen}\, and from Tresse's prized thesis~\cite{
tr1896} may be found in~\cite{ ol1995}.

\subsection*{1.3.~Systems}
Let $x\in \K$, \underline{let $m\geqslant 2$}, let $y:=(y^1, \dots,
y^m) \in \K^m$ and let
\def\theequation{1.4}\begin{equation} 
y_{xx}^1(x) =F^1\left( x, y(x),
y_x(x)\right), \dots \dots, y_{xx}^m(x) =F^m \left(
x, y(x), y_x(x)\right)
\end{equation} 
be a collection of $m$ analytic second order ordinary differential
equations, possibly nonlinear, of the most general form.  

Motivated by physical applications and by
geometrical questions, some authors (Leach~\cite{
le1980}, Grissom-Thompson--Wilkens~\cite{ gtw1989},
Gonz\'alez-L\'opez~\cite{ gl1988}, Fels~\cite{ fe1995},
Crampin-Mart\'{\i}nez-Sarlet~\cite{ cms1996}, 
Doubrov~\cite{ do2000}, 
Grossman~\cite{ gr2000},
Mahomed-Soh~\cite{ ms2001}, and others) have been interested in a
least characterizing those that have the Lie symmetry group of maximal
dimension. In~\cite{ le1980}, based on the belief that the equivalence
between {\bf (1)}, {\bf (2)} and {\bf (3)} of Theorem~1.2 would
persist in the case $m\geqslant 2$, it was conjectured that the
symmetry algebra of {\it every}\, linear system
\def\theequation{1.5}\begin{equation}
y_{ xx}^j 
=
G_0^j(x)+ \sum_{l=1}^m\, y^l\,
G_{1,l}^j(x)+\sum_{l_1=1}^m\, y_x^{l_1}\, H_{l_1}^j(x),
\ \ \ \ \ 
j=1, \dots, m,
\end{equation}
of $m\geqslant 2$ second order differential equations has a Lie symmetry
group locally isomorphic to the full transformation group ${\rm PGL}
(m+1, \K)$ of the projective space $P_{m+1} (\K)$. However,
Gonz\'alez-L\'opez~\cite{ gl1988} infirmed this expectation, and
produced a necessary and sufficient condition ({\it see} Corollary~1.8
below) for local equivalence of such linear systems to the free
particle system
\def\theequation{1.6}\begin{equation}
Y_{ XX}^j = 0, 
\ \ \ \ \ \ \ \ \
j=1, \dots, m.
\end{equation}
Applying Lie's algorithm it is easily seen (\cite{ gg1983}) that the
free particle system has a local symmetry algebra of dimension $\leq
m^2 + 4\, m + 3$, the bound being attained by $Y_{ XX}^j = 0$, $j=1,
\dots, m$, with group ${\rm PGL} (m+1, \K)$.

In 1939, for fiber-preserving transformations only, Chern~\cite{
ch1939} conducted the \'E.~Cartan algorithm through absorptions of
torsion, normalizations and prolongations up to the reduction to an
$\{e\}$-structure. Remarkably, in 1995, Fels~\cite{ fe1995} conducted
the \'E.~Cartan algorithm for general systems $y_{ xx}^j = F^j (x, y,
y_x)$ and for general point transformations. As a byproduct of the
uniqueness of the obtained $\{e\}$-structure for which all invariant
tensors vanish, Fels deduced in~\cite{ fe1995} that the flat system
$Y_{XX}^j=0$, $j=1, \dots, m$, is, up to equivalence, the only system
of second order possessing a symmetry group of maximal dimension.
Alas, the $\{ e\}$-structures obtained by Chern and by Fels are not
parametric, so that the counterpart to {\bf (4)} of Theorem~1.2 was
lacking as soon as $m\geqslant 2$. The extraordinary heaviness of
\'E.~Cartan's method is a well known obstacle.

Recently, Neut~\cite{ n2003} ({\it see} also \cite{ dnp2005}) wrote a
Maple program to compute (among other $\{ e\}$-structures) Fels'
tensors in a parametric, explicit way.  The digital computations
succeeded in the case $m=2$, yielding Theorem~3 of~\cite{ dnp2005}, a
statement that may be checked to be equivalent to Theorem~1.7~{\bf
(3)} just below in the case $m=2$. In the physically most meaningful
case $m=3$, a present-day computer is stuck (\cite{ n2003}, \cite{
dnp2005}).

Extending Lie's less heavy computations
\footnote{ Throughout the article, we do not adopt the summation
convention, because in several subsequent equations, some repeated
indices shall appear that will not be summed. Also, we always put
commas between the indices. For instance $L_{l_1, l_3,y^{ l_2}}^j$
denotes $\partial L_{ l_1, l_3}^j/ \partial y^{ l_2}$ shortly.  As
usual, $\delta_i^j$ is the Kronecker symbol.}, we present here a
complete solution to the characterization of $Y_{ XX}^j = 0$ for
arbitrary $m\geqslant 2$.

\def\thetheorem{1.7}\begin{theorem}
\underline{\rm Suppose $m\geqslant 2$}. 
The following three conditions are equivalent{\rm :}

\begin{itemize}

\smallskip\item[{\bf (1)}]
the system $y_{ xx}^j = F^j (x, y, y_x)$, $j=1, \dots, m$, is
equivalent, under a local point transformation $(x, y^j) \mapsto (X,
Y^j)$ to the free particle equation $Y_{ XX}^j = 0${\rm ;}

\smallskip\item[{\bf (2)}]
the local Lie symmetry group
of $y_{ xx}^j = F^j (x, y, y_x)$ is $(m^2+ 4\, m + 3)$-dimensional and
locally isomorphic to ${\rm PGL} (m+1, \K)${\rm ;}

\smallskip\item[{\bf (3)}]
the right hand sides $F^j(x,y, y_x)$ are of a special form, 
described as follows.

\medskip

\begin{itemize}
\item[{\bf (i)}]
There exist local $\K$-analytic functions $G^j$, $H_{l_1}^j$, $L_{l_1,
l_2}^j$ and $M_{l_1,l_2}$, where $j, l_1, l_2=1,\dots, m$, enjoying
the symmetries $L_{l_1, l_2}^j =L_{l_2, l_1}^j$ and $M_{l_1,l_2}=
M_{l_2, l_1}$ and depending only on $(x,y)$ such that 
$F^j(x,y, y_x)$ may be written as the following specific cubic
polynomial with respect to $y_x${\rm :}
$$
y_{xx}^j=
G^j+
\sum_{l_1=1}^m\, y_x^{l_1}\, H_{l_1}^j+
\sum_{l_1=1}^m\, \sum_{l_2=1}^m\, 
y_x^{l_1} \, y_x^{l_2}\, 
L_{l_1,l_2}^j+y_x^j\cdot
\sum_{l_1=1}^m\, \sum_{l_2=1}^m\, 
y_x^{l_1}\, y_x^{l_2}\, 
M_{l_1,l_2}.
$$
\item[{\bf (ii)}]
The functions $G^j$, $H_{l_1}^j$, $L_{l_1,l_2}^j$ and $M_{l_1,l_2}$
satisfy the following \underline{\rm explicit} system of {\rm four
families} of \underline{first order} partial differential
equations{\rm :}
\def\theequation{I}\begin{equation}
\left\{
\aligned
0= 
&\
-2\, G_{y^{l_1}}^j+
2\, \delta_{l_1}^j\, G_{y^{l_2}}^{l_2}+
H_{l_1,x}^j
-\delta_{l_1}^j\, H_{l_2, x}^{l_2}+ \\
& \
+
2\, \sum_{k=1}^m \, G^k\, L_{l_1,k}^j-2\, 
\delta_{l_1}^j\, \sum_{k=1}^m \, 
G^k\, L_{l_2,k}^{l_2}+ \\
& \
+
\frac{ 1}{2}\, 
\delta_{l_1}^j\, 
\sum_{k=1}^m \, H_{l_2}^k\, H_k^{l_2}-
\frac{ 1}{2}\, 
\sum_{k=1}^m \, H_{l_1}^k\, H_k^j,
\endaligned\right.
\end{equation}
where the indices $j, l_1$ vary in $\{1,2,\dots,m\}${\rm ;} 
\def\theequation{II}\begin{equation}
\left\{
\aligned
0 =
& \
-
\frac{ 1}{2}\, H_{l_1,y^{l_2}}^j
+
\frac{ 1}{6}\, \delta_{l_1}^j\, H_{l_2,y^{l_2}}^{l_2}
+
\frac{ 1}{3}\, \delta_{l_2}^j\, H_{l_1,y^{l_1}}^{l_1}
+ \\
& \
+
L_{l_1,l_2,x}^j-
\frac{ 1}{3}\, \delta_{l_1}^j\, L_{l_2,l_2,x}^{l_2}-
\frac{ 2}{3}\, 
\delta_{l_2}^j\, L_{l_1,l_1,x}^{l_1}+ \\
& \
+ 
G^j\, M_{l_1,l_2}- 
\frac{ 1}{3}\, \delta_{l_1}^j\, G^{l_2}\, M_{l_2,l_2}-
\frac{ 2}{3}\, \delta_{l_2}^j\, G^{l_1}\, M_{l_1,l_1}+ \\
& \
+\frac{ 1}{3}\, \delta_{l_1}^j\,\sum_{k=1}^m 
\, G^k\, M_{l_2,k}-
\frac{ 1}{3}\, \delta_{l_2}^j\, \sum_{k=1}^m \, 
G^k\, M_{l_1,k}- \\
& \
-
\frac{ 1}{2}\, \sum_{k=1}^m \, H_k^j\, L_{l_1,l_2}^k\, 
+\frac{ 1}{2}\, \sum_{k=1}^m \, H_{l_1}^k\, 
L_{l_2, k}^j+ \\
& \
+
\delta_{l_1}^j \left(
\frac{ 1}{6}\, 
\sum_{k=1}^m \, H_k^{l_2} L_{l_2, l_2}^k- 
\frac{ 1}{6}\, \sum_{k=1}^m \, H_{l_2}^k\, L_{l_2, k}^{l_2}
\right) + \\
& \
+
\delta_{l_2}^j\left(
\frac{ 1}{3}\, 
\sum_{k=1}^m \, H_k^{l_1} \, L_{l_1,l_1}^k
-\frac{ 1}{3}\, 
\sum_{k=1}^m \, H_{l_1}^k\, L_{l_1, k}^{l_1}
\right), 
\endaligned\right.
\end{equation}
where the indices $j, l_1, l_2$ vary in 
$\{1,2,\dots,m\}${\rm ;} 
\def\theequation{III}\begin{equation}
\left\{
\aligned
0 =
& \
L_{l_1,l_2, y^{l_3}}^j- 
L_{l_1,l_3, y^{ l_2}}^j+
\delta_{l_3}^j\, M_{l_1,l_2,x}- 
\delta_{l_2}^j\, M_{l_1, l_3,x}+ \\
& \
+\frac{ 1}{2}\, H_{l_3}^j\, M_{l_1, l_2}- 
\frac{ 1}{2}\, H_{l_2}^j\, M_{l_1, l_3}+ \\
& \
+
\frac{ 1}{2}\, \delta_{l_1}^j\, \sum_{k=1}^m \, 
H_{l_3}^k\, M_{l_2, k}- 
\frac{ 1}{2}\, \delta_{l_1}^j\, \sum_{k=1}^m \, 
H_{l_2}^k\, M_{l_3, k}+ \\
& \
+
\frac{ 1}{2}\,
\delta_{l_3}^j\, \sum_{k=1}^m \, 
H_{l_1}^k\, M_{l_2, k}- 
\frac{ 1}{2}\,
\delta_{l_2}^j\, \sum_{ k=1}^m \, H_{l_1}^k\, M_{l_3, k}+ \\
& \
+
\sum_{k= 1}^m \, L_{l_1, l_3}^k\, L_{l_2, k}^j- 
\sum_{k=1}^m \, L_{l_1, l_2}^k\, L_{l_3, k}^j,
\endaligned\right.
\end{equation}
where the indices $j, l_1, l_2, l_3$ vary in 
$\{1, \dots m\}${\rm ;} and 
\def\theequation{IV}\begin{equation}
\left\{
\aligned
0=
M_{l_1,l_2, y^{l_3}}- M_{l_1,l_3, y^{l_2}}-
\sum_{k=1}^m\, L_{l_1,l_2}^k\, M_{l_3, k}+
\sum_{k=1}^m \, L_{l_1,l_3}^k \, M_{l_2,k},
\endaligned\right.
\end{equation}
where the indices $l_1,l_2,l_3$ vary in $\{1,\dots, m\}$.
\end{itemize}
\end{itemize}
\end{theorem}

Let us provide commentaries and explanations. The form of the
right-hand side of {\bf (3)(i)} is the analog of the form of the
right-hand side $F$ in {\bf (4)} of Lie's Theorem~1.2. However, we
notice that the right-hand side of~{\bf (3)(i)} is not the most
general degree three polynomial in the variables $y_x^j$,
$j=1, \dots, m$: some coefficients of the cubic terms vanish. 

Very strikingly, the differential system~{\bf (I)}, {\bf (II)}, {\bf
(III)}, {\bf (IV)} satisfied by the functions $G^j$, $H_{l_1}^j$,
$L_{l_1,l_2}^j$, $M_{l_1, l_2}$ is of {\it first order} for $m\geqslant 2$,
whereas the system {\bf (4)} satisfied by $G$, $H$, $L$, $M$ in Lie's
Theorem~1.2 is of {\it second order} for $m=1$. This confirms the main
theorem of Gonz\'alez-L\'opez~\cite{ gl1988}, that is recovered here
in the special linear case~\thetag{ 1.5}.

\def\thecorollary{1.8}\begin{corollary}
{\rm (\cite{ gl1988})}
For $m\geqslant 2$, a linear {\rm (}nonhomogeneous{\rm )} system $y_{ xx}^j
= G_0^j(x)+ \sum_{l=1}^m\, y^l\, G_{1,l}^j(x)+\sum_{l_1=1}^m\,
y_x^{l_1}\, H_{l_1}^j(x)$ is equivalent to $Y_{ XX}^j = 0$ if and only
if there exists a function $B(x)$ such that the $m\times m$ matrix
$G_1$ may be written under the specific form{\rm :}
\def\theequation{1.9}\begin{equation}
G_{1,l}^j
=
\frac{1}{2}\,H_l^j
-
\frac{1}{4}\,\sum_{k=1}^m\,H_k^j\,H_l^k
+
\delta_l^j\,B.
\end{equation}
\end{corollary}

We also mention that Fels obtained a characterization of equivalence
to the free particle $Y_{ XX}^j = 0$, $j=1, \dots, m$, by the
vanishing of two (nonexplicit) tensors $\widetilde{ S}_{ ikl}^j$ and
$\widetilde{ P}_i^j$ (Corollary~5.1 in~\cite{ fe1995}).  In this
reference, some parametric computations are achieved after restricting
the initial $G$-structure, together with its subsequent prolongations,
to the identity element of the group; explicit expressions of
$\widetilde{ S}_{ ikl}^j \big\vert_{ \rm Id}$ and of $\widetilde{
P}_i^j \big\vert_{ \rm Id}$ are then obtained, through already hard
computations. The vanishing of the two tensors $\widetilde{ P }_i^j$
and $\widetilde{ S }_{ ikl}^j$ {\it at the identity}\, of the
structure group, namely, as computed in Lemma~4.1 of~\cite{ fe1995},
yields (translating into our notation)
\def\theequation{1.10}\begin{equation}
\left\{
\aligned
0 = 
& \
\left.
(\widetilde{ S}_{ikl}^j)\right\vert_{\rm Id}=
F_{y_x^i y_x^ky_x^l}^j- 
\frac{ 1}{n+2}\, 
\sum_{l_1=1}^m\, 
\sum_{\sigma\in \mathfrak{ S}_3}\, 
\delta_{\sigma(l)}^j\, 
F_{y_x^{l_1}y_x^{\sigma(i)}y_x^{\sigma(k)}}^{l_1}, \\
0 = 
& \
\left.
(\widetilde{ P}_i^j)\right\vert_{\rm Id}=
\frac{ 1}{2}\, 
\, D\left(
F_{y_x^i}^j
\right)- 
F_{y^i}^j-\frac{ 1}{4}\, 
\sum_{k=1}^m\, 
F_{y_x^k}^j\, F_{y_x^i}^k- \\
& \
\ \ \ \ \ \ \ \ \ \ \ \ \ \ \ \ \ \ 
- 
\frac{ 1}{m}\, \delta_i^j\left[
\frac{ 1}{2}\,
D \left(
\sum_{k=1}^m\, F_{y_x^k}^k
\right)- 
\sum_{k=1}^m\, 
F_{y^k}^k- 
\frac{ 1}{4}\, 
\sum_{k=1}^m\, 
\sum_{l=1}^m\, 
F_{y_x^k}^l\, 
F_{y_x^l}^k
\right],
\endaligned\right.
\end{equation}
where $D$ is the total differentiation operator $\frac{ \partial
}{\partial x}+ \sum_{l=1}^m\, y_x^l\, \frac{ \partial }{\partial
y^l}+\sum_{l=1}^m\, F^l\, \frac{ \partial}{\partial y_x^l }$, and
where $i,\, j,\, k,\, l= 1,\dots, m$.  Strikingly, one may check that
the first equation is equivalent to {\bf (3)(i)} of Theorem~1.2 and
then that the second equation yields the (complicated) four families
of first order partial differential equations (I), (II), (III) and
(IV).  So, in Corollary~5.1 of~\cite{ fe1995}, one may replace the
vanishing of $\widetilde{ S}_{ ikl}^j$ and of $\widetilde{ P}_i^j$ by
the vanishing of $\widetilde{ S}_{ ikl}^j \big\vert_{ \rm Id}$ and of
$\widetilde{ P}_i^j \big\vert_{ \rm Id}$, which were {\it explicitly
computed}\, there.

This phenomenon could be explained as follows: as soon as the tensors
$\widetilde{ S}_{ikl}^j$ vanish, the system enjoys a projective
connection (appendix of~\cite{ fe1995}); with such a connection, the
tensors $\widetilde{ P}_i^j$ then transform according to a specific
rule via tensorial rotation formulas and their general expression may
be deduced from their expression at the identity\footnote{ Similar
rotation formulas are known in the much simpler case of (pseudo-)
Riemannian metrics, {\it see} Chapter~12 of Olver~\cite{ ol1995}.  It
would be interesting to write a program, finer and more efficient
than~\cite{ n2003}, which would systematically recognize such rotation
formulas in any application of \'E.~Cartan's equivalence method. }. We
have checked this, but as we try to avoid the method of equivalence,
details will not be reproduced here.  Similar observations appear in
Hachtroudi~\cite{ ha1937}.

Even if the expressions~\thetag{ 1.10} are more compact than the
(equivalent) conditions in Theorem~1.7~{\bf (3)}, we prefer the
complete expressions of Theorem~1.7~{\bf (3)}, since they are more
explicit and ready-made for checking whether a physically given
nonlinear system is equivalent to a free particle. If the reader
prefers compact expressions and ``short'' theorems, (s)he may replace
the conditions of Theorem~1.7~{\bf (3)} by~\thetag{ 1.10}.

\def\theopenproblem{1.11}\begin{openproblem}
Characterize explicitly the linearizability of a Newtonian system in
$m \geqslant 2$ degrees of freedom, {\it i.e.} local equivalence
to~{\rm :}
\def\theequation{1.12}\begin{equation}
Y_{ XX}^j = G_0^j(X)+ \sum_{l=1}^m\, Y^l\,
G_{1,l}^j(X)+\sum_{l_1=1}^m\, Y_X^{l_1}\, H_{l_1}^j(X).
\end{equation}
\end{openproblem}

\subsection*{1.13.~Organization, avertissement and
acknowledgment} Section~2 is devoted to a thorough restitution of
Lie's original proof of the equivalence between {\bf (1)} and {\bf
(4)} in Theorem~1.2. Section~3 is devoted to the formulation of
combinatorial formulas yielding the general form of a system
equivalent to $Y_{XX}^j= 0$, $j= 1, \dots, m$, under a local
$\K$-analytic point transformation $(x, y^j) \mapsto (X, Y^j)$, for
general $m\geqslant 2$~; the proof of the main technical Lemma~3.32 is
exposed in Section~5. Section~4 is devoted to the final proof of the
equivalence between {\bf (1)} and {\bf (3)} in Theorem~1.7, the
equivalence between {\bf (1)} and {\bf (2)} being already proved by
Fels~\cite{ fe1995}.

Some word about style and intentions. We wanted the proof of the
equivalence between {\bf (1)} and {\bf (3)} be totally complete, every
tiny detail being rigorously and patiently checked. This is why we
decided to {\it carefully detail each intermediate computational
step}, seeking first the combinatorics of the formal calculations in
the case $m=1$ and devising then the underlying combinatorics for the
case $m\geqslant 2$. Actually, the size of differential expressions is
relatively impressive, as will become soon evident. Thus, no
intermediate symbolic computation will be hidden, hence essentially no
checking work is left to the reader, as would have been the case if we
did not have typed all the computations.

We also would like to point out that except in every specific standard
situations, as with the much studied Riemann and Ricci tensors,
present-day computer programs are not yet powerful enough to apply the
\'E.~Cartan equivalence method when the number of some collection of
variables is a general integer.  All the formulas obtained in
Sections~2, 3 and 4 were first treated completely by hand and then,
some of them were confirmed afterwards with the help of MAPLE release
6 in the cases $m=2$ and $m=3$. The author is indebted to Sylvain Neut
and to Michel Petitot, from the University of Lille~1, for their help
in computer machine confirmations.

\section*{\S2.~Proof of Lie's theorem}

\subsection*{2.1.~Argument}
This preliminary section contains a detailed exposition of Lie's
original proof of the equivalence between {\bf (1)} and {\bf (4)} in
Theorem~1.2. Since our goal is to guess the combinatorics of
computations in several variables, it will be a crucial point for us
to explain thoroughly and patiently each step of Lie's
computation. Without such an intuitive control, it would be hopeless
to conduct any generalization to several variables. Hence we shall
respect a fundamental principle: always explain clearly and completely
what sort of computation is achieved at each step.  Also, we shall
many times introduce some appropriate new notation.

\subsection*{2.2.~Combinatorics of the second order prolongation of
a point transformation} Let $\K=\R$ or $\C$. Let $(x,y)\mapsto
(X(x,y),Y(x,y))$ be a local $\K$-analytic invertible transformation,
defined in a neighborhood of the origin in $\K^2$, which maps the
second order differential equation $y_{xx}=F(x,y,y_x)$ to the flat
equation $Y_{XX}=0$. By assumption, the Jacobian determinant
\def\theequation{2.3}\begin{equation}
\Delta(x\vert y):= 
\left\vert
\begin{array}{cc}
X_x & X_ y \\
Y_x & Y_y
\end{array}
\right\vert
\end{equation}
is nowhere vanishing. Since the equation $Y_{XX}=0$ is left unchanged
by any affine transformation in the $(X,Y)$ space, we can (and we
shall) assume that the transformation is tangent to the identity at
the origin, namely the above Jacobian matrix equals the identity
matrix at $(x,y)=(0,0)$.

The computation how the differential equation in the $(X,Y)$
coordinates is related to the differential equation in the
$(x,y)$-coordinates is classical, {\it cf.}~\cite{lie1883},
\cite{tr1896}, \cite{bk1989}, \cite{ib1992}: let us
remind it. A local graph
$\{y=y(x)\}$ being transformed to a local graph $\{Y=Y(X)\}$, we have
a direct formula for the first derivative $Y_X$:
\def\theequation{2.4}\begin{equation}
Y_X:=\frac{ dY}{dX}= 
\frac{dx \cdot \partial Y(x,y(x))/\partial x}{
dx \cdot \partial X(x,y(x))/\partial x}=
\frac{Y_x+y_x Y_y}{X_x+y_xX_y}.
\end{equation}
This yields the prolongation of the transformation 
to the first order jet space. For the second order prolongation,
introducing the second order total differentiation operator
(which geometrically corresponds to differentiation along graphs
$\{(x,y(x))\}$)
defined by
\def\theequation{2.5}\begin{equation}
D:= \frac{\partial }{\partial x}+
y_x\, \frac{ \partial }{\partial y}+
y_{xx}\, \frac{ \partial }{\partial y_x},
\end{equation}
we may compute, simplify and reorder the expression of
the second order derivative in the $(X,Y)$-coordinates:
\def\theequation{2.6}\begin{equation}
\left\{
\aligned
Y_{XX}:=\frac{d^2 Y}{d X^2}\equiv
\frac{D Y_X}{DX}=
& \
\frac{D\left[(Y_x+y_xY_y)(X_x+y_xX_x)^{-1}\right]}{
X_x+y_xX_y}= \\
= 
& \
\frac{ 1}{[X_x+y_xX_y]^3}\left\{
y_{xx}\left[
X_xY_y-Y_xX_y
\right]+
X_xY_{xx}-Y_xX_{xx}+ \right. \\
+
& \
\left.
y_x\left[
2(X_xY_{xy}-Y_xX_{xy})-(X_{xx}Y_y-Y_{xx}X_y)
\right]+ \right. \\
+
& \
\left.
y_xy_x\left[
X_xY_{yy}-Y_xX_{yy}-2(X_{xy}Y_y-Y_{xy}X_y)
\right]+ \right. \\
+
& \
\left.
y_xy_xy_x\left[
-(X_{yy}Y_y-Y_{yy}X_y)
\right]
\right\}.
\endaligned\right.
\end{equation}
Even if not too complicated, the internal combinatorics of this
expression has to be analyzed and expressed thoroughly. First of all,
as $Y_{XX}=0$ by assumption, we may erase the cubic factor
$[X_x+y_xX_y]^{-3}$. Next, as the factor of $y_{xx}$ in the right-hand
side of~\thetag{ 2.6}, we just recognize the Jacobian $\Delta(x\vert
y)$ expressed in \thetag{ 2.3} above. Also, all the other factors are
modifications of the Jacobian $\Delta(x\vert y)$, whose combinatorics
may be understood as follows.

There exist exactly three possible distinct second order derivatives:
$xx$, $xy$ and $yy$. There are also exactly two columns in~\thetag{
2.3}. By replacing each of the two columns of first order derivative
in $\Delta(x\vert y)$ by any column of second order derivative
(leaving $X$ and $Y$ unchanged), we may build exactly six new
determinants
\def\theequation{2.7}\begin{equation}
\left\{
\aligned
\Delta(xx \vert y) \ \ \ \ \ & 
\Delta(xy \vert y) \ \ \ \ \ & 
\Delta(yy \vert y) \ \ \ \ \ \\
\Delta(x \vert xx) \ \ \ \ \ & 
\Delta(x \vert xy) \ \ \ \ \ & 
\Delta(x \vert yy) \ \ \ \ \ \\
\endaligned\right.
\end{equation}
where for instance
\def\theequation{2.8}\begin{equation}
\left\{
\aligned
\Delta(\underline{xx}\vert y) :=
\left\vert
\begin{array}{cc}
X_{\underline{xx}} & X_y \\
Y_{\underline{xx}} & Y_y 
\end{array}
\right\vert 
\ \ \ \ \ {\rm and} \ \ \ \ \
\Delta(x\vert \underline{xy}) :=
\left\vert
\begin{array}{cc}
X_x & X_{\underline{xy}} \\
Y_x & Y_{\underline{xy}} 
\end{array}
\right\vert .
\endaligned\right.
\end{equation}
Hence, by rewriting~\thetag{ 2.6}, we see that the equation
$y_{xx}=F(x,y,y_x)$ equivalent to $Y_{XX}=0$ may be written under the
general explicit form, involving determinants
\def\theequation{2.9}\begin{equation}
\left\{
\aligned
0= 
& \
y_{xx} \cdot 
\left\vert
\begin{array}{cc}
X_x & X_y \\
Y_x & Y_y 
\end{array}
\right\vert
+
\left\vert
\begin{array}{cc}
X_x & X_{xx} \\
Y_x & Y_{xx} 
\end{array}
\right\vert
+
y_x\cdot
\left\{
2
\left\vert
\begin{array}{cc}
X_x & X_{xy} \\
Y_x & Y_{xy} 
\end{array}
\right\vert-
\left\vert
\begin{array}{cc}
X_{xx} & X_y \\
Y_{xx} & Y_y 
\end{array}
\right\vert
\right\}+ \\
& \
+ y_xy_x\cdot \left\{
\left\vert
\begin{array}{cc}
X_x & X_{yy} \\
Y_x & Y_{yy} 
\end{array}
\right\vert
-2
\left\vert
\begin{array}{cc}
X_{xy} & X_y \\
Y_{xy} & Y_y
\end{array}
\right\vert
\right\}+
y_xy_xy_x\cdot 
\left\{
-
\left\vert
\begin{array}{cc}
X_{yy} & X_y \\
Y_{yy} & Y_y
\end{array}
\right\vert
\right\}
\endaligned\right.
\end{equation}
or equivalently, after solving in $y_{ xx}$, {\it i.e.} after dividing
by the Jacobian $\Delta(x \vert y)$: 
\def\theequation{2.10}\begin{equation}
\left\{
\aligned
y_{xx}=
& \
- \frac{ 
\Delta(x\vert xx)}{\Delta(x\vert y)}+
y_x \cdot \left\{
-2
\frac{ \Delta(x \vert xy)}{\Delta(x\vert y)}+
\frac{\Delta(xx\vert y)}{\Delta(x\vert y)}
\right\}+ \\
& \
+(y_x)^2 \cdot \left\{
-\frac{\Delta(x\vert yy)}{\Delta(x\vert y)}+2
\frac{\Delta(xy\vert y)}{\Delta(x\vert y)}\right\}+
(y_x)^3\cdot 
\left\{
\frac{ \Delta(yy\vert y)}{\Delta(x\vert y)}
\right\}.
\endaligned\right.
\end{equation}
At this point, it will be convenient to slightly contract 
the notation by introducing a new family of 
{\sl square functions} as follows. We first index 
the coordinates $(x,y)$ as $(y^0, y^1)$, namely we
introduce the two notational equivalences 
\def\theequation{2.11}\begin{equation}
\fbox{$y^0\equiv x, \ \ \ \ \
y^1\equiv y$}\,,
\end{equation}
which will be very convenient in the sequel, especially to write down
general combinatorial formulas anticipating our treatment of the case
of $m \geqslant 2$ dependent variables $(y^1,\dots, y^m)$, to be
achieved in Sections~3, 4 and~5 below. With this convention at hand,
our six square functions $\square_{y^{j_1 } y^{ j_2}}^{ k_1}$,
symmetric with respect to the lower indices, where $0\leqslant j_1,
j_2, k_1 \leqslant 1$, are defined by
\def\theequation{2.12}\begin{equation}
\left\{
\aligned
\square_{xx}^0:=\frac{ \Delta(xx\vert y)}{\Delta(x\vert y)}, 
\ \ \ \ \ & 
\square_{xy}^0:=\frac{ \Delta(xy\vert y)}{\Delta(x\vert y)},
\ \ \ \ \ & 
\square_{yy}^0:=\frac{ \Delta(yy\vert y)}{\Delta(x\vert y)}, 
\ \ \ \ \ \\
\square_{xx}^1:=\frac{ \Delta(x\vert xx)}{\Delta(x\vert y)}, 
\ \ \ \ \ & 
\square_{xy}^1:=\frac{ \Delta(x\vert xy)}{\Delta(x\vert y)}, 
\ \ \ \ \ & 
\square_{yy}^1:=\frac{ \Delta(x\vert yy)}{\Delta(x\vert y)}.
\ \ \ \ \ \\
\endaligned\right.
\end{equation}
Here of course, the upper index designates the column upon which
the second order derivative appears, itself being encoded by the two
lower indices. Even if this is hidden in the notation, we shall
remember that the square functions are explicit rational expressions
in terms of the second order jet of the transformation $(x,y)\mapsto
(X(x,y), Y(x,y))$. However, we shall be aware of not confusing
the index in the square functions with a second order partial 
derivative of some function ``$\square^j$'', 
denoted by the square symbol: indeed, 
the partial derivatives are hidden in some determinant. 

At this point, we may summarize what we have
established so far.

\def\thelemma{2.13}\begin{lemma}
The equation $y_{xx}=F(x,y,y_x)$ is equivalent to the flat equation
$Y_{XX}=0$ if and only if there exist two local $\K$-analytic
functions $X(x,y)$ and $Y(x,y)$ such that it may be written under the
form
\def\theequation{2.14}\begin{equation}
y_{xx}=
-\square_{xx}^1+
y_x\cdot \left(
-2\, \square_{xy}^1+\square_{xx}^0
\right)+
(y_x)^2\cdot \left(
-\square_{yy}^1+2\, \square_{xy}^0
\right)+
(y_x)^3\cdot \square_{yy}^0.
\end{equation}
\end{lemma}
At this point, for heuristic reasons, it may be useful to compare the
right-hand side of~\thetag{ 2.14} with the classical expression of the
prolongation to the second order jet space of a general vector field
of the form $L:=X(x,y) \frac{\partial }{ \partial x}+ Y(x,y)
\frac{\partial }{\partial y}$, which is given, according to
\cite{ lie1883}, \cite{ ol1986}, \cite{ bk1989}, by
\def\theequation{2.15}\begin{equation}
\left\{
\aligned
L^{(2)}=
& \
X\, \frac{ \partial }{\partial x}+
Y\, \frac{\partial }{\partial y}+
\left[
Y_x+y_x\cdot(Y_y-X_x)+
(y_x)^2\cdot(-X_y)
\right]\, \frac{\partial }{\partial y_x}+ \\ 
+
& \
\left[
Y_{xx}+y_x\cdot(2Y_{xy}-X_{xx})+
(y_x)^2\cdot(Y_{yy}-2X_{xy})+
(y_x)^3\cdot(-X_{yy})+ \right. \\
+ 
& \
\left.
y_{xx}\cdot(Y_y-2\, X_x)+
y_x y_{xx} \cdot(-3\, X_y)
\right]\, \frac{\partial }{\partial y_{xx}}.
\endaligned\right.
\end{equation}
We immediately see that (up to an overall minus sign) the right-hand
side of~\thetag{ 2.14} is formally analogous to the second line
of~\thetag{ 2.15}~: the letter $X$ corresponds to the symbol
$\square^0$ and the letter $Y$ corresponds to the symbol
$\square^1$. This analogy is no mystery, just because the formula for
$L^{(2)}$ is classically obtained by differentiating at $\varepsilon
=0$ the second order prolongation $[\exp(\varepsilon L)(\cdot)]^{(2)}$
of the flow of $L$~!

In fact, as we assumed that the transformation $(x,y)\mapsto (X(x,y),
Y(x,y))$ is tangent to the identity at the origin, we may think that
$X_x\cong 1$, $X_y\cong 0$, $Y_x\cong 0$ and $Y_y \cong 1$, whence the
Jacobian $\Delta(x\vert y)\cong 1$ and moreover
\def\theequation{2.16}\begin{equation}
\left\{
\aligned
\square_{xx}^0 \cong X_{xx}, \ \ \ \ \ & 
\square_{xy}^0 \cong X_{xy}, \ \ \ \ \ & 
\square_{yy}^0 \cong X_{yy}, \ \ \ \ \ \\
\square_{xx}^1 \cong Y_{xx}, \ \ \ \ \ & 
\square_{xy}^1 \cong Y_{xy}, \ \ \ \ \ & 
\square_{yy}^1 \cong Y_{yy}. \ \ \ \ \ \\
\endaligned\right.
\end{equation}
By means of this (abusive) notational correspondence, we see that, up
to an overall minus sign, the right-hand side of~\thetag{ 2.14}
transforms precisely to the second line of~\thetag{ 2.15}. This
analogy will be useful in devising combinatorial formulas for the
generalization of Lemma~2.13 to the case of $m\geqslant 2$ variables
$(y^1, \dots, y^m)$, {\it see}\, Lemmas~3.22 and~3.32 below.

\subsection*{2.17.~Continuation}
Clearly, since the right-hand side of~\thetag{ 2.14} is a polynomial
of degree three in $y_x$, the first condition of Theorem~1.2 {\bf (4)}
immediately holds. We are therefore led to establish that the second
condition is necessary and sufficient in order that there exist two
local $\K$-analytic functions $X(x,y)$ and $Y(x,y)$ which solve the
following system of nonlinear second order partial differential
equations (remind that the second order jet of $(X,Y)$ is hidden in
the square functions):
\def\theequation{2.18}\begin{equation}
\left\{
\aligned
G = 
& \
-\square_{xx}^1, \\
H =
& \
-2\, \square_{xy}^1+\square_{xx}^0, \\
L= 
& \
-\square_{yy}^1+2\, \square_{xy}^0, \\
M = 
& \
\square_{yy}^0.
\endaligned\right.
\end{equation}
In the remainder of this section, following~\cite{ lie1883}, p.~364,
we shall study this second order system by introducing two auxiliary
systems of partial differential equations which are {\sl complete}, 
and we shall see in \S2.38 below
that the compatibility conditions (insuring
involutivity, hence complete integrability) of the second
auxiliary system exactly provide the two 
partial differential equations appearing in
Theorem~1.2 {\bf (4)}.

\subsection*{2.19.~First auxiliary system}
We notice that in~\thetag{ 2.18}, there are two more square functions
$\square_{xx}^0$, $\square_{xy}^0$, $\square_{yy}^0$,
$\square_{xx}^1$, $\square_{xy}^1$, $\square_{yy}^1$, than functions
$G$, $H$, $L$ and $M$. Hence, as a trick, let us introduce six new
independent functions $\Pi_{j_1,j_2 }^{ k_1}$ of $(x,y)$, symmetric
with respect to the lower indices, for $0\leqslant j_1, j_2, k_1 \leqslant 1$
and let us seek necessary and sufficient conditions in order that
there exist solutions $(X,Y)$ to the {\sl first auxiliary system}:
\def\theequation{2.20}\begin{equation}
\left\{
\aligned
\square_{xx}^0 = \Pi_{0,0}^0, \ \ \ \ \ & 
\square_{xy}^0 = \Pi_{0,1}^0, \ \ \ \ \ & 
\square_{yy}^0 = \Pi_{1,1}^0, \ \ \ \ \ \\
\square_{xx}^1 = \Pi_{0,0}^1, \ \ \ \ \ & 
\square_{xy}^1 = \Pi_{0,1}^1, \ \ \ \ \ & 
\square_{yy}^1 = \Pi_{1,1}^1. \ \ \ \ \ \\ 
\endaligned\right.
\end{equation}
According to the (aprooximate) identities~\thetag{ 2.16}, this system
looks like a complete second order system of partial differential
equations in two variables $(x,y)$ and in two unknowns $(X,Y)$. More
rigorously, by means of elementary algebraic operations, taking
account of the fact that $X_x\cong 1$, $X_y\cong 0$, $Y_x\cong 0$ and
$Y_y \cong 1$, one may transform this sytem in a true second order
{\sl complete}\, system, solved with respect to the top order
derivatives, namely of the form
\def\theequation{2.21}\begin{equation}
\left\{
\aligned
X_{xx} = \Lambda_{0,0}^0, \ \ \ \ \ & 
X_{xy} = \Lambda_{0,1}^0, \ \ \ \ \ & 
X_{yy} = \Lambda_{1,1}^0, \ \ \ \ \ \\
Y_{xx} = \Lambda_{0,0}^1, \ \ \ \ \ & 
Y_{xy} = \Lambda_{0,1}^1, \ \ \ \ \ & 
Y_{yy} = \Lambda_{1,1}^1, \ \ \ \ \ \\ 
\endaligned\right.
\end{equation}
where the $\Lambda_{ j_1, j_2}^{ k_1}$ are 
local $\K$-analytic functions of
$(x, y, X, Y, X_x, X_y, Y_x, Y_y)$. For such a system, the
compatibility conditions [which are necessary and sufficient for the
existence of a solution $(X,Y)$] are easily formulated:
\def\theequation{2.22}\begin{equation}
\left\{
\aligned
(\Lambda_{0,0}^0)_y=(\Lambda_{0,1}^0)_x, \ \ \ \ \ & 
(\Lambda_{0,1}^0)_y=(\Lambda_{1,1}^0)_x, \ \ \ \ \ \\ 
(\Lambda_{0,0}^1)_y=(\Lambda_{0,1}^1)_x, \ \ \ \ \ & 
(\Lambda_{0,1}^1)_y=(\Lambda_{1,1}^1)_x. \ \ \ \ \ & 
\endaligned\right.
\end{equation}
Equivalently, we may express the compatibility conditions directly
with the system~\thetag{ 2.20}, without transforming it to the
form~\thetag{ 2.21}. This direct strategy will be more appropriate.

\subsection*{2.23.~Compatibility conditions for the first
auxiliary system}
Indeed, to begin with, let us remind that the $\Delta(\cdot \vert
\cdot)$ are determinant, hence we have the
skew-symmetry relation $\Delta(x^a
y^b\vert x^cy^d)= -\Delta(x^c y^d\vert x^a y^b)$ and the following two
formulas for partial differentiation
\def\theequation{2.24}\begin{equation}
\left\{
\aligned
{}
\left[\Delta(x^a y^b \vert x^c y^d)\right]_x=
& \ 
\Delta(x^{a+1} y^b \vert x^cy^d)+
\Delta(x^a y^b \vert x^{c+1} y^d), \\
\left[\Delta(x^a y^b \vert x^c y^d)\right]_y=
& \
\Delta(x^a y^{b+1} \vert x^c y^d) +
\Delta(x^a y^b \vert x^c y^{d+1}).
\endaligned\right.
\end{equation}
With these formal rules at hand, as an exercise, let us compute for
instance the following cross differentiation 
(remember that the lower index in the square
functions is {\it not}\, a partial derivative):
\def\theequation{2.25}\begin{equation}
\left\{
\aligned
{}
& \ \
\left(\square_{xx}^0\right)_y-
\left(\square_{xy}^0\right)_x=
\frac{\partial}{\partial y}\left(
\frac{\Delta(xx\vert y)}{\Delta(x\vert y)}
\right)-\frac{\partial}{\partial x}\left(
\frac{\Delta(xy\vert y)}{\Delta(x\vert y)}
\right) = \\
& \
=\frac{ 1}{[\Delta(x\vert y)]^2}\left\{
\underline{
\Delta(xxy\vert y)\cdot \Delta(x\vert y)}_{ 
\octagon \! \! \! \! \! \tiny{\sf a}}+
\Delta(xx \vert yy)\cdot \Delta(x\vert y)- \right. \\
& \ \ \ \ \ \ \ \ \ \ \ \ \ \ \ \ \ \ \ \ \ \ \ \ \
\left.
-\underline{\Delta(xy\vert y)\cdot
\Delta(xx \vert y)}_{ 
\octagon \! \! \! \! \! \tiny{\sf b}}
-
\Delta(x \vert yy ) \cdot \Delta( xx \vert y)- \right. \\
& \ \ \ \ \ \ \ \ \ \ \ \ \ \ \ \ \ \ \ \ \ \ \ \ \
\left.
-\underline{
\Delta(xxy\vert y)\cdot \Delta(x\vert y)}_{
\octagon \! \! \! \! \! \tiny{\sf a}}
-
\underline{
\Delta( xy \vert xy)\cdot \Delta(x\vert y)}_{
\octagon \! \! \! \! \! \tiny{\sf c}}
+ \right. \\
& \ \ \ \ \ \ \ \ \ \ \ \ \ \ \ \ \ \ \ \ \ \ \ \ \
\left.
+
\underline{\Delta(xy\vert y)\cdot \Delta(xx \vert y)}_{
\octagon \! \! \! \! \! \tiny{\sf b}}
+
\Delta(xy \vert y) \cdot \Delta( x\vert xy)
\right\} = \\
& \
=\frac{ 1}{[\Delta(x\vert y)]^2}\left\{
\Delta( xx \vert yy) \cdot \Delta( x \vert y) - 
\Delta( x \vert yy) \cdot \Delta( xx \vert y) + \right. \\
& \ \ \ \ \ \ \ \ \ \ \ \ \ \ \ \ \ \ \ \ \ \ \ \ \
\left.
+
\Delta( xy \vert y)\cdot \Delta( x \vert xy)
\right\}.
\endaligned\right.
\end{equation}
Crucially, we observe that the third order derivatives kill each other
and disappear, {\it see}\, the underlined terms with $\octagon 
\! \! \! \! \tiny{\sf a}$ 
appended. Also, two products of two determinants
$\Delta ( \cdot \vert \cdot)$ 
involving a second order derivative upon one
column of each determinant kill each other: they are underlined with 
$\octagon 
\! \! \! \! \tiny{\sf b}$ 
appended. Finally, by antisymmetry of
determinants, the term $\Delta(xy\vert xy)\cdot
\Delta( x\vert y)$ vanishes gratuitously: it
is underlined with $\octagon 
\! \! \! \! \tiny{\sf c}$ appended.

However, there still remains one term involving second order
derivatives upon the {\it two}\, columns of
a determinant: it is $\Delta ( xx \vert yy)$.

We must transform this unpleasant term $\Delta ( xx\vert yy) \cdot
\Delta ( x \vert y)$ and express it as a product of two determinants,
each involving a second order derivative only in one column. To this
aim, we have:

\def\thelemma{2.26}\begin{lemma}
The following three relations between the differential determinants
$\Delta( \cdot \vert \cdot)$ hold true{\rm :}
\def\theequation{2.27}\begin{equation}
\left\{
\aligned
\Delta( xx \vert xy )\cdot \Delta( x\vert y) =
& \
\Delta( xx \vert y) \cdot \Delta( x \vert xy) -
\Delta( xy \vert y ) \cdot \Delta( x\vert xx), \\
\Delta( xx \vert yy) \cdot \Delta( x\vert y) =
& \
\Delta( xx \vert y) \cdot \Delta( x \vert yy) - 
\Delta( yy \vert y) \cdot \Delta( x \vert xx), \\
\Delta( xy \vert yy) \cdot \Delta( x\vert y) = 
& \
\Delta( xy \vert y)\cdot \Delta( x \vert yy) -
\Delta( yy\vert y )\cdot \Delta( x \vert xy).
\endaligned\right.
\end{equation}
\end{lemma}

\proof
Each of these three formal identities is an immediate direct
consequence of the following Pl\"ucker type identity, easily verified by 
developing all the determinants~:
\def\theequation{2.28}\begin{equation}
\left\vert
\begin{array}{cc}
A_1 & B_1 \\
A_2 & B_2 
\end{array}
\right\vert 
\cdot
\left\vert
\begin{array}{cc}
C_1 & D_1 \\
C_2 & D_2 
\end{array}
\right\vert =
\left\vert
\begin{array}{cc}
A_1 & D_1 \\
A_2 & D_2 
\end{array}
\right\vert 
\cdot
\left\vert
\begin{array}{cc}
C_1 & B_1 \\
C_2 & B_2 
\end{array}
\right\vert-
\left\vert
\begin{array}{cc}
B_1 & D_1 \\
B_2 & D_2 
\end{array}
\right\vert 
\cdot
\left\vert
\begin{array}{cc}
C_1 & A_1 \\
C_2 & A_2 
\end{array}
\right\vert,
\end{equation}
where the variables $A_1, \, A_2, \, B_1, \, B_2, \, C_1, \, C_2, \,
D_1, \, D_2\in \K$ are arbitrary. \endproof

Thanks to the second identity~\thetag{2.27}, we may therefore
transform the result left above in the last two lines of~\thetag{
2.25}; as desired, it will remain determinants having only one second
order derivative per column, so that after division by $[\Delta(x\vert
y)]^2$, we discover {\it a quadratic expression involving only the
square functions themselves}:
\def\theequation{2.29}\begin{equation}
\left\{
\aligned
{}
& \ \
\left(\square_{xx}^0\right)_y-
\left(\square_{xy}^0\right)_x
=
\frac{ 1}{[\Delta(x\vert y)]^2}\left\{
\Delta( xx \vert y)\cdot \Delta( x \vert yy) -
\Delta( yy \vert y) \cdot \Delta( x \vert xx) - \right. \\
& \ \ \ \ \ \ \ \ \ \ \ \ \ \ \ \ \ \ \ \ \ \ \ \ \ 
\ \ \ \ \ \ \ \ \ \ \ \ \ \ \ \ \ \ \ \ \ \ \ \ \ \
\left.
-
\Delta(x \vert yy) \cdot \Delta( xx \vert y)+
\Delta( xy \vert y)\cdot \Delta( x \vert xy)
\right\} \\
& \ \ \ \ \ \ \ \ \ \ \ \ \ \ \ \ \ \ \ \ \ \ \ \ \
\ \ \ \ \
=
\frac{ 1}{[\Delta(x\vert y)]^2}\left\{
-\Delta( yy \vert y)\cdot \Delta( x \vert xx) +
\Delta( xy \vert y)\cdot \Delta( x \vert xy)
\right\} \\
& \ \ \ \ \ \ \ \ \ \ \ \ \ \ \ \ \ \ \ \ \ \ \ \ \
\ \ \ \ \
=
-
\square_{yy}^0 \cdot \square_{xx}^1 +
\square_{xy}^0 \cdot \square_{xy}^1.
\endaligned\right.
\end{equation}
In sum, the result of the cross differentiation $(\square_{xx}^0)_y-
(\square_{xy}^0)_x$ is a quadratic expression in terms of the square
functions themselves! Following the same recipe (with no surprise),
one may establish the following relations, listing all the compatibility
conditions (the first one is nothing
else than~\thetag{ 2.29}):
\def\theequation{2.30}\begin{equation}
\left\{
\aligned
\left(\square_{xx}^0\right)_y-
\left(\square_{xy}^0\right)_x 
=
& \ 
-\square_{xx}^1 \cdot \square_{yy}^0+
\square_{xy}^1\cdot \square_{xy}^0, 
\\
\left(\square_{xy}^0\right)_y-
\left(\square_{yy}^0\right)_x 
=
& \
-\square_{xy}^0 \cdot \square_{yx}^0 -\square_{xy}^1 \cdot
\square_{yy}^0 +
\square_{yy}^0 \cdot \square_{xx}^0+
\square_{yy}^1 \cdot \square_{xy}^0, 
\\
\left(\square_{xx}^1\right)_y-
\left(\square_{xy}^1\right)_x 
=
& \ 
-\square_{xx}^0 \cdot \square_{yx}^1-
\square_{xx}^1\cdot \square_{yy}^1
+\square_{xy}^0\cdot \square_{xx}^1 +
\square_{xy}^1\cdot \square_{xy}^1,
\\
\left(\square_{xy}^1\right)_y-
\left(\square_{yy}^1\right)_x 
=
& \
-\square_{xy}^0 \cdot \square_{yx}^1 +
\square_{yy}^0 \cdot \square_{xx}^1. 
\\
\endaligned\right.
\end{equation}
Instead of checking patiently each of the remaining three
cross above cross differentiation identities, it is better
to establish directly the following 
general relation.

\def\thelemma{2.31}\begin{lemma}
Remind from~\thetag{ 2.11} that we identify $y^0$ with $x$ and $y^1$
with $y$ and let $0 \leqslant j_1,j_2, j_3, k_1 \leqslant 1$. Then
\def\theequation{2.32}\begin{equation}
\left(
\square_{y^{j_1}y^{j_2}}^{k_1}
\right)_{y^{j_3}}-
\left(
\square_{y^{j_1}y^{j_3}}^{k_1}
\right)_{y^{j_2}}= -
\sum_{k_2=0}^1 \,
\square_{y^{j_1}y^{j_2}}^{k_2}\cdot
\square_{y^{j_3}y^{k_2}}^{k_1}+
\sum_{k_2=0}^1 \,
\square_{y^{j_1}y^{j_3}}^{k_2}\cdot
\square_{y^{j_2}y^{k_2}}^{k_1}.
\end{equation}
\end{lemma}

This lemma is left to the reader; anyway, we shall complete the proof
of a generalization of Lemma~2.31 to the case of $m\geqslant 1$ dependent
variables $(y^1,\dots, y^m)$ in Section~2 below (Lemma~3.40).

Coming back to the first auxiliary system~\thetag{ 2.20}, we therefore
have obtained a necessary and sufficient condition for the existence
of $(X,Y)$: the functions $\Pi_{j^1,j^2}^{k_1}$ should satisfy the
following system of first order partial differential equations, just
obtained from~\thetag{ 2.30} by replacing the square functions by the
Pi functions:
\def\theequation{2.33}\begin{equation}
\left\{
\aligned
\left(\Pi_{0,0}^0\right)_y-
\left(\Pi_{0,1}^0\right)_x 
=
& \ 
-\Pi_{0,0}^1 \cdot \Pi_{1,1}^0+
\Pi_{0,1}^1\cdot \Pi_{0,1}^0, 
\\
\left(\Pi_{0,1}^0\right)_y-
\left(\Pi_{1,1}^0\right)_x
=
& \
-\Pi_{0,1}^0 \cdot 
\Pi_{0,1}^0 -\Pi_{0,1}^1 \cdot
\Pi_{1,1}^0 +
\Pi_{1,1}^0 \cdot \Pi_{0,0}^0+
\Pi_{1,1}^1 \cdot \Pi_{0,1}^0, 
\\
\left(\Pi_{0,0}^1\right)_y-
\left(\Pi_{0,1}^1\right)_x 
=
& \ 
-\Pi_{0,0}^0 \cdot \Pi_{0,1}^1-
\Pi_{0,0}^1\cdot \Pi_{1,1}^1
+\Pi_{0,1}^0\cdot \Pi_{1,1}^1 +
\Pi_{0,1}^1\cdot \Pi_{0,1}^1,
\\
\left(\Pi_{0,1}^1\right)_y-
\left(\Pi_{1,1}^1\right)_x 
=
& \
-\Pi_{0,1}^0 \cdot \Pi_{0,1}^1 +
\Pi_{1,1}^0 \cdot \Pi_{0,0}^1. 
\endaligned\right.
\end{equation}

\subsection*{2.34.~Second auxiliary system}
It is now time to come back to the functions $G$, $H$, $L$ and $M$ and
to get rid of the auxiliary ``Pi'' functions. Unfortunately, we cannot
invert directly the linear system~\thetag{ 2.18}, hence we must choose
two specific square functions as {\sl principal unknowns}, and the
best, from a combinatorial point of view, is to choose
$\square_{xx}^0$ and $\square_{yy}^1$. Remind that by~\thetag{ 2.20},
we have $\square_{xx}^0=\Pi_{0,0}^0$ and $\square_{yy}^1=\Pi_{1,1}^1$.
For clarity, it will be useful to adopt the notational equivalences
\def\theequation{2.35}\begin{equation}
\Theta^0\equiv \Pi_{0,0}^0 \ \ \ \ \ \ \ 
{\rm and} \ \ \ \ \ \ \
\Theta^1 \equiv \Pi_{1,1}^1.
\end{equation}
We may therefore quasi-inverse the linear system~\thetag{ 2.18},
obtaining that the four functions $\Pi_{0,0}^1$, $\Pi_{0,1}^1$,
$\Pi_{0,1}^0$ and $\Pi_{1,1}^0$ may be expressed in terms of the
functions $G$, $H$, $L$ and $M$ and in terms of the remaining two
principal unknowns~\thetag{ 2.35}, which yields:
\def\theequation{2.36}\begin{equation}
\left\{
\aligned
\Pi_{0,0}^1=
& \
\square_{xx}^1 =
-G \\
\Pi_{0,1}^1 =
& \
\square_{xy}^1= 
-\frac{ 1}{2} \, H+
\frac{ 1}{2} \, \Theta^0, \\
\Pi_{0,1}^0 =
& \
\square_{xy}^0 =
\frac{ 1}{2}\, L + \frac{ 1}{2} \, 
\Theta^1 , \\
\Pi_{1,1}^0 =
& \ 
\square_{yy}^0 =
M.
\endaligned\right.
\end{equation}
Replacing now each of these four expressions in the compatibility
conditions of the first auxiliary system~\thetag{ 2.33}, solving the
four equations with respect to $\Theta_y^1$, $\Theta_x^0$,
$\Theta_x^1$ and $\Theta_y^0$, we get after hygienic simplifications
what we shall call the {\sl second auxiliary system}, which is a
complete system of first order partial derivatives in the remaining
two principal unknowns $\Theta^0$ and $\Theta^1$:
\def\theequation{2.37}\begin{equation}
\left\{
\aligned
\Theta_y^1 =
& \
-L_y +2M_x 
+HM
-\frac{ 1}{2}\, L^2 +
M\, \Theta^0+
\frac{ 1}{2} \, \left(\Theta^1\right)^2, \\
\Theta_x^0= 
& \
-2\, G_y +H_x+
G\, L -\frac{ 1}{2}\, H^2
-G\, \Theta^1 +
\frac{ 1}{2}\, \left(\Theta^0\right)^2, \\
\Theta_x^1 =
& \
-\frac{ 2}{3}\, 
H_y+
\frac{ 1}{3}\, L_x +2\, G\, M-\frac{ 1}{2}\, 
H\, L-
\frac{ 1}{2}\, H\, \Theta^1+
\frac{ 1}{2}\, L\, \Theta^0 +
\frac{ 1}{2}\, \Theta^0 \, \Theta^1, \\
\Theta_y^0 =
& \
-\frac{ 1}{3}\, H_y +
\frac{ 2}{3}\, L_x+2\, G\, M 
-\frac{ 1}{2}\, H\, L -
\frac{ 1}{2}\, H\, \Theta^1+
\frac{ 1}{2}\, L\, \Theta^0 +
\frac{ 1}{2}\, \Theta^0 \, \Theta^1.
\endaligned\right.
\end{equation}
We do not comment the intermediate computations, since they
offer no new combinatorial discovery.


\subsection*{2.38.~Precise lexicographic rules}
We group first order derivatives before zeroth order derivatives; in
each group, we respect the lexicographic order of appearance given by
the sequence $G$, $H$, $L$, $M$, $\Theta^0$, $\Theta^1$; we always put
rational coefficient of every differential monomial in its left;
consequently, we accept a minus sign just after an equality sign, as
for instance in ${\rm \thetag{ 2.36}}_1$ and in ${\rm \thetag{
2.37}}_2$; for clarity, we prefer to write a complicated differential
equation as $0=\Phi$, with $0$ on the left, instead of $\Phi=0$, since
$\Phi$ may incoporate 10, 20 and up to 150 monomials, as will happen
for instance in the next sections below.

\subsection*{2.39.~Compatibility conditions for the second 
auxiliary system} Clearly, the necessary and sufficient condition for
the existence of solutions $(\Theta^0, \Theta^1)$ to the second
auxiliary system~\thetag{ 2.37} is that the two
cross differentiations vanish:
\def\theequation{2.40}\begin{equation}
\left\{
\aligned
0 = 
& \
\left(\Theta_x^0\right)_y -
\left(\Theta_y^0\right)_x, \\
0 =
& \
\left(\Theta_x^1\right)_y -
\left(\Theta_y^1\right)_x.
\endaligned\right. 
\end{equation}
Using~\thetag{ 2.37}, we shall see that we exactly obtain the two
second order partial differential equations written 
in Theorem~1.2 {\bf (4)}. For completeness, we shall perform
completely the computation of the first compatibility
condition~\thetag{ 2.40} and leave the second as
an (easy) exercise.

First of all, inserting~\thetag{ 2.37} and using the rule of Leibniz
for the differentiation of a product, let us write the crude result,
performing neither any simplification nor any reordering:
\def\theequation{2.41}\begin{equation}
\left\{
\aligned
0 = 
& \
\left(\Theta_x^0\right)_y -
\left(\Theta_y^0\right)_x \\
= 
& \
-2\, G_{yy}+H_{xy}+
G_y\, L+ G\, L_y-
H\, H_y-G_y\, \Theta^1- 
G\, \Theta_y^1+
\Theta^0\,
\Theta_y^0+ \\
& \
+
\frac{ 1}{3}\, H_{xy}- 
\frac{ 2}{3}\, L_{xx}-
2\, G_x\, M- 
2\, G\, M_x+
\frac{ 1}{2}\, H_x\, L+
\frac{ 1}{2}\, H\, L_x+ \\
& \
+
\frac{ 1}{2}\, H_x\, \Theta^1+
\frac{ 1}{2}\, H\, \Theta_x^1-
\frac{ 1}{2}\, L_x\, \Theta^0-
\frac{ 1}{2}\, L\, \Theta_x^0- 
\frac{ 1}{2}\, \Theta_x^0\, \Theta^1- 
\frac{ 1}{2}\, \Theta^0 \, \Theta_x^1. 
\endaligned\right.
\end{equation}
Next, replacing each first order derivative $\Theta_x^0$,
$\Theta_y^0$, $\Theta_x^1$ and $\Theta_y^1$ occuring in~\thetag{ 2.41}
by its expression given in~\thetag{ 2.37}, we obtain (suffering a
little) as a brute result, before any simplification (except that we
put all second order derivatives in the beginning):
\def\theequation{2.42}\begin{equation}
\left\{
\aligned
0 = 
& \
-2\, G_{yy}\, +\frac{ 4}{3}\, 
H_{xy}- \frac{ 2}{3}\, 
L_{xx}+ \\
& \
+ G_y\, L+
G\, L_y-H\, H_y-
G_y\, \Theta^1+
G\, L_y-2\, G\, M_x- 
G\, H\, M+ \\
& \
+
\frac{ 1}{2}\, G\, (L)^2- 
G\, M\, \Theta^0- 
\frac{ 1}{2}\, G\, (\Theta^1)^2- 
\frac{ 1}{3}\, H_y\, \Theta^0+
\frac{ 2}{3}\, L_x\, \Theta^0+ \\
& \
+
2\, G\, M\, \Theta^0-
\frac{ 1}{2}\, H\, L\, \Theta^0- 
\frac{ 1}{2}\, H\, \Theta^0 \, \Theta^1+
\frac{ 1}{2}\, L\, (\Theta^0)^2+
\frac{ 1}{2}\, (\Theta^0)^2\, \Theta^1- \\
& \
- 2\, G_x\, M- 2\, G\, M_x+
\frac{ 1}{2}\, H_x\, L+
\frac{ 1}{2}\, H\, L_x+\frac{ 1}{2}\, 
H_x\, \Theta^1+\\
& \
+
\frac{ 1}{6}\, H\, L_x -
\frac{ 1}{3}\, H\, 
H_y +\, G\, H\, M-
\frac{ 1}{4}\, 
(H)^2\, L-
\frac{ 1}{4}\, (H)^2\, \Theta^1+ \\
& \
+
\frac{ 1}{4}\, H\, L\, \Theta^0 +
\frac{ 1}{4}\, H\, \Theta^0 \, \Theta^1-
\frac{ 1}{2}\, L_x \, \Theta^0
+\, G_y\, L -\frac{ 1}{2}\, H_x\, L-\frac{ 1}{2}\, 
G\, (L)^2 + \\
& \
+
\frac{ 1}{4}\, H^2\, L+\frac{ 1}{2}\, 
G\, L\, \Theta^1 -
\frac{ 1}{4}\, L\, \left(\Theta^0\right)^2
+\, G_y\, \Theta^1 -\frac{ 1}{2}\, H_x\, \Theta^1- \\
& \
-\frac{ 1}{2}\, 
G\, L\, \Theta^1 +\frac{ 1}{4}\, H^2\, \Theta^1+
\frac{ 1}{2}\, G\, (\Theta^1)^2 -
\frac{ 1}{4}\, \left(\Theta^0\right)^2\, \Theta^1-
\frac{ 1}{6}\, L_x\, \Theta^0 +\\
& \ 
+\frac{ 1}{3}\, 
H_y \, \Theta^0
-\, G\, M\, \Theta^0+
\frac{ 1}{4}\, 
H\, L\, \Theta^0+
\frac{ 1}{4}\, H\, \Theta^0\, \Theta^1-
\frac{ 1}{4}\, L\, (\Theta^0)^2 - \\
& \
-
\frac{ 1}{4}\, (\Theta^0)^2 \, \Theta^1.
\endaligned\right.
\end{equation}
Now, we can simplify this brute expression by chasing every 
couple (or triple, or quadruple) of terms killing each other. 
After (patient) simplification and lexicographic ordering, we 
obtain the equation
\def\theequation{2.43}\begin{equation}
\left\{
\aligned
0 = 
-2\, G_{yy} + 
& \
\frac{ 4}{3}\, H_{xy} 
- \frac{ 2}{3}\, L_{xx}+ \\
& \
+ 2(G\,
L)_y - 2\, G_x \, M - 4\, G\, M_x 
+ \frac{ 2}{3}\, H\, L_x - \frac{
4}{3}\, H\, H_y,
\endaligned\right.
\end{equation}
which is exactly the first equation of {\bf (4)}
of Theorem~1.2. The treatment
of the second one is totally similar.  This completes the proof of the
equivalence between {\bf (1)} and {\bf (4)} in Theorem~1.2.
\qed

\subsection*{ 2.44.~Interlude: about hand-computed
formulas} In Section~4 below, when dealing with several dependent
variables $y^1, \dots, y^m$, many simplifications of identities which
are much more massive than~\thetag{ 2.42} will occur several times. It
is therefore welcome to explain how we manage to achieve such
computations, without mistakes at the end and strictly by hand. One of
the trick is to use colors, which, unfortunately, cannot be restituted
in this printed document. Another trick is to {\it underline and to
number the terms which disappear together}, by pair, by triple, by
quadruple, {\it etc.} This trick is illustrated in the detailed
identity~\thetag{ 2.45} below, extracted from our manuscript, which is
a copy of~\thetag{ 2.42} together with the designation of all the
terms which vanish together. Hence, {\it we keep a written track of
each intermediate step of every computation and of every
simplification}. Checking the correctness of a computation simply by
reading is then the easiest way, both for the writer and for the
reader, although of course, it takes time, anyway.

On the contrary, when relying upon a digital computer, most
intermediate steps are invisible; the chase of mistakes is by reading
the program and by testing it on several instances.  Alas, all the
finest intuitions which may awake in the extreme inside of a long
computation are essentially absent, the mind believing that the
machine is stronger for such tasks.  This last belief is in part true,
in the case where straightforward known computations are concerned,
and in part untrue, in the case where some new hidden mathematical
reality is concerned. 

For us, {\it the challenge is to control everything in a sea of
signs}. Computations are to be organized like a living giant coral
tree, all part of which should be clearly visible in a transparent
fluid of thought, and permanently subject to corrections.

Indeed, it often happens that going through a problem involving
massive formal computations, some disharmony or some incoherency is
discovered. Then one has to inspect every living atom in the preceding
branches of the growing coral tree of computations until some very
tiny or ridiculous mistake is found. In addition to making easy the
reading, {\it a perfectly rigorous way of writing the formal
identities which respects a large amount of virtual conventions
facilitates to reorganize rapidly the coral tree after a mistake has
been found}. The accumulation of new virtual conventions, all of which
we cannot speak, constitutes another coral meta-tree and another
profound collection of trick. Finally, we use a blank fluid corrector
to avoid copying to much.

\smallskip

Extracted from our manuscript, here is the identity~\thetag{ 2.42}
with the underlining-numbering of all the vanishing terms (without
the original colours) until we get the final equation~\thetag{ 2.43}:
\[
\small
\aligned
0 = 
& \
-2\, G_{yy}\, 
+
\frac{ 4}{3}\,H_{xy}
-
\frac{ 2}{3}\, L_{xx}
+ \\
& \
+ 
G_y\, L
+
G\, L_y-H\, H_y
-
\underline{ 
G_y\, \Theta^1 }_{
\octagon \! \! \! \! \! \tiny{\sf o}}
+
G\, L_y-2\, G\, M_x
- 
\underline{ 
G\, H\, M }_{
\octagon \! \! \! \! \! \tiny{\sf k}}
+ \\
& \
+
\underline{ 
\frac{ 1}{2}\, G\, (L)^2}_{
\octagon \! \! \! \! \! \tiny{\sf a}}
-
\underline{ 
G\, M\, \Theta^0 }_{
\octagon \! \! \! \! \! \tiny{\sf b}} 
- 
\underline{ 
\frac{ 1}{2}\, G\, (\Theta^1)^2 }_{
\octagon \! \! \! \! \! \tiny{\sf c}} 
- 
\underline{ 
\frac{ 1}{3}\, H_y\, \Theta^0 }_{
\octagon \! \! \! \! \! \tiny{\sf d}} 
+
\underline{ 
\frac{ 2}{3}\, L_x\, \Theta^0 }_{
\octagon \! \! \! \! \! \tiny{\sf e}} 
+ \\
& \
+
\underline{ 
2\, G\, M\, \Theta^0}_{
\octagon \! \! \! \! \! \tiny{\sf b}} 
-
\underline{ 
\frac{ 1}{2}\, H\, L\, \Theta^0 }_{
\octagon \! \! \! \! \! \tiny{\sf f}} 
- 
\underline{ 
\frac{ 1}{2}\, H\, \Theta^0 \, \Theta^1 }_{
\octagon \! \! \! \! \! \tiny{\sf g}} 
+
\underline{ 
\frac{ 1}{2}\, L\, (\Theta^0)^2 }_{
\octagon \! \! \! \! \! \tiny{\sf h}} 
+
\underline{ 
\frac{ 1}{2}\, (\Theta^0)^2\, \Theta^1 }_{
\octagon \! \! \! \! \! \tiny{\sf i}} 
- \\
& \
- 
2\, G_x\, M- 2\, G\, M_x
+
\frac{ 1}{2}\, H_x\, L
+
\frac{ 1}{2}\, H\, L_x
+
\underline{ 
\frac{ 1}{2}\, H_x\, \Theta^1 }_{
\octagon \! \! \! \! \! \tiny{\sf j}} 
+ \\
& \
+
\frac{ 1}{6}\, H\, L_x 
-
\frac{ 1}{3}\, H\, 
H_y 
+
\underline{ 
G\, H\, M }_{\octagon \! \! \! \! \! \tiny{\sf k}} 
-
\underline{ 
\frac{ 1}{4}\, (H)^2\, L }_{
\octagon \! \! \! \! \! \tiny{\sf l}} 
-
\underline{ 
\frac{ 1}{4}\, (H)^2\, \Theta^1 }_{
\octagon \! \! \! \! \! \! \tiny{\sf m}} 
+
\endaligned
\]
\def\theequation{2.45}\begin{equation}
\aligned
& \
+
\underline{ 
\frac{ 1}{4}\, H\, L\, \Theta^0 }_{
\octagon \! \! \! \! \! \tiny{\sf f}} 
+
\underline{ 
\frac{ 1}{4}\, H\, \Theta^0 \, \Theta^1 }_{
\octagon \! \! \! \! \! \tiny{\sf g}} 
-
\underline{ 
\frac{ 1}{2}\, L_x \, \Theta^0}_{
\octagon \! \! \! \! \! \tiny{\sf e}}
+
G_y\, L 
-
\frac{ 1}{2}\, H_x\, L
-
\underline{
\frac{ 1}{2}\, G\, (L)^2 }_{
\octagon \! \! \! \! \! \tiny{\sf a}} 
+ \\
& \
+
\underline{ 
\frac{ 1}{4}\, H^2\, L }_{
\octagon \! \! \! \! \! \tiny{\sf l}} 
+
\underline{ 
\frac{ 1}{2}\, G\, L\, \Theta^1 }_{
\octagon \! \! \! \! \! \tiny{\sf n}} 
-
\underline{ 
\frac{ 1}{4}\, L\, \left(\Theta^0\right)^2 }_{
\octagon \! \! \! \! \! \tiny{\sf h}} 
+
\underline{ 
G_y\, \Theta^1 }_{
\octagon \! \! \! \! \! \tiny{\sf o}} 
-
\underline{ 
\frac{ 1}{2}\, H_x\, \Theta^1 }_{
\octagon \! \! \! \! \! \tiny{\sf j}} 
- \\
& \
-
\underline{ 
\frac{ 1}{2}\, G\, L\, \Theta^1 }_{
\octagon \! \! \! \! \! \tiny{\sf n}} 
+
\underline{ 
\frac{ 1}{4}\, H^2\, \Theta^1 }_{
\octagon \! \! \! \! \! \! \tiny{\sf m}} 
+
\underline{ 
\frac{ 1}{2}\, G\, (\Theta^1)^2 }_{
\octagon \! \! \! \! \! \tiny{\sf c}} 
-
\underline{ 
\frac{ 1}{4}\, \left(\Theta^0\right)^2\, \Theta^1 }_{
\octagon \! \! \! \! \! \tiny{\sf i}} 
-
\underline{ 
\frac{ 1}{6}\, L_x\, \Theta^0 }_{
\octagon \! \! \! \! \! \tiny{\sf e}} 
+ \\ 
\endaligned
\end{equation}
$$
\aligned
& \ 
+
\underline{
\frac{ 1}{3}\, H_y \, \Theta^0 }_{
\octagon \! \! \! \! \! \tiny{\sf d}}
-
\underline{ G\, M\, \Theta^0 }_{
\octagon \! \! \! \! \! \tiny{\sf b}}
+
\underline{ 
\frac{ 1}{4}\, H\, L\, \Theta^0 }_{
\octagon \! \! \! \! \! \tiny{\sf f}} 
+
\underline{ 
\frac{ 1}{4}\, H\, \Theta^0\, \Theta^1 }_{
\octagon \! \! \! \! \! \tiny{\sf g}} 
-
\underline{ 
\frac{ 1}{4}\, L\, (\Theta^0)^2 }_{
\octagon \! \! \! \! \! \tiny{\sf h}} 
- \\
& \
-
\underline{ 
\frac{ 1}{4}\, (\Theta^0)^2 \, \Theta^1 }_{
\octagon \! \! \! \! \! \tiny{\sf i}}.
\endaligned
$$
As may be observed, the order in which we discover the terms which
vanish is governed by chance. After some terms are underlined, they
are automatically disregarded by the eyes, which lightens the chasing
of other terms to be simplified. To collect the remaining terms in
order to obtain the final expression~\thetag{ 2.43}, our method is
similar: we underline the terms which may be summed together.
However, whereas we use the red pencil to underline the vanishing
terms, we use the green pencil to underline the remaining terms. This
small trick is to avoid as much as possible to copy several times some
long formal expressions. Finally, we reorder everything
lexicographically, so as to get the conclusion~\thetag{ 2.43}. In
order to obtain the final equation~\thetag{ 2.43} as efficiently as
possible, we read the remaining terms, picking them directly in
lexicographic order. If, by lack of luck, one or two terms are
forgotten by the eyes and not written in the right place, we copy once
more the very final result in the right order, or we use the blank
corrector.

Of course, such a refined methodology could seem to be essentially
superfluous for such relatively accessible computations. However, when
passing to several dependent variables, the current expressions will
be approximatively five times more massive.  We may really ascertain
that a clever methodology of hand computations is helpful in this
category.

\section*{\S3.~Systems of second order ordinary
differential equations equivalent to free particles}

\subsection*{3.1.~Combinatorics of the second order
prolongation of a point transformation} In this section, we endeavour
to explain how Lie's theorem and proof may be generalized to the
case of several dependent variables. As in the statement of
Theorem~1.7~{\bf (1)}, let us assume that the system $y_{xx}^j= F^j(x,
y, y_x)$, $j=1,\dots, m$, is equivalent under an invertible point
transformation $(x,y)\mapsto (X(x,y), Y(x, y))$ to the free particle
system $Y_{XX}^j=0$, $j=1, \dots, m$. By assumption, 
the Jacobian determinant
\def\theequation{3.2}\begin{equation}
\Delta( x\vert y^1 \vert \dots \vert y^m):=
\left\vert
\begin{array}{cccc}
X_x & X_{y^1} & \dots & X_{y^m} \\
Y_x^1 & Y_{y^1}^1 & \dots & Y_{y^m}^1 \\
\dots & \dots & \dots & \dots \\
Y_x^m & Y_{y^1}^m & \dots & Y_{y^m}^m
\end{array}
\right\vert
\end{equation}
does not vanish at the origin. As in the case $m=1$, since the flat
system $Y_{XX}^j=0$ is left unchanged by any affine transformation, we
can (and we shall) assume that the transformation is tangent to the
identity at the origin, so that the above Jacobian matrix equals the
identity matrix at $(x,y)=(0,0)$, whence in a neighborhood of the
origin it is close to the identity matrix, namely
\def\theequation{3.3}\begin{equation}
X_x \cong 1, \ \ \ \ \
X_{y^j}\cong 0, \ \ \ \ \
Y_x^j\cong 0, \ \ \ \ \
Y_{y^{j_1}}^{j_2}\cong \delta_{j_1}^{j_2}.
\end{equation}
Inductive formulas for the computation how the differential equation
in the $(X, Y)$ coordinates is related to the differential equation in
the $(x,y)$ coordinates may be found in~\cite{ bk1989}, \cite{
ol1995}; the explicit formulas are not achieved in
these references. Let us recall the inductive formulas, just on the
computational level (differential-geometric conceptional background
about graph transformations may be found in~\cite{ ol1986}, Ch.~2).

First of all, we seek how the $Y_X^j:= \frac{ dY^j}{dX}$ are
explicitely related to the $y_x^l$. It suffices to replace, 
in the identity
\def\theequation{3.4}\begin{equation}
Y_X^j\cdot\left(
X_x\, dx+
\sum_{l=1}^m\, X_{y^l}\, dy^l
\right)=
Y_X^j\, dX= 
dY^j=
Y_x^j\, dx+
\sum_{l=1}^m \, 
Y_{y^l}^j\, dy^l
\end{equation}
the differentials $dy^l$ by $y_x^l\, dx$ and then to identify the
coefficient of $dx$ on both sides, which rapidly yields the formulas
\def\theequation{3.5}\begin{equation}
Y_X^j=
\frac{ Y_x^j+\sum_{l=1}^m\, 
y_x^l\, Y_{y^l}^j}{X_x+
\sum_{l=1}^m\, 
y_x^l\, X_{y^l}},
\end{equation}
for $j=1,\dots,m$. 

Next, we seek how the $Y_{XX}^j:= \frac{ d^2Y^j}{dX^2}= \frac{
dY_X^j}{dX}$ are related to the $y_x^{l_1}$, $y_{xx}^{l_2}$. It
suffices to again replace each $dy_l$ by $y_x^l\, dx$ and each
$dy_x^l$ by $y_{xx}^l\, dx$ in the identity
\def\theequation{3.6}\begin{equation}
\left\{
\aligned
Y_{XX}^j\cdot
\left(
X_x\, dx+
\sum_{l=1}^m\, X_{y^l}\, dy^l
\right)
& \
=
Y_{XX}^j\cdot dX= 
dY_X^j \\
& \
=
\frac{ \partial Y_X^j}{\partial x}\, dx+
\sum_{l=1}^m\, 
\frac{ \partial Y_X^j}{\partial y^l}\, dy^l+
\sum_{l=1}^m\, 
\frac{ \partial Y_X^j}{\partial y_x^l}\, 
dy_x^l \\
& \
=
\left(
\frac{ \partial Y_X^j}{\partial x}+
\sum_{l=1}^m\, 
\frac{ \partial Y_X^j}{\partial y^l}\, 
y_x^l+
\sum_{l=1}^m\, 
\frac{ \partial Y_X^j}{\partial y_x^l}\, 
y_{xx}^l
\right)\cdot dx.
\endaligned\right.
\end{equation}
Before entering the precise combinatorics of the explicit expression
of $Y_{XX}^j$, let us observe that the last term of~\thetag{ 3.6}
simply writes $D(Y_X^j)\, dx$, where $D$ denotes the {\sl total
differentiation operator}\, (of order two) defined by
\def\theequation{3.7}\begin{equation}
D:=\frac{\partial }{\partial x}+
\sum_{l=1}^m\, 
y_x^l\, 
\frac{ \partial }{\partial y^l}+
\sum_{l=1}^m\, 
y_{xx}^l \, 
\frac{\partial }{\partial y_x^l}. 
\end{equation}
Since $dX\equiv DX$ after replacing each $dy^l$ by $y_x^l\, dx$, it
follows that we may compactly rewrite~\thetag{ 3.6} as
\def\theequation{3.8}\begin{equation}
Y_{XX}^j \, DX \cdot dx= D \left(Y_X^j\right)\cdot dx
\end{equation}
Consequently, the expressions of $Y_X^j$ (obtained in~\thetag{3.5})
and of $Y_{XX}^j$ are
\def\theequation{3.9}\begin{equation}
Y_X^j= \frac{ DY^j}{DX} \ \ \ \ \ \ \ 
{\rm and} \ \ \ \ \ \ \
Y_{XX}^j = \frac{ D\left(Y_X^j\right)}{DX}= 
\frac{DD Y^j\cdot DX- DD X \cdot
DY^j}{[DX]^3}.
\end{equation}
As, by assumption, the system $y_{ xx}^j= F^j(x, y, y_x)$ transforms
to the flat system $Y_{ X X}^j=0$, after erasing the denominator
of~\thetag{ 3.7}, we come to the equations
\def\theequation{3.10}\begin{equation}
0 = DD Y^j\cdot DX- 
DD X \cdot DY^j,
\end{equation}
for $j=1,\dots, m$. However, this too simple and too compact
expression of the system $y_{xx}^j= F^j(x, y, y_x)$ is of no use and
we must develope 
(patiently!) the explicit expressions of $DDY^j$, of $DX$, of
$DDX$ and of $DY^j$, using the complete expression of $D$ defined
in~\thetag{ 3.6}.

At this point, we would like to stress that {\it it constitutes
already a nontrivial computational and combinatorial task to obtain a
complete explicit formula for the system $y_{xx}^j=F^j(x, y, y_x)$
hidden in the compact form~\thetag{ 3.10}}, which would be the
generalization of the nice formula~\thetag{ 2.9} involving
modifications of the Jacobian determinant. For general $m\geqslant 2$, the
complete proofs are postponed to Section~5 below.

Since it would be intuitively unsatisfactory 
to provide directly the final
simplified expression of the development of~\thetag{ 3.10} in the
general case $m\geqslant 2$, let us firstly describe step by step how one
may guess what is the generalization of~\thetag{ 2.9}.

For instance, in the case $m=2$, by a direct and relatively short
computation which consists in 
developing plainly~\thetag{ 3.10}, we obtain for $j=1,2$:
\def\theequation{3.11}\begin{equation}
\left\{
\aligned
0 = 
& \
-X_{xx}\, Y_x^j+Y_{xx}^j\, X_x+ \\
& \
+
y_x^1\cdot \left[
-X_{xx}\, Y_{y^1}^j+
Y_{xx}^j\, X_{y^1}-
2\, X_{xy^1}\, Y_x^j+
2\, Y_{xy^1}^j\, X_x
\right]+ \\
& \
+
y_x^2\cdot \left[
-X_{xx}\, Y_{y^2}^j+Y_{xx}^j\, X_{y^2}-
2\, X_{xy^2}\, Y_x^j+
2Y_{xy^2}^j \, X_x
\right]+ \\
& \
+
y_x^1\, y_x^1 \cdot
\left[
-2\, X_{xy^1}\, Y_{y^1}^j+2\, Y_{xy^1}^j\, X_{y^1}-
X_{y^1y^1}\, Y_x^j+
Y_{y^1y^1}^j\, X_x
\right]+ \\
& \
+
y_x^1\, y_x^2\cdot
\left[
-2\, X_{xy^1}\, Y_{y^2}^j+
2\, Y_{xy^1}^j\, X_{y^2}-
2\, X_{xy^2}\, Y_{y^1}^j+
2\, Y_{xy^2}^j\, X_{y^1}- \right. \\
& \ 
\left.
\ \ \ \ \ \ \ \ \ \ \ \ \ \ \ \ \ \ \ \ \ \ 
\ \ \ \ \ \ \ \ \ \ \ \ \ \ \ \ \ \ \ \ \ 
\ \ \ \ \ \ \ \ \ \ \ \ \ \ \ \ \ \ \ \ \
-
2\, X_{y^1y^2}\, Y_x^j+
2\, Y_{y^1y^2}^j\, X_x
\right]+ \\
& \
+
y_x^2\, y_x^2\cdot
\left[
-2\, X_{xy^2}\, Y_{y^2}^j+
2\, Y_{xy^2}^j\, X_{y^2}-
X_{y^2y^2}\, Y_x^j+
Y_{y^2y^2}^j\, X_x
\right]+ \\
& \
+
y_x^1\, y_x^1\, y_x^1\cdot
\left[
-X_{y^1y^1}\, Y_{y^1}^j+
Y_{y^1y^1}^j\, X_{y^1}
\right]+ \\
& \
+
y_x^1\, y_x^1 \, y_x^2\cdot
\left[
-X_{y^1y^1}\, Y_{y^2}^j+
Y_{y^1y^1}^j\, X_{y^2}-
2\, X_{y^1y^2}\, Y_{y^1}^j+
2\, Y_{y^1y^2}^j\, X_{y^1}
\right]+ \\
& \
+y_x^1\, y_x^2\, y_x^2\cdot
\left[
-X_{y^2y^2}\, Y_{y^1}^j+
Y_{y^2y^2}^j\, X_{y^1}-
2\, X_{y^1y^2}\, Y_{y^2}^j+
2\, Y_{y^1y^2}^j\, X_{y^2}
\right]+ \\
& \
+
y_x^2\, y_x^2\, y_x^2\cdot
\left[
-X_{y^2y^2}\, Y_{y^2}^j+
Y_{y^2y^2}^j\, X_{y^2}
\right]+ \\
& \
+
y_{xx}^1\cdot 
\left[ -X_{y^1}\, Y_x^j+Y_{y^1}^j\, X_x+
y_x^2\cdot \left\{
-X_{y^1}\, Y_{y^2}^j+Y_{y^1}^j\, X_{y^2}
\right\}
\right]+ \\
& \
+
y_{xx}^2\cdot 
\left[ -X_{y^2}\, Y_x^j+Y_{y^2}^j\, X_x+
y_x^1\cdot \left\{
-X_{y^2}\, Y_{y^1}^j+Y_{y^2}^j\, X_{y^1}
\right\}
\right].
\endaligned\right.
\end{equation}
Unfortunately, {\it the above two equations are not solved 
with respect to $y_{xx}^1$ and to $y_{xx}^2$}. Consequently,
if we abbreviate them as a linear system of the form 
\def\theequation{3.12}\begin{equation}
\left\{
\aligned
0 = 
& \
A^1+ y_{xx}^1 \cdot B_1^1 + 
y_{xx}^2\cdot B_2^1, \\
0 = 
& \
A^2+ y_{xx}^1 \cdot B_1^2 + 
y_{xx}^2\cdot B_2^2,
\endaligned\right.
\end{equation}
we have to solve for $y_{xx}^1$ and for $y_{xx}^2$ by means of the
classical rule of Cramer. Here, it is rather quick to check 
manually that the
determinant of this system has the following nice expression:
\def\theequation{3.13}\begin{equation}
\left\{
\aligned
\left\vert
\begin{array}{cc}
B_1^1 & B_2^1 \\
B_1^2 & B_2^2
\end{array}
\right\vert
& \
=
\Delta(x \vert y^1 \vert y^2) \cdot
\left\{
X_x+ y_x^1\, X_{y^1}+
y_x^2 \, X_{y^2}
\right\} \\
& \
=\Delta(x \vert y^1 \vert y^2) \cdot DX. 
\endaligned\right.
\end{equation}
However, the complete solving for $y_{xx}^1$ and for $y_{xx}^2$
requires some more time. After a direct and rather long hand
computation (or alternately, using Maple or Mathematica) one obtains
formulas involving hidden $3\times 3$ determinants, which have to be
guessed by the intuition; the first equation that
we obtain, namely for $y_{xx}^1$
is as follows:
\def\theequation{3.14}\begin{equation}
\aligned
0 = 
& \
y_{xx}^1 \cdot 
\left\vert
\begin{array}{ccc}
X_{x} & X_{y^1} & X_{y^2} \\
Y_{x}^1 & Y_{y^1}^1 & Y_{y^2}^1 \\
Y_{x}^2 & Y_{y^1}^2 & Y_{y^2}^2
\end{array}
\right\vert+
\left\vert
\begin{array}{ccc}
X_{x} & X_{xx} & X_{y^2} \\
Y_{x}^1 & Y_{xx}^1 & Y_{y^2}^1 \\
Y_{x}^2 & Y_{xx}^2 & Y_{y^2}^2
\end{array}
\right\vert+ \\
& \
+ y_x^1\cdot 
\left\{
2\,
\left\vert
\begin{array}{ccc}
X_{x} & X_{xy^1} & X_{y^2} \\
Y_{x}^1 & Y_{xy^1}^1 & Y_{y^2}^1 \\
Y_{x}^2 & Y_{xy^1}^2 & Y_{y^2}^2
\end{array}
\right\vert
-
\left\vert
\begin{array}{ccc}
X_{xx} & X_{y^1} & X_{y^2} \\
Y_{xx}^1 & Y_{y^1}^1 & Y_{y^2}^1 \\
Y_{xx}^2 & Y_{y^1}^2 & Y_{y^2}^2
\end{array}
\right\vert
\right\}+ \\
& \
+y_x^2\cdot 
\left\{
2\, 
\left\vert
\begin{array}{ccc}
X_{x} & X_{xy^2} & X_{y^2} \\
Y_{x}^1 & Y_{xy^2}^1 & Y_{y^2}^1 \\
Y_{x}^2 & Y_{xy^2}^2 & Y_{y^2}^2
\end{array}
\right\vert
\right\}+ \\
\endaligned
\end{equation}
$$
\aligned
{}
& \
+
y_x^1 \, y_x^1\cdot
\left\{
\left\vert
\begin{array}{ccc}
X_{x} & X_{y^1y^1} & X_{y^2} \\
Y_{x}^1 & Y_{y^1y^1}^1 & Y_{y^2}^1 \\
Y_{x}^2 & Y_{y^1y^1}^2 & Y_{y^2}^2
\end{array}
\right\vert
- 2\, 
\left\vert
\begin{array}{ccc}
X_{xy^1} & X_{y^1} & X_{y^2} \\
Y_{xy^1}^1 & Y_{y^1}^1 & Y_{y^2}^1 \\
Y_{xy^1}^2 & Y_{y^1}^2 & Y_{y^2}^2
\end{array}
\right\vert
\right\}+ \\
& \
+
y_x^1 \, y_x^2 \, 
\left\{
2\, 
\left\vert
\begin{array}{ccc}
X_{x} & X_{y^1y^2} & X_{y^2} \\
Y_{x}^1 & Y_{y^1y^2}^1 & Y_{y^2}^1 \\
Y_{x}^2 & Y_{y^1y^2}^2 & Y_{y^2}^2
\end{array}
\right\vert
- 2\, 
\left\vert
\begin{array}{ccc}
X_{xy^2} & X_{y^1} & X_{y^2} \\
Y_{xy^2}^1 & Y_{y^1}^1 & Y_{y^2}^1 \\
Y_{xy^2}^2 & Y_{y^1}^2 & Y_{y^2}^2
\end{array}
\right\vert
\right\}+ \\
& \
+
y_x^2\, y_x^2 \cdot
\left\{
\left\vert
\begin{array}{ccc}
X_{x} & X_{y^2y^2} & X_{y^2} \\
Y_{x}^1 & Y_{y^2y^2}^1 & Y_{y^2}^1 \\
Y_{x}^2 & Y_{y^2y^2}^2 & Y_{y^2}^2
\end{array}
\right\vert
\right\}+
\endaligned
$$
$$
\aligned
{}
& \
+
y_x^1\, y_x^1 \, y_x^1 \cdot 
\left\{
-
\left\vert
\begin{array}{ccc}
X_{y^1y^1} & X_{y^1} & X_{y^2} \\
Y_{y^1y^1}^1 & Y_{y^1}^1 & Y_{y^2}^1 \\
Y_{y^1y^1}^2 & Y_{y^1}^2 & Y_{y^2}^2
\end{array}
\right\vert
\right\}
+ 
y_x^1\, y_x^1 \, y_x^2 \cdot
\left\{
-2 \, 
\left\vert
\begin{array}{ccc}
X_{y^1y^2} & X_{y^1} & X_{y^2} \\
Y_{y^1y^2}^1 & Y_{y^1}^1 & Y_{y^2}^1 \\
Y_{y^1y^2}^2 & Y_{y^1}^2 & Y_{y^2}^2
\end{array}
\right\vert 
\right\}+ \\
& \
+
y_x^1\, y_x^2 \, y_x^2 \cdot
\left\{
- \, 
\left\vert
\begin{array}{ccc}
X_{y^2y^2} & X_{y^1} & X_{y^2} \\
Y_{y^2y^2}^1 & Y_{y^1}^1 & Y_{y^2}^1 \\
Y_{y^2y^2}^2 & Y_{y^1}^2 & Y_{y^2}^2
\end{array}
\right\vert
\right\}.
\endaligned
$$
This formula and the next have been checked by Sylvain Neut
and Michel Petitot with the help of Maple.
We notice that the size is not negligible, but fortunately, there
appears some combinatorics, much more visible than in~\thetag{ 3.11}.
The second equation that we obtain, namely for $y_{xx}^2$, is as
follows: 
\def\theequation{3.15}\begin{equation}
\aligned
0 = 
& \
y_{xx}^2 \cdot 
\left\vert
\begin{array}{ccc}
X_{x} & X_{y^1} & X_{y^2} \\
Y_{x}^1 & Y_{y^1}^1 & Y_{y^2}^1 \\
Y_{x}^2 & Y_{y^1}^2 & Y_{y^2}^2
\end{array}
\right\vert+
\left\vert
\begin{array}{ccc}
X_{x} & X_{y^1} & X_{xx} \\
Y_{x}^1 & Y_{y^1}^1 & Y_{xx}^1 \\
Y_{x}^2 & Y_{y^2}^2 & Y_{xx}^2
\end{array}
\right\vert+ \\
& \
+ y_x^1\cdot 
\left\{
2\,
\left\vert
\begin{array}{ccc}
X_{x} & X_{y^1} & X_{xy^1} \\
Y_{x}^1 & Y_{y^1}^1 & Y_{xy^1}^1 \\
Y_{x}^2 & Y_{y^1}^2 & Y_{xy^1}^2
\end{array}
\right\vert
\right\}+ \\
& \
+y_x^2\cdot 
\left\{
2\, 
\left\vert
\begin{array}{ccc}
X_{x} & X_{y^1} & X_{xy^2} \\
Y_{x}^1 & Y_{y^1}^1 & Y_{xy^2}^1 \\
Y_{x}^2 & Y_{y^1}^2 & Y_{xy^2}^2
\end{array}
\right\vert
- 
\left\vert
\begin{array}{ccc}
X_{xx} & X_{y^1} & X_{y^2} \\
Y_{xx}^1 & Y_{y^1}^1 & Y_{y^2}^1 \\
Y_{xx}^2 & Y_{y^1}^2 & Y_{y^2}^2
\end{array}
\right\vert
\right\}+ \\
\endaligned
\end{equation}
$$
\aligned
{}
& \
+
y_x^1 \, y_x^1\cdot
\left\{
\left\vert
\begin{array}{ccc}
X_{x} & X_{y^1} & X_{y^1y^1} \\
Y_{x}^1 & Y_{y^1}^1 & Y_{y^1y^1}^1 \\
Y_{x}^2 & Y_{y^1}^2 & Y_{y^1y^1}^2
\end{array}
\right\vert
\right\}+ \\
& \
+
y_x^1 \, y_x^2 \, 
\left\{
2\, 
\left\vert
\begin{array}{ccc}
X_{x} & X_{y^1} & X_{y^1y^2} \\
Y_{x}^1 & Y_{y^1}^1 & Y_{y^1y^2}^1 \\
Y_{x}^2 & Y_{y^1}^2 & Y_{y^1y^2}^2
\end{array}
\right\vert
- 2\, 
\left\vert
\begin{array}{ccc}
X_{xy^1} & X_{y^1} & X_{y^2} \\
Y_{xy^1}^1 & Y_{y^1}^1 & Y_{y^2}^1 \\
Y_{xy^1}^2 & Y_{y^1}^2 & Y_{y^2}^2
\end{array}
\right\vert
\right\}+ \\
& \
+
y_x^2\, y_x^2 \cdot
\left\{
\left\vert
\begin{array}{ccc}
X_{x} & X_{y^1} & X_{y^2y^2} \\
Y_{x}^1 & Y_{y^1}^1 & Y_{y^2y^2}^1 \\
Y_{x}^2 & Y_{y^1}^2 & Y_{y^2y^2}^2
\end{array}
\right\vert
-
2\, 
\left\vert
\begin{array}{ccc}
X_{xy^2} & X_{y^1} & X_{y^2} \\
Y_{xy^2}^1 & Y_{y^1}^1 & Y_{y^2}^1 \\
Y_{xy^2}^2 & Y_{y^1}^2 & Y_{y^2}^2
\end{array}
\right\vert
\right\}
+ \\
\endaligned
$$
$$
\aligned
{}
&
+
y_x^1\, y_x^1 \, y_x^2 \cdot 
\left\{
-
\left\vert
\begin{array}{ccc}
X_{y^1y^1} & X_{y^1} & X_{y^2} \\
Y_{y^1y^1}^1 & Y_{y^1}^1 & Y_{y^2}^1 \\
Y_{y^1y^1}^2 & Y_{y^1}^2 & Y_{y^2}^2
\end{array}
\right\vert
\right\}
+ 
y_x^1\, y_x^2 \, y_x^2 \cdot
\left\{
-2 \, 
\left\vert
\begin{array}{ccc}
X_{y^1y^2} & X_{y^1} & X_{y^2} \\
Y_{y^1y^2}^1 & Y_{y^1}^1 & Y_{y^2}^1 \\
Y_{y^1y^2}^2 & Y_{y^1}^2 & Y_{y^2}^2
\end{array}
\right\vert
\right\}+ \\
& \
+
y_x^2\, y_x^2 \, y_x^2 \cdot
\left\{
- \, 
\left\vert
\begin{array}{ccc}
X_{y^2y^2} & X_{y^1} & X_{y^2} \\
Y_{y^2y^2}^1 & Y_{y^1}^1 & Y_{y^2}^1 \\
Y_{y^2y^2}^2 & Y_{y^1}^2 & Y_{y^2}^2
\end{array}
\right\vert
\right\}.
\endaligned
$$
Importantly, the obtained formulas seem to be analogous to the
formula~\thetag{ 2.9}, since we observe that the coefficients of the
degree three polynomial in the $y_x^l$ are modifications of the
Jacobian determinant $\Delta(x \vert y^1 \vert y^2)$.

To describe the underlying combinatorics, let us observe that there
exist exactly six possible distinct second order derivatives: $xx$,
$xy^1$, $xy^2$, $y^1y^1$, $y^1y^2$ and $y^2y^2$. There are also
exactly three columns in the Jacobian determinant~\thetag{ 3.2}. By
replacing each of the three columns of first order derivatives by a
column of second order detivatives (leaving $X$, $Y^1$ and $Y^2$
unchanged), we may build exactly eighteen new determinants
\def\theequation{3.16}\begin{equation}
\left\{
\aligned
& \
\Delta(xx \vert y^1 \vert y^2) \ \ \ \ \
&
\Delta(x \vert xx \vert y^2) \ \ \ \ \ 
&
\Delta(x \vert y^1 \vert xx) \ \ \ \ \ \\
& \
\Delta(xy^1 \vert y^1 \vert y^2) \ \ \ \ \
&
\Delta(x \vert xy^1 \vert y^2) \ \ \ \ \ 
&
\Delta(x \vert y^1 \vert xy^1) \ \ \ \ \ \\
& \
\Delta(xy^2 \vert y^1 \vert y^2) \ \ \ \ \
&
\Delta(x \vert xy^2 \vert y^2) \ \ \ \ \ 
&
\Delta(x \vert y^1 \vert xy^2) \ \ \ \ \ \\
& \
\Delta(y^1y^1 \vert y^1 \vert y^2) \ \ \ \ \
&
\Delta(x \vert y^1y^1 \vert y^2) \ \ \ \ \ 
&
\Delta(x \vert y^1 \vert y^1y^1) \ \ \ \ \ \\
& \
\Delta(y^1y^2 \vert y^1 \vert y^2) \ \ \ \ \
&
\Delta(x \vert y^1y^2 \vert y^2) \ \ \ \ \ 
&
\Delta(x \vert y^1 \vert y^1y^2) \ \ \ \ \ \\
& \
\Delta(y^2y^2 \vert y^1 \vert y^2) \ \ \ \ \
&
\Delta(x \vert y^2y^2 \vert y^2) \ \ \ \ \ 
&
\Delta(x \vert y^1 \vert y^2y^2), \ \ \ \ \ \\
\endaligned\right.
\end{equation}
where for instance
\def\theequation{3.17}\begin{equation}
\left\{
\aligned
\Delta(\underline{y^1y^2}\vert y^1 \vert y^2):= 
& \
\left\vert
\begin{array}{ccc}
X_{\underline{y^1y^2}} & X_{y^1} & X_{y^2} \\
Y_{\underline{y^1y^2}}^1 & Y_{y^1}^1 & Y_{y^2}^1 \\
Y_{\underline{y^1y^2}}^2 & Y_{y^1}^2 & Y_{y^2}^2 \\
\end{array}
\right\vert \ \ \ \ \ {\rm and} \\
\Delta(x\vert y^1 \vert \underline{xy^2}):= 
& \
\left\vert
\begin{array}{ccc}
X_x & X_{y^1} & X_{\underline{xy^2}} \\
Y_x^1 & Y_{y^1}^1 & Y_{\underline{xy^2}}^1 \\
Y_x^2 & Y_{y^1}^2 & Y_{\underline{xy^2}}^2 \\
\end{array}
\right\vert.
\endaligned\right.
\end{equation}
Hence, using the $\Delta$-notation, we may rewrite the two
equation~\thetag{ 3.14} and~\thetag{ 3.15} 
under a more compact form; after division by 
the Jacobian determinant $\Delta ( x \vert y^1 \vert
y^2)$, the
first equation becomes:
\def\theequation{3.18}\begin{equation}
\left\{
\aligned
0 =
& \
y_{xx}^1 +
\frac{\Delta( x \vert xx \vert y^2)}{\Delta(
x \vert y^1 \vert y^2)}+
y_x^1 \cdot
\left\{
2\, 
\frac{\Delta( x \vert xy^1 \vert y^2)}{
\Delta(x \vert y^1 \vert y^2)}
-
\frac{\Delta( xx \vert y^1 \vert y^2)}{
\Delta(x \vert y^1 \vert y^2)}
\right\}+ \\
& \
+
y_x^2 \cdot
\left\{
2\, 
\frac{\Delta( x \vert xy^2 \vert y^2)}{
\Delta(x \vert y^1 \vert y^2)}
\right\}+
y_x^1\, y_x^1\cdot 
\left\{
\frac{\Delta( x \vert y^1y^1 \vert y^2)}{
\Delta(x \vert y^1 \vert y^2)}-
2\, 
\frac{\Delta( xy^1 \vert y^1 \vert y^2)}{
\Delta(x \vert y^1 \vert y^2)}
\right\} + \\
& \
+
y_x^1 \, y_x^2\cdot
\left\{
2\, 
\frac{\Delta( x \vert y^1y^2 \vert y^2)}{
\Delta(x \vert y^1 \vert y^2)}
- 2\, 
\frac{\Delta( xy^2 \vert y^1 \vert y^2)}{
\Delta(x \vert y^1 \vert y^2)}
\right\}+
y_x^2\, y_x^2 \cdot
\left\{
\frac{\Delta( x \vert y^2y^2 \vert y^2)}{
\Delta(x \vert y^1 \vert y^2)}
\right\}+ \\
& \
+ y_x^1 \, y_x^1 \, y_x^1 \cdot 
\left\{
- 
\frac{\Delta( y^1 y^1 \vert y^1 \vert y^2)}{
\Delta(x \vert y^1 \vert y^2)}
\right\}
+ y_x^1 \, y_x^1 \, y_x^2 \cdot 
\left\{
- 2\, 
\frac{\Delta( y^1 y^2 \vert y^1 \vert y^2)}{
\Delta(x \vert y^1 \vert y^2)}
\right\} + \\
& \
+ y_x^1 \, y_x^2 \, y_x^2 \cdot
\left\{
-
\frac{\Delta( y^2 y^2 \vert y^1 \vert y^2)}{
\Delta(x \vert y^1 \vert y^2)}
\right\}.
\endaligned\right.
\end{equation}
Similarly,
the second equation takes the form:
\def\theequation{3.19}\begin{equation}
\left\{
\aligned
0 =
& \
y_{xx}^2 +
\frac{\Delta( x \vert y^1 \vert xx)}{\Delta(
x \vert y^1 \vert y^2)}+
y_x^1 \cdot
\left\{
2\, 
\frac{\Delta( x \vert y^1 \vert xy^1)}{
\Delta(x \vert y^1 \vert y^2)}
\right\}+ \\
& \
+
y_x^2 \cdot
\left\{
2\, 
\frac{\Delta( x \vert y^1 \vert xy^2)}{
\Delta(x \vert y^1 \vert y^2)}
-
\frac{\Delta( xx \vert y^1 \vert y^2)}{
\Delta(x \vert y^1 \vert y^2)}
\right\}
+
y_x^1\, y_x^1\cdot 
\left\{
\frac{\Delta( x \vert y^1 \vert y^1y^1)}{
\Delta(x \vert y^1 \vert y^2)}
\right\} + \\
& \
+ 
y_x^1 \, y_x^2\cdot
\left\{
2\, 
\frac{\Delta( x \vert y^1 \vert y^1y^2)}{
\Delta(x \vert y^1 \vert y^2)}
- 2\, 
\frac{\Delta( xy^1 \vert y^1 \vert y^2)}{
\Delta(x \vert y^1 \vert y^2)}
\right\}
+ \\
& \
+
y_x^2\, y_x^2 \cdot
\left\{
\frac{\Delta( x \vert y^2 \vert y^2y^2)}{
\Delta(x \vert y^1 \vert y^2)}-
2\, 
\frac{ \Delta( xy^2 \vert y^1 \vert y^2)}{
\Delta(x \vert y^1 \vert y^2)}
\right\}+ \\
& \
+ y_x^1 \, y_x^1 \, y_x^2 \cdot
\left\{
- 
\frac{\Delta( y^1 y^1 \vert y^1 \vert y^2)}{
\Delta(x \vert y^1 \vert y^2)}
\right\}
+ y_x^1 \, y_x^2 \, y_x^2 \cdot 
\left\{
- 2\, 
\frac{\Delta( y^1 y^2 \vert y^1 \vert y^2)}{
\Delta(x \vert y^1 \vert y^2)}
\right\} + \\
& \
+ y_x^2 \, y_x^2 \, y_x^2 \cdot 
\left\{
-
\frac{\Delta( y^2 y^2 \vert y^1 \vert y^2)}{
\Delta(x \vert y^1 \vert y^2)}
\right\}.
\endaligned\right.
\end{equation}
Since the formulas are still of a consequent size, analogously to what
was achieved in Section~2, we shall introduce a new family of {\sl
square functions} as follows. We first index the coordinates $(x,
y^1,\dots, y^m)$ as $(y^0, y^1, \dots, y^m)$, namely we introduce the
notational equivalence
\def\theequation{3.20}\begin{equation}
\fbox{$y^0 \equiv x$},
\end{equation}
which will be very convenient in the sequel, especially in order to
write general formulas. With this convention at hand, our eighteen
square functions $\square_{ y^{l_1} y^{l_2 }}^{ k_1}$, defined for $0
\leqslant j_1, j_2, k_1 \leqslant 2$ are defined by
\def\theequation{3.21}\begin{equation}
\small
\left\{
\aligned
\square_{xx}^0:=&\ \frac{ \Delta(xx\vert y^1 \vert y^2)}{
\Delta(x\vert y^1\vert y^2)} 
\ \ \ \ \ & 
\square_{xy^1}^0:=&\ \frac{ \Delta(xy^1\vert y^1\vert y^2)}{
\Delta(x\vert y^1 \vert y^2)} 
\ \ \ \ \ & 
\square_{xy^2}^0:=&\ \frac{ \Delta(xy^2\vert y^1 \vert y^2)}{
\Delta(x\vert y^1 \vert y^2)} 
\ \ \ \ \ \\
\square_{y^1y^1}^0:=&\ \frac{ \Delta(y^1y^1\vert y^1 \vert y^2)}{
\Delta(x\vert y^1\vert y^2)} 
\ \ \ \ \ & 
\square_{y^1y^2}^0:=&\ \frac{ \Delta(y^1y^2\vert y^1\vert y^2)}{
\Delta(x\vert y^1 \vert y^2)} 
\ \ \ \ \ & 
\square_{y^2y^2}^0:=&\ \frac{ \Delta(y^2y^2\vert y^1 \vert y^2)}{
\Delta(x\vert y^1 \vert y^2)} 
\ \ \ \ \ \\
\square_{xx}^1:=&\ \frac{ \Delta(x\vert xx \vert y^2)}{
\Delta(x\vert y^1\vert y^2)} 
\ \ \ \ \ & 
\square_{xy^1}^1:=&\ \frac{ \Delta(x\vert xy^1\vert y^2)}{
\Delta(x\vert y^1 \vert y^2)} 
\ \ \ \ \ & 
\square_{xy^2}^1:=&\ \frac{ \Delta(x\vert xy^2 \vert y^2)}{
\Delta(x\vert y^1 \vert y^2)} 
\ \ \ \ \ \\
\square_{y^1y^1}^1:=&\ \frac{ \Delta(x\vert y^1y^1 \vert y^2)}{
\Delta(x\vert y^1\vert y^2)} 
\ \ \ \ \ & 
\square_{y^1y^2}^1:=&\ \frac{ \Delta(x\vert y^1y^2\vert y^2)}{
\Delta(x\vert y^1 \vert y^2)} 
\ \ \ \ \ & 
\square_{y^2y^2}^1:=&\ \frac{ \Delta(x\vert y^2y^2 \vert y^2)}{
\Delta(x\vert y^1 \vert y^2)} 
\ \ \ \ \ \\
\square_{xx}^2:=&\ \frac{ \Delta(x\vert y^1 \vert xx)}{
\Delta(x\vert y^1\vert y^2)} 
\ \ \ \ \ & 
\square_{xy^1}^2:=&\ \frac{ \Delta(x\vert y^1\vert xy^1)}{
\Delta(x\vert y^1 \vert y^2)} 
\ \ \ \ \ & 
\square_{xy^2}^2:=&\ \frac{ \Delta(x\vert y^1 \vert xy^2)}{
\Delta(x\vert y^1 \vert y^2)} 
\ \ \ \ \ \\
\square_{y^1y^1}^2:=&\ \frac{ \Delta(x\vert y^1 \vert y^1y^1)}{
\Delta(x\vert y^1\vert y^2)} 
\ \ \ \ \ & 
\square_{y^1y^2}^2:=&\ \frac{ \Delta(x\vert y^1\vert y^1y^2)}{
\Delta(x\vert y^1 \vert y^2)} 
\ \ \ \ \ & 
\square_{y^2y^2}^2:=&\ \frac{ \Delta(x\vert y^1 \vert y^2y^2)}{
\Delta(x\vert y^1 \vert y^2)} 
\ \ \ \ \ \\
\endaligned\right.
\end{equation} 
Obviously, the square functions are symmetric with respect to the
lower indices: $\square_{ y^{ l_1} y^{ l_2}}^{ k_1} = \square_{ y^{
l_2} y^{ l_1}}^{ k_1}$. Here, the upper index designates the column
upon which the second order derivative appears, itself being encoded
by the two lower indices. Even if this is hidden in the notation, we
shall remember that {\it the square functions are explicit rational
expressions in terms of the second order jet of the transformation
$(x,y)\mapsto (X(x,y), Y(x,y))$}. At this point, we may summarize what
we have established so far.

\def\thelemma{3.22}\begin{lemma}
The system of {\rm two} second order ordinary differential equations
$y_{xx}^1= F^1(x, y, y_x)$ and $y_{xx}^2= F^2(x, y, y_x)$ is
equivalent, under a point transformation, to the flat system
$Y_{XX}^1=0$ and $Y_{XX}^2=0$ if and only if there exist three local
$\K$-analytic functions $X(x, y)$, $Y^1(x,y)$ and $Y^2(x,y)$ such that
it may be written under the form
\def\theequation{3.23}\begin{equation}
\left\{
\aligned
0 =
& \
y_{xx}^1+
\square_{xx}^1+
y_x^1\cdot \left(
2\, \square_{xy^1}^1- 
\square_{xx}^0\right)+
y_x^2\cdot \left( 2\, 
\square_{xy^2}^1\right)+\\
& \
+ 
y_x^1 \, y_x^1 \cdot 
\left(
\square_{y^1y^1}^1- 2\, \square_{xy^1}^0
\right)+
y_x^1 y_x^2\cdot 
\left(
2\, \square_{y^1y^2}^1 -
2\, \square_{xy^2}^0
\right)+
y_x^2 \, y_x^2 \cdot
\left(
\square_{y^2y^2}^1
\right)+ \\
& \
+
y_x^1 \, y_x^1 \, y_x^1 \cdot
\left(
-\square_{y^1y^1}^0
\right)
+y_x^1 \, y_x^1 \, y_x^2 \cdot
\left(
-2\,\square_{y^1y^2}^0
\right)
+y_x^1 \, y_x^2 \, y_x^2 \cdot
\left(
-\square_{y^2y^2}^0
\right), \\
0 =
& \
y_{xx}^2+
\square_{xx}^2+
y_x^1\cdot \left(
2\, \square_{xy^1}^2\right)+
y_x^2\cdot \left( 2\, \square_{xy^2}^2-
\square_{xx}^0 \right)+\\
& \
+ 
y_x^1 \, y_x^1 \cdot 
\left(
\square_{y^1y^1}^2
\right)+
y_x^1 y_x^2\cdot
\left(
2\, \square_{y^1y^2}^2 -
2\, \square_{xy^1}^0
\right)+
y_x^2 \, y_x^2 \cdot
\left(
\square_{y^2y^2}^2-2\, \square_{xy^2}^0
\right)+ \\
& \
+
y_x^1 \, y_x^1 \, y_x^2 \cdot
\left(
-\square_{y^1y^1}^0
\right)
+y_x^1 \, y_x^2 \, y_x^2 \cdot
\left(
-2\,\square_{y^1y^2}^0
\right)
+y_x^2 \, y_x^2 \, y_x^2 \cdot
\left(
-\square_{y^2y^2}^0
\right).
\endaligned\right.
\end{equation}
\end{lemma}

\subsection*{3.24.~Second Lie prolongation of 
a vector field} At this point, instead of proceeding further with the
case $m=2$, it is now time to pass to the general case $m\geqslant 2$.
First of all, we would like to remind from~\cite{ gm2003} the complete
explicit expression of the point prolongation to the second order jet
space of a general vector field of the form $L= X\, \frac{ \partial
}{\partial x} + \sum_{j=1}^m\, Y^j\, \frac{ \partial }{\partial y^j}$:
it is a vector field of the form
\def\theequation{3.25}\begin{equation}
L^{(2)}= 
X\, \frac{ \partial }{\partial x}+
\sum_{j=1}^m\, 
Y^j\, \frac{\partial }{\partial y^j}+
\sum_{j=1}^m\, {\bf R}_1^j \, 
\frac{ \partial }{\partial y_x^j}+
\sum_{j=1}^m\,
{\bf R}_2^j\, 
\frac{ \partial }{\partial y_{xx}^j}, 
\end{equation}
where the coefficients ${\bf R}_1^j$ and ${\bf R}_2^j$ are
polynomials in the jet space variables having as
coefficients certain specific linear combinations
of first and second order derivatives of 
$X$ and of the $Y^j$:
\def\theequation{3.26}\begin{equation}
\left\{
\aligned
{\bf R}_1^j =
& \
Y_x^j+\sum_{l_1=1}^m\, 
y_x^{l_1}\cdot\left[
Y_{y^{l_1}}^j-\delta_{l_1}^j\, X_x
\right]+
\sum_{l_1=1}^m\, \sum_{l_2=1}^m\, 
y_x^{l_1}\, y_x^{l_2}\, \cdot 
\left[
-\delta_{l_1}^j\, X_{y^{l_2}}
\right], 
\\
{\bf R}_2^j =
& \
Y_{xx}^j+
\sum_{l_1=1}^m\, 
y_x^{l_1}\cdot
\left[
2\, Y_{xy^{l_1}}^j- 
\delta_{l_1}^j\, X_{xx}
\right]+ \\
& \
+
\sum_{l_1=1}^m\, 
\sum_{l_2=1}^m\, 
y_x^{l_1} \, y_x^{l_2}\cdot
\left[
Y_{y^{l_1}y^{l_2}}^j- 
\delta_{l_1}^j\, X_{xy^{l_2}}-
\delta_{l_2}^j\, X_{xy^{l_1}}
\right]+ \\
& \
+
\sum_{l_1=1}^m\, \sum_{l_2=1}^m\, 
\sum_{l_3=1}^m\, 
y_x^{l_1} \, y_x^{l_2} \, y_x^{l_3}\cdot
\left[
-\delta_{l_1}^j\, X_{y^{l_2} y^{l_3}}
\right]+ \\
& \
+\sum_{l_1=1}^m\, 
y_{xx}^{l_1}\cdot
\left[
Y_{y^{l_1}}^j- 2\, \delta_{l_1}^j\, X_x
\right] + \\
& \
+\sum_{l_1=1}^m\, \sum_{l_2=1}^m\, 
y_x^{l_1}\, 
y_{xx}^{l_2}\cdot
\left[
-\delta_{l_1}^j\, X_{y^{l_2}}- 2\, 
\delta_{l_2}^j\, X_{y^{l_1}}
\right].
\endaligned\right.
\end{equation} 
However, since the notations in~\cite{ gm2003} are different and since
the general case of $n\geqslant 1$ independent variables and $m\geqslant 1$
dependent variables is considered there, it is certainly easier to
reconstiture formulas~\thetag{ 3.26} directly by means of the
inductive formulas described in~\cite{ ol1986}, \cite{ bk1989}). 

Analogously to the observation made in Section~2, we guess that there
exists a formal correspondence between the terms of ${\bf R }_2^j$ not
involving $y_{x x}^l$ and the explicit form of the equation $y_{xx }^j
=F^j(x, y, y_x)$ equivalent to $Y_{XX}^j=0$. In the case $m=2$, we
claim that this formal correspondence also holds true. Indeed, it
suffices to write formula~\thetag{ 3.26} for ${\bf R}_2^j$ modulo the
$y_{x x}^l$, which yields two expressions in total analogy with the
two explicit polynomials appearing in the right-hand side of~\thetag{
3.23}:
\def\theequation{3.27}\begin{equation}
\left\{
\aligned
{\bf R}_2^1 \ ({\rm mod} \ y_{xx}^l) \equiv
& \
Y_{xx}^1 +
y_x^1\cdot \left\{
2\, Y_{xy^1}^1- X_{xx}
\right\}+
y_x^2\cdot \left\{
2\, Y_{xy^2}^1
\right\}+
y_x^1\, y_x^1 \cdot \left\{
Y_{y^1y^1}^1- 2\, X_{xy^1}
\right\}+ \\
& \
+ y_x^1 y_x^2 \cdot \left\{
2\, Y_{y^1y^2}^1 - 2\, X_{xy^2}
\right\}+ 
y_x^2 \, y_x^2 \cdot 
\left\{
Y_{y^2y^2}^1\right\}+ \\
& \
+y_x^1\, y_x^1 \, y_x^1 \cdot
\left\{
-X_{y^1y^1}\right\}+
y_x^1\, y_x^1 \, y_x^2 \cdot
\left\{
-2\, X_{y^1y^2}
\right\}+
y_x^1 \, y_x^2\, y_x^2\cdot
\left\{-X_{y^2y^2}\right\}, \\
{\bf R}_2^2 \ ({\rm mod} \ y_{xx}^l) \equiv
& \
Y_{xx}^2 +
y_x^1\cdot \left\{
2\, Y_{xy^1}^2
\right\}+
y_x^2\cdot \left\{
2\, Y_{xy^2}^2 - X_{xx}
\right\}+
y_x^1 \, y_x^1 \cdot \left\{
Y_{y^1y^1}^2
\right\}+ \\
& \
+ y_x^1 y_x^2 \cdot \left\{
2\, Y_{y^1y^2}^2 - 2\, X_{xy^1}
\right\}+ 
y_x^2 \, y_x^2 \cdot 
\left\{
Y_{y^2y^2}^2-2\, X_{xy^2}\right\}+ \\
& \
+y_x^1\, y_x^1 \, y_x^2 \cdot
\left\{
-X_{y^1y^1}\right\}+
y_x^1\, y_x^2 \, y_x^2 \cdot
\left\{
-2\, X_{y^1y^2}
\right\}+
y_x^2 \, y_x^2\, y_x^2\cdot
\left\{-X_{y^2y^2}\right\}, \\
\endaligned\right.
\end{equation}
Except for inductive inspiration ({\it see} the formulation of
Lemma~3.32 below), this observation will not be used further. At this
stage, it helps at least to maintain a strong intuitive control of the
correctness of the underlying combinatorics.

\subsection*{3.28.~System equivalent to the flat system}
By induction, we therefore guess that the analogy holds for general $m
\geqslant 2$, namely we guess the following combinatorics, which requires
some preliminaries.

As in the beginning of \S3.1, let $x \in \K$, let $y= (y^1, \dots,
y^m)\in \K^m$, let $(x,y) \mapsto (X(x,y), Y(x, y))$ be a local
$\K$-analytic transformation defined in a neighborhood of the origin
in $\K^{m+1}$ and assume that the system $y_{xx}^j= F^j(x, y, y_x)$, 
$j= 1, \dots, m$, is
equivalent to the flat system $Y_{X X}^j=0$, $j=1,\dots,m$. By
assumption, the Jacobian matrix of the equivalence equals the identity
matrix at the origin. Remind that we identify $x$ with $y^0$. For all
$k_1, l_1, l_2=0, \dots, m$, we define a modification
\def\theequation{3.29}\begin{equation}
\Delta( x \vert \dots \vert^{k_1} y^{l_1} y^{l_2} \vert
\dots \vert y^m)
\end{equation}
of the Jacobian determinant as follows. We replace the $k_1$-th
column of the determinant~\thetag{ 3.2}, which consists of first order
derivatives $\cdot_{y^{k_1}}$, by a column which consists of second
order derivatives $\cdot_{y^{l_1}y^{l_2}}$. In~\thetag{ 3.29}, the
notation $\vert^{k_1}$ designates the $k_1$-th column, the first one
being labelled by $k_1=0$ and the last one by $k_1 = m$. With this
notation at hand, we may define the square functions
\def\theequation{3.30}\begin{equation}
\square_{y^{l_1}y^{l_2}}^{k_1}:=
\frac{
\Delta( x \vert \dots 
\vert^{k_1} y^{l_1} y^{l_2} \vert
\dots \vert y^m)}{
\Delta( x \vert \dots 
\vert^{k_1} y^{k_1} \vert \dots
y^m)},
\end{equation}
which are rational expressions in the second order jet of the
transformation $(x,y)\mapsto (X(x, y), Y(x, y))$. As before, the
denominator is the Jacobian determinant of the change of coordinates.

Since, according to~\thetag{ 3.26}, 
the expression of ${\bf R}_2^j \ ({\rm mod} \ y_{xx}^l)$ is
\def\theequation{3.31}\begin{equation}
\left\{
\aligned
{\bf R}_2^j \ ({\rm mod} \ y_{xx}^l)=
Y_{xx}^j
& \
+
\sum_{l_1=1}^m\, 
y_x^{l_1}\cdot
\left[
2\, Y_{xy^{l_1}}^j- 
\delta_{l_1}^j\, X_{xx}
\right]+ \\
& \
+
\sum_{l_1=1}^m\, 
\sum_{l_2=1}^m\, 
y_x^{l_1} \, y_x^{l_2}\cdot
\left[
Y_{y^{l_1}y^{l_2}}^j- 
\delta_{l_1}^j\, X_{xy^{l_2}}-
\delta_{l_2}^j\, X_{xy^{l_1}}
\right]+ \\
& \
+
\sum_{l_1=1}^m\, \sum_{l_2=1}^m\, 
\sum_{l_3=1}^m\, 
y_x^{l_1} \, y_x^{l_2} \, y_x^{l_3}\cdot
\left[
-\delta_{l_1}^j\, X_{y^{l_2} y^{l_3}}
\right], \\
\endaligned\right.
\end{equation}
and since, in the cases $m=1$ and $m=2$, we have already observed
strong analogies between~\thetag{ 3.31} and the complete explicit
expression of the system $y_{xx}^j=F^j(x, y, y_x)$ equivalent to the
flat system $Y_{XX}^j = 0$, we guess that the following lemma 
is formally true.

\def\thelemma{3.32}\begin{lemma}
The system $y_{x x}^j= F^j(x, y, y_x)$, $j=1,\dots,m$, is equivalent
to the flat system $Y_{X X}^j=0$, $j=1,\dots,m$, if and only if there
exist local $\K$-analytic functions $X(x, y)$ and $Y^j(x, y)$,
$j=1,\dots,m$, such that it may be written under the specific form
\def\theequation{3.33}\begin{equation}
\left\{
\aligned
0 = 
y_{xx}^j+
\square_{xx}^j
& \
+
\sum_{l_1=1}^m\, 
y_x^{l_1}\cdot
\left[2\, \square_{xy^{l_1}}^j-
\delta_{l_1}^j\, \square_{xx}^0\right]+ \\
& \
+
\sum_{l_1=1}^m\, \sum_{l_2=1}^m\, 
y_x^{l_1}\, y_x^{l_2}\cdot
\left[
\square_{y^{l_1}y^{l_2}}^j-
\delta_{l_1}^j\, \square_{xy^{l_2}}^0 -
\delta_{l_2}^j\, \square_{xy^{l_1}}^0
\right]+ \\
& \
+
y_x^j \, 
\sum_{l_1=1}^m\, \sum_{l_2=1}^m\, 
y_x^{l_1}\, y_x^{l_2}\cdot
\left[
-\square_{
y^{l_1}y^{l_2}}^0
\right].
\endaligned\right.
\end{equation}
\end{lemma}

The complete proof of this lemma involves only linear algebra
considerations, although with rather massive terms. This makes it
rather lengthy. Consequently, we postpone it to the final Section~5
below.

\subsection*{3.34.~First auxiliary system}
Clearly, if we set
\def\theequation{3.35}\begin{equation}
\left\{
\aligned
G^j 
& \
:=
-\square_{xx}^j, \\
H_{l_1}^j 
& \
:=
-2\, \square_{xy^{l_1}}^j+
\delta_{l_1}^j\, \square_{xx}^0, \\
L_{l_1,l_2}^j
& \
:= 
-\square_{y^{l_1}y^{l_2}}^j+
\delta_{l_1}^j\, 
\square_{xy^{l_2}}^0+
\delta_{l_2}^j\, 
\square_{xy^{l_1}}^0, \\
M_{l_1, l_2}
& \
:=
\square_{y^{l_1}y^{l_2}}^0, 
\endaligned\right.
\end{equation}
we immediately see that the first condition of Theorem~1.7 holds
true. Moreover, we claim that there are $m+1$ more square functions
than functions $G^j$, $H_{l_1}^j$, $L_{l_1,l_2}^j$ and $M_{l_1,
l_2}$. Indeed, taking account of the symmetries, we enumerate:
\def\theequation{3.36}\begin{equation}
\left\{
\aligned
{}
& \
\# \{\square_{xx}^j\}
=
m,
& \
\# \{\square_{xx}^0 \}
=
1, \\
& \
\# \{\square_{xy^{l_1}}^j\}
=
m^2, 
& \
\# \{\square_{xy^{l_1}}^0 \}
= 
m, \\
& \
\# \{\square_{y^{l_1}y^{l_2}}^j\}
= 
\frac{ m^2(m+1)}{2}, 
\ \ \ \ \ \ \ \ \ \ \ 
\ \ \ \ \ \ \ \ \ \ \ 
& \
\# \{\square_{y^{l_1}y^{l_2}}^0 \}
= 
\frac{ m(m+1)}{2}, \\
\endaligned\right.
\end{equation}
whereas
\def\theequation{3.37}\begin{equation}
\left\{
\aligned
{}
& \
\# \{G^j\}
= 
m, 
& \ 
\#
\{H_{l_1}^j \}
=
m^2, \\
& \
\# \{L_{l_1,l_2}^j \}
=
\frac{ m^2(m+1)}{2}, 
& \ 
\ \ \ \ \ \ \ \ \ \ \ 
\ \ \ \ \ \ \ \ \ \ \
\# \{M_{l_1, l_2}\} 
=
\frac{ m(m+1)}{2}.
\endaligned\right.
\end{equation}
Similarly as in Section~2, for $j,l_1,l_2=0,1,\dots,m$, let us
introduce functions $\Pi_{l_1,l_2}^j$ of $(x,y^1, \dots, y^m)$,
symmetric with respect to the lower indices, and let us seek necessary
and sufficient conditions in order that there exist solutions $(X,Y)$
to the {\sl first auxiliary system} defined precisely by:
\def\theequation{3.38}\begin{equation}
\left\{
\aligned
\square_{xx}^0 = \Pi_{0,0}^0, \ \ \ \ \ & 
\square_{xy^{l_1}}^0 = \Pi_{0,l_1}^0, \ \ \ \ \ & 
\square_{y^{l_1}y^{l_2}}^0 = \Pi_{l_1,l_2}^0, \ \ \ \ \ \\
\square_{xx}^j = \Pi_{0,0}^j, \ \ \ \ \ & 
\square_{xy^{l_1}}^j = \Pi_{0,l_1}^j, \ \ \ \ \ & 
\square_{y^{l_1}y^{l_2}}^j = \Pi_{l_1,l_2}^j. \ \ \ \ \ \\ 
\endaligned\right.
\end{equation}

\subsection*{3.39.~Compatibility conditions 
for the first auxiliary system} As in Section~2, the compatibility
conditions for this system will simply be obtained by computing the
cross differentiations. The following statement generalizes
Lemma~2.31 and also provides a proof of
it, in the case $m= 1$.

\def\thelemma{3.40}\begin{lemma} 
For all $j, l_1, l_2, l_3= 0,1, \dots, m$, we have the cross
differentiation relations
\def\theequation{3.41}\begin{equation}
\left(
\square_{y^{l_1}y^{l_2}}^j
\right)_{y^{l_3}}-
\left(
\square_{y^{l_1}y^{l_3}}^j
\right)_{y^{l_2}}= -
\sum_{k=0}^m \,
\square_{y^{l_1}y^{l_2}}^k\cdot \,
\square_{y^{l_3}y^k}^j+
\sum_{k=0}^m \,
\square_{y^{l_1}y^{l_3}}^k\cdot \,
\square_{y^{l_2}y^k}^j.
\end{equation}
\end{lemma}

\proof
To begin with, as a preliminary, let us generalize the Pl\"ucker
identity~\thetag{ 2.28}. Let $C_1, C_2, \dots, C_m, D, E$ be $(m+2)$
column vectors in $\K^m$ and introduce the following notation for the
$m\times (m+2)$ matrix consisting of these vectors:
\def\theequation{3.42}\begin{equation}
\left[
C_1 \vert C_2 \vert \cdots \vert
C_m \vert D \vert E 
\right].
\end{equation}
Extracting columns
from this matrix, we shall construct $m\times m$ determinants which
are modification of the following 
``fundamental'' determinant 
\def\theequation{3.43}\begin{equation}
\left\vert \!
\left\vert
C_1 \vert \cdots \vert
C_m
\right\vert \!
\right\vert\equiv
\left\vert \!
\left\vert
C_1 \vert \cdots \vert^{j_1}C_{j_1} \vert
\cdots \vert^{j_2} C_{j_2} \vert \cdots \vert
C_m
\right\vert \!
\right\vert.
\end{equation}
Here and in the sequel, we use a double vertical line in the beginning
and in the end to denote a determinant. Also, we emphasize two
distinct columns, the $j_1$-th and the $j_2$-th, where $j_2 > j_1$,
since we will modify them. For instance in this matrix, let us replace
these two columns by the column $D$ and by the column $E$, which
yields the determinant
\def\theequation{3.44}\begin{equation}
\left\vert \!
\left\vert
C_1 \vert \cdots \vert^{j_1} D \vert
\cdots \vert^{j_2} E \vert \cdots \vert
C_m
\right\vert \!
\right\vert.
\end{equation}
In this notation, one should understand that {\it only}\, the $j_1$-th
and the $j_2$-th columns are distinct from the columns of the
fundamental $m\times m$ determinant~\thetag{ 3.43}. With this notation
at hand, we can now formulate and prove a preliminary lemma that will
be useful later.

\def\thelemma{3.45}\begin{lemma}
The following quadratic identity between determinants holds
true{\rm :}
\def\theequation{3.46}\begin{equation}
\left\{
\aligned
{} 
& \
\left\vert \! 
\left\vert
C_1 \vert \cdots \vert^{j_1} D 
\vert \cdots \vert^{j_2}
E \vert \cdots \vert C_n
\right\vert \! 
\right\vert
\cdot
\left\vert \! 
\left\vert
C_1 \vert \cdots \vert^{j_1} C_{j_1}
\vert \cdots \vert^{j_2}
C_{j_2} \vert \cdots \vert C_n
\right\vert \! 
\right\vert = \\
& \
=
\left\vert \! 
\left\vert
C_1 \vert \cdots \vert^{j_1} D 
\vert \cdots \vert^{j_2}
C_{j_2} \vert \cdots \vert C_n
\right\vert \! 
\right\vert
\cdot
\left\vert \! 
\left\vert
C_1 \vert \cdots \vert^{j_1} C_{j_1}
\vert \cdots \vert^{j_2}
E \vert \cdots \vert C_n
\right\vert \! 
\right\vert- \\
& \
-
\left\vert \! 
\left\vert
C_1 \vert \cdots \vert^{j_1} E 
\vert \cdots \vert^{j_2}
C^{j_2} \vert \cdots \vert C_n
\right\vert \! 
\right\vert
\cdot
\left\vert \! 
\left\vert
C_1 \vert \cdots \vert^{j_1} C_{j_1}
\vert \cdots \vert^{j_2}
D \vert \cdots \vert C_n
\right\vert \! 
\right\vert.
\endaligned\right.
\end{equation}
\end{lemma}

\proof
After some permutations of columns, this identity 
amounts to
\def\theequation{3.47}\begin{equation}
\left\{
\aligned
{}
& \
\left\vert \! 
\left\vert
C_1 \vert \cdots \vert C_{m-2} \vert D \vert E 
\right\vert \! 
\right\vert 
\cdot
\left\vert \! 
\left\vert
C_1 \vert \cdots \vert C_{m-2} \vert C_{m-1} \vert C_m
\right\vert \! 
\right\vert = \\
& \
=
\left\vert \! 
\left\vert
C_1 \vert \cdots \vert C_{m-2} \vert D \vert C_m
\right\vert \! 
\right\vert 
\cdot
\left\vert \! 
\left\vert
C_1 \vert \cdots \vert C_{m-2} \vert C_{m-1} \vert E
\right\vert \! 
\right\vert - \\
& \
-
\left\vert \! 
\left\vert
C_1 \vert \cdots \vert C_{m-2} \vert E \vert C_m
\right\vert \! 
\right\vert 
\cdot
\left\vert \! 
\left\vert
C_1 \vert \cdots \vert C_{m-2} \vert C_{m-1} \vert D
\right\vert \! 
\right\vert.
\endaligned\right.
\end{equation} 
To establish this identity, we introduce some notation. If $A$ and $B$
are vertical vectors in $\K^m$ and if $i_1,i_2=1,\dots,m$ with $i_1 <
i_2$, we denote
\def\theequation{3.48}\begin{equation}
\Delta_{i_1,i_2}^2(A\vert B):= 
\left\vert
\begin{array}{cc}
A_{i_1} & B_{i_1} \\
A_{i_2} & B_{i_2}
\end{array}
\right\vert.
\end{equation}
If $\left\vert \! \left \vert A_1 \vert A_2 \vert A_3 \vert \cdots
\vert A_m \right\vert \! \right \vert$ is a $m \times m$ determinant,
and if $i_1, i_2=1, \dots, m$ with $i_1 < i_2$, we denote by $M_{i_1,
i_2}^{m-2} (A_3 \vert \cdots \vert A_m)$ the $(m-2) \times ( m-2)$
determinant obtained from the matrix $\left[ A_3 \vert \cdots \vert
A_m \right]$ by erasing the $i_1$-th line and the $i_2$-th line. 
Without proof, we recall an elementary classical formula
\def\theequation{3.49}\begin{equation}
\left\{
\aligned
{}
& \
\left\vert \! 
\left\vert
A_1 \vert A_2 \vert A_3 \vert \cdots
\vert A_m
\right\vert \! 
\right\vert = \\
& \ \ \ \ \ \
=
\sum_{1\leqslant i_1 < i_2 \leqslant m}\, 
(-1)^{i_1+i_2-1} \, 
\Delta_{i_1, i_2}^2(A_1 \vert A_2) \cdot
M_{i_1, i_2}^{m-2} (A_3 \vert
\cdots \vert A_m),
\endaligned\right.
\end{equation}
which may be established by developing the determinant $\left\vert \!
\left\vert A_1 \vert A_2 \vert A_3 \vert \cdots \vert A_m \right\vert
\! \right\vert$ with respect to its first column, and then
re-developing all the obtained $(m-1)\times (m-1)$ determinants with
respect to their first columns. To establish~\thetag{3.47}, we start
with an equivalent version of the identity~\thetag{ 2.28}:
\def\theequation{3.50}\begin{equation}
\left\{
\aligned
\Delta_{i_1,i_2}^2(D\vert E) \cdot
\Delta_{i_3, i_4}^2(C_1\vert C_2)=
& \
\Delta_{i_1, i_2}^2(D\vert C_2) \cdot
\Delta_{i_3, i_4}^2(C_1\vert E) - \\
& \
-
\Delta_{i_1, i_2}^2(E\vert C_2) \cdot
\Delta_{i_3, i_4}^2(C_1\vert D),
\endaligned\right.
\end{equation}
where $1 \leqslant i_1 < i_2 \leqslant m$ and $1 \leqslant i_3 < i_4 \leqslant m$.
Multiplying by $(-1)^{ i_1 + i_2+ i_3+ i_4 -2}$, multiplying by
$M_{i_1,i_2}^{m-2}(C_3 \vert \cdots \vert C_m)$, and multiplying by
$M_{i_3,i_4}^{m-2}(C_3 \vert \cdots \vert C_m)$ applying the double
summation $\sum_{ 1\leqslant i_1 < i_2 \leqslant m} \, \sum_{ 1\leqslant i_3 < i_4
\leqslant m}$, we get
\def\theequation{3.51}\begin{equation}
\left\{
\aligned
{}
& \
\sum_{1\leqslant i_1 < i_2 \leqslant m}\, 
\sum_{1 \leqslant i_3 < i_4 \leqslant m} \, 
(-1)^{i_1+i_2+i_3+i_4-2}\, 
\Delta_{i_1,i_2}^2(D\vert E) \cdot
\Delta_{i_3, i_4}^2(C_1\vert C_2)\cdot \\
& \
\ \ \ \ \ \ \ \ \ \ \ \ \ 
\cdot
M_{i_1,i_2}^{m-2}(C_3 \vert \cdots \vert C_m) \cdot
M_{i_3, i_4}^{m-2}(C_3 \vert \cdots \vert C_m) = \\
& \
=
\sum_{1\leqslant i_1 < i_2 \leqslant m}\, 
\sum_{1 \leqslant i_3 < i_4 \leqslant m} \,
(-1)^{i_1+ i_2 -1} \, (-1)^{ i_3+ i_4 -1} \, 
\left[
\Delta_{i_1, i_2}^2(D\vert C_2) \cdot
\Delta_{i_3, i_4}^2(C_1\vert E) - \right. \\
& \
\ \ \ \ \ 
\left.
-\Delta_{i_1, i_2}^2(E\vert C_2) \cdot
\Delta_{i_3, i_4}^2(C_1\vert D)\right] \cdot
M_{i_1,i_2}^{m-2}(C_3 \vert \cdots \vert C_m) \cdot
M_{i_3, i_4}^{m-2}(C_3 \vert \cdots \vert C_m).
\endaligned\right.
\end{equation}
Thanks to the relation~\thetag{ 3.49}, this last
identity coincides exactly with the desired
identity~\thetag{ 3.47}. The proof
is complete. 
\endproof

We can now establish Lemma~3.40. As a preliminary observation,
by the Leibniz rule for the differentiation of a determinant, 
we must differentiate every column:
\def\theequation{3.52}\begin{equation}
\left\{
\aligned
\left[
\Delta\left(
y^{l_1}y^{k_1} \vert \cdots
\vert y^{l_m}y^{k_m}
\right)\right]_{y^j} =
& \
\Delta\left(
y^j y^{l_1}y^{k_1} \vert \cdots
\vert y^{l_m}y^{k_m}
\right)+ \cdots + \\
& \
+ 
\Delta\left(
y^{l_1}y^{k_1} \vert \cdots
\vert y^j y^{l_m}y^{k_m}
\right).
\endaligned\right.
\end{equation}
Using also the rule for the differentiation of a quotient, we may
endeavour to compute the cross differentiations $\left(
\square_{y^{l_1}y^{l_2}}^j \right)_{y^{l_3}} - \left(
\square_{y^{l_1}y^{l_3}}^j \right)_{y^{l_2}}$ of the left-hand side
of~\thetag{ 3.41}. This will generalize~\thetag{ 2.25}. Sometimes in
the computation, we shall abbreviate the Jacobian determinant
$\Delta(y^0\vert \cdots \vert y^m)$ using the shorter notation
$\Delta$; as before, a product between two elements of $\K$ will often
be denoted by the sign ``$\cdot$'', for clarity.
Here is the computation:
\def\theequation{3.53}\begin{equation}
\left\{
\aligned
{}
& \
\left(
\square_{y^{l_1}y^{l_2}}^j 
\right)_{y^{l_3}}
-
\left(
\square_{y^{l_1}y^{l_3}}^j 
\right)_{y^{l_2}} = \\
& \
=
\frac{ \partial }{\partial y^{l_3}} \left(
\frac{ \Delta\left(
y^0 \vert \cdots \vert^j y^{ l_1} y^{l_2}\vert
\cdots \vert y^m
\right)}{
\Delta\left(
y^0 \vert \cdots \vert y^m
\right)}
\right)
-
\frac{ \partial }{\partial y^{l_2}} \left(
\frac{ \Delta\left(
y^0 \vert \cdots \vert^j y^{ l_1} y^{l_3}\vert
\cdots \vert y^m
\right)}{
\Delta\left(
y^0 \vert \cdots \vert y^m
\right)}
\right) \\
& \
=
\frac{ 1}{[\Delta]^2}
\left[
\aligned
{}
& \
\Delta\left(
y^0 y^{l_3}\vert \cdots \vert^j y^{l_1}y^{l_2} \vert
\cdots \vert y^m
\right)\cdot \Delta + \cdots + \\
& \
+
\underline{
\Delta\left(y^0\vert \cdots
\vert^j y^{l_3}y^{l_1}y^{l_2}\vert
\cdots \vert y^m\right) \cdot \Delta}_{
\octagon \! \! \! \! \! \tiny{\sf a}}
+ 
\cdots + \\
& \
+
\Delta\left(
y^0\vert \cdots \vert^jy^{l_1} y^{l_2} \vert
\cdots \vert y^{l_3} y^m
\right)\cdot \Delta - \\
& \
-
\Delta\left(
y^0 \vert \cdots \vert^j y^{l_1}y^{l_2} \vert
\cdots \vert y^m\right) \cdot
\left[
\Delta\left(
y^0 y^{l_3} \vert \cdots \vert y^m\right)+\cdots+\right. \\
& \ \ \ \ \ \ \ \ \ \ \ \ \ \ \ \ 
\ \ \ \ \ \ \ \ \ \ \ \ \ \ \ \ 
\ \ \ \ \ \ \ \ \ \ \ \ \ \ \ \ 
\ \ \ \ \ \ \ \ \ \ \ \ \ \ \ \ 
\left.
+
\Delta\left(
y^0 \vert \cdots \vert y^{l_3} y^m
\right)
\right]
\endaligned
\right]- \\
& \ \ \ \ \ \ \
-
\frac{ 1}{[\Delta]^2}
\left[
\aligned
{}
& \
\Delta\left(
y^0 y^{l_2}\vert \cdots \vert^j y^{l_1}y^{l_3} \vert
\cdots \vert y^m
\right)\cdot \Delta + \cdots + \\
& \ 
+
\underline{
\Delta\left(y^0\vert \cdots
\vert^j y^{l_2}y^{l_1}y^{l_3}\vert
\cdots \vert y^m\right) \cdot \Delta}_{
\octagon \! \! \! \! \! \tiny{\sf a}}
+ \cdots + \\
& \
+
\Delta\left(
y^0\vert \cdots \vert^jy^{l_1} y^{l_3} \vert
\cdots \vert y^{l_2} y^m
\right)\cdot \Delta - \\
& \
-
\Delta\left(
y^0 \vert \cdots \vert^j y^{l_1}y^{l_3} \vert
\cdots \vert y^m\right) \cdot
\left[
\Delta\left(
y^0 y^{l_2} \vert \cdots \vert y^m\right)+\cdots+ \right. \\
& \ \ \ \ \ \ \ \ \ \ \ \ \ \ \ \ 
\ \ \ \ \ \ \ \ \ \ \ \ \ \ \ \ 
\ \ \ \ \ \ \ \ \ \ \ \ \ \ \ \ 
\ \ \ \ \ \ \ \ \ \ \ \ \ \ \ \ 
\left.
+
\Delta\left(
y^0 \vert \cdots \vert y^{l_2} y^m
\right)
\right]
\endaligned
\right].
\endaligned\right.
\end{equation}
Crucially, we observe that all the determinants involving a third
order derivative upon one of their columns kill each other and
disappear: we have underlined them with 
$\octagon \! \! \! \! \tiny{\sf a}$
appended. However, it still remains plenty of determinants involving a
second order derivative upon two different columns. We must transform
all of them and express them in terms of determinants involving a second
order derivative upon only one column. To this aim, as an application
of our preliminary Lemma~3.45, we have the following relations, valid
for $j_1, j_2, l_1, l_2, l_3, l_4 =0, \dots, m$ and $j_1 < j_2$:
\def\theequation{3.54}\begin{equation}
\left\{
\aligned
{}
& \
\Delta\left(
y^0 \vert \cdots \vert^{j_1} y^{l_1} y^{l_2} \vert \cdots
\vert^{j_2} y^{l_3}y^{l_4} \vert \cdots \vert y^m
\right)\cdot 
\Delta\left(
y^0 \vert \cdots \vert^{j_1} y^{j_1} \vert \cdots \vert^{j_2}
y^{j_2} \vert \cdots \vert y^m
\right) = \\
& \
= 
\Delta\left(
y^0 \vert \cdots \vert^{j_1} y^{l_1} y^{l_2} \vert \cdots
\vert^{j_2} y^{j_2} \vert \cdots \vert y^m
\right)\cdot 
\Delta\left(
y^0 \vert \cdots \vert^{j_1} y^{j_1} \vert \cdots \vert^{j_2}
y^{l_3} y^{l_4} \vert \cdots \vert y^m
\right) \\
& \
-
\Delta\left(
y^0 \vert \cdots \vert^{j_1} y^{l_3} y^{l_4} \vert \cdots
\vert^{j_2} y^{j_2} \vert \cdots \vert y^m
\right)\cdot 
\Delta\left(
y^0 \vert \cdots \vert^{j_1} y^{j_1} \vert \cdots \vert^{j_2}
y^{l_1} y^{l_2} \vert \cdots \vert y^m
\right).
\endaligned\right.
\end{equation}
With these formulas, we may transform the lines number 3, 4, 5 and
8, 9, 10 of~\thetag{ 3.53}. Also, we observe that the
lines 6, 7 and 11, 12 of~\thetag{ 3.53} involve determinants
having a single second order derivative. Taking account 
of the $\frac{ 1}{[\Delta]^2}$ factor, we deduce that the
lines 6, 7 and 11, 12 of~\thetag{ 3.53} may already be
expressed as sums of square functions. 
Achieving all these transformations, we may 
rewrite~\thetag{ 3.53} as follows
\def\theequation{3.55}\begin{equation}
\left\{
\aligned
{}
& \
\left(
\square_{y^{l_1}y^{l_2}}^j 
\right)_{y^{l_3}}
-
\left(
\square_{y^{l_1}y^{l_3}}^j 
\right)_{y^{l_2}} = \\
& \
=
\frac{ 1}{[\Delta]^2}
\left[
\aligned
{}
& \
\Delta\left(
y^0 y^{l_3} \vert \cdots \vert^j y^j \vert \cdots \vert y^m
\right)\cdot \Delta
\left(
y^0 \vert \cdots \vert^j 
y^{l_1}y^{l_2} \vert \cdots \vert y^m
\right)- \\
& \
-
\Delta\left(
y^{l_1} y^{l_2} \vert \cdots \vert^j y^j \vert \cdots \vert y^m
\right)\cdot \Delta
\left(
y^0 \vert \cdots \vert^j 
y^0y^{l_3} \vert \cdots \vert y^m
\right) + \\
& \
+\cdots+ \\
& \
+
\Delta\left(
y^0 \vert \cdots \vert^j y^{l_1} y^{l_2} \vert
\cdots \vert y^m
\right)\cdot
\Delta\left(
y^0 \vert \cdots \vert^j y^j \vert
\cdots \vert y^{l_3}y^m
\right) - \\
& \
-
\Delta\left(
y^0 \vert \cdots \vert^j y^{l_3} y^m \vert
\cdots \vert y^m
\right)\cdot
\Delta\left(
y^0 \vert \cdots \vert^j y^j \vert
\cdots \vert y^{l_1}y^{l_2}
\right)
\endaligned
\right] - \\
& \
\ \ \ \ \ \ 
-
\sum_{k=0}^m\, \square_{y^{l_1}y^{l_2}}^j\, 
\square_{y^{l_3} y^k}^k- \\
& \ \ \ \ \ \
-
\frac{ 1}{[\Delta]^2}
\left[
\aligned
{}
& \
\Delta\left(
y^0 y^{l_2} \vert \cdots \vert^j y^j \vert \cdots \vert y^m
\right)\cdot \Delta
\left(
y^0 \vert \cdots \vert^j 
y^{l_1}y^{l_3} \vert \cdots \vert y^m
\right)- \\
& \
-
\Delta\left(
y^{l_1} y^{l_3} \vert \cdots \vert^j y^j \vert \cdots \vert y^m
\right)\cdot \Delta
\left(
y^0 \vert \cdots \vert^j 
y^0y^{l_2} \vert \cdots \vert y^m
\right) + \\
& \
+\cdots+ \\
& \
+
\Delta\left(
y^0 \vert \cdots \vert^j y^{l_1} y^{l_3} \vert
\cdots \vert y^m
\right)\cdot
\Delta\left(
y^0 \vert \cdots \vert^j y^j \vert
\cdots \vert y^{l_2}y^m
\right) - \\
& \
-
\Delta\left(
y^0 \vert \cdots \vert^j y^{l_2} y^m \vert
\cdots \vert y^m
\right)\cdot
\Delta\left(
y^0 \vert \cdots \vert^j y^j \vert
\cdots \vert y^{l_1}y^{l_3}
\right)
\endaligned
\right] + \\
& \
\ \ \ \ \ \
+
\sum_{k=0}^m\, \square_{y^{l_1}y^{l_3}}^j\, 
\square_{y^{l_2} y^k}^k.
\endaligned\right.
\end{equation} 
Notice that two pairs of ``{\tt cdots}'' terms $+\cdots+$ appearing in
the lines 3, 4 and 8, 9 of~\thetag{ 3.53} are replaced by a single
``{\tt cdots}'' term $+\cdots+$ in the lines 4 and 9 of~\thetag{
3.55}. Importantly, we point that in the ``middle'' of the two ``{\tt
cdots}'' terms $+\cdots+$ appearing in the lines 4 and 10 of~\thetag{
3.55} just above, there are two terms which do not occur: they simply
correspond to the two underlined terms having
$\octagon \! \! \! \! \tiny{\sf a}$
appended appearing in the lines 4 and 9~\thetag{ 3.53}.

Now, taking account of the factor $\frac{ 1}{[\Delta]^2}$, we can
re-express all the terms of~\thetag{ 3.55} as sums of square
functions:
\def\theequation{3.56}\begin{equation}
\left\{
\aligned
{}
& \
\left(
\square_{y^{l_1}y^{l_2}}^j 
\right)_{y^{l_3}}
-
\left(
\square_{y^{l_1}y^{l_3}}^j 
\right)_{y^{l_2}} = \\
& \ \ \ \ \ \ \
=\sum_{k=0; \, k\neq j}^m\, 
\square_{y^{l_3}y^k}^k \, 
\square_{y^{l_1}y^{l_2}}^j- 
\sum_{k=0; \, k\neq j}^m\, 
\square_{y^{l_1}y^{l_2}}^k\, \square_{y^{l_3}y^k}^j- \\
& \ \ \ \ \ \ \
-\sum_{k=0}^m\, \square_{y^{l_1}y^{l_2}}^j\, 
\square_{y^{l_3}y^k}^k- \\
& \ \ \ \ \ \ \
-\sum_{k=0; \, k\neq j}^m\, 
\square_{y^{l_2}y^k}^k \, 
\square_{y^{l_1}y^{l_3}}^j+
\sum_{k=0; \, k\neq j}^m\, 
\square_{y^{l_1}y^{l_3}}^k\, \square_{y^{l_2}y^k}^j+ \\
& \ \ \ \ \ \ \
+
\sum_{k=0}^m\, 
\square_{y^{l_1}y^{l_3}}^j\, 
\square_{y^{l_2}y^k}^k.
\endaligned\right.
\end{equation} 
Finally, we observe that in the two pairs of sums having $k\neq j$
appearing in the lines 2 and 4 just above, we can include the term
$k=j$ in each pair, because these two terms are immediately killed
inside the corresponding pair. In conclusion, after a final obvious
killing of four (among six) complete sums in this modification
of~\thetag{ 3.56}, we obtain the desired formula~\thetag{ 3.41}, with
two sums. This completes the proof of Lemma~3.40 and also at the same
occasion, the proof of Lemma~2.31.
\endproof

\subsection*{ 3.57.~Compatibility conditions for the
first auxiliary system} According to the (approximate)
identities~\thetag{ 3.3}, taking account of
the explicit definitions~\thetag{ 3.30} of
the square functions, we have 
\def\theequation{3.58}\begin{equation}
\left\{
\aligned
\square_{ xx}^0 \cong X_{xx}, \ \ \ \ \ & 
\square_{xy^{ l_1}}^0 \cong X_{ xy^{l_1}}, \ \ \ \ \ & 
\square_{ y^{l_1}y^{ l_2}}^0 \cong X_{ y^{ l_1} y^{ l_2}}, 
\ \ \ \ \ \\
\square_{ xx}^j \cong Y_{xx}^j, \ \ \ \ \ & 
\square_{xy^{ l_1}}^j \cong Y_{ xy^{l_1}}^j, \ \ \ \ \ & 
\square_{ y^{l_1}y^{ l_2}}^j \cong Y_{ y^{ l_1} y^{ l_2}}^j.
\ \ \ \ \ \\
\endaligned\right.
\end{equation}
Consequently, the first auxiliary system~\thetag{ 3.38} looks
approximatively like a complete second order system of partial
differential equations in the $(m+1)$ independent variables $(x, y)$
and in the $(m+1)$ dependent variables $(X, Y)$. By means of
elementary algebraic operations, one may transform this system in a
true second order {\it complete}\, system, solved with respect to the
top order derivatives, namely of the form
\def\theequation{3.59}\begin{equation}
\left\{
\aligned
X_{xx} = 
\Lambda_{0,0}^0, \ \ \ \ \ & 
X_{xy^{ l_1}} =
\Lambda_{0,l_1}^0, \ \ \ \ \ & 
X_{y^{ l_1}y^{ l_2}} = 
\Lambda_{l_1,l_2}^0, \ \ \ \ \ \\
Y_{xx}^j = 
\Lambda_{0,0}^j, \ \ \ \ \ & 
Y_{xy^{ l_1}}^j = 
\Lambda_{0,l_1}^j, \ \ \ \ \ & 
Y_{y^{ l_1}y^{ l_2}} = 
\Lambda_{l_1,l_2}^j, \ \ \ \ \ \\ 
\endaligned\right.
\end{equation}
where the $\Lambda_{ j_1, j_2}^{ k_1}$ are local $\K$-analytic
functions of $(x, y^{l_1}, X, Y^j, X_x, X_{y^{ l_1}}, Y_x^j, Y_{ y^{
l_1}}^j)$. For such a system, the compatibility conditions [which are
necessary and sufficient for the existence of a solution $(X,Y)$]
follow by obvious cross differentiation. Coming back to the system{
3.38}, these compatibility conditions amount to the quadratic-like
compatibility conditions expressed in Lemma~3.40. In conclusion, we
have proved the following intermediate statement.

\def\theproposition{3.60}\begin{proposition} 
There exist functions $X$, $Y^j$ solving the first auxiliary
system~\thetag{ 3.38} of nonlinear second order partial differential
equations {\rm if and only if} the right-hand side functions $\Pi_{
l_1,\, l_2}^j(x,\, y)$ satisfy the quadratic compatibility conditions
\def\theequation{3.61}\begin{equation}
\frac{\partial \Pi_{ l_1, l_2}^j}{\partial y^{ l_3}} - 
\frac{\partial \Pi_{ l_1, l_3}^j}{\partial y^{ l_2}} = 
-\sum_{ k= 0}^m \, 
\Pi_{ l_1, l_2}^k \cdot
\Pi_{ l_3, k}^j + 
\sum_{ k = 0}^m \, 
\Pi_{ l_1, l_3}^k \cdot 
\Pi_{ l_2, k}^j,
\end{equation}
for $j, l_1, l_2, l_3 = 0, 1, \dots, m$.
\end{proposition}

\subsection*{ 3.62.~Principal unknowns}
As there are $(m+1)$ more square (or Pi) functions than the functions
$G^j$, $H_{l_1}^j$, $L_{ l_1, l_2}^j$ and $M_{ l_1, l_2}$ defined
by~\thetag{ 3.35}, we cannot invert directly the linear
system~\thetag{ 3.35} (which is of maximal rank). Hence
we must
choose $(m+1)$ specific square functions, calling them {\sl principal
unknowns}, and similarly as in \S2.34, the
best choice is to choose
$\square_{ xx}^0$ and $\square_{ xx}^j$, for
$j= 1,\dots, m$.
For clarity, it will be useful to 
adopt the notational equivalences
\def\theequation{3.63}\begin{equation}
\Theta^0 \equiv \Pi_{ 0, 0}^0 
\ \ \ \ \ \ \
{\rm and}
\ \ \ \ \ \ \
\Theta^j \equiv \Pi_{ j, j}^j.
\end{equation}
Then we may quasi-inverse the system~\thetag{ 3.35}, 
which yields~:
\def\theequation{3.64}\begin{equation}
\left\{
\aligned
\Pi_{ 0, 0}^j 
&
= 
\square_{ xx}^j 
= 
- G^j, \\
\Pi_{ 0, l_1}^j 
& 
=
\square_{ xy^{ l_1}}^j
= 
-\frac{ 1}{ 2} \,
H_{ l_1}^j + 
\frac{ 1}{ 2} \, \delta_{ l_1}^j \, 
\Theta^0, \\
\Pi_{ l_1, l_2}^j 
&
=
\square_{ y^{ l_1} y^{ l_2}}^j 
=
-L_{ l_1 , l_2}^j +
\frac{ 1}{ 2} \, 
\delta_{ l_1}^j \, 
L_{ l_2, l_2}^{ l_2} + 
\frac{ 1}{ 2} \, \delta_{ l_2}^j \, 
L_{ l_1,l_1}^{ l_1} + 
\frac{ 1}{ 2} \, \delta_{ l_1}^j \, 
\Theta^{ l_2}
+
\frac{ 1}{ 2} \, \delta_{ l_2}^j \,
\Theta^{ l_1}, \\
\Pi_{ 0, l_1}^0 
&
=
\square_{ xy^{ l_1}}^0 
=
\frac{ 1}{ 2} \, L_{ l_1, l_1}^{ l_1} + 
\frac{ 1}{ 2} \, \Theta^{ l_1}, \\
\Pi_{ l_1, l_2}^0 
&
=
\square_{ y^{ l_1} y^{ l_2}}^0 
= 
M_{ l_1,l_2}.
\endaligned\right.
\end{equation}
Before replacing these new expressions of the functions $\Pi_{ 0,
0}^j$, $\Pi_{ 0, l_1}^j$, $\Pi_{ l_1, l_2}^j$, $\Pi_{ l_1, l_2}^0$ and
$\Pi_{ 0, l_1}^{ l_1}$ into the compatibility conditions~\thetag{
3.61}, it is necessary to expound first~\thetag{ 3.61}, taking account
of the original splitting of the indices in the two sets $\{ 0\}$ and
$\{ 1, 2, \dots, m\}$. This yields six families of compatibility
conditions, totally equivalent
to the compact identities~\thetag{ 3.61}:
\def\theequation{3.65}\begin{equation}
\small
\left\{
\aligned
\left( \Pi_{ 0, 0}^j \right)_{ y^{ l_1}}
-
\left( \Pi_{ 0, l_1}^j \right)_x
&
= 
- 
\Pi_{ 0, 0}^0 \, \Pi_{ l_1, 0}^j 
- 
\sum_{ k= 1}^m \, \Pi_{ 0, 0}^k \, \Pi_{ l_1, k}^j
+
\Pi_{0, l_1}^0 \, \Pi_{ 0, 0}^j 
+
\sum_{ k= 1}^m \, \Pi_{ 0, l_1}^k \, \Pi_{ 0, k}^j, \\
\left( \Pi_{ l_1, l_2}^j \right)_x
-
\left( \Pi_{ l_1, 0}^j \right)_{ y^{ l_2}}
&
=
- 
\Pi_{ l_1, l_2}^0 \, \Pi_{ 0, 0 }^j 
- 
\sum_{ k= 1}^m \, \Pi_{ l_1, l_2}^k \, \Pi_{ 0, k}^j 
+ 
\Pi_{ l_1, 0}^0 \, \Pi_{ l_2, 0}^j 
+
\sum_{ k= 1}^m \, \Pi_{ l_1, 0}^k \, \Pi_{ l_2, k}^j, \\
\left(\Pi_{ l_1, l_2}^j\right)_{ y^{ l_3}}
-
\left(\Pi_{ l_1, l_3}^j\right)_{ y^{ l_2}}
&
=
- 
\Pi_{ l_1, l_2}^0 \, \Pi_{ l_3, 0 }^j 
- 
\sum_{ k= 1}^m \, \Pi_{ l_1, l_2}^k \, \Pi_{ l_3, k}^j
+ 
\Pi_{ l_1, l_3}^0 \, \Pi_{ l_2, 0}^j 
+
\sum_{ k= 1}^m \, \Pi_{ l_1, l_3}^k \, \Pi_{ l_2, k}^j, \\
\left(\Pi_{ 0, 0}^0 \right)_{ y^{ l_1}}
-
\left(\Pi_{ 0, l_1}^0 \right)_x
&
= 
-
\underline{
\Pi_{ 0, 0}^0 \, \Pi_{ l_1, 0}^0 }_{ 
\octagon \! \! \! \! \! \tiny{\sf a}} 
- 
\sum_{ k= 1}^m \, \Pi_{ 0, 0}^k \, \Pi_{ l_1, k}^0 
+
\underline{
\Pi_{ 0, l_1}^0 \, \Pi_{ 0, 0}^0 }_{ 
\octagon \! \! \! \! \! \tiny{\sf a}} 
+
\sum_{ k= 1}^m \, \Pi_{ 0, l_1}^k \, \Pi_{ 0, k}^0, \\
\left(\Pi_{ l_1, l_2}^0 \right)_x
-
\left(\Pi_{ l_1, 0}^0 \right)_{ y^{ l_2}}
&
=
- 
\Pi_{ l_1, l_2}^0 \, \Pi_{ 0, 0 }^0 - 
\sum_{ k= 1}^m \, 
\Pi_{ l_1, l_2}^k \, \Pi_{ 0, k}^0 
+ 
\Pi_{ l_1, 0}^0 \, \Pi_{ l_2, 0}^0 
+
\sum_{ k= 1}^m \, \Pi_{ l_1, 0}^k \, \Pi_{ l_2, k}^0, \\
\left(\Pi_{ l_1, l_2}^0 \right)_{ y^{ l_3}}
-
\left(
\Pi_{ l_1, l_3}^0 
\right)_{ y^{ l_2}}
&
=
- 
\Pi_{ l_1, l_2}^0 \, \Pi_{ l_3, 0 }^0 - 
\sum_{ k= 1}^m \, \Pi_{ l_1, l_2}^k \, \Pi_{ l_3, k}^0 
+ 
\Pi_{ l_1, l_3}^0 \, \Pi_{ l_2, 0}^0 
+
\sum_{ k= 1}^m \, \Pi_{ l_1, l_3}^k \, \Pi_{ l_2, k}^0.
\endaligned\right.
\end{equation}

\subsection*{ 3.66.~Convention about sums}
Up to the end of Section~4, we shall abbreviate any sum $\sum_{ k=
1}^m$ or $\sum_{ p = 1}^m$ as $\sum_k$ or $\sum_p$. Such sums will
appear very frequently. For all other sums, we shall precisely write
down the domain of variation of the summation index.

\subsection*{ 3.67.~Continuation} 
Thus, we have to replace~\thetag{ 3.64} in the six identities~\thetag{
3.65}. Firstly, let us expose all the intermediate steps in dealing
with the first identity $(3.65)_1$. Replacing plainly~\thetag{ 3.64}
in~$(3.65)_1$, we get:
\def\theequation{3.68}\begin{equation}
\small
\left\{
\aligned
& 
\left(
\Pi_{0, 0}^j
\right)_{ y^{ l_1}}
- 
\left(
\Pi_{0, l_1}^j
\right)_x
=
G_{ y^{ l_1}}^j 
+ 
\frac{ 1}{ 2} \, H_{ l_1, x}^j 
- 
\frac{ 1}{ 2} \, \delta_{l_1}^j\, \Theta_x^0 
= \\
& 
= 
\frac{ 1}{ 2} \, \Theta^0 \, H_{ l_1}^j 
-
\frac{ 1}{ 2} \, \delta_{ l_1}^j \, \Theta^0 \, \Theta^0 
- \\
&
\ \ \ \ \
-
\sum_k \, \left(-G^k\right) \left(
- 
L_{ l_1, k}^j + \frac{ 1}{ 2}\, \delta_{ l_1}^j \, L_{ k,k}^k 
+ 
\frac{ 1}{ 2} \, \delta_k^j \, L_{ l_1, l_1}^{ l_1} 
+ 
\frac{ 1}{ 2} \, \delta_{l_1}^j \, \Theta^k 
+
\frac{ 1}{ 2} \, \delta_k^j \, \Theta^k
\right) 
+ \\
&
\ \ \ \ \
+ 
\left(
\frac{ 1}{ 2} \, L_{ l_1, l_1}^{ l_1} 
+ 
\frac{ 1}{ 2} \, \Theta^{ l_1}
\right)
\left( -G^j\right)
+ \\
& 
\ \ \ \ \
+
\sum_k \left(
- 
\frac{ 1}{ 2} \, H_{ l_1}^k 
+ 
\frac{ 1}{ 2} \, \delta_{ l_1}^k \, \Theta^0 
\right)
\left(
- 
\frac{ 1}{ 2} \, H_k^j 
+
\frac{ 1}{ 2} \, \delta_k^j \, \Theta^0
\right) = \\
& 
=
\underline{
\frac{ 1}{ 2} \, H_{ l_1}^j \, \Theta^0 
}_{ 
\octagon \! \! \! \! \! \tiny{\sf a}}
- 
\frac{ 1}{ 2} \, \delta_{ l_1}^j \, \Theta^0 \, \Theta^0 
-
\sum_k \, G^k \, L_{ l_1, k}^j 
+ 
\frac{ 1}{ 2} \, \delta_{ l_1}^j \, \sum_k \, G^k \, L_{ k, k}^k 
+
\underline{
\frac{ 1}{ 2} \, G^j \, L_{ l_1, l_1}^{ l_1} 
}_{ 
\octagon \! \! \! \! \! \tiny{\sf b}}
+ \\
& 
\ \ \ \ \
+ 
\frac{ 1}{ 2} \, \delta_{ l_1}^j \, \sum_k \, G^k \, \Theta^k 
+
\underline{
\frac{ 1}{ 2} \, G^j \, \Theta^{ l_1} 
}_{ 
\octagon \! \! \! \! \! \tiny{\sf c}}
- 
\underline{
\frac{ 1}{ 2} \, G^j \, L_{ l_1, l_1}^{ l_1} 
}_{ 
\octagon \! \! \! \! \! \tiny{\sf b}}
- 
\underline{
\frac{ 1}{ 2} \, G^j \, \Theta^{ l_1} 
}_{
\octagon \! \! \! \! \! \tiny{\sf c}}
+ 
\frac{ 1}{ 4} \, \sum_k \, H_{ l_1}^k \, H_k^j 
- \\
& 
\ \ \ \ \
- 
\underline{ 
\frac{ 1}{ 4} \, H_{ l_1}^j \, \Theta^0 
}_{
\octagon \! \! \! \! \! \tiny{\sf a}}
- 
\underline{
\frac{ 1}{ 4} \, H_{ l_1}^j \, \Theta^0 
}_{ 
\octagon \! \! \! \! \! \tiny{\sf a}}
+
\frac{ 1}{ 4} \, \delta_{ l_1}^j \, \Theta^0 \, \Theta^0.
\endaligned\right.
\end{equation}
Eliminating the underlined vanishing terms with the
letters $a$, $b$, $c$ and $d$ appended, multiplying by $-2$ and
reorganizing the identity so as to put the term $\delta_{ l_1}^j\,
\Theta_x^0$ solely in the left-hand side, we obtain
the relation
\def\theequation{3.69}\begin{equation}
\left\{
\aligned
\delta_{ l_1}^j \, \Theta_x^0 
& 
=
-
2\, G_{ y^{ l_1}}^j 
+ 
H_{ l_1, x}^j
+
2 \, \sum_k \, G^k \, L_{ l_1, k}^j 
- 
\delta_{ l_1}^j \, \sum_k \, G^k \, L_{ k, k}^k 
- \\
& 
\ \ \ \ \ \ \ \ \ \ \ \
\ \ \ \ \ \
-
\frac{ 1}{ 2} \, \sum_k \, H_{ l_1}^k \, H_k^j
-
\delta_{ l_1}^j \, \sum_k \, G^k \, \Theta^k
+ 
\frac{ 1}{ 2} \, \delta_{ l_1}^j\, \Theta^0 \, \Theta^0.
\endaligned\right.
\end{equation}

\subsection*{ 3.70.~Conventions for simplifications of
formal expressions} Before proceeding further, let us explain how we
will organize the computations with the formal expressions we shall
encounter until the end of Section~4. {\it Our main goal is to devise
a methodology of writing formal computations which enables to check
every computation visually}, {\sf without being forced to rebuild any
intermediate step}. In fact, it would be unsatisfactory to just claim
that Theorem~1.7~{\bf (3)} 
follows by hidden massive formal computations, so
that we have to guide the rigorous and demanding reader until the very
extremal branches of our coral tree of formal computations.

As an example, suppose that we have to simplify the equation $0 = A_x
+ B_y + A - B+ 2 \, C- \frac{ 1}{ 3} \, D - \frac{ 2}{ 3} \, A +
\frac{ 1}{ 6} \, D + E + B - 2\, C$. In the beginning, the terms $A_x$
and $B_y$ are differentiated once and they 
do not simplify with other terms. To distinguish
them, we underline them plainly and we copy the nine remaining terms
afterwards:
\def\theequation{3.71}\begin{equation}
\small
\left\{
\aligned
0 = 
& \
\underline{ A_x + B_y} + \\
& \
+
\underline{
A
}_{ \fbox{\tiny 1}} 
- 
\underline{ 
B
}_{ \octagon \! \! \! \! \! \tiny{\sf a}} 
+ 
\underline{ 
2 \, C
}_{ \octagon \! \! \! \! \! \tiny{\sf b }}
- 
\underline{ 
\frac{ 1}{ 3} \, D
}_{\fbox{\tiny 2}} 
- 
\underline{ 
\frac{ 2}{ 3} \, A
}_{ \fbox{\tiny 1}}
+
\underline{ 
\frac{ 1}{ 6} \, D
}_{ \fbox{\tiny 2}}
+ 
\underline{ 
E
}_{ \fbox{\tiny 3}}
+
\underline{ 
B
}_{ \octagon \! \! \! \! \! \tiny{\sf a}}
- 
\underline{ 
2\, C
}_{ \octagon \! \! \! \! \! \tiny{\sf b}} \ . 
\endaligned\right.
\end{equation}
Here, each remaining term is also underlined, with a number or with a
letter appended. For reasons of typographical readability, we never
underline the sign, $+$ or $-$ of each term; however, it 
should be understood
that {\it every term always includes its} (not underlined) 
{\it sign}. Until the
end of Section~4, we shall use the roman alphabetic letters $a$, $b$,
$c$, {\it etc.} inside an octagon $\octagon$ to exhibit the vanishing
terms. As readily checked by the eyes, we indeed have $- \underline{ B
}_{ \octagon \! \! \! \! \! \! \tiny{\sf a}} + \underline{ B }_{ \octagon
\! \! \! \! \! \! \tiny{\sf a}} = 0$ and $\underline{ 2\, C}_{ \octagon
\! \! \! \! \! \! \tiny{\sf b }} - \underline{2\, C }_{ \octagon \! \! \!
\! \! \! \tiny{\sf b }} = 0$. Also, until the end of Section~4, we shall
use the numbers $1$, $2$, $3$, {\it etc.} inside a square $\square$ to
exhibit the remaining terms, collected
in a certain order. The numbers have the following
signification: after the simplifications, the equation~\thetag{ 3.71}
may be written
\def\theequation{3.72}\begin{equation}
\small
\left\{
\aligned
0 = 
& \
\underline{ A_x + B_y} + \\
& \
+
\frac{ 1}{ 3}\, A - 
\frac{ 1}{ 6} \, D +
E.
\endaligned\right.
\end{equation}
Here, the plainly underlined terms $\underline{ A_x + B_y}$ do not
count in the numbering (their number is zero, for instance) and {\it
the first term of the second line $\frac{ 1}{ 3} \, A$ correspond to
the addition of all terms $\fbox{\tiny 1}$ in~\thetag{ 3.69}}.
Analogously, the second term $-\frac{ 1}{ 6} \, D$ correspond to the
addition of all terms $\fbox{\tiny 2}$ in~\thetag{ 3.69}. Again, this
guiding facilitates the checking of the correctness of the
computation, using simply the eyes. No hidden delicate computational
step is ``left to the reader'' for the convenience of
the writer.

This principle will be constantly used until the end of Section~4; it
has been systematically used in~\cite{ m2004} and it
could be applied in various other contexts. Again, the advantage is
that it enables to check the correctness of all the formal
computations just by reading, without having to write anything
more. This is also useful for the author.

\subsection*{ 3.73.~Choice of an ordering}
Until the end of Section~4, we shall have to deal with terms $G$, $H$,
$L$, $M$, $\Theta$ together with indices and partial derivatives up to
order two. In order to organize the formal expressions in a way which
provides an easier deciphering, it is convenient to introduce an order
between these differential monomials. In a symbolic index-free
notation, we choose:
\def\theequation{3.74}\begin{equation}
G < H < L < M < \Theta.
\end{equation}
It follows for instance that $G < GH < GL < HHL < HLM\Theta$.
Also, if a sum appears, we choose: $GM < \sum \, GM$.

Here, we have only considered terms of order zero,
without partial differentiation. The first order
partial differentiations are $(\cdot)_x$ and
$(\cdot)_y$, again in symbolic notation, 
dropping the indices. We choose:
\def\theequation{3.75}\begin{equation}
G_x < G_y < H_x < H_y < L_x < L_y < M_x < M_y < \Theta_x < \Theta_y 
< G < \cdots.
\end{equation}
For second order derivatives, we choose:
\def\theequation{3.76}\begin{equation}
G_{ xx} < G_{ xy} < G_{ yy} < H_{ xx} < \cdots < \Theta_{ yy} 
< G_x < \cdots.
\end{equation}
As a final general example including 
indices we have the inequalities
\def\theequation{3.77}\begin{equation}
H_{ l_1, y^{ l_2}}^j <
L_{ l_1, l_2, x}^j < 
G^j M_{ l_1 , l_2} <
\sum_{ k= 1}^m\, 
G^k \, M_{ l_2, k} < 
\sum_{ k= 1}^m \, H_k^{ l_2}\,H_{l_2,l_2}^k,
\end{equation} 
extracted from (II) of Theorem~1.7~{\bf (3)}.

In the sequel, we shall call 
\begin{itemize}
\item[$\bullet$]
terms {\sl of order $0$} monomials like 
$G$, $H$, $L\,M$, $G\,H\,M$;
\item[$\bullet$]
terms {\sl of order $1$} monomials like 
$G_x$, $G_x \, M$, $L_x\,\Theta$;
\item[$\bullet$]
terms {\sl of order $2$} monomials like
$G_{ xy}$, $L_{ yy}$, $M_{ xx}$ (our terms 
of order two will always be linear),
\end{itemize}
according to the top order partial derivatives.

\subsection*{ 3.78.~A mean of checking intuitively the
validity of partial differential relations}
Before replacing~\thetag{ 3.64} in the five remaining identities
$(3.65)_2$, $(3.65)_3$, $(3.65)_4$, $(3.65)_5$ and $(3.65)_6$, let us
observe that if we assume that 
$j\neq l_1$ in~\thetag{ 3.69}, then all the terms
involving $\Theta$ vanish, so that we obtain the {\it nontrivial
partial differential equations}:
\def\theequation{3.79}\begin{equation}
0 = \underline{ 
- 2\, G_{ y^{ l_1}}^{ j} + 
H_{ l_1, x}^{ j}} + 
2\, \sum_k \, G^k \, L_{ l_1, k}^{ j} 
- 
\frac{ 1}{2} \, \sum_k \, H_{ l_1}^k \, H_k^{j},
\end{equation}
for $j\neq l_1$. Here, we have underlined the first order terms
plainly, in order to distinguish them from the terms of order zero.
These equations coincides with (I) of Theorem~1.7~{\bf (3)}, 
again specialized
with $j\neq l_1$. Importantly, we notice that the choice of indices $j
\neq l_1$ is possible only if $m\geqslant 2$. Thus, we have derived {\it a
subpart of (I) as a necessary condition for the point equivalence to
$Y_{XX}^j= 0$, $j= 1, \dots, m\geqslant 2$}. These first order equations
show at once that there is a strong difference with the case $m= 1$.

How can we confirm (at least informally) that the functions $G^j$,
$H_{ l_1}^j$ and $L_{ l_1, l_2}^j$ given by~\thetag{ 3.35} in terms of
$X$ and $Y^j$ do indeed satisfy these equations for $j \neq l_1$?
Dropping the zero order terms in~\thetag{ 3.79} 
above, we obtain an approximated
equation
\def\theequation{3.80}\begin{equation}
0 \equiv 
- 2\, G_{ y^{ l_1}}^j + H_{ l_1, x}^j.
\end{equation}
Here, the sign $\equiv$ precisely means: ``{\sl modulo zero order
terms}''. We claim that this approximated
equation is a consequence of the existence of
$X$, $Y^j$. 

Indeed, according to the approximation~\thetag{ 3.58},
together with the definition~\thetag{ 3.35} of the functions $G^j$ and
$H_{ l_1}^j$, we have
\def\theequation{3.81}\begin{equation}
\left\{
\aligned
G^j 
& \
= - \square_{ xx}^j \cong - Y_{ xx}^j, \\
H_{ l_1}^j 
& \
= - 2\, \square_{ xy^{ l_1}}^j \cong - 2\, 
Y_{ xy^{ l_1}}^j.
\endaligned\right.
\end{equation}
Differentiation of the first line with respect to $y^{ l_1}$ and of
the second line with respect to $x$ yields:
\def\theequation{3.82}\begin{equation}
G_{ y^{ l_1}}^j \cong 
-
Y_{ xxy^{ l_1}}^j 
\ \ \ \ \ \
{\rm and}
\ \ \ \ \ \ 
H_{ l_1, x}^j \cong -2\, Y_{ xy^{ l_1}x}^j, 
\end{equation}
so that we indeed have $0 \equiv - 2\, G_{ y^{ l_1}}^j + H_{ l_1,
x}^j$, approximatively and modulo the derivatives of order $0$, $1$
and $2$ of the functions $X$, $Y^j$.

Similar verifications have been effected constantly in our
manuscript in order to control 
the truth of the formal computations that we shall
expose until the end of Section~4.

\subsection*{ 3.83.~Continuation}
From now on and up to the end of Section~4, the hardest computational
core of the proof may\,\,---\,\,at last\,\,---\,\,be developed.
Further amazing computational obstacles will be encountered.

Replacing plainly~\thetag{ 3.64}
in~$(3.65)_2$, we get:
\def\theequation{3.84}\begin{equation}
\small
\aligned
& 
\left(
\Pi_{l_1, l_2}^j
\right)_x
- 
\left(
\Pi_{l_1, 0}^j
\right)_{ y^{ l_2}}
=
- 
L_{l_1, l_2, x}^j 
+ 
\frac{ 1}{ 2} \, \delta_{ l_1}^j \, L_{ l_2, l_2, x}^{l_2}
+ 
\frac{ 1}{ 2} \, \delta_{ l_2}^j \, L_{ l_1, l_1, x}^{l_1}
+
\frac{ 1}{ 2} \, \delta_{ l_1}^j \, \Theta_x^{ l_2}
+ \\
& \
\ \ \ \ \ \ \ \ \ \ \
\ \ \ \ \ \ \ \ \ \ \
\ \ \ \ \ \ \ \ \ \ \
\ \ \ \ \ \ \ \ \ \ \
+
\frac{ 1}{ 2} \, \delta_{ l_2}^j \, \Theta_x^{ l_1}
+
\frac{ 1}{ 2} \, H_{ l_1, y^{ l_2}}^j 
- 
\frac{ 1}{ 2} \, \delta_{ l_1}^j \, \Theta_{ y^{ l_2}}^0 
= \\
& \
= 
- 
\Pi_{ l_1, l_2}^0 \cdot \Pi_{ 0, 0}^j 
- 
\sum_k \, \Pi_{ l_1, l_2}^k \cdot \Pi_{ 0, k}^j 
+
\Pi_{ l_1, 0}^0 \cdot \Pi_{ l_2, 0}^j 
+ 
\sum_k \, \Pi_{ l_1, 0}^k \cdot \Pi_{ l_2, k}^j 
= \\
& \
= 
M_{ l_1, l_2} \, G^j
- 
\sum_k\, 
\left(
- 
L_{ l_1, l_2}^k 
+
\frac{ 1}{ 2} \, \delta_{ l_1}^k \, L_{ l_2, l_2}^{ l_2}
+ 
\frac{ 1}{ 2} \, \delta_{ l_2}^k \, L_{ l_1, l_1}^{ l_1}
+
\frac{ 1}{ 2} \, \delta_{ l_1}^k \, \Theta^{ l_2}
+
\frac{ 1}{ 2} \, \delta_{ l_2}^k \, \Theta^{ l_1}
\right)
\cdot \\
& \
\ \ \ \ \ \ \ \ \ \ \
\ \ \ \ \ \ \ \ \ \ \
\ \ \ \ \ \ \ \ \ \ \
\ \ \ \ \ \ \ \ \ \ \
\ \ \ \ \ \ \ \ \ \ \
\ \ \ \ \ \ \ \ \ \ \
\cdot
\left(
- 
\frac{ 1}{ 2} \, H_k^j 
+ 
\frac{ 1}{ 2} \, \delta_k^j \, \Theta^0
\right) 
+ \\
& \
\ \ \ \ \
+ 
\left(
\frac{ 1}{ 2} \, L_{ l_1, l_1}^{ l_1} 
+
\frac{ 1}{ 2} \, \Theta^{ l_1}
\right)
\cdot
\left(
- 
\frac{ 1}{ 2} \, H_{ l_2}^j
+ 
\frac{ 1}{ 2} \, \delta_{ l_2}^j \, \Theta^0
\right)
+ \\
& \
\ \ \ \ \
+
\sum_k \, 
\left(
-
\frac{ 1}{ 2} \, H_{ l_1}^k 
+ 
\frac{ 1}{ 2} \, \delta_{ l_1}^k \, \Theta^0 
\right) 
\cdot \\
& \
\ \ \ \ \ \ \ \ \ \ \
\cdot
\left(
- 
L_{ l_2, k}^j 
+
\frac{ 1}{ 2} \, \delta_{ l_2}^j \, L_{ k, k}^k 
+
\frac{ 1}{ 2} \, \delta_k^j \, L_{ l_2, l_2}^{ l_2}
+ 
\frac{ 1}{ 2} \, \delta_{ l_2}^j \, \Theta^k 
+
\frac{ 1}{ 2} \, \delta_k^j \, \Theta^{ l_2}
\right).
\endaligned
\end{equation}
Developing the products and ordering each monomial, we get:
\def\theequation{3.85}\begin{equation}
\small
\aligned
& \
=
\underline{ 
G^j \, M_{ l_1, l_2}
}_{ \fbox{\tiny 1}}
- 
\underline{ 
\frac{ 1}{ 2} \, \sum_k \, H_k^j \, L_{ l_1, l_2}^k 
}_{ \fbox{\tiny 2}}
+
\underline{ 
\frac{ 1}{ 4} \, H_{ l_1}^j \, L_{ l_2, l_2}^{ l_2}
}_{ \octagon \! \! \! \! \! \tiny{\sf a}} 
+
\underline{ 
\frac{ 1}{ 4} \, H_{ l_2}^j \, L_{ l_1, l_1}^{ l_1}
}_{ \octagon \! \! \! \! \! \tiny{\sf b}} 
+ \\
& \
\ \ \ \ \
+
\underline{ 
\frac{ 1}{ 4} \, H_{ l_1}^j \, \Theta^{ l_2}
}_{ \octagon \! \! \! \! \! \tiny{\sf c}} 
+
\underline{ 
\frac{ 1}{ 4} \, H_{ l_2}^j \, \Theta^{ l_1}
}_{ \octagon \! \! \! \! \! \tiny{\sf d}} 
+
\underline{ 
\frac{ 1}{ 2} \, L_{ l_1, l_2}^j \, \Theta^0 
}_{ \octagon \! \! \! \! \! \tiny{\sf e}} 
- 
\underline{ 
\frac{ 1}{ 4} \, \delta_{ l_1}^j\, L_{ l_2, l_2}^{l_2} \, \Theta^0 
}_{ \octagon \! \! \! \! \! \tiny{\sf f}}
- \\ 
& \
\ \ \ \ \
-
\underline{ 
\frac{ 1}{ 4} \, \delta_{ l_2}^j\, L_{ l_1, l_1}^{l_1} \, \Theta^0 
}_{ \octagon \! \! \! \! \! \tiny{\sf g}}
- 
\underline{ 
\frac{ 1}{ 4} \, \delta_{ l_1}^j \, \Theta^0 \, \Theta^{ l_2}
}_{ \octagon \! \! \! \! \! \tiny{\sf h}} 
- 
\underline{ 
\frac{ 1}{ 4} \, \delta_{ l_2}^j \, \Theta^0 \, \Theta^{ l_1}
}_{ \octagon \! \! \! \! \! \tiny{\sf i}} 
-
\underline{ 
\frac{ 1}{ 4} \, H_{ l_2}^j \, L_{ l_1, l_1}^{ l_1}
}_{ \octagon \! \! \! \! \! \tiny{\sf b}}
+ \\ 
& \
\ \ \ \ \
+
\underline{ 
\frac{ 1}{ 4} \, \delta_{ l_2}^j \, L_{ l_1, l_1}^{ l_1} \, \Theta^0
}_{ \fbox{\tiny 6}}
- 
\underline{ 
\frac{ 1}{ 4} \, H_{ l_2}^j \, \Theta^{ l_1}
}_{ \octagon \! \! \! \! \! \tiny{\sf d}} 
+
\underline{ 
\frac{ 1}{ 4} \, \delta_{ l_2}^j \, \Theta^0 \, \Theta^{ l_1}
}_{ \octagon \! \! \! \! \! \tiny{\sf i}} 
+
\underline{ 
\frac{ 1}{ 2} \, \sum_k \, H_{ l_1}^k \, L_{ l_2, k}^j 
}_{ \fbox{\tiny 3}}
- \\
& \
\ \ \ \ \
- 
\underline{ 
\frac{ 1}{ 4} \, \sum_k\, H_{ l_1}^k \, L_{ k, k}^k 
}_{ \fbox{\tiny 4}}
- 
\underline{ 
\frac{ 1}{ 4} \, H_{l_1}^j \, L_{ l_2, l_2}^{ l_2}
}_{ \octagon \! \! \! \! \! \tiny{\sf a}}
- 
\underline{ 
\frac{ 1}{ 4} \, \delta_{ l_2}^j \, \sum_k\, H_{l_1}^k \, \Theta^k
}_{ \fbox{\tiny 5}} 
- 
\underline{ 
\frac{ 1}{ 4} \, H_{ l_1}^j \, \Theta^{ l_2}
}_{ \octagon \! \! \! \! \! \tiny{\sf c}} 
- \\
& \
\ \ \ \ \
-
\underline{ 
\frac{ 1}{ 2} \, L_{ l_2, l_1}^j \, \Theta^0 
}_{ \octagon \! \! \! \! \! \tiny{\sf e}}
+ 
\underline{ 
\frac{ 1}{ 4} \, \delta_{ l_2}^j \, L_{ l_1, l_1}^{ l_1} \, \Theta^0 
}_{ \octagon \! \! \! \! \! \tiny{\sf g}} 
+ 
\underline{ 
\frac{ 1}{ 4} \, \delta_{ l_1}^j \, L_{ l_2, l_2}^{ l_2} \, \Theta^0
}_{ \octagon \! \! \! \! \! \tiny{\sf f}}
+ 
\underline{ 
\frac{ 1}{ 4} \, \delta_{ l_2}^j \, \Theta^0 \, \Theta^{ l_1} 
}_{ \fbox{\tiny 7}} 
+ \\
& \
\ \ \ \ \
+
\underline{ 
\frac{ 1}{ 4} \, \delta_{ l_1}^j \, \Theta^0 \, \Theta^{ l_2}
}_{ \octagon \! \! \! \! \! \tiny{\sf h}}. 
\endaligned
\end{equation}
We simplify according to our general principles and we reorganize the
equality between the first two lines of~\thetag{ 3.84} and~\thetag{
3.85} so as to put all terms $\Theta_x$ in the left-hand side of the
equality and to put all remaining terms in the right-hand side,
respecting the order of \S3.73. We get:
\def\theequation{3.86}\begin{equation}
\aligned
& \
\frac{ 1}{ 2} \, \delta_{ l_1}^j \, \Theta_x^{ l_2} 
+ 
\frac{ 1}{ 2} \, \delta_{ l_2}^j \, \Theta_x^{ l_1}
-
\frac{ 1}{ 2} \, \delta_{ l_1}^j \, \Theta_{ y^{ l_2}}^0 
= \\
& \
\ \ \ \ \
= 
-
\underline{
\frac{ 1}{ 2} \, H_{ l_1, y^{ l_2}}^j 
+ 
L_{ l_1, l_2, x}^j 
- 
\frac{ 1}{ 2} \, \delta_{ l_1}^j L_{ l_2, l_2, x}^{ l_2}
- 
\frac{ 1}{ 2} \, \delta_{ l_2}^j L_{ l_1, l_1, x}^{ l_1}
}
+ \\
& \
\ \ \ \ \ \ \ \ \ \
+ 
G^j \, M_{ l_1, l_2} 
- 
\frac{ 1}{ 2} \, \sum_k \, H_k^j \, L_{ l_1, l_2}^k 
+ 
\frac{ 1}{ 2} \, \sum_k \, H_{ l_1}^k \, L_{ l_2, k}^j 
- \\
& \
\ \ \ \ \ \ \ \ \ \
- 
\frac{ 1}{ 4} \, \delta_{ l_2}^j \, \sum_k \, H_{ l_1}^k \, L_{ k, k}^k 
- 
\frac{ 1}{ 4} \, \delta_{ l_2}^j \sum_k \, H_{ l_1}^k \, \Theta^k 
+
\frac{ 1}{ 4} \, \delta_{ l_2}^j \, L_{l_1,l_1}^{l_1} \, \Theta^0 
+\\
& \
\ \ \ \ \ \ \ \ \ \
+
\frac{ 1}{ 4} \, \delta_{ l_2}^j \, \Theta^0 \, \Theta^{ l_1}.
\endaligned
\end{equation}
Here, we have underlined plainly the four first order terms 
appearing in the second line.

Next, replacing plainly~\thetag{ 3.64} in $(3.65)_3$, we
get:
\def\theequation{3.87}\begin{equation}
\small
\aligned
& 
\left(
\Pi_{ l_1, l_2}^j
\right)_{ y^{ l_3}}
-
\left(
\Pi_{ l_1, l_3}^j
\right)_{ y^{ l_2}}
= \\
& \
\ \ \ \ \
= 
-
L_{ l_1, l_2, y^{ l_3}}^j 
+ 
\frac{ 1}{ 2} \, \delta_{ l_1}^j \, L_{ l_2,l_2, y^{ l_3}}^{l_2}
+
\frac{ 1}{ 2} \, \delta_{ l_2}^j \, L_{ l_1,l_1, y^{ l_3}}^{l_1}
+
\frac{ 1}{ 2} \, \delta_{ l_1}^j \, \Theta_{y^{l_3}}^{l_2}
+
\frac{ 1}{ 2} \, \delta_{ l_2}^j \, \Theta_{y^{l_3}}^{l_1}
+ \\
& \ 
\ \ \ \ \ \ \ \ \ \
+
L_{ l_1, l_3, y^{ l_2}}^j 
-
\frac{ 1}{ 2} \, \delta_{ l_1}^j \, L_{ l_3,l_3, y^{ l_2}}^{l_3}
-
\frac{ 1}{ 2} \, \delta_{ l_3}^j \, L_{ l_1,l_1, y^{ l_2}}^{l_1}
-
\frac{ 1}{ 2} \, \delta_{ l_1}^j \, \Theta_{y^{l_2}}^{l_3}
-
\frac{ 1}{ 2} \, \delta_{ l_3}^j \, \Theta_{y^{l_2}}^{l_1}
= \\
& \
\ \ \ \ \
=
-
\Pi_{l_1,l_2}^0 \cdot \Pi_{l_3,0}^j
- 
\sum_k\, \Pi_{l_1,l_2}^k \cdot \Pi_{ l_3, k}^j 
+ 
\Pi_{ l_1, l_3}^0 \cdot \Pi_{ l_2, 0}^j 
+ 
\sum_k\, \Pi_{l_1, l_3}^k \cdot \Pi_{ l_2, k}^j = \\
& \
\ \ \ \ \
=
- 
M_{ l_1, l_2} 
\cdot
\left(
-
\frac{ 1}{ 2} \, H_{l_3}^j + 
\frac{ 1}{ 2} \, \delta_{ l_3}^j \, \Theta^0
\right)
- 
\sum_k \, 
\left(
- L_{l_1, l_2}^k 
+ 
\frac{ 1}{ 2} \, \delta_{ l_1}^k \, L_{ l_2, l_2}^{ l_2}
+
\right. \\
& \
\ \ \ \ \ \ \ \ \ \
\left.
+ 
\frac{ 1}{ 2} \, \delta_{ l_2}^k \, L_{ l_1, l_1}^{ l_1}
+
\frac{ 1}{ 2} \, \delta_{ l_1}^k \, \Theta^{ l_2}
+
\frac{ 1}{ 2} \, \delta_{ l_2}^k \, \Theta^{ l_1}
\right)
\cdot
\left(
- 
L_{ l_3, k}^j 
+
\frac{ 1}{ 2} \, \delta_{ l_3}^j \, L_{ k, k}^k 
+ 
\right. \\
& \
\ \ \ \ \ \ \ \ \ \
\left.
+
\frac{ 1}{ 2} \, \delta_k^j \, L_{ l_3, l_3}^{ l_3}
+ 
\frac{ 1}{ 2} \, \delta_{ l_3}^j \, \Theta^k 
+
\frac{ 1}{ 2} \, \delta_k^j \, \Theta^{ l_3}
\right)
+ \\
& \
\ \ \ \ \ \ \ \ \ \
+
M_{ l_1, l_3} 
\cdot
\left(
-
\frac{ 1}{ 2} \, H_{l_2}^j + 
\frac{ 1}{ 2} \, \delta_{ l_2}^j \, \Theta^0
\right)
+
\sum_k \, 
\left(
- L_{l_1, l_3}^k 
+ 
\frac{ 1}{ 2} \, \delta_{ l_1}^k \, L_{ l_3, l_3}^{ l_3}
+
\right. \\
& \
\ \ \ \ \ \ \ \ \ \
\left.
+ 
\frac{ 1}{ 2} \, \delta_{ l_3}^k \, L_{ l_1, l_1}^{ l_1}
+
\frac{ 1}{ 2} \, \delta_{ l_1}^k \, \Theta^{ l_3}
+
\frac{ 1}{ 2} \, \delta_{ l_3}^k \, \Theta^{ l_1}
\right)
\cdot
\left(
- 
L_{ l_2, k}^j 
+
\frac{ 1}{ 2} \, \delta_{ l_2}^j \, L_{ k, k}^k 
+ 
\right. \\
& \
\ \ \ \ \ \ \ \ \ \
\left.
+
\frac{ 1}{ 2} \, \delta_k^j \, L_{ l_2, l_2}^{ l_2}
+ 
\frac{ 1}{ 2} \, \delta_{ l_2}^j \, \Theta^k 
+
\frac{ 1}{ 2} \, \delta_k^j \, \Theta^{ l_2}
\right). \\
\endaligned
\end{equation}
Developing the products and ordering each monomial, we get:
\def\theequation{3.88}\begin{equation}
\small
\aligned
& \
= 
\underline{ 
\frac{ 1}{ 2} \, H_{ l_3}^j \, M_{ l_1, l_2} 
}_{ \fbox{\tiny 1}}
- 
\underline{ 
\frac{ 1}{ 2} \, \delta_{ l_3}^j \, M_{l_1,l_2} \, \Theta^0 
}_{ \fbox{\tiny 16}}
- 
\underline{ 
\sum_k \, L_{ l_1,l_2}^k \, L_{l_3,k}^j
}_{ \fbox{\tiny 6}}
+
\underline{ 
\frac{ 1}{ 2} \, \delta_{ l_3}^j \,\sum_k\,L_{l_1,l_2}^k\,L_{k,k}^k
}_{ \fbox{\tiny 7}}
+ \\
& \
\ \ \ \ \
+
\underline{ 
\frac{ 1}{ 2} \, L_{l_1,l_2}^j \, L_{l_3,l_3}^{l_3}
}_{ \octagon \! \! \! \! \! \tiny{\sf a}} 
+
\underline{ 
\frac{ 1}{ 2} \, \delta_{ l_3}^j \,\sum_k\,L_{l_1,l_2}^k\,\Theta^k
}_{ \fbox{\tiny 13}}
+
\underline{ 
\frac{ 1}{ 2} \, L_{l_1,l_2}^j\, \Theta^{l_3}
}_{ \octagon \! \! \! \! \! \tiny{\sf b}} 
+
\underline{ 
\frac{ 1}{ 2} \, L_{ l_2, l_2}^{ l_2} \, L_{l_3,l_1}^j
}_{ \octagon \! \! \! \! \! \tiny{\sf c}} 
- \\
& \
\ \ \ \ \
-
\underline{ 
\frac{ 1}{ 4}\,\delta_{l_3}^j\,L_{l_2,l_2}^{l_2}\,L_{l_1,l_1}^{l_1}
}_{ \fbox{\tiny 4}}
- 
\underline{ 
\frac{ 1}{ 4} \, \delta_{ l_1}^j \,L_{l_2,l_2}^{l_2}\,L_{l_3,l_3}^{l_3}
}_{ \octagon \! \! \! \! \! \tiny{\sf d}} 
- 
\underline{ 
\frac{ 1}{ 4} \,\delta_{ l_3}^j \,L_{l_2,l_2}^{l_2}\,\Theta^{l_1}
}_{ \octagon \! \! \! \! \! \tiny{\sf e}} 
-
\underline{ 
\frac{ 1}{ 4} \,\delta_{ l_1}^j \,L_{l_2,l_2}^{l_2}\,\Theta^{l_3}
}_{ \octagon \! \! \! \! \! \tiny{\sf f}} 
+ \\
& \
\ \ \ \ \
+
\underline{ 
\frac{ 1}{ 2} \, L_{l_1,l_1}^{l_1} \, L_{l_3,l_2}^j 
}_{ \octagon \! \! \! \! \! \tiny{\sf g}} 
- 
\underline{ 
\frac{ 1}{ 4}\,\delta_{l_3}^j\,L_{l_1,l_1}^{l_1}\,L_{l_2,l_2}^{l_2}
}_{ \octagon \! \! \! \! \! \tiny{\sf h}} 
- 
\underline{ 
\frac{ 1}{ 4}\,\delta_{l_2}^j\,L_{l_1,l_1}^{l_1}\,L_{l_3,l_3}^{l_3}
}_{ \octagon \! \! \! \! \! \tiny{\sf i}} 
-
\underline{ 
\frac{ 1}{ 4} \, \delta_{ l_3}^j \,L_{l_1,l_1}^{l_1}\,\Theta^{l_2} 
}_{ \octagon \! \! \! \! \! \tiny{\sf j}} 
- \\
& \
\ \ \ \ \
-
\underline{ 
\frac{ 1}{ 4} \, \delta_{ l_2}^j \,L_{l_1,l_1}^{l_1}\,\Theta^{l_3} 
}_{ \octagon \! \! \! \! \! \tiny{\sf k}} 
+
\underline{ 
\frac{ 1}{ 2} \, L_{l_3,l_1}^j \, \Theta^{l_2}
}_{ \octagon \! \! \! \! \! \tiny{\sf l}} 
- 
\underline{ 
\frac{ 1}{ 4} \, \delta_{ l_3}^j \,L_{l_1,l_1}^{l_1}\,\Theta^{l_2}
}_{ \fbox{\tiny 10}}
- 
\underline{ 
\frac{ 1}{ 4} \, \delta_{ l_1}^j \,L_{l_3,l_3}^{l_3}\,\Theta^{l_2}
}_{ \octagon \! \! \! \! \! \! \tiny{\sf m}}
- \\
& \
\ \ \ \ \
- 
\underline{ 
\frac{ 1}{ 4} \, \delta_{ l_3}^j \,\Theta^{l_2}\,\Theta^{l_1}
}_{ \octagon \! \! \! \! \! \tiny{\sf n}} 
- 
\underline{ 
\frac{ 1}{ 4} \, \delta_{ l_1}^j \,\Theta^{l_2}\,\Theta^{l_3}
}_{ \octagon \! \! \! \! \! \tiny{\sf o}}
+
\underline{ 
\frac{ 1}{ 2} \,L_{l_3,l_2}^j\,\Theta^{l_1}
}_{ \octagon \! \! \! \! \! \tiny{\sf p}} 
- 
\underline{ 
\frac{ 1}{ 4} \, \delta_{ l_3}^j \,L_{l_2,l_2}^{l_2}\,\Theta^{l_1}
}_{ \fbox{\tiny 12}}
- \\
& \
\ \ \ \ \
-
\underline{ 
\frac{ 1}{ 4} \, \delta_{ l_2}^j \,L_{l_3,l_3}^{l_3}\,\Theta^{l_1}
}_{ \octagon \! \! \! \! \! \tiny{\sf q}} 
- 
\underline{ 
\frac{ 1}{ 4} \, \delta_{ l_3}^j \,\Theta^{l_1}\,\Theta^{l_2}
}_{ \fbox{\tiny 18}}
- 
\underline{ 
\frac{ 1}{ 4} \,\delta_{ l_2}^j \,\Theta^{l_1}\,\Theta^{l_3}
}_{ \octagon \! \! \! \! \! \tiny{\sf r}} 
-
\endaligned
\end{equation}
$$
\aligned
& \
\ \ \ \ \
-
\underline{ 
\frac{ 1}{ 2} \, H_{ l_2}^j \, M_{ l_1, l_3} 
}_{ \fbox{\tiny 2}}
+
\underline{ 
\frac{ 1}{ 2} \, \delta_{ l_2}^j \, M_{l_1,l_3} \, \Theta^0 
}_{ \fbox{\tiny 15}}
+
\underline{ 
\sum_k \, L_{ l_1,l_3}^k \, L_{l_2,k}^j
}_{ \fbox{\tiny 5}}
-
\underline{ 
\frac{ 1}{ 2} \, \delta_{ l_2}^j \,\sum_k\,L_{l_1,l_3}^k\,L_{k,k}^k
}_{ \fbox{\tiny 8}}
- \\
& \
\ \ \ \ \
-
\underline{
\frac{ 1}{ 2} \, L_{l_1,l_3}^j \, L_{l_2,l_2}^{l_2}
}_{ \octagon \! \! \! \! \! \tiny{\sf c}} 
-
\underline{ 
\frac{ 1}{ 2} \, \delta_{ l_2}^j \,\sum_k\,L_{l_1,l_3}^k\,\Theta^k
}_{ \fbox{\tiny 14}}
-
\underline{ 
\frac{ 1}{ 2} \, L_{l_1,l_3}^j\, \Theta^{l_2}
}_{ \octagon \! \! \! \! \! \tiny{\sf l}} 
-
\underline{ 
\frac{ 1}{ 2} \, L_{ l_3, l_3}^{ l_3} \, L_{l_2,l_1}^j
}_{ \octagon \! \! \! \! \! \tiny{\sf a}} 
+ \\
& \
\ \ \ \ \
+
\underline{ 
\frac{ 1}{ 4}\,\delta_{l_2}^j\,L_{l_3,l_3}^{l_3}\,L_{l_1,l_1}^{l_1}
}_{ \fbox{\tiny 3}}
+ 
\underline{ 
\frac{ 1}{ 4} \, \delta_{ l_1}^j \,L_{l_3,l_3}^{l_3}\,L_{l_2,l_2}^{l_2}
}_{ \octagon \! \! \! \! \! \tiny{\sf d}} 
+ 
\underline{ 
\frac{ 1}{ 4} \,\delta_{ l_2}^j \,L_{l_3,l_3}^{l_3}\,\Theta^{l_1}
}_{ \octagon \! \! \! \! \! \tiny{\sf q}} 
+
\underline{ 
\frac{ 1}{ 4} \,\delta_{ l_1}^j \,L_{l_3,l_3}^{l_3}\,\Theta^{l_2}
}_{ \octagon \! \! \! \! \! \! \tiny{\sf m}} 
- \\
& \
\ \ \ \ \
-
\underline{ 
\frac{ 1}{ 2} \, L_{l_1,l_1}^{l_1} \, L_{l_2,l_3}^j 
}_{ \octagon \! \! \! \! \! \tiny{\sf g}} 
+ 
\underline{ 
\frac{ 1}{ 4}\,\delta_{l_2}^j\,L_{l_1,l_1}^{l_1}\,L_{l_3,l_3}^{l_3}
}_{ \octagon \! \! \! \! \! \tiny{\sf i}} 
+ 
\underline{ 
\frac{ 1}{ 4}\,\delta_{l_3}^j\,L_{l_1,l_1}^{l_1}\,L_{l_2,l_2}^{l_2}
}_{ \octagon \! \! \! \! \! \tiny{\sf h}} 
+
\underline{ 
\frac{ 1}{ 4} \, \delta_{ l_2}^j \,L_{l_1,l_1}^{l_1}\,\Theta^{l_3} 
}_{ \fbox{\tiny 9}}
+ \\
& \
\ \ \ \ \
+
\underline{ 
\frac{ 1}{ 4} \, \delta_{ l_3}^j \,L_{l_1,l_1}^{l_1}\,\Theta^{l_2} 
}_{ \octagon \! \! \! \! \! \tiny{\sf j}} 
-
\underline{ 
\frac{ 1}{ 2} \, L_{l_2,l_1}^j \, \Theta^{l_3}
}_{ \octagon \! \! \! \! \! \tiny{\sf b}} 
+ 
\underline{ 
\frac{ 1}{ 4} \, \delta_{ l_2}^j \,L_{l_1,l_1}^{l_1}\,\Theta^{l_3}
}_{ \octagon \! \! \! \! \! \tiny{\sf k}} 
+ 
\underline{ 
\frac{ 1}{ 4} \, \delta_{ l_1}^j \,L_{l_2,l_2}^{l_2}\,\Theta^{l_3}
}_{ \octagon \! \! \! \! \! \! \tiny{\sf f}}
+ \\
& \
\ \ \ \ \
+ 
\underline{ 
\frac{ 1}{ 4} \, \delta_{ l_2}^j \,\Theta^{l_3}\,\Theta^{l_1}
}_{ \fbox{\tiny 17}}
+ 
\underline{ 
\frac{ 1}{ 4} \, \delta_{ l_1}^j \,\Theta^{l_3}\,\Theta^{l_2}
}_{ \octagon \! \! \! \! \! \tiny{\sf o}}
-
\underline{ 
\frac{ 1}{ 2} \,L_{l_2,l_3}^j\,\Theta^{l_1}
}_{ \octagon \! \! \! \! \! \tiny{\sf p}} 
+ 
\underline{ 
\frac{ 1}{ 4} \, \delta_{ l_2}^j \,L_{l_3,l_3}^{l_3}\,\Theta^{l_1}
}_{ \fbox{\tiny 11}}
+ \\
& \
\ \ \ \ \
+
\underline{ 
\frac{ 1}{ 4} \, \delta_{ l_3}^j \,L_{l_2,l_2}^{l_2}\,\Theta^{l_1}
}_{ \octagon \! \! \! \! \! \tiny{\sf e}} 
+ 
\underline{ 
\frac{ 1}{ 4} \, \delta_{ l_2}^j \,\Theta^{l_1}\,\Theta^{l_3}
}_{ \octagon \! \! \! \! \! \tiny{\sf r}} 
+ 
\underline{ 
\frac{ 1}{ 4} \,\delta_{ l_3}^j \,\Theta^{l_1}\,\Theta^{l_2}
}_{ \octagon \! \! \! \! \! \tiny{\sf n}}.
\endaligned
$$
We simplify and we reorganize the equality between the second and
third lines of~\thetag{ 3.88} and~\thetag{ 3.85} so as to put all
terms $\Theta_y$ in the left-hand side of the equality and to put all
remaining terms in the right-hand side, respecting
the order of \S3.73. We get:
\def\theequation{3.89}\begin{equation}
\aligned
&
\frac{ 1}{ 2} \, \delta_{ l_1}^j \, \Theta_{y^{l_3}}^{l_2}
-
\frac{ 1}{ 2} \, \delta_{ l_1}^j \, \Theta_{y^{l_2}}^{l_3}
+
\frac{ 1}{ 2} \, \delta_{ l_2}^j \, \Theta_{y^{l_3}}^{l_1}
-
\frac{ 1}{ 2} \, \delta_{ l_3}^j \, \Theta_{y^{l_2}}^{l_1}
= \\
& \
\ \ \ \ \
= 
\underline{
L_{l_1,l_2,y^{l_3}}^j 
- 
L_{l_1,l_3,y^{l_2}}^j 
+
\frac{ 1}{ 2} \, \delta_{ l_1}^j \, L_{l_3,l_3,y^{l_2}}^{l_3}
- 
\frac{ 1}{ 2} \, \delta_{ l_1}^j \, L_{l_2,l_2,y^{l_3}}^{l_2}
}
+ \\
& \
\ \ \ \ \ \ \ \ \ \ 
+
\underline{
\frac{ 1}{ 2} \, \delta_{ l_3}^j \, L_{l_1,l_1,y^{l_2}}^{l_1}
-
\frac{ 1}{ 2} \, \delta_{ l_2}^j \, L_{l_1,l_1,y^{l_3}}^{l_1}
}
+ \\
& \
\ \ \ \ \ \ \ \ \ \ 
+
\frac{ 1}{ 2} \, H_{l_3}^j \, M_{ l_1, l_2}
- 
\frac{ 1}{ 2} \, H_{l_2}^j \, M_{ l_1, l_3}
+
\frac{ 1}{ 4}\,\delta_{l_2}^j\,L_{l_3,l_3}^{l_3}\,L_{l_1,l_1}^{l_1}
-
\frac{ 1}{ 4}\,\delta_{l_3}^j\,L_{l_2,l_2}^{l_2}\,L_{l_1,l_1}^{l_1}
+ \\
& \
\ \ \ \ \ \ \ \ \ \
+
\sum_k\, L_{l_1, l_3}^k\, L_{l_2,k}^j 
- 
\sum_k\, L_{l_1, l_2}^k\, L_{l_3,k}^j
+ \\
& \
\ \ \ \ \ \ \ \ \ \
+
\frac{ 1}{ 2} \,\delta_{l_3}^j \,\sum_k \, L_{l_1,l_2}^k\,L_{k,k}^k
- 
\frac{ 1}{ 2} \,\delta_{l_2}^j \,\sum_k \, L_{l_1,l_3}^k\,L_{k,k}^k
+ \\
& \
\ \ \ \ \ \ \ \ \ \
+
\frac{ 1}{ 4} \, \delta_{ l_2}^j \,L_{l_1,l_1}^{l_1}\,\Theta^{l_3} 
-
\frac{ 1}{ 4} \, \delta_{ l_3}^j \,L_{l_1,l_1}^{l_1}\,\Theta^{l_2}
+
\frac{ 1}{ 4} \, \delta_{ l_2}^j \,L_{l_3,l_3}^{l_3}\,\Theta^{l_1} 
-
\frac{ 1}{ 4} \, \delta_{ l_3}^j \,L_{l_2,l_2}^{l_2}\,\Theta^{l_1}
+ \\
& \
\ \ \ \ \ \ \ \ \ \
+
\frac{ 1}{ 2} \, \delta_{ l_3}^j \,\sum_k\,L_{l_1,l_2}^k\,\Theta^k
-
\frac{ 1}{ 2} \, \delta_{ l_2}^j \,\sum_k\,L_{l_1,l_3}^k\,\Theta^k
+ \\
& \
\ \ \ \ \ \ \ \ \ \
+
\frac{ 1}{ 2} \, \delta_{ l_2}^j \,M_{l_1,l_3}\,\Theta^0 
-
\frac{ 1}{ 2} \, \delta_{ l_3}^j \,M_{l_1,l_2}\,\Theta^0 
+
\frac{ 1}{ 4} \, \delta_{ l_2}^j \,\Theta^{l_3}\,\Theta^{l_1}
-
\frac{ 1}{ 4} \, \delta_{ l_3}^j \,\Theta^{l_2}\,\Theta^{l_1}.
\endaligned
\end{equation}

Next, replacing plainly~\thetag{ 3.64} in $(3.65)_4$, we
get:
\def\theequation{3.90}\begin{equation}
\small
\aligned
&
\left(
\Pi_{0,0}^0
\right)_{y^{l_1}}
-
\left(
\Pi_{0,l_1}^0
\right)_x
= \\
& \
\ \ \ \ \
=
\Theta_{y^{ l_1}}^0 
- 
\frac{ 1}{ 2} \, L_{ l_1,l_1,x}^{l_1}
- 
\frac{ 1}{ 2} \, \Theta_x^{l_1}
= \\
& \
\ \ \ \ \
= 
-
\underline{
\Pi_{0,0}^0 \cdot \Pi_{l_1, 0}^0 
}_{ \octagon \! \! \! \! \! \tiny{\sf a}} 
- 
\sum_k\, \Pi_{0, 0}^k \cdot \Pi_{l_1, k}^0
+
\underline{
\Pi_{ 0, l_1}^0 \cdot \Pi_{0,0}^0 
}_{ \octagon \! \! \! \! \! \tiny{\sf a}} 
+ 
\sum_k\, \Pi_{0, l_1}^k \cdot \Pi_{ 0, k}^0 
= \\
& \
\ \ \ \ \
= 
- 
\sum_k \, \left( -G^k \right) 
\cdot
\left(
M_{ l_1, k}
\right)
+
\sum_k \, 
\left(
- 
\frac{ 1}{ 2} \, H_{l_1}^k 
+ 
\frac{ 1}{ 2} \, \delta_{ l_1}^k \, \Theta^0
\right) 
\cdot \\
& \
\ \ \ \ \ \ \ \ \ \
\cdot
\left(
\frac{ 1}{ 2} \, L_{ k,k}^k 
+ 
\frac{ 1}{ 2} \, \Theta^k
\right) 
= \\
& \
\ \ \ \ \
= 
\sum_k \, G^k \, M_{l_1,k} 
- 
\frac{ 1}{ 4} \, \sum_k \, H_{ l_1}^k \, L_{ k,k}^k 
- 
\frac{ 1}{ 4} \, \sum_k \, H_{l_1}^k \, \Theta^k 
+ 
\frac{ 1}{ 4} \, L_{l_1,l_1}^{l_1} \, \Theta^0 
+ 
\frac{ 1}{ 4} \, \Theta^0 \, \Theta^{l_1}.
\endaligned
\end{equation}
Reorganizing the equality so as to put the terms $\Theta_x$ and
$\Theta_y$ alone in the left-hand side, we get:
\def\theequation{3.91}\begin{equation}
\aligned
- 
\frac{ 1}{ 2} \, \Theta_x^{l_1}
+ 
\Theta_{ y^{l_1}}^0 
& \
= 
\underline{
\frac{ 1}{ 2} \, L_{l_1,l_1,x}^{l_1}
} 
+ \\
& \
\ \ \ \ \
+
\sum_k\, G^k \, M_{l_1,k}
- 
\frac{ 1}{ 4} \, \sum_k \, H_{ l_1}^k \, L_{ k,k}^k 
- 
\frac{ 1}{ 4} \, \sum_k \, H_{l_1}^k \, \Theta^k 
+ \\
& \
\ \ \ \ \
+
\frac{ 1}{ 4} \, L_{ l_1, l_1}^{l_1} \, \Theta^0 
+
\frac{ 1}{ 4} \, \Theta^0 \, \Theta^{l_1}.
\endaligned
\end{equation}

Next, replacing plainly~\thetag{ 3.64} in $(3.65)_5$, we
get:
\def\theequation{3.92}\begin{equation}
\small
\aligned
& 
\left(
\Pi_{l_1,l_2}^0
\right)_x
- 
\left(
\Pi_{l_1, 0}^0 
\right)_{y^{l_2}}
= \\
& \
\ \ \ \ \
= 
M_{l_1,l_2,x}
- 
\frac{ 1}{ 2} \, L_{ l_1, l_1, y^{l_2}}^{l_1}
- 
\frac{ 1}{ 2} \, \Theta_{y^{l_2}}^{l_1}
= \\
& \
\ \ \ \ \
=
- 
\Pi_{l_1, l_2} \cdot \Pi_{0, 0}^0
- 
\sum_k \, \Pi_{l_1, l_2}^k \cdot \Pi_{0, k}^0
+ 
\Pi_{ l_1, 0}^0 \cdot \Pi_{l_2, 0}^0 
+ 
\sum_k\, \Pi_{l_1,0}^k \cdot \Pi_{l_2, k}^0
= \\
& \
\ \ \ \ \
= 
- 
M_{l_1, l_2} \, \Theta^0 
- 
\sum_k\, 
\left(
- 
L_{l_1,l_2}^k 
+
\frac{ 1}{ 2} \, \delta_{ l_1}^k \, L_{l_2, l_2}^{l_2}
+
\frac{ 1}{ 2} \, \delta_{ l_2}^k \, L_{l_1, l_1}^{l_1} +
\right. \\
& \
\ \ \ \ \ \ \ \ \ \
\left.
+
\frac{ 1}{ 2} \, \delta_{ l_1}^k \, \Theta^{l_2}
+
\frac{ 1}{ 2} \, \delta_{ l_2}^k \, \Theta^{l_1}
\right)
\cdot
\left(
\frac{ 1}{ 2} \, L_{k,k}^k 
+
\frac{ 1}{ 2} \, \Theta^k
\right)
+ \\
& \
\ \ \ \ \ \ \ \ \ \
+
\left(
\frac{ 1}{ 2} \, L_{ l_1, l_1}^{l_1} 
+
\frac{ 1}{ 2} \, \Theta^{l_1}
\right)
\cdot
\left(
\frac{ 1}{ 2} \, L_{l_2,l_2}^{l_2} 
+
\frac{ 1}{ 2} \, \Theta^{l_2}
\right)
+ \\
& \
\ \ \ \ \ \ \ \ \ \
+
\sum_k\, 
\left(
-
\frac{ 1}{ 2} \, H_{l_1}^k 
+
\frac{ 1}{ 2} \, \delta_{ l_1}^k \, \Theta^0 
\right)
\cdot
M_{l_2,k} =
\endaligned
\end{equation}
$$
\small
\aligned
& \
\ \ \ \ \
=
- 
\underline{ 
M_{l_1,l_2} \, \Theta^0 
}_{ \fbox{\tiny 7}}
+ 
\underline{ 
\frac{ 1}{ 2} \, \sum_k\,L_{l_1,l_2}^k\, L_{k,k}^k
}_{ \fbox{\tiny 3}}
+
\underline{ 
\frac{ 1}{ 2} \, \sum_k\, L_{ l_1, l_2}^k \, \Theta^k
}_{ \fbox{\tiny 6}}
- 
\underline{ 
\frac{ 1}{ 4} \, L_{l_2,l_2}^{l_2}\,L_{l_1,l_1}^{l_1}
}_{ \fbox{\tiny 2}}
- \\
& \
\ \ \ \ \ \ \ \ \ \
- 
\underline{ 
\frac{ 1}{ 4} \, L_{ l_2,l_2}^{l_2} \, \Theta^{l_1}
}_{ \octagon \! \! \! \! \! \tiny{\sf a}} 
- 
\underline{ 
\frac{ 1}{ 4} \, L_{l_1,l_1}^{l_1} \, L_{l_2,l_2}^{l_2}
}_{ \octagon \! \! \! \! \! \tiny{\sf b}} 
-
\underline{ 
\frac{ 1}{ 4} \, L_{ l_1, l_1}^{l_1} \, \Theta^{l_2}
}_{ \octagon \! \! \! \! \! \tiny{\sf c}} 
- 
\underline{ 
\frac{ 1}{ 4} \, L_{ l_1, l_1}^{l_1} \, \Theta^{l_2}
}_{ \fbox{\tiny 4}} 
- \\
& \
\ \ \ \ \ \ \ \ \ \
- 
\underline{ 
\frac{ 1}{ 4} \, \Theta^{l_1} \, \Theta^{l_2}
}_{ \octagon \! \! \! \! \! \tiny{\sf d}} 
-
\underline{ 
\frac{ 1}{ 4} \, L_{ l_2, l_2}^{l_2} \, \Theta^{l_1}
}_{ \fbox{\tiny 5}}
-
\underline{ 
\frac{ 1}{ 4} \, \Theta^{l_1} \, \Theta^{l_2}
}_{ \fbox{\tiny 8}}
+
\underline{ 
\frac{ 1}{ 4} \, L_{l_1,l_1}^{l_1} \, L_{l_2,l_2}^{l_2}
}_{ \octagon \! \! \! \! \! \tiny{\sf b}}
+ \\
& \
\ \ \ \ \ \ \ \ \ \
+
\underline{ 
\frac{ 1}{ 4} \, L_{ l_1, l_1}^{l_1} \, \Theta^{l_2}
}_{ \octagon \! \! \! \! \! \tiny{\sf c}} 
+
\underline{ 
\frac{ 1}{ 4} \, L_{ l_2, l_2}^{l_2} \, \Theta^{l_1}
}_{ \octagon \! \! \! \! \! \tiny{\sf a}}
+
\underline{ 
\frac{ 1}{ 4} \, \Theta^{l_1} \, \Theta^{l_2}
}_{ \octagon \! \! \! \! \! \tiny{\sf d}}
-
\underline{
\frac{ 1}{ 2} \, \sum_k \, H_{ l_1}^k \, M_{l_2, k}
}_{ \fbox{\tiny 1}} 
+ \\
& \
\ \ \ \ \ \ \ \ \ \
+
\underline{ 
\frac{ 1}{ 2} \, M_{ l_2, l_1} \, \Theta^0
}_{ \fbox{\tiny 7}}.
\endaligned
$$
Multiplying by $-2$ and reorganizing the
equality, we get:
\def\theequation{3.93}\begin{equation}
\aligned
\Theta_{ y^{l_2}}^{l_1}
& 
=
-
\underline{
L_{l_1, l_1, y^{l_2}}^{l_1}
+ 
2\, M_{l_1, l_2, x}}
+ \\
& \
\ \ \ \ \
+
\sum_k\, H_{ l_1}^k \, M_{l_2,k}
+ 
\frac{ 1}{ 2} \, L_{l_1,l_1}^{l_1} \, L_{l_2,l_2}^{l_2}
- 
\sum_k\, L_{l_1,l_2}^k\, L_{k,k}^k 
+ \\
& \
\ \ \ \ \
+
\frac{ 1}{ 2} \, L_{l_1, l_1}^{l_1} \, \Theta^{l_2}
+
\frac{ 1}{ 2} \, L_{l_2, l_2}^{l_2} \, \Theta^{l_1}
- 
\sum_k\, L_{l_1, l_2}^k \, \Theta^k 
+ \\
& \
\ \ \ \ \
+
M_{ l_1, l_2} \, \Theta^0 
+ 
\frac{ 1}{ 2} \,
\Theta^{l_1} \, \Theta^{l_2}.
\endaligned
\end{equation}

Next, replacing plainly~\thetag{ 3.64} in $(3.65)_4$, we
get:
\def\theequation{3.94}\begin{equation}
\small
\aligned
&
\left(
\Pi_{l_1,l_2}^0 
\right)_{y^{l_3}}
-
\left(
\Pi_{l_1,l_3}^0 
\right)_{y^{l_2}}
= \\
& \
\ \ \ \ \
= 
M_{ l_1, l_2, y^{l_3}} 
- 
M_{ l_1, l_3, y^{l_2}} 
= \\ 
& \
\ \ \ \ \
= 
- 
\Pi_{l_1,l_2}^0 \cdot \Pi_{l_3, 0}^0 
- 
\sum_k\, \Pi_{l_1, l_2}^k \cdot \Pi_{l_3, k}^0 
+ 
\Pi_{l_1, l_3}^0 \cdot \Pi_{l_2, 0}^0 
+ 
\sum_k \, \Pi_{l_1, l_3}^k \cdot \Pi_{ l_2, k}^0 
= \\
& \
\ \ \ \ \
= 
- 
M_{l_1, l_2} 
\left(
\frac{ 1}{ 2} \, L_{l_3, l_3}^{l_3} 
+
\frac{ 1}{ 2} \, \Theta^{l_3}
\right)
- 
\sum_k\, M_{l_3, k}
\left(
- 
L_{l_1, l_2}^k 
+ 
\frac{ 1}{ 2} \, \delta_{l_1}^k \, L_{l_2, l_2}^{l_2}
+
\right. \\
& \
\ \ \ \ \ \ \ \ \ \
\left.
+
\frac{ 1}{ 2} \, \delta_{ l_2}^k \, L_{l_1, l_1}^{l_1}
+ 
\frac{ 1}{ 2} \, \delta_{ l_1}^k \, \Theta^{l_2}
+
\frac{ 1}{ 2} \, \delta_{ l_2}^k \, \Theta^{l_1}
\right) 
+ \\
& \
\ \ \ \ \ \ \ \ \ \
+
M_{l_1, l_3} 
\left(
\frac{ 1}{ 2} \, L_{l_2, l_2}^{l_2} 
+
\frac{ 1}{ 2} \, \Theta^{l_2}
\right)
+
\sum_k\, M_{l_2, k}
\left(
-
L_{l_1, l_3}^k 
+
\frac{ 1}{ 2} \, \delta_{l_1}^k \, L_{l_3, l_3}^{l_3}
+
\right. \\
& \
\ \ \ \ \ \ \ \ \ \
\left.
+
\frac{ 1}{ 2} \, \delta_{ l_3}^k \, L_{l_1, l_1}^{l_1}
+
\frac{ 1}{ 2} \, \delta_{ l_1}^k \, \Theta^{l_3}
+
\frac{ 1}{ 2} \, \delta_{ l_3}^k \, \Theta^{l_1}
\right).
\endaligned
\end{equation}
Developing the products and
ordering each monomial, we get:
\def\theequation{3.95}\begin{equation}
\small
\aligned
& 
= 
- 
\underline{ 
\frac{ 1}{ 2} \, L_{l_3, l_3}^{ l_3} \, M_{ l_1, l_2}
}_{ \octagon \! \! \! \! \! \tiny{\sf a}} 
- 
\underline{ 
\frac{ 1}{ 2} \, M_{ l_1, l_2} \, \Theta^{ l_3}
}_{ \octagon \! \! \! \! \! \tiny{\sf b}}
+ 
\underline{ 
\sum_k\, L_{l_1, l_2}^k \, M_{ l_3, k}
}_{ \fbox{\tiny 1}}
- 
\underline{ 
\frac{ 1}{ 2} \, L_{ l_2, l_2}^{ l_2} \, M_{ l_3, l_1}
}_{ \octagon \! \! \! \! \! \tiny{\sf c}} 
- \\
& \
\ \ \ \ \
- 
\underline{ 
\frac{ 1}{ 2} \, L_{ l_1, l_1}^{ l_1} \, M_{ l_3, l_2}
}_{ \octagon \! \! \! \! \! \tiny{\sf d}} 
-
\underline{ 
\frac{ 1}{ 2} \, M_{ l_3, l_1} \, \Theta^{ l_2}
}_{ \octagon \! \! \! \! \! \tiny{\sf e}}
- 
\underline{ 
\frac{ 1}{ 2} \, M_{ l_3, l_2} \, \Theta^{ l_1}
}_{ \octagon \! \! \! \! \! \tiny{\sf f}} 
+ 
\underline{ 
\frac{ 1}{ 2} \, L_{ l_2, l_2}^{ l_2} \, M_{ l_3, l_1}
}_{ \octagon \! \! \! \! \! \tiny{\sf c}} 
+ \\ 
& \
\ \ \ \ \
+
\underline{ 
\frac{ 1}{ 2} \, M_{ l_1, l_3} \, \Theta^{ l_2}
}_{ \octagon \! \! \! \! \! \tiny{\sf e}} 
- 
\underline{ 
\sum_k \, L_{ l_1, l_3}^k \, M_{ l_2, k}
}_{ \fbox{\tiny 2}}
+ 
\underline{ 
\frac{ 1}{ 2} \, L_{l_3, l_3}^{ l_3} \, M_{ l_2, l_1}
}_{ \octagon \! \! \! \! \! \tiny{\sf a}} 
+ 
\underline{ 
\frac{ 1}{ 2} \, L_{ l_1, l_1}^{ l_1} \, M_{ l_2, l_3}
}_{ \octagon \! \! \! \! \! \tiny{\sf d}}
+ \\
& \
\ \ \ \ \
+
\underline{ 
\frac{ 1}{ 2} \, M_{ l_2, l_1} \, \Theta^{ l_3}
}_{ \octagon \! \! \! \! \! \tiny{\sf b}}
+ 
\underline{ 
\frac{ 1}{ 2} \, M_{ l_2, l_3} \, \Theta^{ l_1}
}_{ \octagon \! \! \! \! \! \tiny{\sf f}}.
\endaligned
\end{equation}
Simplifying, we obtain the family (IV) in the
statement of Theorem~1.7~{\bf (3)}:
\def\theequation{3.96}\begin{equation}
0 = 
\underline{
M_{l_1, l_2, y^{ l_3}} 
- 
M_{l_1, l_3, y^{ l_2}}
}
- 
\sum_k \, L_{ l_1, l_2}^k \, M_{ l_3, k}
+ 
\sum_k \, L_{ l_1, l_3}^k \, M_{ l_2, k}.
\end{equation}

\subsection*{ 3.97.~Solving $\Theta_x^0$, $\Theta_{ y^{ l_1}}^0$, 
$\Theta_x^{ l_1}$ and $\Theta_{ y^{ l_2}}^{l_1}$} It is now easy to
solve all first order partial derivatives of the functions $\Theta^0$
and $\Theta^l$.
Equation~\thetag{ 3.93} already 
provides the solution for $\Theta_{ y^{ l_2}}^{ l_1}$.
We state the result as an independent proposition.

\def\theproposition{3.98}\begin{proposition}
As a consequence of the six families of equations~\thetag{ 3.69},
\thetag{ 3.86}, \thetag{ 3.89}, \thetag{ 3.91}, \thetag{ 3.93} and
\thetag{ 3.96} the first order derivatives $\Theta_x^0$, $\Theta_{ y^{
l_1}}^0$, $\Theta_x^{ l_1}$ and $\Theta_{ y^{ l_2}}^{l_1}$ of the
principal unknowns are given by{\rm :}
\def\theequation{3.99}\begin{equation}
\left\{
\aligned
\Theta_x^0 
& 
=
-
\underline{ 
2 \, G_{ y^{ l_1}}^{ l_1} 
+
H_{ l_1, x}^{ l_1}
} 
+ \\
& \
\ \ \ \ \
+ 
2\, \sum_k\, G^k\, L_{ l_1,k}^{l_1}
- 
\sum_k\, G^k\, L_{k,k}^k
- 
\frac{ 1}{ 2} \,
\sum_k\, H_{ l_1}^k\, H_k^{l_1}
- \\
& \
\ \ \ \ \
-
\sum_k\, G^k \, \Theta^k 
+ 
\frac{ 1}{ 2} \,
\Theta^0 \, \Theta^0.
\endaligned\right.
\end{equation}
\def\theequation{3.100}\begin{equation}
\left\{
\aligned
\Theta_{y^{ l_1}}^0
& 
=
\underline{
\frac{ 2}{ 3} \, L_{ l_1, l_1, x}^{ l_1}
- 
\frac{ 1}{ 3} \, H_{ l_1, y^{ l_1}}^{l_1}
} 
+ \\
& \
\ \ \ \ \
+
\frac{ 2}{ 3} \, G^{ l_1} \, M_{ l_1, l_1}
+ 
\frac{ 4}{ 3} \, \sum_k\, G^k\, M_{ l_1, k}
- 
\frac{ 1}{ 3} \, \sum_k\, H_k^{l_1} \, L_{ l_1, l_1}^k
+ \\
& \
\ \ \ \ \
+ 
\frac{ 1}{ 3} \, \sum_k\, H_{ l_1}^k \, L_{ l_1, k}^{l_1}
- 
\frac{ 1}{ 2} \, \sum_k \, H_{l_1}^k L_{k,k}^k
- 
\frac{ 1}{ 2} \,\sum_k \, H_{l_1}^k\, \Theta^k 
+ \\
& \
\ \ \ \ \
+ 
\frac{ 1}{ 2} \, L_{ l_1, l_1}^{ l_1} \, \Theta^0
+ 
\frac{ 1}{ 2} \,
\Theta^0 \, \Theta^{ l_1}.
\endaligned\right.
\end{equation}
\def\theequation{3.101}\begin{equation}
\left\{
\aligned
\Theta_x^{l_1}
& 
=
-
\underline{
\frac{ 2}{ 3} \, H_{ l_1, y^{ l_1}}^{ l_1}
+ 
\frac{ 1}{ 3} \, L_{l_1, l_1, x}^{ l_1}
}
+ \\
& \
\ \ \ \ \
+ 
\frac{ 4}{ 3} \, G^{ l_1} \, M_{ l_1, l_1}
+ 
\frac{ 2}{ 3} \, \sum_k \, G^k \, M_{ l_1,k} 
- 
\frac{ 2}{ 3} \, \sum_k \, H_k^{ l_1} \, L_{ l_1, l_1}^k 
+ \\
& \
\ \ \ \ \
+ 
\frac{ 2}{ 3} \, \sum_k \, H_{ l_1}^k \, L_{ l_1, k}^{ l_1}
- 
\frac{ 1}{ 2} \, \sum_k \, H_{ l_1}^k \, L_{ k,k}^k 
- 
\frac{ 1}{ 2} \, \sum_k \, H_{ l_1}^k \, \Theta^k
+ \\
& \
\ \ \ \ \
+
\frac{ 1}{ 2} \, L_{l_1, l_1}^{ l_1} \, \Theta^0 
+ 
\frac{ 1}{ 2} \, \Theta^0 \, \Theta^{ l_1}.
\endaligned\right.
\end{equation}
\def\theequation{3.102}\begin{equation}
\left\{
\aligned
\Theta_{ y^{l_2}}^{l_1}
& 
=
-
\underline{
L_{l_1, l_1, y^{l_2}}^{l_1}
+ 
2\, M_{l_1, l_2, x}}
+ \\
& \
\ \ \ \ \
+
\sum_k\, H_{ l_1}^k \, M_{l_2,k}
+ 
\frac{ 1}{ 2} \, L_{l_1,l_1}^{l_1} \, L_{l_2,l_2}^{l_2}
- 
\sum_k\, L_{l_1,l_2}^k\, L_{k,k}^k 
+ \\
& \
\ \ \ \ \
+
\frac{ 1}{ 2} \, L_{l_1, l_1}^{l_1} \, \Theta^{l_2}
+
\frac{ 1}{ 2} \, L_{l_2, l_2}^{l_2} \, \Theta^{l_1}
- 
\sum_k\, L_{l_1, l_2}^k \, \Theta^k 
+ \\
& \
\ \ \ \ \
+
M_{ l_1, l_2} \, \Theta^0 
+ 
\frac{ 1}{ 2} \,
\Theta^{l_1} \, \Theta^{l_2}.
\endaligned\right.
\end{equation}
\end{proposition}

We notice that the right-hand side
of~\thetag{ 3.99} should be independent
of $l_1$; this phenomenon will be explained in a while.

\proof
For $\Theta_x^0$ in~\thetag{ 3.99}, it suffices to 
put $j:= l_1$ in~\thetag{ 3.69}. 

To obtain $\Theta_{ y^{ l_1}}^0$, we put 
$j:= l_2$ and $l_2 := l_1$ in~\thetag{ 3.86}, which 
yields:
\def\theequation{3.103}\begin{equation}
\small
\aligned
\Theta_x^{ l_1} 
- 
\frac{ 1}{ 2} \, \Theta_{ y^{ l_1}}^0 
& = 
- 
\underline{
\frac{ 1}{ 2} \, H_{ l_1, y^{ l_1}}^{ l_1}
}
+ \\
& \
\ \ \ \ \
+ 
G^{ l_1}\, M_{ l_1, l_1} 
- 
\frac{ 1}{ 2} \, \sum_k\, H_k^{ l_1} \, L_{ l_1, l_1}^k 
+ \\
& \
\ \ \ \ \
+ 
\frac{ 1}{ 2} \, \sum_k\, H_{ l_1}^k \, L_{ l_1, k}^{ l_1}
- 
\frac{ 1}{ 4} \, \sum_k \, H_{ l_1}^k \, L_{ k,k}^k 
- \\
& \
\ \ \ \ \
- 
\frac{ 1}{ 4} \, \sum_k \, H_{ l_1}^k \, \Theta^k 
+ 
\frac{ 1}{ 4} \, L_{ l_1, l_1}^{ l_1} \, \Theta^0 
+ 
\frac{ 1}{ 4} \, \Theta^0 \, \Theta^{ l_1}.
\endaligned
\end{equation}
We may easily 
solve $\Theta_{ y^{ l_1}}^0$ and $\Theta_x^{ l_1}$ thanks to this
equation~\thetag{ 3.103} and thanks to~\thetag{ 3.91}: indeed, to
obtain~\thetag{ 3.100}, it suffices to compute $\frac{ 4}{ 3} \cdot
(3.91) + \frac{ 2}{ 3} \cdot (3.103)$; to obtain~\thetag{ 3.101}, it
suffices to compute $\frac{ 2}{ 3} \cdot (3.91) + \frac{ 4}{ 3} \cdot
(3.103)$. Finally,
\thetag{ 3.102} is a copy of~\thetag{ 3.93}. 
This completes the proof.
\endproof

\subsection*{3.104.~Appearance 
of the crucial four families of first order partial differential
relations (I), (II), (III) and (IV) of Theorem~1.7~{\bf (3)}} However, in
solving $\Theta_x^0$, $\Theta_{ y^{ l_1}}^0$, $\Theta_x^{ l_1}$ and
$\Theta_{ y^{ l_2}}^{l_1}$ from our six families of equations~\thetag{
3.69}, \thetag{ 3.86}, \thetag{ 3.89}, \thetag{ 3.91}, \thetag{ 3.93}
and \thetag{ 3.96}, only a subpart of these equations has been used.
We notice that the two families of equations~\thetag{ 3.91}
and~\thetag{ 3.93} have been used completely and that the family of
equations~\thetag{ 3.96}, which does not involve $\Theta$, coincides
precisely with the system (IV) of Theorem~1.7~{\bf (3)}. To insure that
$\Theta_x^0$, $\Theta_{ y^{ l_1}}^0$, $\Theta_x^{ l_1}$ and $\Theta_{
y^{ l_2}}^{l_1}$ as written in Proposition~3.98 are true solutions, it
is necessary and sufficient that they satisfy the remaining equations.
Thus, we have to replace these solutions~\thetag{ 3.99}, \thetag{
3.100}, \thetag{ 3.101} and~\thetag{ 3.102} in the three remaining
families~\thetag{ 3.69}, \thetag{ 3.86} and \thetag{ 3.89}.

Firstly, let us insert inside~\thetag{ 3.69} the value
of $\Theta_x^0$ given by the equation~\thetag{ 3.99}, in which
the index $l_1$ is replaced in advance by 
an arbitrary index $l_2$. We get:
\def\theequation{3.105}\begin{equation}
\small
\aligned
0 
& =
- 
\underline{
2\, G_{ y^{ l_1}}^j 
+ 
2\, \delta_{ l_1}^j \, G_{ y^{ l_2}}^{ l_2}
+ 
H_{ l_1, x}^j 
- 
\delta_{ l_1}^j \, H_{ l_2, x}^{ l_2}
}
+ \\
& \
\ \ \ \ \
+ 
2\, \sum_k\, G^k \, L_{ l_1, k}^j
- 
2\, \delta_{ l_1}^j \, \sum_k\, G^k \, L_{ l_2, k}^{ l_2}
- 
\underline{ 
\delta_{ l_1}^j \, \sum_k\, G^k \, L_{ k,k}^k
}_{ \octagon \! \! \! \! \! \tiny{\sf a}} 
+ \\
& \
\ \ \ \ \
+ 
\underline{
\delta_{ l_1}^j \, \sum_k\, G^k \, L_{ k,k}^k 
}_{ \octagon \! \! \! \! \! \tiny{\sf a}} 
- 
\frac{ 1}{ 2} \, \sum_k\, H_{ l_1}^k \, H_k^j 
+ 
\frac{ 1}{ 2} \, \delta_{ l_1}^j \, \sum_k\, H_{ l_2}^k\,H_k^{ l_2}
- \\
& \
\ \ \ \ \
- 
\underline{
\delta_{ l_1}^j \, \sum_k\, G^k \, \Theta^k 
}_{ \octagon \! \! \! \! \! \tiny{\sf b}} 
+ 
\underline{
\delta_{ l_1}^j \, \sum_k\, G^k \, \Theta^k 
}_{ \octagon \! \! \! \! \! \tiny{\sf b}} 
+ 
\underline{
\frac{ 1}{ 2} \, \delta_{ l_1}^j \, \Theta^0 \, \Theta^0
}_{ \octagon \! \! \! \! \! \tiny{\sf c}} 
- 
\underline{
\frac{ 1}{ 2} \, \delta_{ l_1}^j \, \Theta^0 \, \Theta^0
}_{ \octagon \! \! \! \! \! \tiny{\sf c}} \ .
\endaligned
\end{equation}
We simplify, which yields the 
family (I) of partial differential relations of
Theorem~1.7~{\bf (3)}:
\def\theequation{3.106}\begin{equation}
\aligned
0= 
&\
-
\underline{
2\, G_{y^{l_1}}^j
+
2\, \delta_{l_1}^j\, G_{y^{l_2}}^{l_2}
+
H_{l_1,x}^j
-
\delta_{l_1}^j\, H_{l_2, x}^{l_2}
}
+ \\
& \
+
2\, \sum_k \, G^k\, L_{l_1,k}^j
-
2\, \delta_{l_1}^j\, \sum_k \, 
G^k\, L_{l_2,k}^{l_2}
+ \\
& \
-
\frac{ 1}{2}\, \sum_k \, H_{l_1}^k\, H_k^j
+
\frac{ 1}{2}\, \delta_{l_1}^j\, \sum_k \, H_{l_2}^k\, H_k^{l_2}.
\endaligned
\end{equation}

Secondly, let us insert inside~\thetag{ 3.86} the values of
$\Theta_x^{ l_1}$, $\Theta_x^{ l_2}$ given by~\thetag{ 3.101} and the
value of $\Theta_{ y^{ l_1}}^0$ given by~\thetag{ 3.100}. We place all
the terms in the right-hand side of the equality and we place the
first order terms in the beginning (first
three lines just below). We obtain:
\def\theequation{3.107}\begin{equation}
\small
\aligned
0 &
=
- 
\underline{ 
\frac{ 1}{ 2} \, H_{ l_1, y^{ l_2}}^j 
}_{ \fbox{\tiny 1}}
+ 
\underline{ 
L_{ l_1, l_2, x}^j
}_{ \fbox{\tiny 4}}
- 
\underline{ 
\frac{ 1}{ 2} \, \delta_{ l_1}^j \, L_{ l_2, l_2, x}^{ l_2}
}_{ \fbox{\tiny 5}}
- 
\underline{ 
\frac{ 1}{ 2} \, \delta_{ l_1}^j \,L_{ l_1, l_1, x}^{ l_1}
}_{ \fbox{\tiny 6}}
+ \\
& \
\ \ \ \ \
+
\underline{ 
\frac{ 1}{ 3} \, \delta_{ l_1}^j \,H_{l_2, y^{l_2}}^{l_2}
}_{ \fbox{\tiny 2}}
- 
\underline{ 
\frac{ 1}{ 6} \,\delta_{ l_1}^j \,L_{l_2,l_2,x}^{l_2}
}_{ \fbox{\tiny 5}}
+
\underline{ 
\frac{ 1}{ 3} \,\delta_{ l_2}^j \, H_{l_1,y^{l_1}}^{l_1}
}_{ \fbox{\tiny 3}}
- 
\underline{ 
\frac{ 1}{ 6} \, \delta_{ l_2}^j \,L_{l_1,l_1,x}^{l_1}
}_{ \fbox{\tiny 6}}
+ \\
& \
\ \ \ \ \
+ 
\underline{ 
\frac{ 1}{ 3} \, \delta_{ l_1}^j \,L_{l_2,l_2,x}^{l_2}
}_{ \fbox{\tiny 5}}
- 
\underline{ 
\frac{ 1}{ 6} \, \delta_{ l_1}^j \,H_{l_2, y^{l_2}}^{l_2}
}_{ \fbox{\tiny 2}}
+ \\
& \
\ \ \ \ \
+ 
\underline{ 
G^j M_{ l_1, l_2}
}_{ \fbox{\tiny 7}}
- 
\underline{ 
\frac{ 1}{ 2} \, \sum_k\, H_k^j \, L_{ l_1, l_2}^k
}_{ \fbox{\tiny 12}}
+ 
\underline{ 
\frac{ 1}{ 2} \, \sum_k\, H_{l_1}^k \, L_{ l_2, k}^j
}_{ \fbox{\tiny 13}}
- 
\underline{ 
\frac{ 1}{ 4} \, \delta_{ l_2}^j \,\sum_k\,H_{l_1}^k\,L_{k,k}^k
}_{ \octagon \! \! \! \! \! \tiny{\sf a}} 
- \\
& \
\ \ \ \ \
-
\underline{ 
\frac{ 1}{ 4} \, \delta_{ l_2}^j \,\sum_k\,H_{l_1}^k\,\Theta^k
}_{ \octagon \! \! \! \! \! \tiny{\sf b}} 
+
\underline{ 
\frac{ 1}{ 4} \, \delta_{ l_2}^j \, L_{ l_1, l_1}^{ l_1} \, \Theta^0 
}_{ \octagon \! \! \! \! \! \tiny{\sf c}} 
+ 
\underline{ 
\frac{ 1}{ 4} \, \delta_{ l_2}^j \, \Theta^0 \, \Theta^{l_1} 
}_{ \octagon \! \! \! \! \! \tiny{\sf d}}
- \\
& \
\ \ \ \ \
-
\underline{ 
\frac{ 2}{3} \, \delta_{ l_1}^j \, G^{ l_2} \, M_{ l_2, l_2}
}_{ \fbox{\tiny 8}} 
- 
\underline{ 
\frac{ 1}{ 3} \,\delta_{ l_1}^j \,\sum_k\,G^k\,M_{l_2,k}
}_{ \fbox{\tiny 10}}
+ 
\underline{ 
\frac{ 1}{ 3} \,\delta_{ l_1}^j \,\sum_k\,H_k^{l_2}\,L_{l_2,l_2}^k
}_{ \fbox{\tiny 14}}
- 
\underline{ 
\frac{ 1}{ 3} \, \delta_{ l_1}^j \,\sum_k\, H_{l_2}^k\,L_{l_2,k}^{l_2}
}_{ \fbox{\tiny 15}}
+
\endaligned
\end{equation}
$$
\small
\aligned
& \
\ \ \ \ \
+
\underline{ 
\frac{ 1}{ 4} \,\delta_{ l_1}^j \,\sum_k\,H_{l_2}^k\,L_{k,k}^k
}_{ \octagon \! \! \! \! \! \tiny{\sf e}} 
+
\underline{ 
\frac{ 1}{ 4} \,\delta_{ l_1}^j \,\sum_k\,H_{l_2}^k\,\Theta^k
}_{ \octagon \! \! \! \! \! \tiny{\sf f}} 
-
\underline{ 
\frac{ 1}{ 4} \,\delta_{ l_1}^j \,L_{l_2,l_2}^{l_2}\, \Theta^0
}_{ \octagon \! \! \! \! \! \tiny{\sf g}} 
- 
\underline{ 
\frac{ 1}{ 4} \,\delta_{ l_1}^j \,\Theta^0 \,\Theta^{ l_2}
}_{ \octagon \! \! \! \! \! \tiny{\sf h}} 
- \\
& \
\ \ \ \ \ 
-
\underline{ 
\frac{ 2}{ 3} \, \delta_{ l_2}^j \,G^{l_1}\, M_{l_1,l_1}
}_{ \fbox{\tiny 9}}
- 
\underline{ 
\frac{ 1}{ 3} \,\delta_{ l_2}^j \,\sum_k\,G^k\, M_{l_1,k}
}_{ \fbox{\tiny 11}}
+ 
\underline{ 
\frac{ 1}{ 3} \,\delta_{ l_2}^j \,\sum_k\,H_k^{l_1}L_{l_1,l_1}^k
}_{ \fbox{\tiny 16}}
- 
\underline{ 
\frac{ 1}{ 3} \,\delta_{ l_2}^j \,\sum_k\,H_{l_1}^k\,L_{l_1,k}^{l_1}
}_{ \fbox{\tiny 17}}
+ \\
& \
\ \ \ \ \
+
\underline{ 
\frac{ 1}{ 4} \,\delta_{ l_2}^j \,\sum_k\,H_{l_1}^k\,L_{k,k}^k
}_{ \octagon \! \! \! \! \! \tiny{\sf a}} 
+
\underline{ 
\frac{ 1}{ 4} \,\delta_{ l_2}^j \,\sum_k\, H_{l_1}^k\, \Theta^k
}_{ \octagon \! \! \! \! \! \tiny{\sf b}}
- 
\underline{ 
\frac{ 1}{ 4} \, \delta_{ l_2}^j \, L_{ l_1, l_1}^{ l_1} \, \Theta^0 
}_{ \octagon \! \! \! \! \! \tiny{\sf c}}
-
\underline{ 
\frac{ 1}{ 4} \, \delta_{ l_2}^j \, \Theta^0 \, \Theta^0 
}_{ \octagon \! \! \! \! \! \tiny{\sf d}}
+ \\
& \
\ \ \ \ \
+
\underline{ 
\frac{ 1}{ 3} \,\delta_{ l_1}^j \,G^{l_2}\, M_{l_2,l_2}
}_{ \fbox{\tiny 8}}
+ 
\underline{ 
\frac{ 2}{ 3} \,\delta_{ l_1}^j \,\sum_k\,G^k\,M_{l_2,k}
}_{ \fbox{\tiny 10}}
- 
\underline{ 
\frac{ 1}{ 6} \,\delta_{ l_1}^j \,\sum_k\,H_k^{l_2}\,L_{l_2,l_2}^k
}_{ \fbox{\tiny 14}}
+
\underline{ 
\frac{ 1}{ 6} \,\delta_{ l_1}^j \,\sum_k\, H_{l_2}^k\, L_{l_2,k}^{l_2}
}_{ \fbox{\tiny 15}} 
- \\
& \
\ \ \ \ \
-
\underline{ 
\frac{ 1}{ 4} \,\delta_{ l_1}^j \,\sum_k\,H_{l_2}^k\,L_{k,k}^k
}_{ \octagon \! \! \! \! \! \tiny{\sf e}} 
- 
\underline{ 
\frac{ 1}{ 4} \,\delta_{ l_1}^j \,\sum_k\,H_{l_2}^k\,\Theta^k
}_{ \octagon \! \! \! \! \! \tiny{\sf f}} 
+
\underline{ 
\frac{ 1}{ 4} \,\delta_{ l_1}^j \,L_{l_2,l_2}^{l_2}\, \Theta^0
}_{ \octagon \! \! \! \! \! \tiny{\sf g}} 
+
\underline{ 
\frac{ 1}{ 4} \,\delta_{ l_1}^j \,\Theta^0 \,\Theta^{ l_2}
}_{ \octagon \! \! \! \! \! \tiny{\sf h}} \ .
\endaligned
$$
Simplifying and ordering, we obtain the
family (II) of partial differential relations
of Theorem~1.7~{\bf (3)}:
\def\theequation{3.108}\begin{equation}
\aligned
0
&
=
-
\underline{
\frac{ 1}{ 2} \, H_{l_1, y^{l_2}}^j
+ 
\frac{ 1}{ 6} \, \delta_{ l_1}^j \, H_{ l_2, y^{l_2}}^{l_2}
+
\frac{ 1}{ 3} \, \delta_{ l_2}^j \, H_{l_1, y^{l_1}}^{l_1}
}
+ \\
& \
\ \ \ \ \
+
\underline{
L_{l_1,l_2,x}^j 
- 
\frac{ 1}{ 3} \, \delta_{ l_1}^j \,L_{l_2, l_2, x}^{l_2}
- 
\frac{ 2}{ 3} \,\delta_{ l_2}^j \, L_{l_1, l_1, x}^{l_1}
}
+ \\
& \
\ \ \ \ \ 
+
G^j\, M_{l_1, l_2}
- 
\frac{ 1}{ 3} \, \delta_{ l_1}^j \, G^{l_2} \, M_{ l_2, l_2}
- 
\frac{ 2}{ 3} \,\delta_{ l_2}^j \, G^{l_1} \, M_{l_1, l_1}
+ 
\frac{ 1}{ 3} \,\delta_{ l_1}^j \,\sum_k\,G^k\,M_{l_2,k}
- \\
& \
\ \ \ \ \
- 
\frac{ 1}{ 3} \,\delta_{ l_2}^j \,\sum_k\, G^k\, M_{l_1,k} 
- 
\frac{ 1}{ 2} \, \sum_k\, H_k^j \, L_{l_1, l_2}^k 
+ 
\frac{ 1}{ 2} \, \sum_k\, H_{ l_1}^k \, L_{ l_2, k}^j 
+ \\
& \
\ \ \ \ \
+
\frac{ 1}{ 6} \,\delta_{ l_1}^j \,\sum_k\,H_k^{l_2}\, L_{l_2,l_2}^k
- 
\frac{ 1}{ 6}\,\delta_{l_1}^j\,\sum_k\, H_{l_2}^k\,L_{l_2, k}^{l_2}
+ \\
& \
\ \ \ \ \
+
\frac{ 1}{ 3} \,\delta_{ l_2}^j \,\sum_k\,H_k^{l_1}\,L_{l_1,l_1}^k 
- 
\frac{ 1}{ 3}\,\delta_{l_2}^j \,\sum_k\,H_{l_1}^k\,L_{l_1,k}^{l_1}.
\endaligned
\end{equation}

Thirdly, let us insert inside~\thetag{ 3.89} the values of $\Theta_{
y^{l_3}}^{l_2}$, of $\Theta_{ y^{l_2}}^{l_3}$, of $\Theta_{
y^{l_3}}^{l_1}$ and of $\Theta_{ y^{l_2}}^{l_1}$ given by~\thetag{
3.102}. We place all the terms in the right-hand side of the equality
and we place the first order terms in the beginning (first four lines
just below). We obtain:
\def\theequation{3.109}\begin{equation}
\small
\aligned
0 
&
=
\underline{ 
\frac{ 1}{ 2} \,\delta_{ l_1}^j \, L_{l_2, l_2, y^{l_3}}^{l_2}
}_{ \octagon \! \! \! \! \! \tiny{\sf a}} 
- 
\underline{ 
\delta_{ l_1}^j \, M_{l_2,l_3,x}
}_{ \octagon \! \! \! \! \! \tiny{\sf b}}
+ 
\underline{ 
\frac{ 1}{ 2} \,\delta_{ l_2}^j \, L_{l_1,l_1,y^{l_3}}^{l_1}
}_{ \octagon \! \! \! \! \! \tiny{\sf c}} 
-
\underline{ 
\delta_{l_2}^j \, M_{l_1, l_3, x}
}_{ \fbox{\tiny 4}}
- \\
& \
\ \ \ \ \
- 
\underline{ 
\frac{ 1}{ 2} \,\delta_{ l_1}^j \,L_{l_3,l_3,y^{l_2}}^{l_3}
}_{ \octagon \! \! \! \! \! \tiny{\sf d}} 
+
\underline{ 
\delta_{ l_1}^j \,M_{l_3,l_2,x}
}_{ \octagon \! \! \! \! \! \tiny{\sf b}} 
-
\underline{
\frac{ 1}{ 3} \,\delta_{ l_3}^j \, L_{l_1,l_1,y^{l_2}}^{l_1}
}_{ \octagon \! \! \! \! \! \tiny{\sf e}} 
+
\underline{ 
\delta_{l_3}^j \, M_{ l_1, l_2,x}
}_{ \fbox{\tiny 3}}
+ \\
& \
\ \ \ \ \
+
\underline{ 
L_{l_1,l_2,y^{l_3}}^j
}_{ \fbox{\tiny 1}}
-
\underline{ 
L_{l_1,l_3,y^{l_2}}^j 
}_{ \fbox{\tiny 2}} 
+ 
\underline{ 
\frac{ 1}{ 2} \,\delta_{ l_1}^j \,L_{l_3,l_3,y^{l_2}}^{l_3}
}_{ \octagon \! \! \! \! \! \tiny{\sf d}} 
- 
\underline{ 
\frac{ 1}{ 2} \,\delta_{ l_1}^j \, L_{l_2, l_2, y^{l_3}}^{l_2}
}_{ \octagon \! \! \! \! \! \tiny{\sf a}} 
+ \\
& \
\ \ \ \ \
+
\underline{
\frac{ 1}{ 3} \,\delta_{ l_3}^j \, L_{l_1,l_1,y^{l_2}}^{l_1}
}_{ \octagon \! \! \! \! \! \tiny{\sf e}}
-
\underline{ 
\frac{ 1}{ 2} \,\delta_{ l_2}^j \, L_{l_1,l_1,y^{l_3}}^{l_1}
}_{ \octagon \! \! \! \! \! \tiny{\sf c}} 
+ \\
\endaligned
\end{equation}
$$
\small
\aligned
& \
\ \ \ \
+
\underline{ 
\frac{ 1}{ 2} \, H_{ l_3}^j \, M_{ l_1, l_2}
}_{ \fbox{\tiny 5}}
- 
\underline{
\frac{ 1}{ 2} \, H_{ l_2}^j \, M_{ l_1, l_3} 
}_{ \fbox{\tiny 6}}
+
\underline{ 
\frac{ 1}{ 4}\,\delta_{ l_2}^j \,L_{l_3,l_3}^{l_3}\,L_{l_1,l_1}^{l_1}
}_{ \octagon \! \! \! \! \! \tiny{\sf f}}
- 
\underline{ 
\frac{ 1}{ 4}\,\delta_{ l_3}^j \,L_{l_2,l_2}^{l_2}\,L_{l_1,l_1}^{l_1}
}_{ \octagon \! \! \! \! \! \tiny{\sf g}}
+ \\
& \
\ \ \ \ \
+
\underline{ 
\sum_k\, L_{l_1, l_3}^k\, L_{l_2, k}^j 
}_{ \fbox{\tiny 11}}
-
\underline{ 
\sum_k\, L_{l_1, l_2}^k\, L_{l_3, k}^j 
}_{ \fbox{\tiny 12}}
+ \\
& \
\ \ \ \ \
+
\underline{ 
\frac{ 1}{ 2} \,\delta_{ l_3}^j \,\sum_k\,L_{l_1,l_2}^k\,L_{k,k}^k
}_{ \octagon \! \! \! \! \! \tiny{\sf h}}
- 
\underline{ 
\frac{ 1}{ 2} \,\delta_{ l_2}^j \,\sum_k\,L_{l_1,l_3}^k\,L_{k,k}^k
}_{ \octagon \! \! \! \! \! \tiny{\sf i}}
+ \\
& \
\ \ \ \ \
+
\underline{ 
\frac{ 1}{ 4} \,\delta_{ l_2}^j \,L_{l_1,l_1}^{l_1}\,\Theta^{l_3}
}_{ \octagon \! \! \! \! \! \tiny{\sf j}} 
-
\underline{ 
\frac{ 1}{ 4} \,\delta_{ l_3}^j \,L_{l_1,l_1}^{l_1}\,\Theta^{l_2}
}_{ \octagon \! \! \! \! \! \tiny{\sf k}}
+
\underline{ 
\frac{ 1}{ 4} \,\delta_{ l_2}^j \,L_{l_3,l_3}^{l_3}\,\Theta^{l_1}
}_{ \octagon \! \! \! \! \! \tiny{\sf l}}
- 
\underline{ 
\frac{ 1}{ 4} \,\delta_{ l_3}^j \,L_{l_2,l_2}^{l_2}\,\Theta^{l_1}
}_{ \octagon \! \! \! \! \! \! \tiny{\sf m}} 
+ \\
& \
\ \ \ \ \
+ 
\underline{ 
\frac{ 1}{ 2} \,\delta_{ l_3}^j \,\sum_k\,L_{l_1,l_2}^k\, \Theta^k
}_{ \octagon \! \! \! \! \! \tiny{\sf n}} 
-
\underline{ 
\frac{ 1}{ 2} \,\delta_{ l_2}^j \,\sum_k\,L_{l_1,l_3}^k\, \Theta^k
}_{ \octagon \! \! \! \! \! \tiny{\sf o}}
+ \\
& \
\ \ \ \ \
+
\underline{
\frac{ 1}{ 2} \,\delta_{ l_2}^j \,M_{l_1,l_3} \, \Theta^0 
}_{ \octagon \! \! \! \! \! \tiny{\sf p}} 
- 
\underline{
\frac{ 1}{ 2} \,\delta_{ l_3}^j \,M_{l_1,l_2} \, \Theta^0 
}_{ \octagon \! \! \! \! \! \tiny{\sf q}} 
+ 
\underline{ 
\frac{ 1}{ 4} \,\delta_{ l_2}^j \,\Theta^{l_3}\,\Theta^{l_1}
}_{ \octagon \! \! \! \! \! \tiny{\sf r}} 
-
\underline{ 
\frac{ 1}{ 4} \,\delta_{ l_3}^j \,\Theta^{l_2}\,\Theta^{l_1}
}_{ \octagon \! \! \! \! \! \tiny{\sf s}}
- \\
\endaligned
$$
$$
\small
\aligned
& \
\ \ \ \ \
-
\underline{ 
\frac{ 1}{ 2} \, \delta_{ l_1}^j \,\sum_k\,H_{l_2}^k\,M_{l_3,k}
}_{ \fbox{\tiny 8}}
-
\underline{ 
\frac{ 1}{ 4} \,\delta_{ l_1}^j \,L_{l_2,l_2}^{l_2}\,L_{l_3,l_3}^{l_3}
}_{ \octagon \! \! \! \! \! \tiny{\sf t}}
+
\underline{ 
\frac{ 1}{ 2} \,\delta_{ l_1}^j \,\sum_k\,L_{l_2,l_3}\,L_{k,k}^k
}_{ \octagon \! \! \! \! \! \tiny{\sf u}} 
- \\
& \
\ \ \ \ \
-
\underline{ 
\frac{ 1}{ 4} \,\delta_{ l_1}^j \, L_{l_2,l_2}^{l_2}\,\Theta^{l_3}
}_{ \octagon \! \! \! \! \! \tiny{\sf v}} 
- 
\underline{ 
\frac{ 1}{ 4} \,\delta_{ l_1}^j \, L_{l_3,l_3}^{l_3}\,\Theta^{l_2}
}_{ \octagon \! \! \! \! \! \! \tiny{\sf w}} 
+ 
\underline{ 
\frac{ 1}{ 2} \,\delta_{ l_1}^j \,\sum_k\, L_{l_2,l_3}^k\,\Theta^k
}_{ \octagon \! \! \! \! \! \tiny{\sf x}} 
- \\
& \
\ \ \ \ \
-
\underline{ 
\frac{ 1}{ 2} \,\delta_{ l_1}^j \,M_{l_2,l_3}\, \Theta^0
}_{ \octagon \! \! \! \! \! \tiny{\sf y}} 
- 
\underline{ 
\frac{ 1}{ 4} \,\delta_{ l_1}^j \,\Theta^{l_2}\, \Theta^{l_3}
}_{ \octagon \! \! \! \! \! \tiny{\sf z}} 
- \\
\endaligned
$$
$$
\small
\aligned
& \
\ \ \ \ \
- 
\underline{ 
\frac{ 1}{ 2} \,\delta_{ l_2}^j \,\sum_k\,H_{l_1}^k\,M_{l_3,k}
}_{ \fbox{\tiny 10}}
-
\underline{ 
\frac{ 1}{ 4}\,\delta_{ l_2}^j \,L_{l_1,l_1}^{l_1}\,L_{l_3,l_3}^{l_3}
}_{ \octagon \! \! \! \! \! \tiny{\sf f}}
+
\underline{ 
\frac{ 1}{ 2} \,\delta_{ l_2}^j \,\sum_k\,L_{l_1,l_3}^k\,L_{k,k}^k
}_{ \octagon \! \! \! \! \! \tiny{\sf i}}
- \\
& \
\ \ \ \ \
- 
\underline{ 
\frac{ 1}{ 4} \,\delta_{ l_2}^j \,L_{l_1,l_1}^{l_1}\,\Theta^{l_3}
}_{ \octagon \! \! \! \! \! \tiny{\sf j}} 
- 
\underline{ 
\frac{ 1}{ 4} \,\delta_{ l_2}^j \,L_{l_3,l_3}^{l_3}\,\Theta^{l_1}
}_{ \octagon \! \! \! \! \! \tiny{\sf l}}
+
\underline{ 
\frac{ 1}{ 2} \,\delta_{ l_2}^j \,\sum_k\,L_{l_1,l_3}^k\, \Theta^k
}_{ \octagon \! \! \! \! \! \tiny{\sf o}}
- \\
& \
\ \ \ \ \
-
\underline{
\frac{ 1}{ 2} \,\delta_{ l_2}^j \,M_{l_1,l_3} \, \Theta^0 
}_{ \octagon \! \! \! \! \! \tiny{\sf p}}
- 
\underline{ 
\frac{ 1}{ 4} \,\delta_{ l_2}^j \,\Theta^{l_3}\,\Theta^{l_1}
}_{ \octagon \! \! \! \! \! \tiny{\sf r}} 
+ \\
\endaligned
$$
$$
\small
\aligned
& \
\ \ \ \ \
+
\underline{ 
\frac{ 1}{ 2} \,\delta_{ l_1}^j \,\sum_k\, H_{l_3}^k \,M_{l_2,k}
}_{ \fbox{\tiny 7}}
+ 
\underline{ 
\frac{ 1}{ 4} \,\delta_{ l_1}^j \,L_{l_3,l_3}^{l_3}\,L_{l_2,l_2}^{l_2}
}_{ \octagon \! \! \! \! \! \tiny{\sf t}}
- 
\underline{ 
\frac{ 1}{ 2} \,\delta_{ l_1}^j \,\sum_k\,L_{l_3,l_2}\,L_{k,k}^k
}_{ \octagon \! \! \! \! \! \tiny{\sf u}}
+ \\
& \
\ \ \ \ \
+
\underline{ 
\frac{ 1}{ 4} \,\delta_{ l_1}^j \, L_{l_3,l_3}^{l_3}\,\Theta^{l_2}
}_{ \octagon \! \! \! \! \! \! \tiny{\sf w}} 
+ 
\underline{ 
\frac{ 1}{ 4} \,\delta_{ l_1}^j \, L_{l_2,l_2}^{l_2}\,\Theta^{l_3}
}_{ \octagon \! \! \! \! \! \tiny{\sf v}} 
- 
\underline{ 
\frac{ 1}{ 2} \,\delta_{ l_1}^j \,\sum_k\, L_{l_3,l_2}^k\,\Theta^k
}_{ \octagon \! \! \! \! \! \tiny{\sf x}} 
+ \\
& \
\ \ \ \ \
+
\underline{ 
\frac{ 1}{ 2} \,\delta_{ l_1}^j \,M_{l_3,l_2}\, \Theta^0
}_{ \octagon \! \! \! \! \! \tiny{\sf y}} 
+
\underline{ 
\frac{ 1}{ 4} \,\delta_{ l_1}^j \,\Theta^{l_3}\, \Theta^{l_2}
}_{ \octagon \! \! \! \! \! \tiny{\sf z}} 
+ \\ 
\endaligned
$$
$$
\small
\aligned
& \
\ \ \ \ \
+
\underline{ 
\frac{ 1}{ 2} \,\delta_{ l_3}^j \, \sum_k\, H_{l_1}^k \, M_{ l_2,k}
}_{ \fbox{\tiny 9}}
+
\underline{ 
\frac{ 1}{ 4}\,\delta_{ l_3}^j \,L_{l_1,l_1}^{l_1}\,L_{l_2,l_2}^{l_2}
}_{ \octagon \! \! \! \! \! \tiny{\sf g}}
- 
\underline{ 
\frac{ 1}{ 2} \,\delta_{ l_3}^j \,\sum_k\,L_{l_1,l_2}^k\,L_{k,k}^k
}_{ \octagon \! \! \! \! \! \tiny{\sf h}}
+ \\
& \
\ \ \ \ \
+
\underline{ 
\frac{ 1}{ 4} \,\delta_{ l_3}^j \,L_{l_1,l_1}^{l_1}\,\Theta^{l_2}
}_{ \octagon \! \! \! \! \! \tiny{\sf k}}
+
\underline{ 
\frac{ 1}{ 4} \,\delta_{ l_3}^j \,L_{l_2,l_2}^{l_2}\,\Theta^{l_1}
}_{ \octagon \! \! \! \! \! \! \tiny{\sf m}} 
- 
\underline{ 
\frac{ 1}{ 2} \,\delta_{ l_3}^j \,\sum_k\,L_{l_1,l_2}^k\, \Theta^k
}_{ \octagon \! \! \! \! \! \tiny{\sf n}} 
+ \\
& \
\ \ \ \ \
+ 
\underline{
\frac{ 1}{ 2} \,\delta_{ l_3}^j \,M_{l_1,l_2} \, \Theta^0 
}_{ \octagon \! \! \! \! \! \tiny{\sf q}}
+ 
\underline{ 
\frac{ 1}{ 4} \,\delta_{ l_3}^j \,\Theta^{l_1}\,\Theta^{l_2}
}_{ \octagon \! \! \! \! \! \tiny{\sf s}} \ .
\endaligned
$$
Simplifying and ordering, we obtain the
family (III) of
partial differential relations of
Theorem~1.7~{\bf (3)}:
\def\theequation{3.110}\begin{equation}
\aligned
0 
& 
=
\underline{
L_{l_1, l_2, y^{l_3}}^j 
- 
L_{l_1, l_3, y^{l_2}}^j 
+ 
\delta_{ l_3}^j \, M_{l_1, l_2, x}
- 
\delta_{ l_2}^j \, M_{l_1, l_3, x}
}
+ \\
& \
\ \ \ \ \
+ 
\frac{ 1}{ 2} \, H_{l_3}^j \, M_{l_1, l_2} 
- 
\frac{ 1}{ 2} \, H_{ l_2}^j \, M_{l_1, l_3}
+ \\
& \
\ \ \ \ \
+ 
\frac{ 1}{ 2} \, \delta_{ l_1}^j \,\sum_k\,H_{l_3}^k\,M_{l_2,k}
- 
\frac{ 1}{ 2} \, \delta_{ l_1}^j \,\sum_k\,H_{l_2}^k\,M_{l_3,k}
+ \\
& \
\ \ \ \ \
+
\frac{ 1}{ 2} \, \delta_{ l_3}^j \,\sum_k\,H_{l_1}^k\,M_{l_2,k}
-
\frac{ 1}{ 2} \, \delta_{ l_2}^j \,\sum_k\,H_{l_1}^k\,M_{l_3,k}
+ \\
& \
\ \ \ \ \
+ 
\sum_k\, L_{l_1, l_3}^k \, L_{l_2, k}^j 
- 
\sum_k\, L_{l_1, l_2}^k \, L_{l_3, k}^j.
\endaligned
\end{equation}

\subsection*{3.111.~Arguments for the proof of Theorem~1.7~{\bf (3)}: 
necessity
and sufficiency of (I), (II), (III), (IV)} Let us summarize the
implications that have been established so far, from the beginning of
Section~3. Recall that $m\geqslant 2$.

\smallskip

\begin{itemize}
\item[$\bullet$]
There exist functions $X$, $Y^j$ of
$(x,y)$ transforming the 
system $y_{xx}^j = F^j (x, y, y_x)$, $j= 1, \dots,
m$, to the free particle system
$Y_{XX}^j = 0$, $j= 1, \dots, m$.
\item[$\Downarrow$]
\item[$\bullet$]
There exist functions $\Pi_{ l_1, l_2}^j$ of
$(x, y)$, $0 \leqslant j, l_1, l_2 \leqslant m$,
satisfying the first auxiliary system~\thetag{ 3.38}
of partial differential equations.
\item[$\Downarrow$]
\item[$\bullet$]
There exist (principal unknowns) functions $\Theta^0$, $\Theta^j$
satisfying the six families of partial differential 
equations~\thetag{ 3.69}, 
\thetag{ 3.86}, \thetag{ 3.89}, \thetag{ 3.91}, 
\thetag{ 3.93} and~\thetag{ 3.96}.
\item[$\Downarrow$]
\item[$\bullet$]
The functions $G^j$, $H_{l_1}^j$, $L_{l_1, l_2}^j$
and $M_{l_1,l_2}$ satisfy the four families of
partial differential equations
(I), (II), (III) and (IV)
of Theorem~1.7~{\bf (3)}.
\end{itemize}

\smallskip

The four families of first order partial differential
equations~\thetag{ 3.99}, \thetag{ 3.100}, \thetag{ 3.101}
and~\thetag{ 3.102} satisfied by the
principal unknowns will be called the {\sl second auxiliary system}.
It is a complete system.

To achieve the proof of Theorem~1.7~{\bf (3)}, 
we have to establish the reverse
implications. More precisely: 

\smallskip


\begin{itemize}
\item[$\bullet$]
Some given functions $G^j$, $H_{l_1}^j$, $L_{l_1, l_2}^j= L_{ l_2,
l_1}^j$ and $M_{l_1,l_2} = M_{ l_2, l_1}$ of $(x,y)$ satisfy the four
families of partial differential equations (I), (II), (III) and (IV)
of Theorem~1.7~{\bf (3)}, or equivalently, the
partial differential equations~\thetag{ 3.106}, 
\thetag{ 3.108}, \thetag{ 3.110} and~\thetag{ 3.96}. 
\item[$\Downarrow$]
\item[$\bullet$]
There exist functions $\Theta^0$, $\Theta^j$ satisfying the second
auxiliary system~\thetag{ 3.99}, \thetag{ 3.100}, \thetag{ 3.101}
and~\thetag{ 3.102}.
\item[$\Downarrow$]
\item[$\bullet$]
These solution functions $\Theta^0$, $\Theta^j$ satisfy the six
families of partial differential equations~\thetag{ 3.69}, \thetag{
3.86}, \thetag{ 3.89}, \thetag{ 3.91}, \thetag{ 3.93} and~\thetag{
3.96}.
\item[$\Downarrow$]
\item[$\bullet$]
There exist functions $\Pi_{ l_1, l_2}^j$ of $(x, y)$, $0 \leqslant j, l_1,
l_2 \leqslant m$, satisfying the first auxiliary system~\thetag{ 3.38} of
partial differential equations.
\item[$\Downarrow$]
\item[$\bullet$]
There exist functions $X$, $Y^j$ of $(x,y)$ transforming the system
$y_{xx}^j = F^j (x, y, y_x)$, $j= 1, \dots, m$, to the free particle
system $Y_{XX}^j = 0$, $j= 1, \dots, m$.
\end{itemize}

\smallskip

The above last three implications have been already implicitely
established in the preceding paragraphs, as may be checked by
inspecting Lemma~3.40 and the formal computations after \S3.62.

Thus, {\it it remains only to establish the first implication in the
above reverse list}. Since the second auxiliary system~\thetag{ 3.99},
\thetag{ 3.100}, \thetag{ 3.101} and~\thetag{ 3.102} is complete and
of first order, a necessary and sufficient condition for the existence
of solutions follows by writing out the following four families of
cross-differentiations: 
\def\theequation{3.112}\begin{equation}
\left\{
\aligned
0 
& 
=
\left(
\Theta_x^0 
\right)_{y^{l_1}}
-
\left(
\Theta_{y^{l_1}}^0
\right)_x, 
\\
0 
& 
=
\left(
\Theta_{y^{l_1}}^0
\right)_{y^{ l_2}}
-
\left(
\Theta_{y^{l_2}}^0
\right)_{y^{ l_1}},
\\
0 
& 
=
\left(
\Theta_x^{l_1} 
\right)_{y^{l_2}}
-
\left(
\Theta_{y^{l_2}}^{l_1}
\right)_x, 
\\
0 
& 
=
\left(
\Theta_{y^{l_2}}^{l_1}
\right)_{y^{ l_3}}
-
\left(
\Theta_{y^{l_3}}^{l_1}
\right)_{y^{ l_2}}.
\endaligned\right.
\end{equation}
In the hardest techical part of this paper (Section~4 below), we
verify that these four families of compatibility conditions are a
consequence of (I), (II), (III) and (IV). For reasons of space, we
shall in fact only study the first family of compatibility conditions,
{\it i.e.} the first line of~\thetag{ 3.112}. In the our
manuscript, we have treated the remaining three families
of compatibility conditions similarly and completely, up to the very
end of every branch of the coral tree of computations. However, we would
like to mention that typesetting the remaining three cases would add
at least fifty pages of Latex to Section~4. Thus, we prefer to expose
thoroughly the treatment of the first family of compatibility
conditions, explaining implicitely how to guess 
the treatment of the remaining three.

\section*{ \S4.~Compatibility conditions for 
the second auxiliary system}

So, we have to develope the first line of~\thetag{ 3.112}: we replace
$\Theta_x^0$ by its expression~\thetag{ 3.99}, we differentiate it
with respect to $y^{l_1}$, we replace $\Theta_{ y^{l_1}}^0$ by its
expression~\thetag{ 3.100}, we differentiate it with respect to $x$
and we substract. We get:
\def\theequation{4.1}\begin{equation}
\small
\aligned
0
&
= 
\left(
\Theta_x^0 
\right)_{y^{l_1}}
-
\left(
\Theta_{y^{l_1}}^0
\right)_x
\\
&
=
-
\underline{\underline{
2\, G_{ y^{l_1}y^{ l_1}}^{l_1} 
+ 
H_{ l_1, xy^{l_1}}^{l_1}
}}
+ \\
& \
\ \ \ \ \
+
2\, \sum_k\, G_{y^{l_1}}^k\,L_{l_1,k}^{l_1}
+
2\,\sum_k\,G^k\,L_{l_1,k,y^{l_1}}^{l_1}
-
\sum_k\,G_{y^{l_1}}^k\,L_{k,k}^k
-
\sum_k\,G^k\,L_{k,k,y^{l_1}}^k
- \\
& \
\ \ \ \ \
- 
\frac{ 1}{ 2} \,\sum_k\,H_{l_1,y^{l_1}}^k\,H_k^{l_1}
-
\frac{ 1}{ 2} \,\sum_k\,H_{l_1}^k\,H_{k,y^{l_1}}^{l_1}
-
\sum_k\,G_{y^{l_1}}^k\,\Theta^k 
-
\sum_k\,G^k\,\underline{\Theta_{y^{l_1}}^k}
+ \\
& \
\ \ \ \ \
+
\Theta^0 \, \underline{\Theta_{y^{l_1}}^0}
- \\
& \
\ \ \ \ \
-
\underline{\underline{
\frac{ 2}{ 3} \, L_{l_1,l_1,xx}^{l_1}
+
\frac{ 1}{ 3} \, H_{l_1, y^{l_1}x}^{l_1}
}} 
- \\
& \
\ \ \ \ \
- 
\frac{ 2}{ 3} \, G_x^{l_1}\,M_{l_1,l_1} 
-
\frac{ 2}{ 3} \, G^{l_1}\,M_{l_1,l_1,x}
-
\frac{ 4}{ 3} \,\sum_k\,G_x^k\,M_{l_1,k}
-
\frac{ 4}{ 3} \,\sum_k\,G^k\,M_{l_1,k,x}
+ \\
& \
\ \ \ \ \
+
\frac{ 1}{ 3} \, \sum_k\, H_{ k,x}^{l_1}\, L_{ l_1,l_1}^k
+ 
\frac{ 1}{ 3} \, \sum_k \, H_k^{l_1} \, L_{ l_1,l_1,x}^k
-
\frac{ 1}{ 3} \, \sum_k\, H_{ l_1,x}^k\, L_{ l_1,k}^{l_1}
-
\frac{ 1}{ 3} \, \sum_k\, H_{ l_1}^k\, L_{ l_1,k,x}^{l_1}
+ \\
& \
\ \ \ \ \
+
\frac{ 1}{ 2} \,\sum_k\, H_{l_1,x}^k\,L_{k,k}^k
+ 
\frac{ 1}{ 2} \,\sum_k\, H_{l_1}^k\,L_{k,k,x}^k
+ 
\frac{ 1}{ 2} \,\sum_k\, H_{l_1,x}^k\, \Theta^k
+
\frac{ 1}{ 2} \,\sum_k\, H_{l_1}^k\, \underline{\Theta_x^k}
- \\
& \
\ \ \ \ \
-
\frac{ 1}{ 2} \, L_{l_1,l_1,x} \, \Theta^0
-
\frac{ 1}{ 2} \, L_{l_1,l_1} \, \underline{\Theta_x^0}
-
\frac{ 1}{ 2} \,\underline{\Theta_x^0} \,\Theta^{l_1}
-
\frac{ 1}{ 2} \,\Theta^0 \,\underline{\Theta_x^{l_1}}.
\endaligned
\end{equation}
Here, we underline twice the second order terms. Also, we have
underlined once the six terms: $\underline{ \Theta_{ y^{ l_1 }}^k}$,
$\underline{ \Theta_{ y^{ l_1 }}^0}$, $\underline{ \Theta_x^k }$,
$\underline{ \Theta_x^0 }$, $\underline{ \Theta_x^0}$ and $\underline{
\Theta_x^{l_1}}$. They must be replaced by their values given
in~\thetag{ 3.99}, \thetag{ 3.100}, \thetag{ 3.101} and~\thetag{
3.102}. In this replacement, some double sums appear. As before, we
use the first index $k= 1,\dots, m$ for 
single summation and then the second
index $p= 1,\dots, m$ for double
summation. Finally, we put all the second
order terms in the first line, not disturbing the
order of appearance of the 73 remaining terms. We get:
\def\theequation{4.2}\begin{equation}
\small
\aligned
0 
& 
=
-
\underline{\underline{
2 \, G_{ y^{l_1} y^{l_1}}^{l_1}
+ 
\frac{ 4}{ 3} \, H_{l_1, xy^{l_1}}^{l_1}
- 
\frac{ 2}{ 3} \, L_{l_1, l_1, xx}^{l_1}
}}
+ \\
& \
\ \ \ \ \
+ 
\underline{ 
2\, \sum_k\, G_{ y^{ l_1}}^k\,L_{l_1,k}^{l_1}
}_{ \fbox{\tiny 4}}
+ 
\underline{ 
2\,\sum_k\,G^k\,L_{l_1,k,y^{l_1}}^{l_1}
}_{ \fbox{\tiny 16}}
-
\underline{ 
\sum_k\,G_{y^{l_1}}^k\, L_{k,k}^k
}_{ \fbox{\tiny 5}}
-
\underline{ 
\sum_k\,G^k\,L_{k,k,y^{l_1}}^k
}_{ \octagon \! \! \! \! \! \tiny{\sf a}} 
- \\
& \
\ \ \ \ \
-
\underline{ 
\frac{ 1}{ 2} \,\sum_k\,H_{l_1,y^{l_1}}^k\,H_k^{l_1}
}_{ \fbox{\tiny 10}}
- 
\underline{ 
\frac{ 1}{ 2} \,\sum_k\,H_{l_1}^k\,H_{k,y^{l_1}}^{l_1}
}_{ \fbox{\tiny 11}}
- 
\underline{ 
\sum_k\,G_{y^{l_1}}^k\,\Theta^k
}_{ \octagon \! \! \! \! \! \tiny{\sf \alpha}} 
+ \\
\endaligned
\end{equation}
$$
\small
\aligned
& \
\ \ \ \ \
+ 
\underline{ 
\sum_k\,G^k\,L_{k,k,y^{l_1}}^k
}_{ \octagon \! \! \! \! \! \tiny{\sf a}} 
- 
\underline{ 
2\, \sum_k\, G^k\,M_{k,l_1,x}
}_{ \fbox{\tiny 18}}
-
\underline{ 
\sum_k\,\sum_p\, G^k\,H_k^p\,M_{l_1,p}
}_{ \fbox{\tiny 19}}
- \\
& \ 
\ \ \ \ \
- 
\underline{ 
\frac{ 1}{ 2} \,\sum_k\,G^k\,L_{k,k}^k\,L_{l_1,l_1}^{l_1}
}_{ \octagon \! \! \! \! \! \tiny{\sf b}} 
+
\underline{ 
\sum_k\,\sum_p\,G^k\,L_{k,l_1}^p\,L_{p,p}^p
}_{ \fbox{\tiny 23}}
- 
\underline{ 
\frac{ 1}{ 2} \,\sum_k\,G^k\,L_{k,k}^k\,\Theta^{l_1}
}_{ \octagon \! \! \! \! \! \tiny{\sf c}} 
- \\
& \
\ \ \ \ \
- 
\underline{ 
\frac{ 1}{ 2} \, \sum_k\,G^k\,L_{l_1,l_1}^{l_1}\,\Theta^k
}_{ \octagon \! \! \! \! \! \tiny{\sf d}} 
+
\underline{ 
\sum_k\,\sum_p\,G^k\,L_{k,l_1}^p\,\Theta^p
}_{ \octagon \! \! \! \! \! \tiny{\sf \varepsilon}} 
-
\underline{ 
\sum_k\,G^k\,M_{k,l_1}\,\Theta^0
}_{ \octagon \! \! \! \! \! \tiny{\sf e}} 
-
\underline{ 
\frac{ 1}{ 2} \,\sum_k\,G^k\,\Theta^k\,\Theta^{l_1}
}_{ \octagon \! \! \! \! \! \tiny{\sf f}} 
+ \\
\endaligned
$$
$$
\small
\aligned
& \
\ \ \ \ \
+
\underline{ 
\frac{ 2}{ 3} \, L_{l_1,l_1,x}^{l_1}\,\Theta^0
}_{ \octagon \! \! \! \! \! \tiny{\sf g}} 
- 
\underline{ 
\frac{ 1}{ 3} \,H_{l_1,y^{l_1}}^{l_1}\,\Theta^0
}_{ \octagon \! \! \! \! \! \tiny{\sf h}} 
+ \\
& \
\ \ \ \ \
+
\underline{ 
\frac{2}{3}\,G^{l_1}\,M_{l_1,l_1}\,\Theta^0
}_{ \octagon \! \! \! \! \! \tiny{\sf i}}
+
\underline{ 
\frac{ 4}{3}\,\sum_k\,G^k\,M_{l_1,k}\,\Theta^0
}_{ \octagon \! \! \! \! \! \tiny{\sf e}} 
-
\underline{ 
\frac{ 1}{ 3} \,\sum_k\,H_k^{l_1}\,L_{l_1,l_1}^{l_1}\,\Theta^0
}_{ \octagon \! \! \! \! \! \tiny{\sf j}} 
+
\underline{ 
\frac{ 1}{ 3} \,\sum_k\,H_{l_1}^k\,L_{l_1,k}^{l_1}\,\Theta^0
}_{ \octagon \! \! \! \! \! \tiny{\sf k}} 
- \\
& \
\ \ \ \ \
-
\underline{ 
\frac{ 1}{ 2} \,\sum_k\,H_{l_1}^k\,L_{k,k}^k\,\Theta^0
}_{ \octagon \! \! \! \! \! \tiny{\sf l}} 
- 
\underline{ 
\frac{ 1}{ 2} \,\sum_k\,H_{l_1}^k\, \Theta^k\,\Theta^0
}_{ \octagon \! \! \! \! \! \! \tiny{\sf m}} 
+
\underline{ 
\frac{ 1}{ 2} \,L_{l_1,l_1}^{l_1}\,\Theta^0\,\Theta^0
}_{ \octagon \! \! \! \! \! \tiny{\sf n}} 
+
\underline{ 
\frac{ 1}{ 2} \,\Theta^0\,\Theta^0\,\Theta^{l_1}
}_{ \octagon \! \! \! \! \! \tiny{\sf o}}
- \\
\endaligned
$$
$$
\small
\aligned
& \
\ \ \ \ \
-
\underline{
\frac{ 2}{3}\, G_x^{l_1}\,M_{l_1,l_1}
}_{ \fbox{\tiny 1}}
- 
\underline{
\frac{ 2}{3}\, G^{l_1}\,M_{l_1,l_1,x}
}_{ \fbox{\tiny 17}}
-
\underline{ 
\frac{4}{ 3} \, \sum_k\,G_x^k\,M_{l_1,k}
}_{ \fbox{\tiny 2}}
-
\underline{ 
\frac{4}{ 3} \, \sum_k\,G^k\,M_{l_1,k,x}
}_{ \fbox{\tiny 18}}
+ \\
& \
\ \ \ \ \
+
\underline{ 
\frac{ 1}{ 3} \,\sum_k\,H_{k,x}^{l_1}\,L_{l_1,l_1}^k
}_{ \fbox{\tiny 9}}
+
\underline{ 
\frac{ 1}{ 3} \,\sum_k\,H_k^{l_1}\,L_{l_1,l_1,x}^k
}_{ \fbox{\tiny 14}}
-
\underline{ 
\frac{ 1}{ 3} \,\sum_k\,H_{l_1,x}^k\,L_{l_1,k}^{l_1}
}_{ \fbox{\tiny 7}}
- 
\underline{ 
\frac{ 1}{ 3} \,\sum_k\,H_{l_1}^k\,L_{l_1,k,x}^{l_1}
}_{ \fbox{\tiny 13}}
+ \\
& \
\ \ \ \ \
+
\underline{ 
\frac{ 1}{ 2} \,\sum_k\,H_{l_1,x}^k\,L_{k,k}^k
}_{ \fbox{\tiny 8}}
+ 
\underline{ 
\frac{ 1}{ 2} \,\sum_k\, H_{l_1}^k\, L_{k,k,x}^k\,
}_{ \fbox{\tiny 15}}
+
\underline{ 
\frac{ 1}{ 2} \, \sum_k\, H_{l_1,x}^k\, \Theta^k
}_{ \octagon \! \! \! \! \! \tiny{\sf \gamma}} 
- \\
\endaligned
$$
$$
\small
\aligned
& \
\ \ \ \ \
-
\underline{ 
\frac{ 1}{ 3} \,\sum_k\, H_{l_1}^k\, H_{k,y^k}^k
}_{ \fbox{\tiny 12}}
+
\underline{ 
\frac{ 1}{ 6} \,\sum_k\,H_{l_1}^k\,L_{k,k,x}^k
}_{ \fbox{\tiny 15}}
+
\underline{ 
\frac{ 2}{ 3} \,\sum_k\,H_{l_1}^k\, G^k\,M_{k,k}
}_{ \fbox{\tiny 20}}
+ \\
& \
\ \ \ \ \
+
\underline{ 
\frac{ 1}{ 3} \,\sum_k\,\sum_p\,H_{l_1}^k\,G^p\,M_{k,p}
}_{ \fbox{\tiny 21}}
-
\underline{ 
\frac{ 1}{ 3} \, \sum_k\,\sum_p\,H_{l_1}^k\,H_p^k\,L_{k,k}^p
}_{ \fbox{\tiny 24}}
+
\underline{ 
\frac{ 1}{ 3} \,\sum_k\,\sum_p\,H_{l_1}^k\,H_k^p\,L_{k,p}^k
}_{ \fbox{\tiny 25}}
- \\
& \
\ \ \ \ \
-
\underline{ 
\frac{ 1}{ 4} \,\sum_k\,\sum_p\,H_{l_1}^k\,H_k^p\,L_{p,p}^p
}_{ \fbox{\tiny 26}}
-
\underline{ 
\frac{ 1}{ 4} \,\sum_k\,\sum_p\,H_{l_1}^k\,H_k^p\,\Theta^p
}_{ \octagon \! \! \! \! \! \tiny{\sf \eta}} 
+
\underline{ 
\frac{ 1}{ 4} \,\sum_k\,H_{l_1}^k\,L_{k,k}^k\,\Theta^0
}_{ \octagon \! \! \! \! \! \tiny{\sf l}} 
+ \\
& \
\ \ \ \ \
+
\underline{ 
\frac{ 1}{ 4} \,\sum_k\, H_{l_1}^k\,\Theta^0\,\Theta^k
}_{ \octagon \! \! \! \! \! \! \tiny{\sf m}} 
-
\underline{ 
\frac{ 1}{ 2} \,L_{l_1,l_1,x}^{l_1}\, \Theta^0
}_{ \octagon \! \! \! \! \! \tiny{\sf g}} 
+ \\
\endaligned
$$
$$
\small
\aligned
& \
\ \ \ \ \
+
\underline{ 
G_{y^{l_1}}^{l_1}\, L_{l_1,l_1}^{l_1}
}_{ \fbox{\tiny 3}}
-
\underline{ 
\frac{ 1}{ 2} \, H_{l_1,x}^{l_1}\, L_{l_1,l_1}^{l_1}
}_{ \fbox{\tiny 6}}
- \\
& \
\ \ \ \ \
-
\underline{ 
\sum_k\,G^k\,L_{l_1,l_1}^{l_1}\,L_{l_1,k}^{l_1}
}_{ \fbox{\tiny 22}}
+
\underline{ 
\frac{ 1}{ 2} \,\sum_k\,G^k\,L_{l_1,l_1}^{l_1}\,L_{k,k}^k
}_{ \octagon \! \! \! \! \! \tiny{\sf b}} 
+
\underline{ 
\frac{ 1}{ 4} \, \sum_k\, H_{l_1}^k\,H_k^{l_1}\,L_{l_1,l_1}^{l_1}
}_{ \fbox{\tiny 27}}
+ 
\\
& \
\ \ \ \ \
+
\underline{ 
\frac{ 1}{ 2} \, \sum_k\, G^k\,L_{l_1,l_1}^{l_1}\,\Theta^k
}_{ \octagon \! \! \! \! \! \tiny{\sf d}} 
-
\underline{ 
\frac{ 1}{ 4} \, L_{l_1,l_1}^{l_1}\,\Theta^0\, \Theta^0
}_{ \octagon \! \! \! \! \! \tiny{\sf n}}
+ \\
\endaligned
$$
$$
\small
\aligned
& \
\ \ \ \ \
+
\underline{ 
G_{y^{l_1}}^{l_1}\,\Theta^{l_1}
}_{ \octagon \! \! \! \! \! \tiny{\sf \beta}} 
-
\underline{ 
\frac{ 1}{ 2} \, H_{l_1,x}^{l_1}\,\Theta^{l_1}
}_{ \octagon \! \! \! \! \! \tiny{\sf \delta}} 
- \\
& \
\ \ \ \ \
-
\underline{ 
\sum_k\, G^k\,L_{l_1,k}^{l_1}\,\Theta^{l_1}
}_{ \octagon \! \! \! \! \! \tiny{\sf \zeta}} 
+
\underline{ 
\frac{ 1}{ 2} \,\sum_k\, G^k\,L_{k,k}^k\,\Theta^{l_1}
}_{ \octagon \! \! \! \! \! \tiny{\sf c}} 
+
\underline{ 
\frac{ 1}{ 4} \, \sum_k\, H_{l_1}^k\,H_k^{l_1}\,\Theta^{l_1}
}_{ \octagon \! \! \! \! \! \tiny{\sf \theta}} 
+
\underline{ 
\frac{ 1}{ 2} \,\sum_k\, G^k\,\Theta^k\,\Theta^{l_1}
}_{ \octagon \! \! \! \! \! \tiny{\sf f}} 
- \\
& \
\ \ \ \ \
-
\underline{ 
\frac{ 1}{ 4} \, \Theta^0 \, \Theta^0\, \Theta^{l_1}
}_{ \octagon \! \! \! \! \! \tiny{\sf o}} 
+ \\
\endaligned
$$
$$
\small
\aligned
& \
\ \ \ \ \
+
\underline{
\frac{ 1}{ 3} \, H_{l_1,y^{l_1}}^{l_1}\,\Theta^0 
}_{ \octagon \! \! \! \! \! \tiny{\sf h}} 
-
\underline{ 
\frac{ 1}{ 6} \, L_{l_1,l_1,x}^{l_1} \,\Theta^0
}_{ \octagon \! \! \! \! \! \tiny{\sf g}} 
- \\
& \
\ \ \ \ \
-
\underline{ 
\frac{ 2}{ 3} \, G^{l_1}\,M_{l_1,l_1}\,\Theta^0
}_{ \octagon \! \! \! \! \! \tiny{\sf i}} 
-
\underline{ 
\frac{1}{3}\,\sum_k\,G^k\,M_{l_1,k}\,\Theta^0
}_{ \octagon \! \! \! \! \! \tiny{\sf e}}
+
\underline{ 
\frac{1}{3}\,\sum_k\, H_k^{l_1}\,L_{l_1,l_1}^k\,\Theta^0
}_{ \octagon \! \! \! \! \! \tiny{\sf j}}
-
\underline{ 
\frac{1}{3}\,\sum_k\,H_{l_1}^k\,L_{l_1,k}^{l_1}\,\Theta^0
}_{ \octagon \! \! \! \! \! \tiny{\sf k}}
+ \\
& \
\ \ \ \ \
+
\underline{ 
\frac{1}{4}\,\sum_k\,H_{l_1}^k\,L_{k,k}^k\,\Theta^0
}_{ \octagon \! \! \! \! \! \tiny{\sf l}}
+
\underline{ 
\frac{1}{4}\,\sum_k\, H_{l_1}^k\,\Theta^0\,\Theta^k
}_{ \octagon \! \! \! \! \! \! \tiny{\sf m}}
-
\underline{ 
\frac{1}{4}\,L_{l_1,l_1}^{l_1}\,\Theta^0\,\Theta^0
}_{ \octagon \! \! \! \! \! \tiny{\sf n}}
-
\underline{ 
\frac{1}{4}\,\Theta^0\,\Theta^0\,\Theta^{l_1}
}_{ \octagon \! \! \! \! \! \tiny{\sf o}} \ .
\endaligned
$$
As usual, all the terms underlined with the 15 roman alphabetic letters
${\sf a}, {\sf b}, \ldots, {\sf n}, {\sf o}$ 
appended vanish evidently. Furthermore, 
we claim that the eight terms underlined
with the 8 Greek alphabetic letters 
$\alpha$, $\beta$, $\gamma$, $\delta$, 
$\varepsilon$, $\zeta$, $\eta$ and $\theta$ also 
vanish:
\def\theequation{4.3}\begin{equation}
\small
\aligned
0 
& 
=
?\!\!
=
-
\sum_k\, G_{ y^{l_1}}^k\, \Theta^k
+ 
G_{y^{l_1}}^{l_1}\, \Theta^{l_1}
+
\frac{1}{2}\,\sum_k\, H_{l_1,x}^k\, \Theta^k
- 
\frac{1}{2}\, H_{l_1,x}^{l_1}\, \Theta^{l_1}
+ \\
& \
\ \ \ \ \ \ \ \ \ \
+
\sum_k\, \sum_p\,G^k \, L_{k,l_1}^p\, \Theta^p
- 
\sum_k\, G^k\,L_{l_1,k}^{l_1}\,\Theta^{l_1}
- 
\frac{1}{4}\,\sum_k\,\sum_p\,H_{l_1}^k\,H_k^p\,\Theta^p
+ 
\frac{1}{4}\,\sum_k\,H_{l_1}^k\,H_k^{l_1}\,\Theta^{l_1}.
\endaligned
\end{equation}
Indeed, it suffices to observe that this 
identity coincides with
\def\theequation{4.4}\begin{equation}
\small
\aligned
0 
= 
\frac{1}{2}\,\sum_k\, \Theta^k
\left(
\left.
(3.106)\right\vert_{j:= k; \ l_1:=l_1; \ l_2:=l_1}
\right).
\endaligned
\end{equation}
Simplifying then~\thetag{ 4.2}, we get the explicit
formulation of the first family of compatibility 
conditions for the second auxiliary system:
\def\theequation{4.5}\begin{equation}
\small
\aligned
0 
& 
=
?\!\!
=
\underline{\underline{
2 \, G_{ y^{l_1} y^{l_1}}^{l_1}
+ 
\frac{ 4}{ 3} \, H_{l_1, xy^{l_1}}^{l_1}
- 
\frac{ 2}{ 3} \, L_{l_1, l_1, xx}^{l_1}
}}
- \\
& \
\ \ \ \ \
-
\frac{2}{3}\,G_x^{l_1}\,M_{l_1,l_1}
-
\frac{4}{3}\,\sum_k\,G_x^k\,M_{l_1,k}
+
G_{y^{l_1}}^{l_1}\,L_{l_1,l_1}^{l_1}
+ \\
& \
\ \ \ \ \
+
2\,\sum_k\, G_{y^{l_1}}^k\, L_{l_1,k}^{l_1}
-
\sum_k\, G_{y^{l_1}}^k\, L_{k,k}^k
- \\
& \
\ \ \ \ \
- 
\frac{1}{2}\,H_{l_1,x}^{l_1}\, L_{l_1,l_1}^{l_1}
- 
\frac{1}{3}\,\sum_k\, H_{l_1,x}^k\, L_{l_1,k}^{l_1}
+
\frac{1}{2}\,\sum_k\, H_{l_1,x}^k\,L_{k,k}^k
+ \\
& \
\ \ \ \ \
+
\frac{1}{3}\,\sum_k\, H_{k,x}^{l_1}\, L_{l_1,l_1}^k
-
\frac{1}{2}\,\sum_k\, H_{l_1, y^{l_1}}^k\, H_k^{l_1}
- 
\frac{1}{2}\,\sum_k\, H_{k, y^{l_1}}^{l_1} \, H_{l_1}^k
- \\
& \
\ \ \ \ \
-
\frac{1}{3}\,\sum_k\, H_{k,y^k}^k\, H_{l_1}^k 
- \\
\endaligned
\end{equation}
$$
\small
\aligned
& \
\ \ \ \ \
-
\frac{1}{3}\,\sum_k\, L_{l_1,k,x}^{l_1}\, H_{l_1}^k
+
\frac{1}{3}\,\sum_k\, L_{l_1,l_1,x}^k\, H_k^{l_1}
+
\frac{2}{3}\,\sum_k\, L_{k,k,x}^k\, H_{l_1}^k
+ \\
& \
\ \ \ \ \
+
2\,\sum_k\, L_{l_1,k,y^{l_1}}^{l_1}\, G^k
- \\
& \
\ \ \ \ \
-
\frac{2}{3}\,M_{l_1,l_1,x}\, G^{l_1}
-
\frac{10}{3}\, \sum_k\, M_{k,l_1,x}\, G^k
- \\
& \
\ \ \ \ \
-
\sum_k\, \sum_p\, G^k\, H_k^p\, M_{l_1,p}
+
\frac{2}{3}\,\sum_k\, G^k\, H_{l_1}^k\, M_{k,k}
+ 
\frac{1}{3}\,\sum_k\, \sum_p G^p\, H_{l_1}^k\, M_{k,p}
- \\
& \ 
\ \ \ \ \
- 
\sum_k\, G^k\, L_{l_1, k}^{l_1}\,L_{l_1, l_1}^{l_1}
+
\sum_k\, \sum_p\, G^k\, L_{k,l_1}^p\, L_{p,p}^p
- \\
& \
\ \ \ \ \
- 
\frac{1}{3}\,\sum_k\, \sum_p\, H_{l_1}^k\, H_p^k \, L_{k,k}^p
+ 
\frac{1}{3}\,\sum_k\, \sum_p\,H_{l_1}^k\,H_k^p\,L_{k,p}^k
- \\
& \
\ \ \ \ \
- 
\frac{1}{4}\,\sum_k\, \sum_p\, H_{l_1}^k\, H_k^p\, L_{p,p}^p
+ 
\frac{1}{4}\, \sum_k\, H_{l_1}^k\, H_k^{l_1}\, L_{l_1,l_1}^{l_1}.
\endaligned
$$
We can now state the main technical lemma of this section and of this
paper.

\def\thelemma{4.6}\begin{lemma}
The second order partial differential relations~\thetag{ 4.5} hold
true for $l=1, \dots, m$, and they are a consequence, by
differentiations and by linear combinations, of the fundamental first
order partial differential equations~\thetag{ 3.106}, \thetag{ 3.108},
\thetag{ 3.110} and \thetag{ 3.96}.
\end{lemma}

\subsection*{4.7.~Reconstitution of the appropriate linear combinations}
The remaining of Section~4 is entirely devoted to the proof of this
statement. From the manual computational point of view, the difficulty
of the task is due to the fact that one has to manipulate
formal expressions having from 10 to 50 terms. So the
real question is: {\it how can we reconstitute the
linear combinations and the differentiations which lead to the
goal~\thetag{ 4.5} from the data~\thetag{ 3.106}, \thetag{ 3.108},
\thetag{ 3.110} and \thetag{ 3.96}}?.

The main trick is to first neglect the first order and the zero order
terms in the goal~\thetag{ 4.5}. Using the symbol ``$\equiv$'' to
denote ``{\sl modulo first order and the zero order terms}'', we
formulate the following sub-goal:
\def\theequation{4.8}\begin{equation} \aligned 0 & \equiv ?\!\!
\equiv - \underline{\underline{ 2 \, G_{ y^{l_1} y^{l_1}}^{l_1} +
\frac{ 4}{ 3} \, H_{l_1, xy^{l_1}}^{l_1} - \frac{ 2}{ 3} \, L_{l_1,
l_1, xx}^{l_1} }}, \endaligned
\end{equation}
for $l_1 = 1, \dots, m$. Before estabilishing that these partial
differential relations are a consequence of the data~\thetag{ 3.106},
\thetag{ 3.108}, \thetag{ 3.110} and \thetag{ 3.96} (written with a
similar sign $\equiv$), let us check that they are a consequence of
the existence of the change of coordinates $(x, y) \mapsto (X,Y)$
(however, recall that, in establishing the reverse implications of
\S3.111, we still do not know that such a change of coordinates really
exists); this will confirm the coherence and the validity of our
computations. Importantly, we have been able to achieve systematic
corrections of our computations by always checking them alongside with
the existence of the change of coordinates $(x, y) \mapsto (X,Y)$.

Coming back to the definition~\thetag{ 3.35} and to the
approximation~\thetag{ 3.58}, we have:
\def\theequation{4.9}\begin{equation} \aligned G^{l_1} & = -
\square_{xx}^{l_1} \cong - Y_{xx}^{l_1}, \\ H_{l_1}^{l_1} & = -
2\square_{xy^{l_1}}^{l_1} + \square_{xx}^0 \cong - 2\,
Y_{xy^{l_1}}^{l_1} + X_{xx}, \\ L_{l_1, l_1}^{l_1} & = -
\square_{y^{l_1}y^{l_1}}^{l_1} + 2\, \square_{xy^{l_1}}^0 \cong -
Y_{y^{l_1}y^{l_1}}^{l_1} + 2\, X_{xy^{l_1}}. \endaligned
\end{equation}
Differentiating the first two lines with respect to $y^{l_1}$ and the
third line with respect to $x$, and replacing the sign $\cong$ by the
sign $\equiv$ (in a non-rigorous way, this corresponds essentially to
neglecting the derivatives of order $0$, $1$, $2$ and $3$ of $X$,
$Y^j$ and to neglecting the difference between the Jacobian matrix of
the transformation and the identity matrix), we get:
\def\theequation{4.10}\begin{equation} \small \aligned
\underline{\underline{ G_{y^{l_1}y^{l_1}}^{l_1} }} & \equiv -
Y_{xxy^{l_1}y^{l_1}}^{l_1}, \\ \underline{\underline{ H_{l_1,
xy^{l_1}}^{l_1} }} & \equiv - 2\, Y_{xy^{l_1}xy^{l_1}}^{l_1} +
X_{xxxy^{l_1}}, \\ \underline{\underline{ L_{l_1, l_1, xx}^{l_1} }} &
\equiv - Y_{y^{l_1}y^{l_1}xx}^{l_1} + 2\, X_{xy^{l_1}xx}. \endaligned
\end{equation}
Hence the linear combination $-2\cdot (4.10)_1+ \frac{ 4}{3} \cdot
(4.10)_2 - \frac{ 2}{3} (4.10)_3$ yields the desired result:
\def\theequation{4.11}\begin{equation} \small \aligned 0 & \equiv
?\!\! \equiv - \underline{\underline{ 2\, G_{y^{l_1}y^{l_1}}^{l_1} +
\frac{4}{3}\,H_{l_1, xy^{l_1}}^{l_1} - \frac{2}{3}\,L_{l_1, l_1,
xx}^{l_1} }} \\ & \ \ \ \ \ \equiv \underline{ 2\, Y_{xx
y^{l_1}y^{l_1}}^{l_1} }_{ \octagon \! \! \! \! \! \tiny{\sf a}} -
\underline{ \frac{8}{3}\,Y_{xy^{l_1}xy^{l_1}}^{l_1} }_{ \octagon \! \!
\! \! \! \tiny{\sf a}} + \underline{ \frac{4}{3}\,X_{xxxy^{l_1}} }_{
\octagon \! \! \! \! \! \tiny{\sf b}} + \underline{
\frac{2}{3}\,Y_{y^{l_1}y^{l_1}xx}^{l_1} }_{ \octagon \! \! \! \! \!
\tiny{\sf a}} - \underline{ \frac{4}{3}\, X_{xy^{l_1} xx} }_{ \octagon
\! \! \! \! \! \tiny{\sf b}} \\ & \ \ \ \ \ \equiv 0, \ \ \ \ \ {\rm
indeed}! \endaligned
\end{equation}
Thanks to this straightforward computation, we guess that the
approximate partial differential relations~\thetag{ 4.8} are a
consequence of the approximate relations~\thetag{ 3.106}, \thetag{
3.108}, \thetag{ 3.110} and~\thetag{ 3.96}, namely:
\def\theequation{4.12}\begin{equation} \aligned (3.106)^{\rm mod} : \
\ \ \ \ 0 & \equiv - 2\, G_{y^{l_1}}^j + 2\, \delta_{l_1}^j\,
G_{y^{l_2}}^{l_2} + H_{l_1, x}^j - \delta_{l_1}^j\, H_{l_2, x}^{l_2},
\\ (3.108)^{\rm mod} : \ \ \ \ \ 0 & \equiv - \frac{1}{2}\, H_{l_1,
y^{l_2}}^j + \frac{1}{6}\,\delta_{l_1}^j\, H_{l_2, y^{l_2}}^{l_2} +
\frac{1}{3}\, \delta_{l_2}^j\, H_{l_1, y^{l_1}}^{l_1} + \\ & \ \ \ \ \
\ + L_{l_1, l_2, x}^j - \frac{1}{3}\,\delta_{l_1}^j\, L_{l_2, l_2,
x}^{l_2} - \frac{2}{3}\, \delta_{l_2}^j\, L_{l_1, l_1, x}^{l_1}, \\
(3.110)^{\rm mod} : \ \ \ \ \ 0 & \equiv L_{l_1, l_2, y^{l_3}}^j -
L_{l_1, l_3, y^{l_2}}^j + \delta_{l_3}^j\, M_{l_1, l_2, x} -
\delta_{l_2}^j\, M_{l_1, l_3, x}, \\ (3.96)^{\rm mod} : \ \ \ \ \ 0 &
\equiv M_{l_1, l_2, y^{l_3}} - M_{l_1, l_3, y^{l_2}}. \\ \endaligned
\end{equation}
Here, the sign $\equiv$ means ``{\sl modulo zero order terms}''.
Before proceding further, recall the correspondence between partial
differential relations: \def\theequation{4.13}\begin{equation} \small
\aligned \text{\rm (I) = (3.106)}, \\ \text{\rm (II) = (3.108)}, \\
\text{\rm (III) = (3.110)}, \\ \text{\rm (IV) = (3.96)}. \endaligned
\end{equation}
However, these couples of equivalent identities are written slightly
differently, as may be read by comparison. To fix ideas and to
facilitate the eye-checking of our subsequent computations, {\it we
shall only use and refer to the exact writing of~\thetag{ 3.106},
of~\thetag{ 3.108}, of~\thetag{ 3.110} and of~\thetag{ 3.96}}.

\subsection*{ 4.14.~Construction of a guide}
So we want to show that the approximate relation~\thetag{ 4.8} is a
consequence, by differentiations and by linear combinations, of the
approximate identities~\thetag{ 4.12}. The interest of working with
approximate identities is that formal computations are lightened
substantially. After having discovered which linear combinations and
which differentiations are appropriate, {\it i.e.} after having
constructed a ``guide'', in \S4.22 below, we shall write down the
complete computations, including all zero order terms, following our
guide.

We shall use two indices $l_1$ and $l_2$ with $1 \leqslant l_1, l_2\leqslant m$
and, crucially, $l_2 \neq l_1$. Again, the assumption $m\geqslant 2$ is
used strongly.

Firstly, put $j:= l_1$ in $(3.106)^{\rm mod}$ with $l_2\neq l_1$ and
differentiate with respect to $y^{l_1}$:
\def\theequation{4.15}\begin{equation}
\small
\aligned
0
\equiv
- 
2\, G_{y^{l_1}y^{l_1}}^{l_1}
+
2\, G_{y^{l_2}y^{l_1}}^{l_2}
+
H_{l_1, xy^{l_1}}^{l_1}
-
H_{l_2, xy^{l_1}}^{l_2}.
\endaligned
\end{equation}
Secondly, put $j:= l_2$ in $(3.106)^{\rm mod}$ with $l_1\neq l_2$ and
differentiate with respect to $y^{l_2}$:
\def\theequation{4.16}\begin{equation}
\small
\aligned
0 
\equiv
- 
2\, G_{y^{l_1}y^{l_2}}^{l_2}
+
H_{l_1, xy^{l_2}}^{l_2}.
\endaligned
\end{equation}
Thirdly, put $j:= l_2$ in $(3.108)^{\rm mod}$ with $l_1 \neq l_2$ and
differentiate with respect to $x$:
\def\theequation{4.17}\begin{equation}
\small
\aligned
0 
\equiv
-\frac{1}{2}\, H_{l_1, y^{l_2}x}^{l_2}
+
\frac{1}{3}\, H_{l_1, y^{l_1}x}^{l_1}
+ 
L_{l_1, l_2, xx}^{l_2}
-
\frac{2}{3}\,
L_{l_1, l_1, xx}^{l_1}.
\endaligned
\end{equation}
Fourthly, put $j:= l_1$ in $(3.108)^{\rm mod}$ with $l_2 \neq l_1$ and
differentiate with respect to $x$:
\def\theequation{4.18}\begin{equation}
\small
\aligned
0 
\equiv
-\frac{1}{2}\,
H_{l_1, y^{l_2}x}^{l_1}
+
\frac{1}{6}\,
H_{l_2, y^{l_2}x}^{l_2}
+
L_{l_1, l_2, xx}^{l_1}
- 
\frac{1}{3}\,
L_{l_2, l_2, xx}^{l_2}.
\endaligned
\end{equation}
Fithly, permute the indices $(l_1, l_2) \mapsto (l_2, l_1)$:
\def\theequation{4.19}\begin{equation}
\small
\aligned
0 
\equiv
-\frac{1}{2}\,
H_{l_2, y^{l_1}x}^{l_2}
+
\frac{1}{6}\,
H_{l_1, y^{l_1}x}^{l_1}
+
L_{l_2, l_1, xx}^{l_2}
- 
\frac{1}{3}\,
L_{l_1, l_1, xx}^{l_1}.
\endaligned
\end{equation}
Finally, compute the linear combination $(4.15)+ (4.16) + 2\cdot (4.17) -
2\cdot (4.19)$:
\def\theequation{4.20}\begin{equation}
\small
\aligned
0
&
\equiv
- 
\underline{ 
2\, G_{y^{l_1}y^{l_1}}^{l_1}
}_{ \fbox{\tiny 1}}
+ 
\underline{ 
2\, G_{ y^{l_1} y^{l_2}}^{l_2}
}_{ \octagon \! \! \! \! \! \tiny{\sf a}}
+ 
\underline{ 
H_{l_1, xy^{l_1}}^{l_1}
}_{ \fbox{\tiny 2}}
- 
\underline{ 
H_{l_2, xy^{l_1}}^{l_2}
}_{ \octagon \! \! \! \! \! \tiny{\sf b}}
- \\
& \
\ \ \ \ \
-
\underline{ 
2\, G_{ y^{l_1} y^{l_2}}^{l_2}
}_{ \octagon \! \! \! \! \! \tiny{\sf a}}
+
\underline{ 
H_{l_1, xy^{l_2}}^{l_2}
}_{ \octagon \! \! \! \! \! \tiny{\sf c}}
- \\
& \
\ \ \ \ \
-
\underline{ 
H_{l_1, xy^{l_2}}^{l_2}
}_{ \octagon \! \! \! \! \! \tiny{\sf c}}
+
\underline{ 
\frac{2}{3}\, H_{l_1, l_1, xy^{l_1}}^{l_1}
}_{ \fbox{\tiny 2}}
+
\underline{ 
2\, L_{l_1, l_2, xx}^{l_2}
}_{ \octagon \! \! \! \! \! \tiny{\sf d}}
-
\underline{ 
\frac{4}{3}\, L_{l_1, l_1, xx}^{l_1}
}_{ \fbox{\tiny 3}}
+ \\
& \
\ \ \ \ \
+
\underline{ 
H_{l_2, xy^{l_1}}^{l_2}
}_{ \octagon \! \! \! \! \! \tiny{\sf b}}
- 
\underline{ 
\frac{1}{3}\, H_{l_1, xy^{l_1}}^{l_1}
}_{ \fbox{\tiny 2}}
-
\underline{ 
2\, L_{l_2, l_1, xx}^{l_2}
}_{ \octagon \! \! \! \! \! \tiny{\sf d}}
+
\underline{ 
\frac{2}{3}\, L_{l_1, l_1, xx}^{l_1}
}_{ \fbox{\tiny 3}}.
\endaligned
\end{equation}
We indeed get the desired approximate identity:
\def\theequation{4.21}\begin{equation}
0
\equiv
-
2 \, G_{ y^{l_1} y^{l_1}}^{l_1}
+
\frac{ 4}{ 3} \, H_{l_1, xy^{l_1}}^{l_1}
- 
\frac{ 2}{ 3} \, L_{l_1, l_1, xx}^{l_1}.
\end{equation}

\subsection*{ 4.22.~Complete computation}
Now that the guide is constructed, we
can achieve the complete computations.

Firstly, put $j:= l_1$ in $(3.106)$ with $l_2\neq l_1$ and
differentiate with respect to $y^{l_1}$:
\def\theequation{4.23}\begin{equation}
\small
\aligned
0
&
\equiv
- 
\underline{\underline{
2\, G_{y^{l_1}y^{l_1}}^{l_1}
+
2\, G_{y^{l_2}y^{l_1}}^{l_2}
+
H_{l_1, xy^{l_1}}^{l_1}
-
H_{l_2, xy^{l_1}}^{l_2} 
}}
+ \\
& \
\ \ \ \ \
+ 
2\, \sum_k\, G_{y^{l_1}}^k\, L_{l_1, k}^{l_1}
+
2\, \sum_k\, G^k\, L_{l_1, k, y^{l_1}}^{l_1}
-
2\, \sum_k\, G_{y^{l_1}}^k\, L_{l_2, k}^{l_2}
-
2\, \sum_k\, G^k\, L_{l_2, k, y^{l_1}}^{l_2}
- \\
& \
\ \ \ \ \
-
\frac{1}{2}\, \sum_k\, H_{l_1, y^{l_1}}^k \, H_k^{l_1}
-
\frac{1}{2}\, \sum_k\, H_{l_1}^k \, H_{k, y^{l_1}}^{l_1}
+
\frac{1}{2}\, \sum_k\, H_{l_2, y^{l_1}}^k\,H_{k}^{l_2}
+
\frac{1}{2}\, \sum_k\, H_{l_2}^k\,H_{k, y^{l_1}}^{l_2}.
\endaligned
\end{equation}
Secondly, put $j:= l_2$ in $(3.106)$ with $l_1\neq l_2$ and
differentiate with respect to $y^{l_2}$:
\def\theequation{4.24}\begin{equation}
\small
\aligned
0 
&
\equiv
- 
\underline{\underline{
2\, G_{y^{l_1}y^{l_2}}^{l_2}
+
H_{l_1, xy^{l_2}}^{l_2}
}}
+ \\
& \
\ \ \ \ \
2\, \sum_k\, G_{y^{l_2}}^k\, L_{l_1,k}^{l_2}
+
2\, \sum_k\, G^k\, L_{l_1,k,y^{l_2}}^{l_2}
-
\frac{1}{2}\,\sum_k\,H_{l_1,y^{l_2}}^k\,H_{k}^{l_2}
-
\frac{1}{2}\,\sum_k\,H_{l_1}^k\,H_{k,y^{l_2}}^{l_2}.
\endaligned
\end{equation}
Thirdly, put $j:= l_2$ in $(3.108)$ with $l_1 \neq l_2$ and
differentiate with respect to $x$:
\def\theequation{4.25}\begin{equation}
\small
\aligned
0 
&
\equiv
-
\underline{\underline{
\frac{1}{2}\, H_{l_1, y^{l_2}x}^{l_2}
+
\frac{1}{3}\, H_{l_1, y^{l_1}x}^{l_1}
+ 
L_{l_1, l_2, xx}^{l_2}
-
\frac{2}{3}\,
L_{l_1, l_1, xx}^{l_1}
}}
+ \\
& \
\ \ \ \ \
+
G_x^{l_2}\, M_{l_1,l_2}
+
G_x^{l_2} \, M_{l_1,l_2,x}
-
\frac{2}{3}\, G_x^{l_1} \, M_{l_1, l_1}
-
\frac{2}{3}\, G^{l_1} \, M_{l_1, l_1,x}
- \\
& \
\ \ \ \ \
-
\frac{1}{3}\, \sum_k\, G_x^k \, M_{l_1, k}
-
\frac{1}{3}\, \sum_k\, G^k \, M_{l_1, k, x}
-
\frac{1}{2}\,\sum_k\, H_{k,x}^{l_2}\, L_{l_1, l_2}^k
-
\frac{1}{2}\,\sum_k\, H_{k}^{l_2}\, L_{l_1, l_2,x}^k 
+ \\
& \
\ \ \ \ \
+
\frac{1}{2}\,\sum_k\, H_{l_1, x}^k \, L_{l_2, k}^{l_2}
+
\frac{1}{2}\,\sum_k\, H_{l_1}^k \, L_{l_2,k , x}^{l_2}
+
\frac{1}{3}\, \sum_k\, H_{k,x}^{l_1}\, L_{l_1, l_1}^k
+
\frac{1}{3}\, \sum_k\, H_{k}^{l_1}\, L_{l_1, l_1,x}^k
- \\
& \
\ \ \ \ \
-
\frac{1}{3}\,\sum_k\, H_{l_1, x}^k\, L_{l_1, k}^{l_1}
-
\frac{1}{3}\,\sum_k\, H_{l_1}^k\, L_{l_1, k, x}^{l_1}.
\endaligned
\end{equation}
Fourthly, put $j:= l_1$ in $(3.108)$ with $l_2 \neq l_1$,
differentiate with respect to $x$ and permute the indices $(l_1, l_2)
\mapsto (l_2, l_1)$:
\def\theequation{4.26}\begin{equation}
\small
\aligned
0
&
\equiv
-
\underline{\underline{
\frac{1}{2}\,
H_{l_2, y^{l_1}x}^{l_2}
+
\frac{1}{6}\,
H_{l_1, y^{l_1}x}^{l_1}
+
L_{l_2, l_1, xx}^{l_2}
- 
\frac{1}{3}\,
L_{l_1, l_1, xx}^{l_1}
}}
+ \\
& \
\ \ \ \ \
+
G_x^{l_2}\, M_{l_2, l_1}
+
G^{l_2}\, M_{l_2, l_1,x}
-
\frac{1}{3}\, G_x^{l_1}\, M_{l_1,l_1}
-
\frac{1}{3}\, G^{l_1}\, M_{l_1,l_1,x}
+ \\
& \
\ \ \ \ \
+
\frac{1}{3}\, \sum_k\, G_x^k\, M_{l_1,k}
+
\frac{1}{3}\, \sum_k\, G^k\, M_{l_1,k,x}
-
\frac{1}{2}\,\sum_k\, H_{k,x}^{l_2}\, L_{l_2, l_1}^k
-
\frac{1}{2}\,\sum_k\, H_k^{l_2}\, L_{l_2, l_1,x}^k
+ \\
& \
\ \ \ \ \ 
+
\frac{1}{2}\, \sum_k\, H_{l_2,x}^k\, L_{l_1, k}^{l_2}
+
\frac{1}{2}\, \sum_k\, H_{l_2}^k\, L_{l_1, k,x}^{l_2}
+
\frac{1}{6}\, \sum_k\, H_{k,x}^{l_1} \, L_{l_1, l_1}^k
+
\frac{1}{6}\, \sum_k\, H_{k}^{l_1} \, L_{l_1, l_1, x}^k
- \\
& \
\ \ \ \ \
- 
\frac{1}{6}\, \sum_k\, H_{l_1, x}^k\, L_{l_1, k}^{l_1}
-
\frac{1}{6}\, \sum_k\, H_{l_1}^k\, L_{l_1, k, x}^{l_1}.
\endaligned
\end{equation}
Finally, compute the linear combination
$(4.23)+ (4.24) + 2\cdot (4.25) - 2\cdot (4.26)$:
\def\theequation{4.27}\begin{equation}
\small
\aligned
0 
& 
= 
-
\underline{\underline{
2 \, G_{ y^{l_1} y^{l_1}}^{l_1}
+
\frac{ 4}{ 3} \, H_{l_1, xy^{l_1}}^{l_1}
- 
\frac{ 2}{ 3} \, L_{l_1, l_1, xx}^{l_1}
}}
- \\
& \
\ \ \ \ \
- 
\frac{2}{3}\, G_x^{l_1} \, M_{l_1, l_1}
-
\frac{4}{3}\, \sum_k\, G_x^k \, M_{l_1, k}
+ \\
& \
\ \ \ \ \
+ 
2 \, \sum_k\, G_{y^{l_1}}^k \, L_{l_1, k}^{l_1}
-
2 \, \sum_k\, G_{y^{l_1}}^k \, L_{l_2, k}^{l_2}
+
2\, \sum_k\, G_{y^{l_2}}^k\, L_{l_1,k}^{l_2} 
+ \\
& \
\ \ \ \ \
+
\sum_k\, H_{l_1, x}^k \, L_{l_2, k}^{l_2}
+
\frac{1}{3}\, \sum_k\, H_{k,x}^{l_1} \, L_{l_1, l_1}^k
-
\frac{1}{3}\, \sum_k\, H_{l_1, x}^k \, L_{l_1, k}^{l_1}
- 
\sum_k\, H_{l_2,x}^k \, L_{l_1, k}^{l_2}
+ \\
& \
\ \ \ \ \
+
\frac{1}{2}\, \sum_k\, H_{l_2, y^{l_1}}^k \, H_k^{l_2}
+ 
\frac{1}{2}\, \sum_k\, H_{k, y^{l_1}}^{l_2}\, H_{l_2}^k
-
\frac{1}{2}\,\sum_k\, H_{l_1, y^{l_1}}^k\, H_k^{l_1}
- \\
& \
\ \ \ \ \
-
\frac{1}{2}\,\sum_k\, H_{l_1}^k\, H_{k, y^{l_1}}^{l_1}
-
\frac{1}{2}\,\sum_k\, H_{l_1, y^{l_2}}^k\, H_{k}^{l_1}
- 
\frac{1}{2}\, \sum_k\, H_{l_1}^k\, H_{k,y^{l_2}}^{l_2}
+ \\
& \
\ \ \ \ \
+ 
\sum_k\, L_{l_2, k, x}^{l_2}\, H_{l_1}^k
+ 
\frac{1}{3}\, \sum_k\, L_{l_1, l_1, x}^k\, H_k^{l_1}
-
\frac{1}{3}\, \sum_k\, L_{l_1, k,x}^{l_1}\, H_{l_1}^k
-
\sum_k\, L_{l_1, k,x}^{l_2}\, H_{l_2}^k
+ \\
& \
\ \ \ \ \
+ 
2\,\sum_k\, L_{l_1, k, y^{l_1}}^{l_1} \, G^k
+ 
2\, \sum_k\, L_{l_1, k, y^{l_2}}^{l_2} \, G^k
- 
2\, \sum_k\, L_{l_2, k, y^{l_1}}^{l_2}\, G^k
- \\
& \
\ \ \ \ \
-
\frac{2}{3}\, M_{l_1, l_1, x} \, G^{l_1}
-
\frac{4}{3}\, \sum_k\, M_{l_1, k, x} \, G^k.
\endaligned
\end{equation}
In this partial differential relation, importantly, the second order
terms are exactly the same as in our goal~\thetag{ 4.5}.
Unfortunately, the first order and the zero order terms are not the
same.

\subsection*{4.28.~Formulation of a new goal}
Thus, in order to get rid of the second order expression $ 2 \, G_{
y^{l_1} y^{l_1}}^{l_1} + \frac{ 4}{ 3} \, H_{l_1, xy^{l_1}}^{l_1} -
\frac{ 2}{ 3} \, L_{l_1, l_1, xx}^{l_1} $, we substract: $(4.5) -
(4.27)$. In the result, we write the first order terms in a certain
way, adapted in advance to our subsequent computations. For this
substraction yielding~\thetag{ 4.29} just below, we have not
underlined the terms in~\thetag{ 4.5} and in~\thetag{ 4.27}. However,
they may be underlined with a pencil to check that the result~\thetag{
4.29} is correct. We get:
\def\theequation{4.29}\begin{equation}
\small
\aligned
0 
& 
=
?\!\!
= 
\\
&
=
- 
\underline{
\sum_k\, G_{y^{l_1}}^k\, L_{k,k}^k
+ 
G_{y^{l_1}}^{l_1} \, L_{l_1, l_1}^{l_1}
+ 
\frac{1}{2}\, \sum_k\, H_{l_1, x}^k \, L_{k,k}^k
-
\frac{1}{2}\,
H_{l_1, x}^{l_1}\, L_{l_1, l_1}^{l_1}
}
+ \\
& \
\ \ \ \ \
+
\underline{
2\, \sum_k\, G_{y^{l_1}}^k\, L_{l_2, k}^{l_2}
- 
\sum_k\, H_{l_1, x}^k \, L_{l_2, k}^{l_2}
}
- \\
& \
\ \ \ \ \
-
\underline{
2\, \sum_k\, G_{y^{l_2}}^k\, L_{l_1, k}^{l_2}
+ 
\sum_k\, H_{l_2, x}^k \, L_{l_1, k}^{l_2}
}
+ \\
& \
\ \ \ \ \ \
+
\underline{
\frac{1}{2}\, \sum_k\, H_{k, y^{l_2}}^{l_2}\, H_{l_1}^k
-
\frac{1}{3}\, \sum_k\, H_{k,y^k}^k\, H_{l_1}^k\, 
- 
\sum_k\, L_{l_2, k,x}^{l_2}\, H_{l_1}^k
+
\frac{2}{3}\,\sum_k\, L_{k,k,x}^k\, H_{l_1}^k
}
- \\
& \
\ \ \ \ \
- 
\underline{
\frac{1}{2}\,\sum_k\, H_{k,y^{l_1}}^{l_2} \, H_{l_2}^k
+ 
\sum_k\, L_{l_1, k,x}^{l_2}\, H_{l_2}^k
}
+ \\
& \
\ \ \ \ \
+
\underline{
\frac{1}{2}\,\sum_k\, H_{l_1, y^{l_2}}^k\, H_k^{l_2}
- 
\frac{1}{2}\, \sum_k\, H_{l_2,y^{l_1}}^k\, H_k^{l_2}
}
+ \\
& \
\ \ \ \ \
+
\underline{
2\, \sum_k\, L_{l_2, k, y^{l_1}}^{l_2} \, G^k
- 
2\, \sum_k\, L_{l_1, k, y^{l_2}}^{l_2} \, G^k
-
2\, \sum_k\, M_{k,l_1, x}\, G^k
}
- \\
\endaligned
\end{equation}
$$
\small
\aligned
& \
\ \ \ \ \
-
\sum_k\, \sum_p\, G^k\, H_k^p\, M_{l_1, p}
+
\frac{2}{3}\, \sum_k\, G^k\, H_{l_1}^k\, M_{k,k}
+
\frac{1}{3}\, \sum_k\, \sum_p\, G^p \,H_{l_1}^k\, M_{k,p}
- \\
& \
\ \ \ \ \
-
\sum_k\, G^k\, L_{l_1, l_1}^{l_1} \, L_{l_1, k}^{l_1}
+ 
\sum_k\, \sum_p\, G^k \, L_{k,l_1}^p\, L_{p,p}^p
- 
\frac{1}{3}\, \sum_k\, \sum_p\, H_{l_1}^k\, H_p^k\, L_{k,k}^p
+ \\
& \
\ \ \ \ \
+
\frac{1}{3}\, \sum_k\, \sum_p\, H_{l_1}^k\, H_k^p\, L_{k,p}^k
-
\frac{1}{4}\, \sum_k\,\sum_p\, H_{l_1}^k\, H_k^p \, L_{p,p}^p
+ 
\frac{1}{4}\, \sum_k\, H_{l_1}^k\, H_k^{l_1}\, L_{l_1, l_1}^{l_1}.
\endaligned
$$
We have underlined plainly the first order terms appearing
in lines $1$, $2$, $3$, $4$, $5$, $6$ and $7$.

\subsection*{ 4.30.~Reconstitution of the subgoal~\thetag{ 4.29} from
~\thetag{ 3.106}, from~\thetag{ 3.108} and from~\thetag{ 3.110}} Now,
it suffices to establish that the first order partial differential
relations~\thetag{ 4.29} for $1\leqslant l_1, l_2 \leqslant m$ and
$l_2 \neq l_1$ (crucial assumption) are a consequence of~\thetag{
3.106}, of~\thetag{ 3.108} and of~\thetag{ 3.110} by linear
combinations. The auxiliary index $l_2$, which is absent in the
goal~\thetag{ 4.5}, will disappear at the end. Differentiations will
not be applied anymore.  Also, the partial differential
relations~\thetag{ 3.96}, which were not used above, will neither be
used in the sequel. However, they are strongly used in the treatment
of the remaining three compatibility conditions $(3.112)_2$,
$(3.112)_3$ and $(3.112)_4$, the detail of which we do not copy in the
typesetted paper. Finally, the construction of a guide for the
subgoal~\thetag{ 4.29} may be guessed similarly as in~\S4.14 above. We
shall provide the final computations directly, without any guide: they
consists of the {\it seven}\, partial differential relations~\thetag{
4.32}, \thetag{ 4.33}, \thetag{ 4.35}, \thetag{ 4.37}, \thetag{ 4.39},
\thetag{ 4.43} and \thetag{ 4.45} below. At the end, we shall make the
addition~\thetag{ 4.47} below, producing the desired subgoal~\thetag{
4.29} := \thetag{4.32} + \thetag{ 4.33} + \thetag{ 4.35} + \thetag{
4.37} + \thetag{ 4.39} + \thetag{ 4.43} + \thetag{ 4.45}, with the
numerotation of terms corresponding to the order of appearance of the
terms of~\thetag{ 4.29}, as usual.

Firstly, put $j:= k$, $l_1 := l_1$ and $l_2 := l_1$ in~\thetag{
3.106}:
\def\theequation{4.31}\begin{equation}
\small
\aligned
0 
& 
= 
-
\underline{
2\, G_{y^{l_1}}^k
+
2\, \delta_{l_1}^k\,G_{y^{l_1}}^{l_1}
+
H_{l_1, x}^k 
- 
\delta_{l_1}^k\, H_{l_1, x}^{l_1}
}
+ \\
& \
\ \ \ \ \
+ 
2\, \sum_p\, G^p \, L_{l_1, p}^k
- 
2\, \delta_{l_1}^k\, \sum_p\, G^p\, L_{l_1,p}^{l_1}
-
\frac{1}{2}\, \sum_p\, H_{l_1}^p\, H_p^k 
+ 
\frac{1}{2}\, \delta_{l_1}^k\, \sum_p\, H_{l_1}^p\, H_p^{l_1}.
\endaligned
\end{equation}
Apply the operator $\frac{1}{2}\, \sum_k\, L_{k,k}^k(\cdot)$
to the preceding equality, namely 
compute $\frac{1}{2}\, \sum_k\, L_{k,k}^k \cdot (4.31)$. This yields:
\def\theequation{4.32}\begin{equation}
\small
\aligned
0 
& 
=
-
\underline{
\sum_k\, G_{y^{l_1}}^k\, L_{k,k}^k 
+ 
G_{y^{l_1}}^{l_1}
+
\frac{1}{2}\,\sum_k\, H_{l_1, x}^k\, L_{k,k}^k
-
\frac{1}{2}\, H_{l_1,x}^{l_1}\, L_{l_1, l_1}^{l_1}
}
+ \\
& \
\ \ \ \ \ 
+
\sum_k\, \sum_p \, G^p\, L_{k,k}^k\, L_{l_1, p}^k
-
\sum_p\, G^p\, L_{l_1, l_1}^{l_1} \, L_{l_1, p}^k
-
\frac{1}{4}\,\sum_k\, \sum_p\, H_{l_1}^p\, H_p^k\, L_{k,k}^k
+ \\
& \
\ \ \ \ \
+
\frac{1}{4}\, \sum_p\, H_{l_1}^p\, H_p^{l_1} \, L_{l_1, l_1}^{l_1}.
\endaligned
\end{equation}
Secondly, apply the operator $-\sum_k\, L_{l_2,k}^{l_2}(\cdot)$
to~\thetag{ 4.32}, namely compute $-\sum_k\, L_{l_2,k}^{l_2} \cdot
(4.32)$. This yields:
\def\theequation{4.33}\begin{equation}
\small
\aligned
0 
& 
= 
\underline{
2\, \sum_k\, G_{y^{l_1}}^k\, L_{l_2, k}^{l_2}
-
2\, G_{y^{l_1}}^{l_1} \, L_{l_2, l_1}^{l_2}
-
\sum_k\, H_{l_1, x}^k\, L_{l_2, k}^{l_2}
+ 
H_{l_1, x}^{l_1}\, L_{l_2, l_1}^{l_2}
}
- \\
& \
\ \ \ \ \ 
-
2\, \sum_k\, \sum_p\, G^p\, L_{l_2, k}^{l_2}\, L_{l_1,p}^k
+
2\, \sum_p\, G^p\, L_{l_2, l_1}^{l_2} \, L_{l_1, p}^{l_1}
+
\frac{1}{2}\,\sum_k\, \sum_p \, H_{l_1}^p\,H_p^k\,L_{l_2, k}^{l_2}
- \\
& \
\ \ \ \ \
-
\frac{1}{2}\,\sum_p\, H_{l_1}^p\, H_p^{l_1}\,L_{l_2, l_1}^{l_2}.
\endaligned
\end{equation}
Thirdly, put $j:= k$, $l_1:= l_2$ and $l_2 := l_1$ with $l_2\neq l_1$
in~\thetag{ 3.106}:
\def\theequation{4.34}\begin{equation}
\small
\aligned
0 
& 
= 
-
\underline{
2\, G_{y^{l_2}}^k
+
2\, \delta_{l_2}^k\, G_{y^{l_1}}^{l_1}
+
H_{l_2, x}^k
-
\delta_{l_2}^k\, H_{l_1, x}^{l_1}
}
+ \\
& \
\ \ \ \ \
+
2\, \sum_p\, G^p\, L_{l_2, p}^k
-
2\, \delta_{l_2}^k\, \sum_p\, G^p\, L_{l_1, p}^{l_1}
-
\frac{1}{2}\,\sum_p\,H_{l_2}^p\,H_p^k
+
\frac{1}{2}\,\delta_{l_2}^k\,\sum_p\, H_{l_1}^p\, H_p^{l_1}.
\endaligned
\end{equation}
Next, apply the operator $\sum_k\, L_{l_1,k}^{l_2}(\cdot)$
to~\thetag{ 4.34}, namely 
compute $\sum_k\, L_{l_1,k}^{l_2} \cdot (4.34)$. This yields:
\def\theequation{4.35}\begin{equation}
\small
\aligned
0
& 
=
-
\underline{
2\, \sum_k\,G_{y^{l_2}}^k\, L_{l_1,k}^{l_2}
+
2\, G_{y^{l_1}}^{l_1}\, L_{l_1, l_2}^{l_2}
+
\sum_k\, H_{l_2, x}^k\, L_{l_1,k}^{l_2}
-
H_{l_1, x}^{l_1}\, L_{l_1, l_2}^{l_2}
}
+ \\
& \
\ \ \ \ \
+ 
2\,\sum_k\,\sum_p\,G^p\,L_{l_2,p}^k\,L_{l_1,k}^{l_2}
-
2\,\sum_p\,G^p\,L_{l_1,p}^{l_1}\,L_{l_1,l_2}^{l_2}
-
\frac{1}{2}\,\sum_k\,\sum_p\,H_{l_2}^p\,H_p^k\,L_{l_1,k}^{l_2}
+ \\
& \
\ \ \ \ \
+
\frac{1}{2}\,\sum_p\,H_{l_1}^p\,H_p^{l_1}\,L_{l_1,l_2}^{l_2}.
\endaligned
\end{equation}
Fourthly, put $j:= l_2$, $l_1 := k$ and $l_2:= l_2$ in~\thetag{
3.108}:
\def\theequation{4.36}\begin{equation}
\small
\aligned
0
& 
=
-
\underline{
\frac{1}{2}\, H_{k,y^{l_2}}^{l_2}
+
\frac{1}{6}\,\delta_k^{l_2}\,H_{l_2,y^{l_2}}^{l_2}
+
\frac{1}{3}\,H_{k,y^k}^k
+
L_{k,l_2,x}^{l_2}
}
- \\
& \
\ \ \ \ \
-
\underline{
\frac{1}{3}\, \delta_{l_2}^k\, L_{l_2, l_2,x}^{l_2}
-
\frac{2}{3}\, L_{k,k,x}^k
}
+ \\
& \
\ \ \ \ \
+ 
G^{l_2}\, M_{k,l_2} 
-
\frac{1}{3}\,\delta_{l_2}^k\, G^{l_2} \, M_{l_2, l_2}
-
\frac{2}{3}\, G^k\, M_{k,k}
+
\frac{1}{3}\, \delta_k^{l_2} \,\sum_p\,G^p\,M_{l_2,p}
- \\
& \
\ \ \ \ \
-
\frac{1}{3}\, \sum_p\, G^p\, M_{k,p}
-
\frac{1}{2}\,\sum_p\, H_p^{l_2}\,L_{k,l_2}^p
+
\frac{1}{2}\,\sum_p\, H_k^p\,L_{l_2,p}^{l_2}
+
\frac{1}{6}\, \delta_k^{l_2}\,\sum_p\,H_p^{l_2}\,L_{l_2,l_2}^p
- \\
& \
\ \ \ \ \
-
\frac{1}{6}\,\delta_k^{l_2}\, \sum_p\,H_{l_2}^p\,L_{l_2,p}^{l_2}
+
\frac{1}{3}\,\sum_p\,H_p^k\,L_{k,k}^p
-
\frac{1}{3}\,\sum_p\,H_k^p\, L_{k,p}^k.
\endaligned
\end{equation}
Next, apply the operator $-\sum_k\, H_{l_1}^k(\cdot)$ to~\thetag{
4.36}, namely compute $-\sum_k\, H_{l_1}^k \cdot (4.36)$. This yields:
\def\theequation{4.37}\begin{equation}
\small
\aligned
0 
& 
=
\underline{
\frac{1}{2}\,\sum_k\, H_{k,y^{l_2}}^{l_2}\,H_{l_1}^k
-
\frac{1}{6}\, H_{l_2, y^{l_2}}^{l_2}\, H_{l_1}^{l_2}
-
\frac{1}{3}\,\sum_k\,H_{k,y^k}^k\,H_{l_1}^k
}
- \\
& \
\ \ \ \ \
-
\underline{
\sum_k\, L_{k,l_2,x}^{l_2}\, H_{l_1}^k
+
\frac{1}{3}\, L_{l_2,l_2,x}^{l_2}\,H_{l_1}^{l_2}
+
\frac{2}{3}\,\sum_k\,L_{k,k,x}^k\,H_{l_1}^k
}
- \\
& \
\ \ \ \ \ 
-
\sum_k\, G^{l_2}\, H_{l_1}^k\, M_{k,l_2}
+
\frac{1}{3}\, G^{l_2}\, H_{l_1}^{l_2}\,M_{l_2,l_2}
+
\frac{2}{3}\, \sum_k\, G^k\,H_{l_1}^k\, M_{k,k}
- \\
& \
\ \ \ \ \
-
\frac{1}{3}\, \sum_p\,G^p\,H_{l_1}^{l_2}\,M_{l_2,p}
+
\frac{1}{3}\,\sum_k\,\sum_p\,G^pH_{l_1}^k\,M_{k,p}
+
\frac{1}{2}\,\sum_k\,\sum_p\,H_{l_1}^k\,H_p^{l_2}\,L_{k,l_2}^p
- \\
& \
\ \ \ \ \
-
\frac{1}{2}\,\sum_k\,\sum_p\, H_{l_1}^k\,H_k^p\,L_{l_2,p}^{l_2}
- 
\frac{1}{6}\,\sum_p\,H_{l_1}^{l_2}\,H_p^{l_2}\,L_{l_2,l_2}^p
+
\frac{1}{6}\,\sum_p\,H_{l_1}^{l_2}\,H_{l_2}^p\,L_{l_2,p}^{l_2}
- \\
& \
\ \ \ \ \
-
\frac{1}{3}\,\sum_k\,\sum_p\,H_{l_1}^k\,H_p^k\,L_{k,k}^p
+
\frac{1}{3}\,\sum_k\,\sum_p\,H_{l_1}^k\,H_k^p\,L_{k,p}^k.
\endaligned
\end{equation}
Fithly, put $j:= l_2$, $l_1 := k$ and
$l_2:= l_1$ in~\thetag{ 3.108}:
\def\theequation{4.38}\begin{equation}
\small
\aligned
0
& 
=
-
\underline{
\frac{1}{2}\,H_{k,y^{l_1}}^{l_2}
+
\frac{1}{6}\,\delta_k^{l_2}\,H_{l_1,y^{l_1}}^{l_1}
+
L_{k,l_1,x}^{l_2}
-
\frac{1}{3}\, \delta_k^{l_2}\,L_{l_1,l_1,x}^{l_1}
}
+ \\
& \
\ \ \ \ \
+
G^{l_2}\, M_{k,l_1}
-
\frac{1}{3}\, \delta_k^{l_2}\, G^{l_1}\, M_{l_1, l_1}
+
\frac{1}{3}\, \delta_k^{l_2} \, \sum_p\,G^p\,M_{l_1,p}
-
\frac{1}{2}\,\sum_p\,H_p^{l_2}\,L_{k,l_1}^p
+ \\
& \
\ \ \ \ \
+
\frac{1}{2}\,\sum_p\, H_k^p\,L_{l_1,p}^{l_2}
+
\frac{1}{6}\,\delta_k^{l_2}\,\sum_p\, H_p^{l_1}\,L_{l_1,l_1}^p
-
\frac{1}{6}\,\delta_k^{l_2}\,\sum_p\,H_{l_1}^p\,L_{l_1,p}^{l_1}.
\endaligned
\end{equation}
Next, apply the operator $\sum_k\, H_{l_2}^k(\cdot)$ to~\thetag{
4.38}, namely compute $\sum_k\, H_{l_2}^k\cdot (4.38)$. This yields:
\def\theequation{4.39}\begin{equation}
\small
\aligned
0 
& 
=
-
\underline{
\frac{1}{2}\,\sum_k\,H_{k,y^{l_1}}^{l_2}\,H_{l_2}^k
+
\frac{1}{6}\,H_{l_1,y^{l_1}}^{l_1}\,H_{l_2}^{l_2}
+
\sum_k\,L_{k,l_1,x}^{l_2}\,H_{l_2}^k
}
- \\
& \
\ \ \ \ \
-
\underline{
\frac{1}{3}\, L_{l_1, l_1,x}^{l_1}\, H_{l_2}^{l_2}
}
+ \\
& \
\ \ \ \ \
+
\sum_k\,G^{l_2}\, H_{l_2}^k\,M_{k,l_1}
-
\frac{1}{3}\,G^{l_1}\,H_{l_2}^{l_2}\, M_{l_1,l_1}
+
\frac{1}{3}\,\sum_p\, G^p\,H_{l_2}^{l_2}\,M_{l_1,p}
- \\
& \
\ \ \ \ \
-
\frac{1}{2}\,\sum_k\,\sum_p\,H_{l_2}^k\,H_p^{l_2}\,L_{k,l_1}^p
+
\frac{1}{2}\,\sum_k\,\sum_p\,H_{l_2}^k\,H_k^p\,L_{l_1,p}^{l_2}
+
\frac{1}{6}\,\sum_p\,H_{l_2}^{l_2}\,H_p^{l_1}\,L_{l_1,l_1}^p
- \\
& \
\ \ \ \ \
-
\frac{1}{6}\,\sum_p\,H_{l_2}^{l_2}\,H_{l_1}^p\,L_{l_1,p}^{l_1}.
\endaligned
\end{equation}
Sixthly, we form the expression:
\def\theequation{4.40}\begin{equation}
\aligned
\left.
(3.108)
\right\vert_{
j:=k; \ l_1:=l_1; \ l_2:=l_2}
-
\left.
(3.108)
\right\vert_{
j:=k; \ l_1:=l_2; \ l_2:=l_1}.
\endaligned
\end{equation}
Writing term by term the substractions, we get:
\def\theequation{4.41}\begin{equation}
\small
\aligned
0
& 
=
-
\underline{ 
\frac{1}{2}\,H_{l_1,y^{l_2}}^k
}_{ \fbox{\tiny 1}}
+
\underline{ 
\frac{1}{2}\,H_{l_2,y^{l_1}}^k
}_{ \fbox{\tiny 2}}
+
\underline{ 
\frac{1}{6}\,\delta_{l_1}^k\,H_{l_2,y^{l_2}}^{l_2}
}_{ \fbox{\tiny 3}}
-
\underline{ 
\frac{1}{6}\,\delta_{l_2}^k\,H_{l_1,y^{l_1}}^{l_1}
}_{ \fbox{\tiny 4}}
+ \\
& \
\ \ \ \ \
+
\underline{ 
\frac{1}{3}\, \delta_{l_2}^k\, H_{l_1, y^{l_1}}^{l_1}
}_{ \fbox{\tiny 4}}
-
\underline{ 
\frac{1}{3}\, \delta_{l_1}^k\, H_{l_2, y^{l_2}}^{l_2}
}_{ \fbox{\tiny 3}}
+
\underline{
L_{l_1,l_2,x}^k 
}_{ \octagon \! \! \! \! \! \tiny{\sf a}}
-
\underline{
L_{l_2,l_1,x}^k
}_{ \octagon \! \! \! \! \! \tiny{\sf a}}
- \\
& \
\ \ \ \ \
-
\underline{
\frac{1}{3}\, \delta_{l_1}^k\, L_{l_2,l_2,x}^{l_2}
}_{ \fbox{\tiny 5}}
+
\underline{
\frac{1}{3}\, \delta_{l_2}^k\,L_{l_1,l_1,x}^{l_1}
}_{ \fbox{\tiny 6}}
-
\underline{
\frac{2}{3}\,\delta_{l_2}^k\,L_{l_1,l_1,x}^{l_1}
}_{ \fbox{\tiny 6}}
+
\underline{
\frac{2}{3}\,\delta_{l_1}^k\,L_{l_2,l_2,x}^{l_2}
}_{ \fbox{\tiny 5}}
+ \\
& \
\ \ \ \ \
+
\underline{
G^k\, M_{l_1,l_2}
}_{ \octagon \! \! \! \! \! \tiny{\sf b}}
-
\underline{
G^k\, M_{l_2,l_1}
}_{ \octagon \! \! \! \! \! \tiny{\sf b}}
-
\underline{
\frac{1}{3}\,\delta_{l_1}^k\,G^{l_2}\,M_{l_2,l_2}
}_{ \fbox{\tiny 7}}
+
\underline{
\frac{1}{3}\,\delta_{l_2}^k\,G^{l_1}\,M_{l_1,l_1}
}_{ \fbox{\tiny 8}}
- \\
& \
\ \ \ \ \
-
\underline{
\frac{2}{3}\, \delta_{l_2}^k\,G^{l_1}\,M_{l_1,l_1}
}_{ \fbox{\tiny 8}}
+
\underline{
\frac{2}{3}\, \delta_{l_1}^k\,G^{l_2}\,M_{l_2,l_2}
}_{ \fbox{\tiny 7}}
+
\underline{
\frac{1}{3}\, \delta_{l_1}^k\, \sum_p\, G^p\, M_{l_2,p}
}_{ \fbox{\tiny 9}}
-
\underline{
\frac{1}{3}\, \delta_{l_2}^k\, \sum_p\, G^p\, M_{l_1,p}
}_{ \fbox{\tiny 10}}
- \\
& \
\ \ \ \ \
-
\underline{
\frac{1}{3}\, \delta_{l_2}^k\, \sum_p\,G^p\,M_{l_1,p}
}_{ \fbox{\tiny 10}}
+
\underline{
\frac{1}{3}\, \delta_{l_1}^k\, \sum_p\,G^p\,M_{l_2,p}
}_{ \fbox{\tiny 9}}
-
\underline{
\frac{1}{2}\,\sum_p\,H_p^k\,L_{l_1,l_2}^p
}_{ \octagon \! \! \! \! \! \tiny{\sf c}}
+
\underline{
\frac{1}{2}\,\sum_p\,H_p^k\,L_{l_2,l_1}^p
}_{ \octagon \! \! \! \! \! \tiny{\sf c}}
+ \\
& \
\ \ \ \ \
+
\underline{
\frac{1}{2}\,\sum_p\, H_{l_1}^p\, L_{l_2,p}^k
}_{ \fbox{\tiny 11}}
-
\underline{
\frac{1}{2}\,\sum_p\, H_{l_2}^p\, L_{l_1,p}^k
}_{ \fbox{\tiny 12}}
+
\underline{
\frac{1}{6}\,\delta_{l_1}^k\, \sum_p\,H_p^{l_2}\,L_{l_2,l_2}^p
}_{ \fbox{\tiny 13}}
-
\underline{
\frac{1}{6}\,\delta_{l_2}^k\, \sum_p\,H_p^{l_1}\,L_{l_1,l_1}^p
}_{ \fbox{\tiny 14}}
- \\
& \
\ \ \ \ \
-
\underline{
\frac{1}{6}\,\delta_{l_1}^k\, \sum_p\, H_{l_2}^p\,L_{l_2,p}^{l_2}
}_{ \fbox{\tiny 15}}
+
\underline{
\frac{1}{6}\,\delta_{l_2}^k\, \sum_p\, H_{l_1}^p\,L_{l_1,p}^{l_1}
}_{ \fbox{\tiny 16}}
+
\underline{
\frac{1}{3}\,\delta_{l_2}^k\, \sum_p\, H_p^{l_1}\,L_{l_1,l_1}^p
}_{ \fbox{\tiny 14}}
- \\
& \
\ \ \ \ \
-
\underline{
\frac{1}{3}\,\delta_{l_1}^k\, \sum_p\, H_p^{l_2}\,L_{l_2,l_2}^p
}_{ \fbox{\tiny 13}}
-
\underline{
\frac{1}{3}\,\delta_{l_2}^k\,\sum_p\,H_{l_1}^p\,L_{l_1,p}^{l_1}
}_{ \fbox{\tiny 16}}
+
\underline{
\frac{1}{3}\,\delta_{l_1}^k\,\sum_p\,H_{l_2}^p\,L_{l_2,p}^{l_2}
}_{ \fbox{\tiny 15}} \ .
\endaligned
\end{equation}
Simplifying, we get:
\def\theequation{4.42}\begin{equation}
\small
\aligned
0 
& 
=
-
\underline{
\frac{1}{2}\,H_{l_1,y^{l_2}}^k
+
\frac{1}{2}\,H_{l_2,y^{l_1}}^k
-
\frac{1}{6}\,\delta_{l_1}^k\, H_{l_2, y^{l_2}}^{l_2}
+
\frac{1}{6}\,\delta_{l_2}^k\, H_{l_1, y^{l_1}}^{l_1}
}
+ \\
& \
\ \ \ \ \
+ 
\underline{
\frac{1}{3}\, \delta_{l_1}^k\, L_{l_2,l_2,x}^{l_2}
-
\frac{1}{3}\, \delta_{l_2}^k\, L_{l_1,l_1,x}^{l_1}
}
+ \\
& \
\ \ \ \ \
+
\frac{1}{3}\, \delta_{l_1}^k\,G^{l_2}\,M_{l_2,l_2}
-
\frac{1}{3}\, \delta_{l_2}^k\,G^{l_1}\,M_{l_1,l_1}
+
\frac{2}{3}\, \delta_{l_1}^k\,\sum_p\,G^p\,M_{l_2,p}
-
\frac{2}{3}\, \delta_{l_2}^k\,\sum_p\,G^p\,M_{l_1,p}
+ \\
& \
\ \ \ \ \
+
\frac{1}{2}\,\sum_p\,H_{l_1}^p\,L_{l_2,p}^k
-
\frac{1}{2}\,\sum_p\,H_{l_2}^p\,L_{l_1,p}^k
-
\frac{1}{6}\,\delta_{l_1}^k\,\sum_p\,H_p^{l_2}\,L_{l_2,l_2}^p
+
\frac{1}{6}\,\delta_{l_2}^k\,\sum_p\,H_p^{l_1}\,L_{l_1,l_1}^p
+ \\
& \
\ \ \ \ \
+
\frac{1}{6}\,\delta_{l_1}^k\,\sum_p\,H_{l_2}^p\,L_{l_2,p}^{l_2}
-
\frac{1}{6}\,\delta_{l_2}^k\,\sum_p\,H_{l_1}^p\,L_{l_1,p}^{l_1}.
\endaligned
\end{equation}
Next, apply the operator $-\sum_k\, H_k^{l_2}(\cdot)$ to~\thetag{
4.42}, namely compute $-\sum_k\, H_k^{l_2} \cdot (4.42)$. This yields:
\def\theequation{4.43}\begin{equation}
\small
\aligned
0 
& 
=
\underline{
\frac{1}{2}\,\sum_k\,H_{l_1,y^{l_2}}^k\,H_k^{l_2}
-
\frac{1}{2}\,\sum_k\,H_{l_2,y^{l_1}}^k\,H_k^{l_2}
+
\frac{1}{6}\,H_{l_2,y^{l_2}}^{l_2}\,H_{l_1}^{l_2}
}
- \\
& \
\ \ \ \ \
-
\underline{
\frac{1}{6}\,H_{l_1,y^{l_1}}^{l_1}\,H_{l_2}^{l_2}
-
\frac{1}{3}\,L_{l_2,l_2,x}^{l_2}\,H_{l_1}^{l_2}
+
\frac{1}{3}\,L_{l_1,l_1,x}^{l_1}\,H_{l_2}^{l_2}
}
- \\
& \
\ \ \ \ \
-
\frac{1}{3}\,G^{l_2}\,H_{l_1}^{l_2}\,M_{l_2,l_2}
+
\frac{1}{3}\,G^{l_1}\,H_{l_2}^{l_2}\,M_{l_1,l_1}
-
\frac{2}{3}\,\sum_p\,G^p\,H_{l_1}^{l_2}\,M_{l_2,p}
+ \\
& \
\ \ \ \ \
+
\frac{2}{3}\,\sum_p\,G^p\,H_{l_2}^{l_2}\,M_{l_1,p}
-
\frac{1}{2}\,\sum_k\,\sum_p\,H_k^{l_2}\,H_{l_1}^p\,L_{l_2,p}^k
+
\frac{1}{2}\,\sum_k\,\sum_p\,H_k^{l_2}\,H_{l_2}^p\,L_{l_1,p}^k
+ \\
& \
\ \ \ \ \
+
\frac{1}{6}\,\sum_p\,H_{l_1}^{l_2}\,H_p^{l_2}\,L_{l_2,l_2}^p
-
\frac{1}{6}\,\sum_p\,H_{l_2}^{l_2}\,H_p^{l_1}\,L_{l_1,l_1}^p
-
\frac{1}{6}\,\sum_p\,H_{l_1}^{l_2}\,H_{l_2}^p\,L_{l_2,p}^{l_2}
+ \\
& \
\ \ \ \ \
+
\frac{1}{6}\,\sum_p\,H_{l_2}^{l_2}\,H_{l_1}^p\,L_{l_1,p}^{l_1}.
\endaligned
\end{equation}
Seventhly, put $j:= l_2$, $l_1:= k$, $l_2:= l2$
and $l_3:= l_1$ in~\thetag{ 3.108}:
\def\theequation{4.44}\begin{equation}
\small
\aligned
0
& 
=
\underline{
L_{k,l_2,y^{l_1}}^{l_2}
-
L_{k,l_1,y^{l_2}}^{l_2}
-
M_{k,l_1,x}
}
+ \\
& \
\ \ \ \ \
+
\frac{1}{2}\,H_{l_1}^{l_2}\,M_{k,l_2}
-
\frac{1}{2}\,H_{l_2}^{l_2}\,M_{k,l_1}
+
\frac{1}{2}\,\delta_k^{l_2}\,\sum_p\,H_{l_1}^p\,M_{l_2,p}
-
\frac{1}{2}\,\delta_k^{l_2}\,\sum_p\,H_{l_2}^p\,M_{l_1,p}
- \\
& \
\ \ \ \ \
-
\frac{1}{2}\,\sum_p\,H_k^p\,M_{l_1,p}
+
\sum_p\,L_{k,l_1}^p\,L_{l_2,p}^{l_2}
-
\sum_p\,L_{k,l_2}^p\,L_{l_1,p}^{l_2},
\endaligned
\end{equation}
and then apply the operator $2\, \sum_k\, G^k(\cdot)$:
\def\theequation{4.45}\begin{equation}
\small
\aligned
0
& 
=
2\, \sum_k\,L_{k,l_2,y^{l_1}}^{l_2}\,G^k
-
2\,\sum_k\,L_{k,l_1,y^{l_2}}^{l_2}\,G^k
-
2\,\sum_k\,M_{k,l_1,x}\,G^k
+ \\
& \
\ \ \ \ \
+
\sum_k\,G^k\,H_{l_1}^{l_2}\,M_{l_1,l_2}
-
\sum_k\,G^k\,H_{l_2}^{l_2}\,M_{k,l_1}
+
\sum_p\,G^{l_2}\,H_{l_1}^p\,M_{l_2,p}
- \\
& \
\ \ \ \ \ 
-
\sum_p\, G^{l_2}\,H_{l_2}^p\,M_{l_1,p}
-
\sum_k\,\sum_p\,G^k\,H_k^p\,M_{l_1,p}
+
2\,\sum_k\,\sum_p\,G^k\,L_{k,l_1}^p\,L_{l_2,p}^{l_2}
- \\
& \
\ \ \ \ \ 
-
2\,\sum_k\,\sum_p\,G^k\,L_{k,l_2}^p\,L_{l_1,p}^{l_2}.
\endaligned
\end{equation}
Finally, achieve the addition
\def\theequation{4.46}\begin{equation}
(4.32)
+
(4.33)
+
(4.35)
+
(4.37)
+
(4.39)
+
(4.43)
+
(4.45).
\end{equation}
We copy these seven formal expression, we underline the vanishing
terms and we number the remaining terms so as to respect the order of
appearance of the terms of the subgoal~\thetag{ 4.29}:
\def\theequation{4.47}\begin{equation}
\small
\aligned
0 
& 
=
-
\underline{
\sum_k\, G_{y^{l_1}}^k\, L_{k,k}^k 
}_{ \fbox{\tiny 1}}
+ 
\underline{
G_{y^{l_1}}^{l_1}
}_{ \fbox{\tiny 2}}
+
\underline{
\frac{1}{2}\,\sum_k\, H_{l_1, x}^k\, L_{k,k}^k
}_{ \fbox{\tiny 3}}
-
\underline{
\frac{1}{2}\, H_{l_1,x}^{l_1}\, L_{l_1, l_1}^{l_1}
}_{ \fbox{\tiny 4}}
+ \\
& \
\ \ \ \ \ 
+
\underline{
\sum_k\, \sum_p \, G^p\, L_{k,k}^k\, L_{l_1, p}^k
}_{ \fbox{\tiny 24}}
-
\underline{
\sum_p\, G^p\, L_{l_1, l_1}^{l_1} \, L_{l_1, p}^k
}_{ \fbox{\tiny 23}}
-
\underline{
\frac{1}{4}\,\sum_k\, \sum_p\, H_{l_1}^p\, H_p^k\, L_{k,k}^k
}_{ \fbox{\tiny 27}}
+ \\
& \
\ \ \ \ \
+
\underline{
\frac{1}{4}\, \sum_p\, H_{l_1}^p\, H_p^{l_1} \, L_{l_1, l_1}^{l_1}
}_{ \fbox{\tiny 28}}
+ \\
\endaligned
\end{equation}
$$
\small
\aligned
& 
+
\underline{
2\, \sum_k\, G_{y^{l_1}}^k\, L_{l_2, k}^{l_2} 
}_{ \fbox{\tiny 5}}
-
\underline{
2\, G_{y^{l_1}}^{l_1} \, L_{l_2, l_1}^{l_2}
}_{ \octagon \! \! \! \! \! \tiny{\sf a}}
-
\underline{
\sum_k\, H_{l_1, x}^k\, L_{l_2, k}^{l_2} 
}_{ \fbox{\tiny 6}}
+ 
\underline{
H_{l_1, x}^{l_1}\, L_{l_2, l_1}^{l_2}
}_{ \octagon \! \! \! \! \! \tiny{\sf b}}
- \\
& \
\ \ \ \ \ 
-
\underline{
2\, \sum_k\, \sum_p\, G^p\, L_{l_2, k}^{l_2}\, L_{l_1,p}^k
}_{ \octagon \! \! \! \! \! \tiny{\sf g}}
+
\underline{
2\, \sum_p\, G^p\, L_{l_2, l_1}^{l_2} \, L_{l_1, p}^{l_1}
}_{ \octagon \! \! \! \! \! \tiny{\sf h}}
+
\underline{
\frac{1}{2}\,\sum_k\, \sum_p \, H_{l_1}^p\,H_p^k\,L_{l_2, k}^{l_2}
}_{ \octagon \! \! \! \! \! \tiny{\sf i}}
- \\
& \
\ \ \ \ \
-
\underline{
\frac{1}{2}\,\sum_p\, H_{l_1}^p\, H_p^{l_1}\,L_{l_2, l_1}^{l_2}
}_{ \octagon \! \! \! \! \! \tiny{\sf j}} 
- \\
\endaligned
$$
$$
\small
\aligned
& 
-
\underline{
2\, \sum_k\,G_{y^{l_2}}^k\, L_{l_1,k}^{l_2}
}_{ \fbox{\tiny 7}}
+
\underline{
2\, G_{y^{l_1}}^{l_1}\, L_{l_1, l_2}^{l_2} 
}_{ \octagon \! \! \! \! \! \tiny{\sf a}}
+
\underline{
\sum_k\, H_{l_2, x}^k\, L_{l_1,k}^{l_2}
}_{ \fbox{\tiny 8}}
-
\underline{
H_{l_1, x}^{l_1}\, L_{l_1, l_2}^{l_2} 
}_{ \octagon \! \! \! \! \! \tiny{\sf b}} 
+ \\
& \
\ \ \ \ \
+ 
\underline{
2\,\sum_k\,\sum_p\,G^p\,L_{l_2,p}^k\,L_{l_1,k}^{l_2}
}_{ \octagon \! \! \! \! \! \tiny{\sf k}}
-
\underline{
2\,\sum_p\,G^p\,L_{l_1,p}^{l_1}\,L_{l_1,l_2}^{l_2}
}_{ \octagon \! \! \! \! \! \tiny{\sf h}}
-
\underline{
\frac{1}{2}\,\sum_k\,\sum_p\,H_{l_2}^p\,H_p^k\,L_{l_1,k}^{l_2}
}_{ \octagon \! \! \! \! \! \tiny{\sf l}}
+ \\
& \
\ \ \ \ \
+
\underline{
\frac{1}{2}\,\sum_p\,H_{l_1}^p\,H_p^{l_1}\,L_{l_1,l_2}^{l_2}
}_{ \octagon \! \! \! \! \! \tiny{\sf j}} 
+ \\
\endaligned
$$
$$
\small
\aligned
& 
+
\underline{
\frac{1}{2}\,\sum_k\, H_{k,y^{l_2}}^{l_2}\,H_{l_1}^k
}_{ \fbox{\tiny 9}}
-
\underline{
\frac{1}{6}\, H_{l_2, y^{l_2}}^{l_2}\, H_{l_1}^{l_2}
}_{ \octagon \! \! \! \! \! \tiny{\sf c}}
-
\underline{
\frac{1}{3}\,\sum_k\,H_{k,y^k}^k\,H_{l_1}^k
}_{ \fbox{\tiny 10}}
- \\
& \
\ \ \ \ \
-
\underline{
\sum_k\, L_{k,l_2,x}^{l_2}\, H_{l_1}^k
}_{ \fbox{\tiny 11}}
+
\underline{
\frac{1}{3}\, L_{l_2,l_2,x}^{l_2}\,H_{l_1}^{l_2}
}_{ \octagon \! \! \! \! \! \tiny{\sf d}}
+
\underline{
\frac{2}{3}\,\sum_k\,L_{k,k,x}^k\,H_{l_1}^k
}_{ \fbox{\tiny 12}}
- \\
& \
\ \ \ \ \ 
-
\underline{
\sum_k\, G^{l_2}\, H_{l_1}^k\, M_{k,l_2}
}_{ \octagon \! \! \! \! \! \! \tiny{\sf m}}
+
\underline{
\frac{1}{3}\, G^{l_2}\, H_{l_1}^{l_2}\,M_{l_2,l_2}
}_{ \octagon \! \! \! \! \! \tiny{\sf n}}
+
\underline{
\frac{2}{3}\, \sum_k\, G^k\,H_{l_1}^k\, M_{k,k}
}_{ \fbox{\tiny 21}}
- \\
& \
\ \ \ \ \
-
\underline{
\frac{1}{3}\, \sum_p\,G^p\,H_{l_1}^{l_2}\,M_{l_2,p}
}_{ \octagon \! \! \! \! \! \tiny{\sf o}}
+
\underline{
\frac{1}{3}\,\sum_k\,\sum_p\,G^pH_{l_1}^k\,M_{k,p}
}_{ \fbox{\tiny 22}}
+
\underline{
\frac{1}{2}\,\sum_k\,\sum_p\,H_{l_1}^k\,H_p^{l_2}\,L_{k,l_2}^p
}_{ \octagon \! \! \! \! \! \tiny{\sf p}}
- \\
& \
\ \ \ \ \
-
\underline{
\frac{1}{2}\,\sum_k\,\sum_p\, H_{l_1}^k\,H_k^p\,L_{l_2,p}^{l_2}
}_{ \octagon \! \! \! \! \! \tiny{\sf i}}
- 
\underline{
\frac{1}{6}\,\sum_p\,H_{l_1}^{l_2}\,H_p^{l_2}\,L_{l_2,l_2}^p
}_{ \octagon \! \! \! \! \! \tiny{\sf q}}
+
\underline{
\frac{1}{6}\,\sum_p\,H_{l_1}^{l_2}\,H_{l_2}^p\,L_{l_2,p}^{l_2}
}_{ \octagon \! \! \! \! \! \tiny{\sf r}}
- \\
& \
\ \ \ \ \
-
\underline{
\frac{1}{3}\,\sum_k\,\sum_p\,H_{l_1}^k\,H_p^k\,L_{k,k}^p
}_{ \fbox{\tiny 25}}
+
\underline{
\frac{1}{3}\,\sum_k\,\sum_p\,H_{l_1}^k\,H_k^p\,L_{k,p}^k
}_{ \fbox{\tiny 26}} 
- \\
\endaligned
$$
$$
\aligned
&
-
\underline{
\frac{1}{2}\,\sum_k\,H_{k,y^{l_1}}^{l_2}\,H_{l_2}^k 
}_{ \fbox{\tiny 13}}
+
\underline{
\frac{1}{6}\,H_{l_1,y^{l_1}}^{l_1}\,H_{l_2}^{l_2}
}_{ \octagon \! \! \! \! \! \tiny{\sf e}}
+
\underline{
\sum_k\,L_{k,l_1,x}^{l_2}\,H_{l_2}^k
}_{ \fbox{\tiny 14}}
- \\
& \
\ \ \ \ \
-
\underline{
\frac{1}{3}\, L_{l_1, l_1,x}^{l_1}\, H_{l_2}^{l_2}
}_{ \octagon \! \! \! \! \! \tiny{\sf f}}
+ \\
& \
\ \ \ \ \
+
\underline{
\sum_k\,G^{l_2}\, H_{l_2}^k\,M_{k,l_1}
}_{ \octagon \! \! \! \! \! \tiny{\sf s}}
-
\underline{
\frac{1}{3}\,G^{l_1}\,H_{l_2}^{l_2}\, M_{l_1,l_1}
}_{ \octagon \! \! \! \! \! \tiny{\sf t}}
+
\underline{
\frac{1}{3}\,\sum_p\, G^p\,H_{l_2}^{l_2}\,M_{l_1,p}
}_{ \octagon \! \! \! \! \! \tiny{\sf u}}
- \\
& \
\ \ \ \ \
-
\underline{
\frac{1}{2}\,\sum_k\,\sum_p\,H_{l_2}^k\,H_p^{l_2}\,L_{k,l_1}^p
}_{ \octagon \! \! \! \! \! \tiny{\sf v}}
+
\underline{
\frac{1}{2}\,\sum_k\,\sum_p\,H_{l_2}^k\,H_k^p\,L_{l_1,p}^{l_2}
}_{ \octagon \! \! \! \! \! \tiny{\sf l}}
+
\underline{
\frac{1}{6}\,\sum_p\,H_{l_2}^{l_2}\,H_p^{l_1}\,L_{l_1,l_1}^p
}_{ \octagon \! \! \! \! \! \! \tiny{\sf w}}
- \\
& \
\ \ \ \ \
-
\underline{
\frac{1}{6}\,\sum_p\,H_{l_2}^{l_2}\,H_{l_1}^p\,L_{l_1,p}^{l_1}
}_{ \octagon \! \! \! \! \! \tiny{\sf x}}
+
\endaligned
$$
$$
\small
\aligned
& 
+
\underline{ 
\frac{1}{2}\,\sum_k\,H_{l_1,y^{l_2}}^k\,H_k^{l_2}
}_{ \fbox{\tiny 15}}
-
\underline{ 
\frac{1}{2}\,\sum_k\,H_{l_2,y^{l_1}}^k\,H_k^{l_2}
}_{ \fbox{\tiny 16}}
+
\underline{ 
\frac{1}{6}\,H_{l_2,y^{l_2}}^{l_2}\,H_{l_1}^{l_2}
}_{ \octagon \! \! \! \! \! \tiny{\sf c}}
- \\
& \
\ \ \ \ \
-
\underline{
\frac{1}{6}\,H_{l_1,y^{l_1}}^{l_1}\,H_{l_2}^{l_2}
}_{ \octagon \! \! \! \! \! \tiny{\sf e}}
-
\underline{ 
\frac{1}{3}\,L_{l_2,l_2,x}^{l_2}\,H_{l_1}^{l_2}
}_{ \octagon \! \! \! \! \! \tiny{\sf d}}
+
\underline{ 
\frac{1}{3}\,L_{l_1,l_1,x}^{l_1}\,H_{l_2}^{l_2}
}_{ \octagon \! \! \! \! \! \tiny{\sf f}}
- \\
& \
\ \ \ \ \
-
\underline{
\frac{1}{3}\,G^{l_2}\,H_{l_1}^{l_2}\,M_{l_2,l_2}
}_{ \octagon \! \! \! \! \! \tiny{\sf n}}
+
\underline{
\frac{1}{3}\,G^{l_1}\,H_{l_2}^{l_2}\,M_{l_1,l_1}
}_{ \octagon \! \! \! \! \! \tiny{\sf t}}
-
\underline{
\frac{2}{3}\,\sum_p\,G^p\,H_{l_1}^{l_2}\,M_{l_2,p}
}_{ \octagon \! \! \! \! \! \tiny{\sf o}}
+ \\
& \
\ \ \ \ \
+
\underline{
\frac{2}{3}\,\sum_p\,G^p\,H_{l_2}^{l_2}\,M_{l_1,p}
}_{ \octagon \! \! \! \! \! \tiny{\sf u}}
-
\underline{
\frac{1}{2}\,\sum_k\,\sum_p\,H_k^{l_2}\,H_{l_1}^p\,L_{l_2,p}^k
}_{ \octagon \! \! \! \! \! \tiny{\sf p}}
+
\underline{
\frac{1}{2}\,\sum_k\,\sum_p\,H_k^{l_2}\,H_{l_2}^p\,L_{l_1,p}^k
}_{ \octagon \! \! \! \! \! \tiny{\sf v}}
+ \\
& \
\ \ \ \ \
+
\underline{
\frac{1}{6}\,\sum_p\,H_{l_1}^{l_2}\,H_p^{l_2}\,L_{l_2,l_2}^p
}_{ \octagon \! \! \! \! \! \tiny{\sf q}}
-
\underline{
\frac{1}{6}\,\sum_p\,H_{l_2}^{l_2}\,H_p^{l_1}\,L_{l_1,l_1}^p
}_{ \octagon \! \! \! \! \! \! \tiny{\sf w}}
-
\underline{
\frac{1}{6}\,\sum_p\,H_{l_1}^{l_2}\,H_{l_2}^p\,L_{l_2,p}^{l_2}
}_{ \octagon \! \! \! \! \! \tiny{\sf r}}
+ \\
& \
\ \ \ \ \
+
\underline{
\frac{1}{6}\,\sum_p\,H_{l_2}^{l_2}\,H_{l_1}^p\,L_{l_1,p}^{l_1}
}_{ \octagon \! \! \! \! \! \tiny{\sf x}} 
+ \\
\endaligned
$$
$$
\small
\aligned
& 
+
\underline{ 
2\, \sum_k\,L_{k,l_2,y^{l_1}}^{l_2}\,G^k
}_{ \fbox{\tiny 17}}
-
\underline{ 
2\,\sum_k\,L_{k,l_1,y^{l_2}}^{l_2}\,G^k
}_{ \fbox{\tiny 18}}
-
\underline{ 
2\,\sum_k\,M_{k,l_1,x}\,G^k
}_{ \fbox{\tiny 19}}
+ \\
& \
\ \ \ \ \
+
\underline{ 
\sum_k\,G^k\,H_{l_1}^{l_2}\,M_{l_1,l_2}
}_{ \octagon \! \! \! \! \! \tiny{\sf o}}
-
\underline{ 
\sum_k\,G^k\,H_{l_2}^{l_2}\,M_{k,l_1}
}_{ \octagon \! \! \! \! \! \tiny{\sf u}}
+
\underline{ 
\sum_p\,G^{l_2}\,H_{l_1}^p\,M_{l_2,p}
}_{ \octagon \! \! \! \! \! \! \tiny{\sf m}}
- \\
& \
\ \ \ \ \ 
-
\underline{ 
\sum_p\, G^{l_2}\,H_{l_2}^p\,M_{l_1,p}
}_{ \octagon \! \! \! \! \! \tiny{\sf s}}
-
\underline{
\sum_k\,\sum_p\,G^k\,H_k^p\,M_{l_1,p}
}_{ \fbox{\tiny 20}}
+
\underline{ 
2\,\sum_k\,\sum_p\,G^k\,L_{k,l_1}^p\,L_{l_2,p}^{l_2}
}_{ \octagon \! \! \! \! \! \tiny{\sf g}}
- \\
& \
\ \ \ \ \ 
-
\underline{ 
2\,\sum_k\,\sum_p\,G^k\,L_{k,l_2}^p\,L_{l_1,p}^{l_2}
}_{ \octagon \! \! \! \! \! \tiny{\sf k}} \ .
\endaligned
$$

\medskip

In conclusion, there is exact coincidence with the subgoal~\thetag{
4.29}. The proof that the first family $(3.112)_1$ of compatibility
conditions of the second auxiliary system~\thetag{ 3.99}, \thetag{
3.100}, \thetag{ 3.101} and~\thetag{ 3.102} are a consequence of (I),
(II), (III) and (IV) of Theorem~1.7~{\bf (3)} 
is complete. Granted that the
treatment of the other three families of compatibility conditions
$(3.112)_2$, $(3.112)_3$ and $(3.112)_4$ is similar (and as well
painful), we consider that the proof of
the equivalence between~{\bf (1)} and~{\bf (3)}
in Theorem~1.7 is complete, now.
\qed

\section*{\S5.~General form of the point transformation \\
of the free particle system}
This section is devoted to the exposition of a complete proof of
Lemma~3.32. To start with, we must develope the fundamental
equations~\thetag{ 3.10}, for $j=1,\dots,m$. Recalling that the total
differentiation operator is given by $D = \frac{\partial }{\partial x}+
\sum_{l_1=1}^m\, y_x^{ l_1} \cdot \frac{\partial }{\partial y^{l_1}}+
\sum_{l_1=1}^m\, y_{ xx}^{ l_1} \cdot \frac{ \partial }{\partial
y_x^{l_1}}$, we compute first
\def\theequation{5.1}\begin{equation}
\left\{
\aligned
DDX
& \
= 
D \, \left[
X_x+\sum_{l_1=1}^m\, 
y_x^{l_1} \cdot X_{y^{l_1}}
\right] \\
& \
=
X_{xx}+ 2\, 
\sum_{l_1=1}^m\, 
y_x^{l_1} \cdot X_{xy^{l_1}}+
\sum_{l_1=1}^m\, \sum_{l_2=1}^m\, 
y_x^{l_1} \, y_x^{l_2} \cdot
X_{y^{l_1}y^{l_2}}+ 
\sum_{l_1=1}^m\, 
y_{xx}^{l_1} \cdot X_{y^{l_1}},
\endaligned\right.
\end{equation}
and
\def\theequation{5.2}\begin{equation}
\left\{
\aligned
DDY^j
& \
= 
D \, \left[
Y_x^j+\sum_{l_1=1}^m\, 
y_x^{l_1} \cdot Y_{y^{l_1}}^j
\right] \\
& \
=
Y_{xx}^j+ 2\, 
\sum_{l_1=1}^m\, 
y_x^{l_1} \cdot Y_{xy^{l_1}}^j+
\sum_{l_1=1}^m\, \sum_{l_2=1}^m\, 
y_x^{l_1} \, y_x^{l_2} \cdot
Y_{y^{l_1}y^{l_2}}^j+ 
\sum_{l_1=1}^m\, 
y_{xx}^{l_1} \cdot Y_{y^{l_1}}^j.
\endaligned\right.
\end{equation}
Now, we can develope the equation $0= -DY^j \cdot DD X+ DX \cdot
DDY^j$, which yields
\def\theequation{5.3}\begin{equation}
\aligned
0 = 
& \
-\left[
Y_x^j +\sum_{l_1=1}^m\, 
y_x^{l_1}\cdot Y_{y^{l_1}}^j
\right]\cdot
\left[
X_{xx}+
2\, \sum_{l_1=1}^m\, 
y_x^{l_1} \cdot X_{xy^{l_1}}+ \right. \\
& \ \ \ \ \ \ \ \ \ \ \ \ \
\left.
+
\sum_{l_1=1}^m\, \sum_{l_2=1}^m\, 
y_x^{l_1} \, y_x^{l_2}\cdot
X_{y^{l_1} y^{l_2}}+
\sum_{l_1=1}^m\, y_{xx}^{l_1}\cdot X_{y^{l_1}}
\right]+ \\
& \
+
\left[
X_x+\sum_{l_1=1}^m\, 
y_x^{l_1} \cdot X_{y^{l_1}}
\right]\cdot
\left[
Y_{xx}^j+
2\, \sum_{l_1=1}^m\, 
y_x^{l_1} \cdot Y_{xy^{l_1}}^j+ 
\right. \\
& \ \ \ \ \ \ \ \ \ \ \ \ \
\left.
+
\sum_{l_1=1}^m\, \sum_{l_2=1}^m\, 
y_x^{l_1} y_x^{l_2} \cdot Y_{y^{l_1}y^{l_2}}^j
+
\sum_{l_1=1}^m\, 
y_{xx}^{l_1} \cdot Y_{y^{l_1}}^j
\right]=
\endaligned
\end{equation}
$$
\aligned
=
& \
-X_{xx}\, Y_x^j+Y_{xx}^j\, X_x+ \\
& \
\ \ \ \ \
+\sum_{l_1=1}^m\,
y_x^{l_1} \cdot
\left[
-2\, X_{xy^{l_1}} \, Y_x^j+
2\, Y_{xy^{l_1}}^j\, X_x- \right. \\
& \
\ \ \ \ \
\ \ \ \ \
\ \ \ \ \
\ \ \ \ \
\ \ \ \ \
\ \ \ \ \
\ \ \ \ \
\ \ \ \ \
\left. 
- 
X_{xx}\, Y_{y^{l_1}}^j+
Y_{xx}^j\, X_{y^{l_1}}
\right] + \\
& \
\ \ \ \ \
+\sum_{l_1=1}^m\,
\sum_{l_2=1}^m\, 
y_x^{l_1} \, y_x^{l_2} \cdot
\left[
-X_{y^{l_1} y^{l_2}}\, Y_x^j+
Y_{y^{l_1}y^{l_2}}^j\, 
X_x- \right. \\
& \
\ \ \ \ \
\ \ \ \ \
\ \ \ \ \
\ \ \ \ \
\ \ \ \ \
\ \ \ \ \
\ \ \ \ \
\ \ \ \ \
\left.
-2\, X_{xy^{l_2}}\, Y_{y^{l_1}}^j+
2\, Y_{xy^{l_2}}^j\, X_{y^{l_1}}
\right]+ \\
& \
\ \ \ \ \
+\sum_{l_1=1}^m\,
\sum_{l_2=1}^m\, 
\sum_{l_3=1}^m\, 
y_x^{l_1} \, y_x^{l_2} y_x^{l_3} \cdot
\left[
-X_{y^{l_2}y^{l_3}}\, Y_{y^{l_1}}^j+
Y_{y^{l_2}y^{l_3}}^j\, X_{y^{l_1}}
\right]+ \\
& \
\ \ \ \ \
+ 
\sum_{l_1=1}^m\, 
y_{xx}^{l_1} \cdot 
\left[
-X_{y^{l_1}}\, Y_x^j+
Y_{y^{l_1}}^j\, X_x
\right]+ \\
& \
\ \ \ \ \
+
\sum_{l_1=1}^m\, 
\sum_{l_2=1}^m\, 
y_{xx}^{l_1}\, y_x^{l_2} \cdot
\left[
-X_{y^{l_1}}\, Y_{y^{l_2}}^j+
Y_{y^{l_1}}^j\, X_{y^{l_2}}
\right].
\endaligned
$$
The goal is to show that after solving these $m$
equations for $j=1,\dots,m$ with respect to 
the $y_{xx}^l$, $l=1,\dots, m$, one obtains the 
expression~\thetag{ 3.33} of Lemma~3.32, or equivalently, 
using the $\Delta$ notation instead of the square
notation, one obtains
\def\theequation{5.4}\begin{equation}
\left\{
\aligned
0 = 
& \
y_{xx}^j\cdot \Delta\left(
x \vert y^1 \vert \cdots \vert y^m
\right)+
\Delta \left(
x \vert y^1 \vert \cdots \vert^j xx \vert \cdots
\vert y^m
\right)+ \\
& \
+\sum_{l_1=1}^m\, 
y_x^{l_1} \cdot
\left[
2\, \Delta
\left(
x\vert y^1 \vert \cdots \vert^j xy^{l_1} \vert
\cdots \vert y^m
\right)- \right. \\
& \
\left. 
\ \ \ \ \ \ \ \ \ \ \ \ \ \ \ \ \ \ \ 
\ \ \ \ \ \ \ \ \ \ \ \ \ \ 
-
\delta_{l_1}^j\, 
\Delta\left(
xx \vert y^1 \vert \cdots \vert y^m
\right)
\right]+ \\
& \
+
\sum_{l_1=1}^m\, 
\sum_{l_2=1}^m\, 
y_x^{l_1} \, y_x^{l_2} \cdot
\left[
\Delta\left( x \vert y^1 \vert \cdots
\vert^j y^{l_1} y^{l_2} \vert \cdots
\vert y^m\right)- \right. \\
& \
\left. 
\ \ \ \ \ \ \ \ \ \ \ \ \ \ \ \ \ \ \ 
\ \ \ \ \ \ \ \ \ \ \ \ \ \ \ \ \ \ \ 
-
2\, \delta_{l_1}^j\, 
\Delta\left(xy^{l_2} \vert y^1 \vert \cdots \vert y^m\right) 
\right]+ \\
& \
+\sum_{l_1=1}^m\, \sum_{l_2=1}^m\, \sum_{l_3=1}^m\, 
y_x^{l_1} \, y_x^{l_2} \, y_x^{l_3}\, 
\left[
-\delta_{l_1}^j\, \Delta\left(
y^{l_2} y^{l_3} \vert y^1 \vert \cdots \vert y^m
\right)
\right].
\endaligned\right.
\end{equation}
Unfortunately, the equations~\thetag{ 5.3} are not
solved with respect to the $y_{xx}^j$, because in its last
line, we notice that the $y_{xx}^{l_2}$ are
mixed with the $y_x^{l_1}$. Consequently, we have
to solve
a linear system of $m$ equations with the
unknowns $y_{xx}^j$ of the form
\def\theequation{5.5}\begin{equation}
\left\{
\aligned
0 = 
& \
A^j+ 
\sum_{l_1=1}^m\, 
y_{xx}^{l_1} \cdot
\left[
-X_{y^{l_1}}\, Y_x^j+
Y_{y^{l_1}}^j\, X_x+
\sum_{l_2=1}^m\, 
y_x^{l_2}\cdot \left[
-X_{y^{l_1}}\, Y_{y^{l_2}}^j+
Y_{y^{l_1}}^j\, X_{y^{l_2}}
\right],
\right],
\endaligned\right.
\end{equation}
for $j=1,\dots,m$, where $A^j$ is an abbreviation for the terms
appearing in the lines 5, 6, 7, 8, 9 and 10 of~\thetag{ 5.3}, or even
more compactly, changing the index $j$ to the index
$k$
\def\theequation{5.6}\begin{equation}
\left\{
\aligned
0 = 
& \
A^k+ 
\sum_{l_1=1}^m\, 
y_{xx}^{l_1} \cdot B_{l_1}^k,
\endaligned\right.
\end{equation}
for $k=1,\dots,m$, where $B_{l_1}^k$ is an abbreviation for the terms
in the brackets in~\thetag{ 5.5}. 

Thanks to the assumption that the determinant~\thetag{ 3.2} is the
identity determinant at $(x,y)=(0,0)$, we deduce that the determinant
of the $m\times m$ matrix $(B_{l_1}^k)_{1\leqslant l_1 \leqslant m}^{1\leqslant k
\leqslant m}$ is also the identity determinant at $(x,y, y_x)= (0,0,0)$.
It follows that the determinant of the $m\times m$ matrix
$(B_{l_1}^k)_{1 \leqslant l_1 \leqslant m}^{ 1\leqslant k \leqslant m}$ is nonvanishing in
a neighborhood of the origin in the first order jet space.
Consequently, we can apply the rule of Cramer to solve the $y_{xx}^j$
explicitely interms of the $A^k$ and of the $B_{ l_1}^k$ as follows
\def\theequation{5.7}\begin{equation}
\left\{
\aligned
y_{xx}^j = 
& \
-
\frac{ 
\left\vert 
\begin{array}{ccccc}
B_1^1 & \cdots & A^1 & \cdots & B_m^1 \\
\cdots & \cdots & \cdots & \cdots & \cdots \\
B_1^m & \cdots & A^m& \cdots & B_m^m \\
\end{array}
\right\vert
}{
\left\vert
\begin{array}{ccccc}
B_1^1 & \cdots & B^1 & \cdots & B_m^1 \\
\cdots & \cdots & \cdots & \cdots & \cdots \\
B_1^m & \cdots & B^m& \cdots & B_m^m \\
\end{array}
\right\vert
}
\endaligned\right.
\end{equation}
where on the numerator, the only modification of the determinant of
the matrix $(B_{l_1}^k)_{ 1 \leqslant l_1 \leqslant m}^{ 1\leqslant k \leqslant m}$ is the
replacement of its $j$-th column by the column vector $A$. We have to
show that after replacing the $A^k$ and the $B_{l_1}^k$ by their
complete expressions, one indeed obtains the desired equation~\thetag{
5.4}. As in~\thetag{ 3.43}, we shall introduce a notation for the two
$m\times m$ determinants appearing in~\thetag{ 5.6}: we write this
quotient under the form
\def\theequation{5.8}\begin{equation}
\left\{
\aligned
y_{xx}^j = 
& \
-
\frac{ 
\left\vert \! 
\left\vert
B_1^k \vert \cdots \vert^j A^k \vert \cdots \vert B_m^k
\right\vert \!
\right\vert
}{ 
\left\vert \! 
\left\vert
B_1^k \vert \cdots \vert^j B_j^k \vert \cdots \vert B_m^k
\right\vert \!
\right\vert
},
\endaligned\right.
\end{equation}
where it is understood that $B_1^k,\dots,B_j^k,\dots,B_m^k$ and $A^k$
are column vectors whose index $k$ (for their lines) varies from $1$
to $m$. This representation of determinants emphasizing only its
columns will be appropriate for later manipulations.

Our first task is to compute the determinant in the denominator
of~\thetag{ 5.8}. Recalling that we make the notational
identification $y^0 \equiv x$, it will be convenient to reexpress the
$B_{l_1}^k$ in a slightly compacter form,
using the total differentiation operator $D$:
\def\theequation{5.9}\begin{equation}
\left\{
\aligned
B_{l_1}^k
& \
= 
-X_{y^{l_1}}\, Y_x^k+
Y_{y^{l_1}}^k\, X_x+
\sum_{l_2=1}^m\, 
y_x^{l_2}\cdot \left[
-X_{y^{l_1}}\, Y_{y^{l_2}}^k +
Y_{y^{l_1}}^k \, X_{y^{l_2}} \right] = \\
& \
=
Y_{y^{l_1}}^k\cdot DX- 
X_{y^{l_1}}\cdot DY^k.
\endaligned\right.
\end{equation}

\def\thelemma{5.10}\begin{lemma}
We have the following expression for the determinant of the matrix
$(B_{l_1}^k)_{ 1 \leqslant l_1 \leqslant m}^{ 1\leqslant k \leqslant m}${\rm :}
\def\theequation{5.11}\begin{equation}
\left\{
\aligned
{}
& \
\left\vert \! 
\left\vert
Y_{y^1}^k\cdot DX- X_{y^1} \cdot DY^k \vert
\cdots \vert
Y_{y^m}^k\cdot DX - X_{y^m} \cdot DY^k
\right\vert \!
\right\vert= \\
& \
=
[DX]^{m-1} \cdot \Delta\left(
x \vert y^1 \vert \cdots \vert y^m
\right).
\endaligned\right.
\end{equation}
\end{lemma}

\proof
By multilinearity, we may develope the determinant written in the
first line of~\thetag{ 5.11}. Since it contains two terms in each
columns, we should obtain a sum of $2^m$ determinants. However, since
the obtained determinants vanish as soon as the column $DY^k$
(multiplied by various factors $X_{y^l}$) appears at least two
different places, it remains only $(m+1)$ nonvanishing determinants,
those for which the column $DY^k$ appears at most once:
\def\theequation{5.12}\begin{equation}
\left\{
\aligned
{}
& \
\left\vert \! 
\left\vert
Y_{y^1}^k\cdot DX- X_{y^1} \cdot DY^k \vert
\cdots \vert
Y_{y^m}^k\cdot DX - X_{y^m} \cdot DY^k
\right\vert \!
\right\vert= \\
& \
=
[DX]^m\cdot
\left\vert \! 
\left\vert 
Y_{y^1}^k \vert \cdots \vert Y_{y^m}^k
\right\vert \!
\right\vert- 
[DX]^{m-1} \, X_{y^1} \cdot
\left\vert \! 
\left\vert 
DY^k \vert
Y_{y^2}^k \vert \cdots \vert Y_{y^m}^k
\right\vert \!
\right\vert-\cdots - \\
& \
\ \ \ \ \ \ 
- [DX]^{m-1} \, X_{y^m} \cdot
\left\vert \! 
\left\vert 
Y_{y^1}^k \vert \cdots \vert Y_{y^{m-1}}^k
\vert DY^k
\right\vert \!
\right\vert.
\endaligned\right.
\end{equation} 
To establish the desired expression appearing in the second line
of~\thetag{ 5.11}, we factor out by $[DX]^{ m-1}$ and we develope all
the remaining total differentiation operators $D$. 
Since $y^0 \equiv x$, we have $y_x^0 = 1$, and this enables
us to contract $X_x + \sum_{ l_1}^m \, 
y_x^{ l_1} \, X_{ y^{ l_1}}$ as
$\sum_{ l_1= 0}^m \, y_x^{ l_1} \, X_{ y^{ l_1}}$. 
So, we achieve the following
computation (further explanations and comments just afterwards):
\def\theequation{5.13}\begin{equation}
\left\{
\aligned
{}
& \
= 
[DX]^{m-1} \cdot
\left\{
\sum_{l_1=0}^m\, 
y_x^{l_1} \, X_{y^{l_1}}\cdot
\left\vert \! 
\left\vert 
Y_{y^1}^k \vert \cdots \vert Y_{y^m}^k
\right\vert \!
\right\vert- \right. \\
& \
\ \ \ \ \ \ \ \ \ \ \ \ \ \ \ \ \ \ \ \ \ \ \ \ \ 
\ \ \
\left.
-
\sum_{l_1=0}^m \, 
y_x^{l_1} \, X_{y^1} \cdot
\left\vert \! 
\left\vert 
Y_{y^{l_1}}^k \vert Y_{y^2}^k\vert 
\cdots \vert Y_{y^m}^k
\right\vert \!
\right\vert - \right. \\
& \
\ \ \ \ \ \ \ \ \ \ \ \ \ \ \ \ \ \ \ \ \ \ \ \ \ 
\ \ \
\left.
- \cdots - 
\sum_{l_1=0}^m\, y_x^{l_1} \, X_{y^m} \cdot
\left\vert \! 
\left\vert 
Y_{y^1}^k\vert 
\cdots \vert Y_{y^{m-1}}^k \vert
Y_{y^{l_1}}^k
\right\vert \!
\right\vert
\right\} \\
& \
=[DX]^{m-1} \cdot
\left\{
\sum_{l_1=0}^m\, 
y_x^{l_1} \, X_{y^{l_1}}\cdot
\left\vert \! 
\left\vert 
Y_{y^1}^k \vert \cdots \vert Y_{y^m}^k
\right\vert \!
\right\vert-
X_{y^1}\cdot
\left\vert \! 
\left\vert 
Y_x^k \vert Y_{y^2}^k \vert
\cdots \vert Y_{y^m}^k
\right\vert \!
\right\vert- 
\right. \\
& \
\ \ \ \ \ \ \ \ \ \ \ \ \ \ \ \ \ \ \ \ \ \ \ \ \ 
\ \ 
\left.
- y_x^1 \, X_{y^1} \cdot
\left\vert \! 
\left\vert 
Y_{y^1}^k \vert Y_{y^2}^k \vert
\cdots \vert Y_{y^m}^k
\right\vert \!
\right\vert- \cdots - \right . \\
& \
\ \ \ \ \ \ \ \ \ \ \ \ \ \ \ \ \ \ \ \ \ \ \ \ \ 
\ \ 
\left.
-X_{y^m}\cdot
\left\vert \! 
\left\vert 
Y_{y^1}^k \vert \cdots \vert Y_{y^{m-1}}^k 
\vert Y_x^k
\right\vert \!
\right\vert- 
y_x^m\, X_{y^m} \cdot
\left\vert \! 
\left\vert 
Y_{y^1}^k \vert \cdots \vert Y_{y^{m-1}}^k 
\vert Y_{y^m}^k
\right\vert \!
\right\vert
\right\} \\
& \
=
[DX]^{m-1} \cdot
\left\{
X_x\cdot 
\left\vert \! 
\left\vert 
Y_{y^1}^k \vert
\cdots \vert Y_{y^m}^k
\right\vert \!
\right\vert -
X_{y^1} \cdot
\left\vert \! 
\left\vert 
Y_x^k \vert
Y_{y^2}^k \vert
\cdots \vert Y_{y^m}^k
\right\vert \!
\right\vert 
- \cdots - \right. \\
& \
\left.
\ \ \ \ \ \ \ \ \ \ \ \ \ \ \ \ \ \ \ \ \ \ \ \ \ \ \ 
- X_{y^m} \cdot 
\left\vert \! 
\left\vert 
Y_{y^1}^k \vert
\cdots \vert Y_{y^{m-1}}^k \vert
Y_x^k
\right\vert \!
\right\vert
\right\} \\
& \
= 
[DX]^{m-1} \cdot \left\{
\Delta(x \vert y^1 \vert \cdots \vert y^m)
\right\}.
\endaligned\right.
\end{equation}
For the passage to the equality of line 4, using the fact that a
determinant having two identical columns vanishes, we observe that in
each of the $m$ sums $\sum_{ l_1 = 0}^m$ appearing in lines 2 and 3
(including the {\tt cdots}), there remains only two non-vanishing
determinants. For the passage to the equality of line 7, we just sum
up all the linear combinations of determinants appearing in lines 4, 5
and 6. Finally, for the passage to the equality of line 9, we
recognize the development of the fundamental Jacobian
determinant~\thetag{ 3.2} along its first line
$(X_x, X_{ y^1}, \dots, X_{ y^m})$, modulo some permutations of
columns in the $m\times m$ minors. The proof is complete.
\endproof

Our second task, similar but computationnally more heavy, is to
compute the determinant in the numerator of~\thetag{ 5.8}. First of
all, we have to re-express the $A^k$ defined implicitely
between~\thetag{ 5.3} and~\thetag{ 5.5} using the total
differentiation operator to contract them
as follows
\def\theequation{5.14}\begin{equation}
\left\{
\aligned
A^k =
& \
DX\cdot Y_{xx}^k- 
DY^k\cdot X_{xx}+
2\, \sum_{l_1=1}^m\, 
y_x^{l_1} \cdot \left[
DX \cdot Y_{xy^{l_1}}^k- 
DY^k\cdot X_{xy^{l_1}}
\right]+ \\
& \
+ 
\sum_{l_1=1}^m\, 
\sum_{l_2=1}^m\, 
y_x^{l_1} \, 
y_x^{l_2}\cdot
\left[
DX \cdot Y_{y^{l_1}y^{l_2}}^k- 
DY^k\cdot X_{y^{l_1}y^{l_2}}
\right].
\endaligned\right.
\end{equation}
Replacing this expression of $A^k$ in~\thetag{ 5.8}, taking account of
the expression of the denominator already obtained in the second line
of~\thetag{ 5.11} and abbreviating $\Delta( x\vert y^1 \vert \cdots
\vert y^m)$ as $\Delta$, we may write~\thetag{ 5.8} in length and then
develope it by linearity as follows
\def\theequation{5.15}\begin{equation}
\left\{
\aligned
y_{xx}^j = 
& \
\frac{ -1}{ [DX]^{m-1} \cdot \Delta} \cdot
\left[
\aligned
{}
& \
\left\vert \! 
\left\vert
Y_{y^1}^k\cdot DX - X_{y^1} \cdot DY^k 
\vert \cdots \right. \right. \\
& \
\left. \left. 
\ \
\cdots 
\vert^j \,
DX\cdot Y_{xx}^k - DY^k \cdot Y_{xx}^k + 
\right. \right. \\ 
& \
\left. \left. 
\ \ \ \ \ \ \ \ \
+2\, \sum_{l_1=1}^m\, 
y_x^{l_1} \cdot
\left[
DX\cdot Y_{xy^{l_1}}^k - 
DY^k \cdot X_{xy^{l_1}}
\right] + 
\right. \right. \\ 
& \
\left. \left. 
\ \ \ \ \ \ \ \ \
+
\sum_{l_1=1}^m \, 
\sum_{l_2=1}^m\, 
y_x^{l_1} \, y_x^{l_2} \, 
\cdot
\left[
DX\cdot Y_{y^{l_1} y^{l_2}}^k - 
DY^k \cdot X_{y^{l_1} y^{l_2}} \right]
\vert \cdots 
\right. \right. \\ 
& \
\left. \left.
\ \ \ \ \ \ \ \ \ \ \ \ \ \ 
\ \ \ \ \ \ \ \ \ \ \ \ \ \ 
\ \ \ \ \ \ \ \ \ \ \ \ \ \ 
\ \ \ \ \ \ \ \ \ \ \ \ \ \ 
\cdots \vert
Y_{y^m}^k \cdot DX - X_{y^m} \cdot DY^k
\right\vert \!
\right\vert
\endaligned
\right] \\
=
& \
\frac{ -1}{[DX]^{m-1} \cdot \Delta} \cdot
\left[
\aligned
{}
& \
\left\vert \! 
\left\vert 
Y_{y^1}^k\cdot DX- X_{y^1}\cdot DY^k 
\vert \cdots
\right. \right. \\ 
& \
\left. \left.
\ \ \ \ \ \
\cdots
\vert^j \,
DX\cdot Y_{xx}^k - DY^k \cdot X_{xx} \vert
\cdots
\right. \right. \\ 
& \
\left. \left.
\ \ \ \ \ \ \ \ \ \ \ \ \ \ 
\cdots
\vert 
Y_{y^m}^k \cdot
DX - X_{y^m} \cdot DY^k
\right\vert \!
\right\vert+ \\
& \
+2\, \sum_{l_1=1}^m\, 
y_x^{l_1} \cdot
\left\vert \! 
\left\vert 
Y_{y^1}^k\cdot DX - X_{y^1} \cdot DY^k \vert
\cdots
\right. \right. \\ 
& \
\left. \left.
\ \ \ \ \ \ \ 
\cdots
\vert^j\, 
DX\cdot Y_{xy^{l_1}}^k - 
DY^k\cdot X_{xy^{l_1}} \vert \cdots
\right. \right. \\ 
& \
\left. \left.
\ \ \ \ \ \ \ \ \ \ \ \ \ \ 
\cdots
\vert 
Y_{y^m}^k \cdot DX - X_{y^m} \cdot DY^k 
\right\vert \!
\right\vert+ \\& \
+
\sum_{l_1=1}^m\,
\sum_{l_2=1}^m\, 
y_x^{l_1}\, y_x^{l_2} \cdot
\left\vert \! 
\left\vert 
Y_{y^1}^k\cdot DX - X_{y^1} \cdot DY^k \vert
\cdots
\right. \right. \\ 
& \
\left. \left.
\ \ \ \ \ \ \ 
\cdots
\vert^j\, 
DX\cdot Y_{y^{l_1}y^{l_2}}^k - 
DY^k\cdot X_{y^{l_1}y^{l_2}} \vert \cdots
\right. \right. \\ 
& \
\left. \left.
\ \ \ \ \ \ \ \ \ \ \ \ \ \ 
\ \ \ \ \ \ \ \ \ \ \ \ \ \ \ \ \ \ 
\cdots
\vert 
Y_{y^m}^k \cdot DX - X_{y^m} \cdot DY^k 
\right\vert \!
\right\vert \\
\endaligned
\right].
\endaligned\right.
\end{equation}
As it is delicate to read, let us say that lines 2, 3 and 4 just
express the $j$-th colum $\vert^j \, A^k \vert$ of the determinant
$\vert \! \vert B_1 k \vert \cdots \vert j \, A^k \vert \cdots \vert
B_m^k \vert \! \vert$, after replacement of $A^k$ by its complete
expression~\thetag{ 5.14}.

In lines 6, 7, 8; in lines 9, 10, 11; and in lines 12, 13, 14, there
are three families of $m\times m$ determinants containing a linear
combination (soustraction) having exactly two terms in each column. As
in the proof of Lemma~5.10, by multilinarity, we have to develope each
such determinant. In principle, for each development, we should get
$2^m$ terms, but since the obtained determinants vanish as soon as the
column $DY^k$ (modulo a multiplication by some factor) appears at
least twice, it remains only $(m+1)$ nonvanishing determinants, those
for which the column $DY^k$ appears at most once. In addition, for
each of the obtained determinant, the factor $[DX]^{m-1}$ appears
(sometimes even the factor $[DX]^m$), so that this factor compensates
the factor $[DX]^{m-1}$ in the numerator. In sum, the continuation of
the huge computation yields:
\def\theequation{5.16}\begin{equation}
y_{ xx}^j
= 
-\frac{1}{\Delta} \cdot
\left[
\aligned
{}
&
DX\cdot
\left\vert \! 
\left\vert
Y_{y^1}^k \vert \cdots \vert^j\, 
Y_{xx}^k \vert \cdots \vert Y_{y^m}^k
\right\vert \!
\right\vert- \\
& \
-
X_{y^1} \cdot
\left\vert \! 
\left\vert
DY^k \vert \cdots \vert^j\, 
Y_{xx}^k\vert \cdots \vert Y_{y^m}^k
\right\vert \!
\right\vert-\cdots - \\
& \
-
X_{xx}\cdot
\left\vert \! 
\left\vert
Y_{y^1}^k \vert \cdots \vert^j\, 
DY^k \vert \cdots \vert Y_{y^m}^k
\right\vert \!
\right\vert - \cdots - \\
& \
-
X_{y^m} \cdot 
\left\vert \! 
\left\vert
Y_{y^1}^k \vert \cdots \vert^j\, 
Y_{xx}^k \vert \cdots \vert DY^k
\right\vert \!
\right\vert + \\
& \
+2\, \sum_{l_1=1}^m\, 
y_x^{l_1} \, DX \cdot
\left\vert \! 
\left\vert
Y_{y^1}^k \vert \cdots \vert^j\, 
Y_{xy^{l_1}}^k \vert \cdots \vert Y_{y^m}^k
\right\vert \!
\right\vert
- \\
& \
-
2\, \sum_{l_1=1}^m\, 
y_x^{l_1} \, X_{y^1} \cdot
\left\vert \! 
\left\vert
DY^k \vert \cdots \vert^j\, 
Y_{xy^{l_1}}^k \vert \cdots \vert Y_{y^m}^k
\right\vert \!
\right\vert - \cdots - \\
& \
-
2\, \sum_{l_1=1}^m\, 
y_x^{l_1} \, X_{xy^{l_1}} \cdot
\left\vert \! 
\left\vert
Y_{y^1}^k \vert \cdots \vert^j\, 
DY^k \vert \cdots \vert Y_{y^m}^k
\right\vert \!
\right\vert - \cdots - \\
& \
-2\, \sum_{l_1=1}^m\, 
y_x^{l_1} \, X_{y^m} \cdot
\left\vert \! 
\left\vert
Y_{y^1}^k \vert \cdots \vert^j\, 
Y_{xy^{l_1}}^k \vert \cdots \vert DY^k
\right\vert \!
\right\vert+ \\
& \
+
\sum_{l_1=1}^m\, 
\sum_{l_2=1}^m\, 
y_x^{l_1} \, 
y_x^{l_2} \, 
DX \cdot
\left\vert \! 
\left\vert
Y_{y^1}^k \vert \cdots \vert^j\, 
Y_{y^{l_1}y^{l_2}}^k \vert \cdots \vert Y_{y^m}^k
\right\vert \!
\right\vert- \\
& \
-
\sum_{l_1=1}^m\, 
\sum_{l_2=1}^m\, 
y_x^{l_1} \, 
y_x^{l_2} \, 
X_{y^1} \cdot
\left\vert \! 
\left\vert
DY^k \vert \cdots \vert^j\, 
Y_{y^{l_1}y^{l_2}}^k \vert \cdots \vert Y_{y^m}^k
\right\vert \!
\right\vert - \cdots - \\
& \
-
\sum_{l_1=1}^m\, 
\sum_{l_2=1}^m\, 
y_x^{l_1} \, 
y_x^{l_2} \, 
X_{y^{l_1}y^{l_2}} \cdot
\left\vert \! 
\left\vert
Y_{y^1}^k \vert \cdots \vert^j\, 
DY^k \vert \cdots \vert Y_{y^m}^k
\right\vert \!
\right\vert - \cdots - \\
& \
-
\sum_{l_1=1}^m\, 
\sum_{l_2=1}^m\, 
y_x^{l_1} \, 
y_x^{l_2} \, 
X_{y^m} \cdot
\left\vert \! 
\left\vert
Y_{y^1}^k \vert \cdots \vert^j\, 
Y_{y^{l_1}y^{l_2}}^k \vert \cdots \vert DY^k
\right\vert \!
\right\vert
\endaligned
\right].
\end{equation}
To establish the desired expression~\thetag{ 5.4}, we must develope
all the total differentiation operators $D$ of the terms $DX$ placed
as factor and of the terms $DY^k$ placed in various columns of
determinants. We notice that in developing $DY^k$, we obtain columns
$Y_{y^l}^k$ (multiplied by the factor $y_x^l$) and for {\it only}
three (or two) values of $l= 0, 1, \dots, m$, this column does not
already appear in the corresponding determinant, so that $(m-1)$
determinants vanish and only $3$ (or $2$) remain nonzero. Taking
account of these simplifications, we have the continuation
\def\theequation{5.17}\begin{equation}
-y_{xx}^j \cdot \Delta= {\rm I}+ {\rm II} + 
{\rm III}, 
\end{equation}
where the term I is the development of lines 1, 2, 3, 4 of~\thetag{
5.16}; the term II is the development of lines 5, 6, 7, 8 of~\thetag{
5.16}; and the term III is the development of lines 9, 10, 11, 12
of~\thetag{ 5.16}. So we get firstly (further explanations follows):
\def\theequation{5.18}\begin{equation}
{\rm I}:=
\aligned
{}
& \
\underline{ \sum_{l=0}^m\, 
y_x^l \, X_{y^l} \cdot
\left\vert \!
\left\vert
Y_{y^1}^k \vert \cdots \vert^j \, 
Y_{xx}^k \vert \cdots \vert
Y_{y^m}^k
\right\vert \! 
\right\vert}_{\tiny \fbox{ 1}}
- \\
& \
-
X_{y^1} \cdot
\left\vert \!
\left\vert
Y_x^k \vert \cdots \vert^j \, 
Y_{xx}^k \vert \cdots \vert
Y_{y^m}^k
\right\vert \! 
\right\vert
- \\
& \
-
\underline{y_x^1\, X_{y^1} \cdot
\left\vert \!
\left\vert
Y_{y^1}^k \vert \cdots \vert^j \, 
Y_{xx}^k \vert \cdots \vert
Y_{y^m}^k
\right\vert \! 
\right\vert}_{\tiny \fbox{ 1}} 
- \\
& \
-
y_x^j\, X_{y^j} \cdot
\left\vert \!
\left\vert
Y_{y^j}^k \vert \cdots \vert^j \, 
Y_{xx}^k \vert \cdots \vert
Y_{y^m}^k
\right\vert \! 
\right\vert -\cdots - \\
& \
-
X_{xx}\cdot
\left\vert \!
\left\vert
Y_{y^1}^k \vert \cdots \vert^j \, 
Y_x^k \vert \cdots \vert
Y_{y^m}^k
\right\vert \! 
\right\vert - \\
& \
-
y_x^j \, X_{xx} \cdot
\left\vert \!
\left\vert
Y_{y^1}^k \vert \cdots \vert^j \, 
Y_{y^j}^k \vert \cdots \vert
Y_{y^m}^k
\right\vert \! 
\right\vert - \cdots - \\
& \
-
X_{y^m}\cdot
\left\vert \!
\left\vert
Y_{y^1}^k \vert \cdots \vert^j \, 
Y_{xx}^k \vert \cdots \vert
Y_x^k
\right\vert \! 
\right\vert - \\
& \
-
\underline{y_x^m \, X_{y^m} \cdot
\left\vert \!
\left\vert
Y_{y^1}^k \vert \cdots \vert^j \, 
Y_{xx}^k \vert \cdots \vert
Y_{y^m}^k
\right\vert \! 
\right\vert}_{\tiny \fbox{1}} - \\
& \
-y_x^j\, X_{y^m} \cdot
\left\vert \!
\left\vert
Y_{y^1}^k \vert \cdots \vert^j \, 
Y_{xx}^k \vert \cdots \vert
Y_{y^j}^k
\right\vert \! 
\right\vert,
\endaligned
\end{equation}
and secondly (we discuss afterwards the 
annihilation of the underlined terms):
\def\theequation{5.19}\begin{equation}
{\rm II}:=
\aligned
{}
& \
\underline{ 2\, \sum_{l_1=1}^m\, 
\sum_{l=0}^m\, 
y_x^{l_1} \, y_x^l\, 
X_{y^l}\cdot
\left\vert \!
\left\vert
Y_{y^1}^k \vert \cdots \vert^j \, 
Y_{xy^{l_1}}^k \vert \cdots \vert
Y_{y^m}^k
\right\vert \! 
\right\vert}_{\tiny \fbox{ 2}} - \\
& \
-
2\, \sum_{l_1=1}^m\, 
y_x^{l_1}\, X_{y^1} \cdot
\left\vert \!
\left\vert
Y_x^k \vert \cdots \vert^j \, 
Y_{xy^{l_1}}^k \vert \cdots \vert
Y_{y^m}^k
\right\vert \! 
\right\vert - \\
& \
- 
\underline{ 2\, \sum_{l_1=1}^m\, 
y_x^{l_1} \, y_x^1 \, X_{y^1} \cdot
\left\vert \!
\left\vert
Y_{y^1}^k \vert \cdots \vert^j \, 
Y_{xy^{l_1}}^k \vert \cdots \vert
Y_{y^m}^k 
\right\vert \! 
\right\vert}_{\tiny \fbox{ 2}} - \\
& \
-
2\, \sum_{l_1=1}^m\, 
y_x^{l_1} \, y_x^j \, X_{y^1} \cdot
\left\vert \!
\left\vert
Y_{y^j}^k \vert \cdots \vert^j \, 
Y_{xy^{l_1}}^k \vert \cdots \vert
Y_{y^m}^k 
\right\vert \! 
\right\vert - \cdots - \\
& \
-
2\, \sum_{l_1=1}^m\, 
y_x^{l_1} \, X_{xy^{l_1}} \cdot
\left\vert \!
\left\vert
Y_{y^1}^k \vert \cdots \vert^j \, 
Y_x^k \vert \cdots \vert
Y_{y^m}^k 
\right\vert \! 
\right\vert - \\
& \
-
2\, \sum_{l_1=1}^m\, 
y_x^{l_1}\, y_x^j \, X_{xy^{l_1}} \cdot
\left\vert \!
\left\vert
Y_{y^1}^k \vert \cdots \vert^j \, 
Y_{y^j}^k \vert \cdots \vert
Y_{y^m}^k 
\right\vert \! 
\right\vert - \cdots - \\
& \
-
2\, \sum_{l_1=1}^m\, 
y_x^{l_1}\, X_{y^m} \cdot
\left\vert \!
\left\vert
Y_{y^1}^k \vert \cdots \vert^j \, 
Y_{xy^{l_1}}^k \vert \cdots \vert
Y_x^k 
\right\vert \! 
\right\vert - \\
& \
-
\underline{ 2\, \sum_{l_1=1}^m\, 
y_x^{l_1}\, y_x^m\, X_{y^m} \cdot
\left\vert \!
\left\vert
Y_{y^1}^k \vert \cdots \vert^j \, 
Y_{xy^{l_1}}^k \vert \cdots \vert
Y_{y^m}^k 
\right\vert \! 
\right\vert}_{\tiny \fbox{ 2}} - \\
& \
-
2\, \sum_{l_1=1}^m\, 
y_x^{l_1}\, y_x^j \, X_{y^m} \cdot
\left\vert \!
\left\vert
Y_{y^1}^k \vert \cdots \vert^j \, 
Y_{xy^{l_1}}^k \vert \cdots \vert
Y_{y^j}^k 
\right\vert \! 
\right\vert,
\endaligned
\end{equation}
and where thirdly (we are nearly the end of the proof):
\def\theequation{5.20}\begin{equation}
{\rm III}:=
\aligned
{}
& \ 
\underline{\sum_{l_1=1}^m\, 
\sum_{l_2=1}^m\, 
\sum_{l=0}^m\, 
y_x^{l_1} \, y_x^{l_2} \, y_x^l\, 
X_{y^l} \cdot
\left\vert \!
\left\vert
Y_{y^1}^k \vert \cdots \vert^j \, 
Y_{y^{l_1}y^{l_2}}^k \vert \cdots \vert
Y_{y^m}^k 
\right\vert \! 
\right\vert}_{\tiny \fbox{ 3}} - \\
& \
-
\sum_{l_1=1}^m\, 
\sum_{l_2=1}^m\, 
y_x^{l_1} \, y_x^{l_2} \,
X_{y^1} \cdot
\left\vert \!
\left\vert
Y_x^k \vert \cdots \vert^j \, 
Y_{y^{l_1}y^{l_2}}^k \vert \cdots \vert
Y_{y^m}^k 
\right\vert \! 
\right\vert - \\
& \
-
\underline{ \sum_{l_1=1}^m\, 
\sum_{l_2=1}^m\, 
y_x^{l_1} \, y_x^{l_2} \, y_x^1 \,
X_{y^1} \cdot
\left\vert \!
\left\vert
Y_{y^1}^k \vert \cdots \vert^j \, 
Y_{y^{l_1}y^{l_2}}^k \vert \cdots \vert
Y_{y^m}^k 
\right\vert \! 
\right\vert}_{\tiny \fbox{ 3}} - \\
& \
-
\sum_{l_1=1}^m\, 
\sum_{l_2=1}^m\, 
y_x^{l_1} \, y_x^{l_2} \, y_x^j \,
X_{y^1} \cdot
\left\vert \!
\left\vert
Y_{y^j}^k \vert \cdots \vert^j \, 
Y_{y^{l_1}y^{l_2}}^k \vert \cdots \vert
Y_{y^m}^k 
\right\vert \! 
\right\vert - \cdots - \\
& \
-
\sum_{l_1=1}^m\, 
\sum_{l_2=1}^m\, 
y_x^{l_1} \, y_x^{l_2}
X_{y^{l_1}y^{l_2}} \cdot
\left\vert \!
\left\vert
Y_{y^1}^k \vert \cdots \vert^j \, 
Y_x^k \vert \cdots \vert
Y_{y^m}^k 
\right\vert \! 
\right\vert - \\
& \
-
\sum_{l_1=1}^m\, 
\sum_{l_2=1}^m\, 
y_x^{l_1} \, y_x^{l_2} \, y_x^j \, 
X_{y^{l_1}y^{l_2}} \cdot
\left\vert \!
\left\vert
Y_{y^1}^k \vert \cdots \vert^j \, 
Y_{y^j}^k \vert \cdots \vert
Y_{y^m}^k 
\right\vert \! 
\right\vert - \cdots - \\
& \
-
\sum_{l_1=1}^m\, 
\sum_{l_2=1}^m\, 
y_x^{l_1} \, y_x^{l_2} 
X_{y^m} \cdot
\left\vert \!
\left\vert
Y_{y^1}^k \vert \cdots \vert^j \, 
Y_{y^{l_1}y^{l_2}}^k \vert \cdots \vert
Y_x^k 
\right\vert \! 
\right\vert - \\
& \
-
\underline{ \sum_{l_1=1}^m\, 
\sum_{l_2=1}^m\, 
y_x^{l_1} \, y_x^{l_2} \, y_x^m 
X_{y^m} \cdot
\left\vert \!
\left\vert
Y_{y^1}^k \vert \cdots \vert^j \, 
Y_{y^{l_1}y^{l_2}}^k \vert \cdots \vert
Y_{y^m}^k 
\right\vert \! 
\right\vert}_{\tiny \fbox{ 3}} - \\
& \
-
\sum_{l_1=1}^m\, 
\sum_{l_2=1}^m\, 
y_x^{l_1} \, y_x^{l_2} \, y_x^j 
X_{y^m} \cdot
\left\vert \!
\left\vert
Y_{y^1}^k \vert \cdots \vert^j \, 
Y_{y^{l_1}y^{l_2}}^k \vert \cdots \vert
Y_{y^j}^k 
\right\vert \! 
\right\vert.
\endaligned
\end{equation}
Now, we explain the annihilation of the underlined terms. Consider I:
in the first sum $\sum_{ l=0}^m$, all the terms except only the two
corresponding to $l=0$ and to $l=j$ are annihilated by the other terms
with ${\tiny \fbox{ 1}}$ appended: indeed, one must take account of
the fact that in the expression of I, we have two sums represented by
some {\tt cdots}, the nature of which was defined without ambiguity in
the passage from~\thetag{ 5.15} to~\thetag{ 5.16}.

Similar simplifications occur for II and for III. Consequently, 
we obtain firstly:
\def\theequation{5.21}\begin{equation}
{\rm I}:=
\aligned
{}
& \
X_x \cdot
\left\vert \!
\left\vert
Y_{y^1}^k \vert \cdots \vert^j \, 
Y_{xx}^k \vert \cdots \vert
Y_{y^m}^k 
\right\vert \! 
\right\vert + \\
& \
+
y_x^j \, X_{y^j} \cdot
\left\vert \!
\left\vert
Y_{y^1}^k \vert \cdots \vert^j \, 
Y_{xx}^k \vert \cdots \vert
Y_{y^m}^k 
\right\vert \! 
\right\vert - \\
& \
-
X_{y^1} \cdot 
\left\vert \!
\left\vert
Y_x^k \vert \cdots \vert^j \, 
Y_{xx}^k \vert \cdots \vert
Y_{y^m}^k 
\right\vert \! 
\right\vert - \\
& \
-
y_x^j \, X_{y^j} \cdot
\left\vert \!
\left\vert
Y_{y^j}^k \vert \cdots \vert^j \, 
Y_{xx}^k \vert \cdots \vert
Y_{y^m}^k 
\right\vert \! 
\right\vert - \cdots - \\
& \
-
X_{xx}\cdot
\left\vert \!
\left\vert
Y_{y^1}^k \vert \cdots \vert^j \, 
Y_x^k \vert \cdots \vert
Y_{y^m}^k 
\right\vert \! 
\right\vert - \\
& \
- 
y_x^j\, X_{xx} \cdot
\left\vert \!
\left\vert
Y_{y^1}^k \vert \cdots \vert^j \, 
Y_{y^j}^k \vert \cdots \vert
Y_{y^m}^k 
\right\vert \! 
\right\vert - \cdots - \\
& \
-
X_{y^m} \cdot
\left\vert \!
\left\vert
Y_{y^1}^k \vert \cdots \vert^j \, 
Y_{xx}^k \vert \cdots \vert
Y_x^k 
\right\vert \! 
\right\vert - \\
& \
- 
y_x^j \, X_{y^m} \cdot
\left\vert \!
\left\vert
Y_{y^1}^k \vert \cdots \vert^j \, 
Y_{xx}^k \vert \cdots \vert
Y_{y^j}^k 
\right\vert \! 
\right\vert;
\endaligned
\end{equation}
just above, the first two lines consist of the
two terms in the sum underlined at the
first line of~\thetag{ 5.18}
which are not annihilated; 
secondly we obtain: 
\def\theequation{5.22}\begin{equation}
{\rm II}:=
\aligned
{}
& \
2\sum_{l_1=1}^m\, y_x^{l_1} \, X_x \cdot
\left\vert \!
\left\vert
Y_{y^1}^k \vert \cdots \vert^j \, 
Y_{xy^{l_1}}^k \vert \cdots \vert
Y_{y^m}^k 
\right\vert \! 
\right\vert+ \\
& \
+
2\sum_{l_1=1}^m\, y_x^{l_1} \, y_x^j \,
\, X_{y^j} \cdot
\left\vert \!
\left\vert
Y_{y^1}^k \vert \cdots \vert^j \, 
Y_{xy^{l_1}}^k \vert \cdots \vert
Y_{y^m}^k 
\right\vert \! 
\right\vert- \\
& \
-
2\sum_{l_1=1}^m\, y_x^{l_1} 
\, X_{y^1} \cdot
\left\vert \!
\left\vert
Y_x^k \vert \cdots \vert^j \, 
Y_{xy^{l_1}}^k \vert \cdots \vert
Y_{y^m}^k 
\right\vert \! 
\right\vert - \\
& \
-
2\sum_{l_1=1}^m\, y_x^{l_1} \, y_x^j \,
X_{y^1} \cdot
\left\vert \!
\left\vert
Y_{y^j}^k \vert \cdots \vert^j \, 
Y_{xy^{l_1}}^k \vert \cdots \vert
Y_{y^m}^k 
\right\vert \! 
\right\vert - \cdots - \\
& \
-
2\sum_{l_1=1}^m\, y_x^{l_1}
X_{xy^{l_1}} \cdot
\left\vert \!
\left\vert
Y_{y^1}^k \vert \cdots \vert^j \, 
Y_x^k \vert \cdots \vert
Y_{y^m}^k 
\right\vert \! 
\right\vert - \\
& \
-
2\sum_{l_1=1}^m\, y_x^{l_1} \, y_x^j \,
X_{xy^{l_1}} \cdot
\left\vert \!
\left\vert
Y_{y^1}^k \vert \cdots \vert^j \, 
Y_x^k \vert \cdots \vert
Y_{y^m}^k 
\right\vert \! 
\right\vert - \cdots - \\
& \
-
2\sum_{l_1=1}^m\, y_x^{l_1} \,
X_{y^m} \cdot
\left\vert \!
\left\vert
Y_{y^1}^k \vert \cdots \vert^j \, 
Y_{xy^{l_1}}^k \vert \cdots \vert
Y_x^k 
\right\vert \! 
\right\vert - \\
& \
-
2\sum_{l_1=1}^m\, y_x^{l_1} \, y_x^j \, 
X_{y^m} \cdot
\left\vert \!
\left\vert
Y_{y^1}^k \vert \cdots \vert^j \, 
Y_{xy^{l_1}}^k \vert \cdots \vert
Y_{y^j}^k 
\right\vert \! 
\right\vert; 
\endaligned
\end{equation}
similarly, the first two lines above consist of the
two terms in the sum underlined at the
first line of~\thetag{ 5.19}
which are not annihilated;
and thirdly we obtain:
\def\theequation{5.23}\begin{equation}
{\rm III}:=
\aligned
{}
& \ 
\sum_{l_1=1}^m\, 
\sum_{l_2=1}^m\, 
y_x^{l_1}\, y_x^{l_2} \, 
X_x \cdot
\left\vert \!
\left\vert
Y_{y^1}^k \vert \cdots \vert^j \, 
Y_{y^{l_1}y^{l_2}}^k \vert \cdots \vert
Y_{y^m}^k 
\right\vert \! 
\right\vert + \\
& \
+
\sum_{l_1=1}^m\, 
\sum_{l_2=1}^m\, 
y_x^{l_1}\, y_x^{l_2} \, y_x^j \,
X_{y^j} \cdot
\left\vert \!
\left\vert
Y_{y^1}^k \vert \cdots \vert^j \, 
Y_{y^{l_1}y^{l_2}}^k \vert \cdots \vert
Y_{y^m}^k 
\right\vert \! 
\right\vert - \\
& \
-
\sum_{l_1=1}^m\, 
\sum_{l_2=1}^m\, 
y_x^{l_1}\, y_x^{l_2} \, 
X_{y^1} \cdot
\left\vert \!
\left\vert
Y_x^k \vert \cdots \vert^j \, 
Y_{y^{l_1}y^{l_2}}^k \vert \cdots \vert
Y_{y^m}^k 
\right\vert \! 
\right\vert - \\
& \
-
\sum_{l_1=1}^m\, 
\sum_{l_2=1}^m\, 
y_x^{l_1}\, y_x^{l_2} \, y_x^j \,
X_{y^1} \cdot
\left\vert \!
\left\vert
Y_{y^j}^k \vert \cdots \vert^j \, 
Y_{y^{l_1}y^{l_2}}^k \vert \cdots \vert
Y_{y^m}^k 
\right\vert \! 
\right\vert - \cdots - \\
& \
-
\sum_{l_1=1}^m\, 
\sum_{l_2=1}^m\, 
y_x^{l_1}\, y_x^{l_2} \,
X_{y^{l_1}y^{l_2}} \cdot
\left\vert \!
\left\vert
Y_{y^1}^k \vert \cdots \vert^j \, 
Y_x^k \vert \cdots \vert
Y_{y^m}^k 
\right\vert \! 
\right\vert - \\
& \
-
\sum_{l_1=1}^m\, 
\sum_{l_2=1}^m\, 
y_x^{l_1}\, y_x^{l_2} \, y_x^j \,
X_{y^{l_1}y^{l_2}} \cdot
\left\vert \!
\left\vert
Y_{y^1}^k \vert \cdots \vert^j \, 
Y_{y^j}^k \vert \cdots \vert
Y_{y^m}^k 
\right\vert \! 
\right\vert - \cdots - \\
\endaligned
\end{equation}
$$
\aligned
& \
-
\sum_{l_1=1}^m\, 
\sum_{l_2=1}^m\, 
y_x^{l_1}\, y_x^{l_2} \,
X_{y^m} \cdot
\left\vert \!
\left\vert
Y_{y^1}^k \vert \cdots \vert^j \, 
Y_{y^{l_1}y^{l_2}}^k \vert \cdots \vert
Y_x^k 
\right\vert \! 
\right\vert - \\
& \
-
\sum_{l_1=1}^m\,
\sum_{l_2=1}^m\, 
y_x^{l_1}\, y_x^{l_2} \, y_x^j \, 
X_{y^m} \cdot
\left\vert \!
\left\vert
Y_{y^1}^k \vert \cdots \vert^j \, 
Y_{y^{l_1}y^{l_2}}^k \vert \cdots \vert
Y_{y^j}^k 
\right\vert \! 
\right\vert.
\endaligned
$$
Collecting the odd lines of~\thetag{ 5.21}, we obtain exactly $(m+1)$
terms which correspond to the development of the determinant $\Delta(x
\vert \cdots \vert^j \, xx \vert \cdots \vert y^m)$ along its first
line, modulo permutations of columns of the associated $m\times m$
minors; collecting the even lines of~\thetag{ 5.21}, we obtain exactly
$(m+1)$ terms which correspond to the development of the determinant
$-y_x^j \cdot \Delta(xx \vert y^1 \vert \cdots \vert y^m)$ along its
first lines, modulo permutations of columns of the associated
$m\times m$ minors. Similar observations hold about II and
III.

In , we may rewrite the final expressions of these three terms:
firstly
\def\theequation{5.24}\begin{equation}
\left\{
\aligned
{\rm I} = 
& \
\Delta( x \vert \cdots \vert^j\, 
xx \vert \cdots \vert y^m) -
y_x^j \cdot 
\Delta( xx \vert y^1 \vert \cdots \vert y^m), \\
{\rm II} = 
& \
2\sum_{l_1=1}^m\, 
y_x^{l_1} \cdot \Delta( x \vert
\cdots \vert^j\, xy^{l_1} \vert \cdots
\vert y^m)- \\
& \
\ \ \ \ \ 
-
2\, y_x^j \, 
\sum_{l_1=1}^m\, 
y_x^{l_1} \cdot \Delta(xy^{l_1} \vert y^1 \vert
\cdots \vert y^m), \\
{\rm III} =
& \
\sum_{l_1=1}^m\, 
\sum_{l_2=1}^m\, 
y_x^{l_1} \, 
y_x^{l_2} \cdot
\Delta( x \vert \cdots \vert^j\, 
y^{l_1} y^{l_2} \vert \cdots \vert y^m)- \\
& \
\ \ \ \ \ \ \ \ \ 
- 
y_x^j \, 
\sum_{l_1=1}^m\, 
\sum_{l_2=1}^m\, 
y_x^{l_1} \, 
y_x^{l_2} \cdot
\Delta( y^{l_1} y^{l_2} \vert y^1 \vert \cdots 
\vert y^m).
\endaligned\right.
\end{equation}
Coming back to~\thetag{ 5.17}, we
obtain the desired expression~\thetag{ 5.4}.

The proof of the\,\,---\,\,technical, though involving only
linear algebra\,\,---\,\,Lemma~3.32 is complete. \qed

\vfill\end{document}